\documentclass{ws-jktr}
\usepackage{graphicx}

\def\hiduke{ }

\newcommand{\byouga}[3]{
\begin{center}
\includegraphics[width=#2cm,clip]{#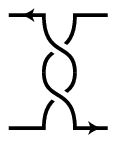}\vspace{0cm}
\end{center}
\begin{center} Fig. #3 
\end{center}
}

\newcommand{\inclf}[2]
{\includegraphics[width=#2cm,clip]{#1.eps}
}

\newcommand{\incskein}[1]{
\raisebox{-8pt}{
\includegraphics[height=20pt,clip]{#1.eps}
}}

\def\figyonichi{\byouga{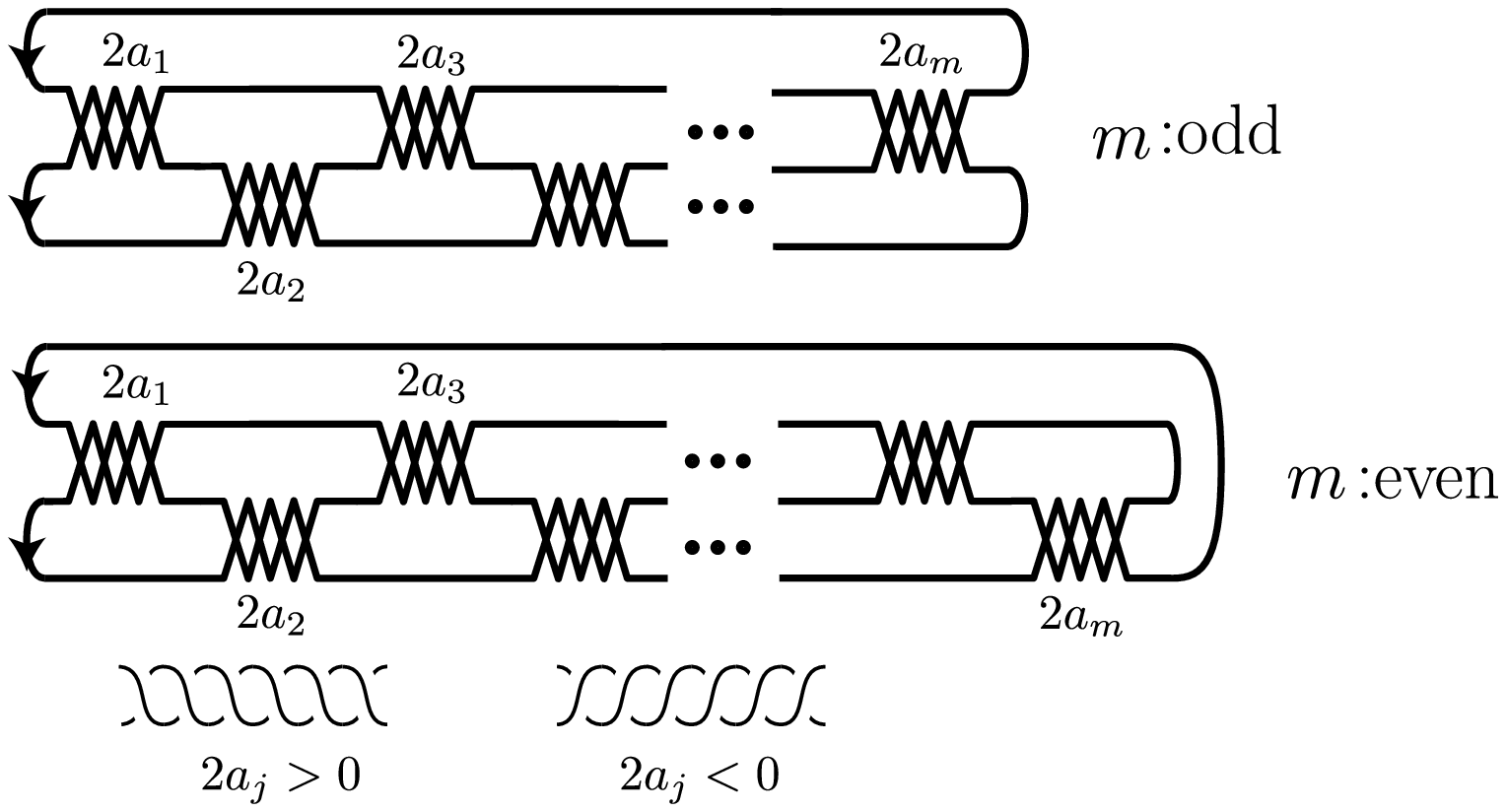}{7}{4.1}}
\def\figyonni{\byouga{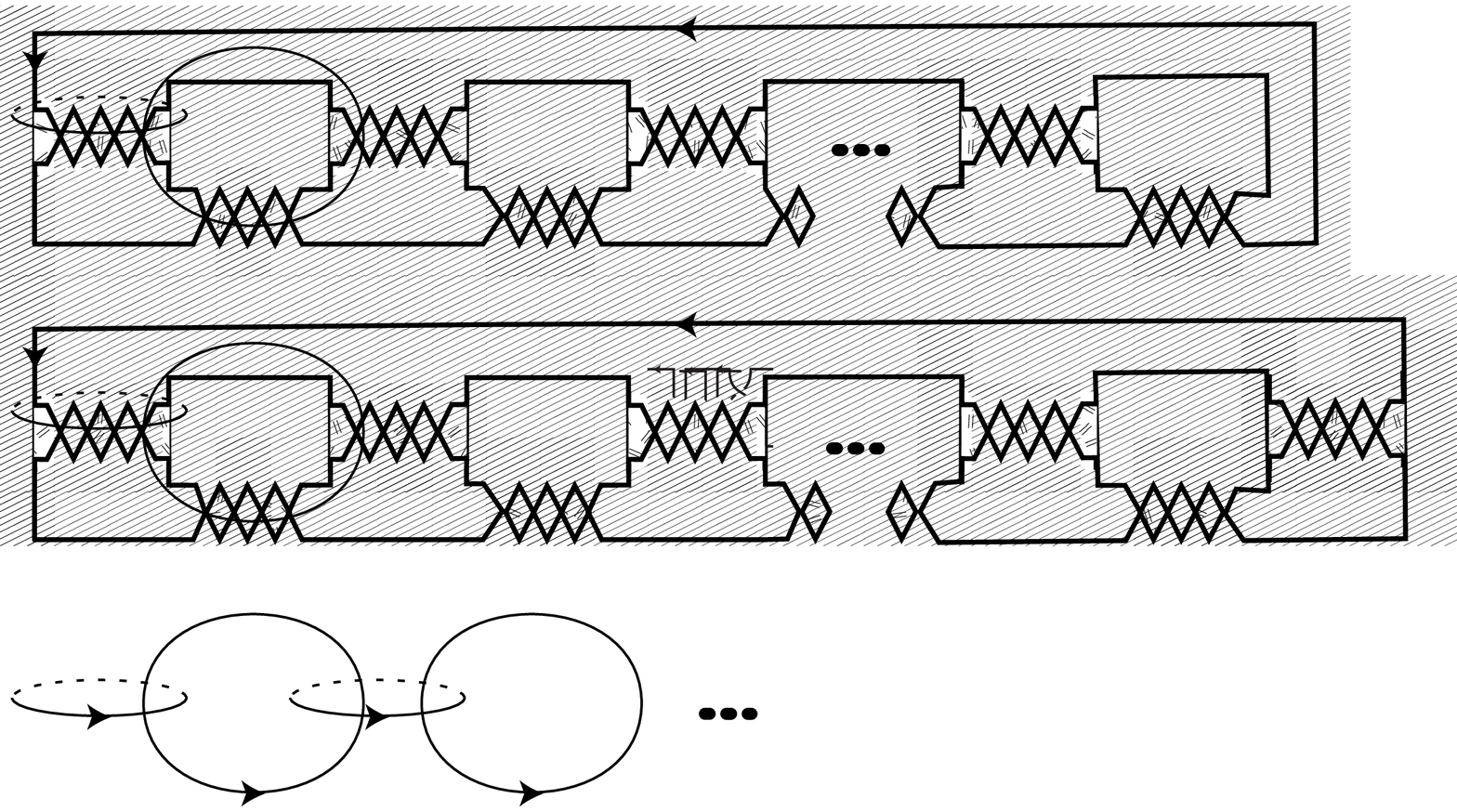}{6.8}{4.2}}
\def\yonsan{\byouga{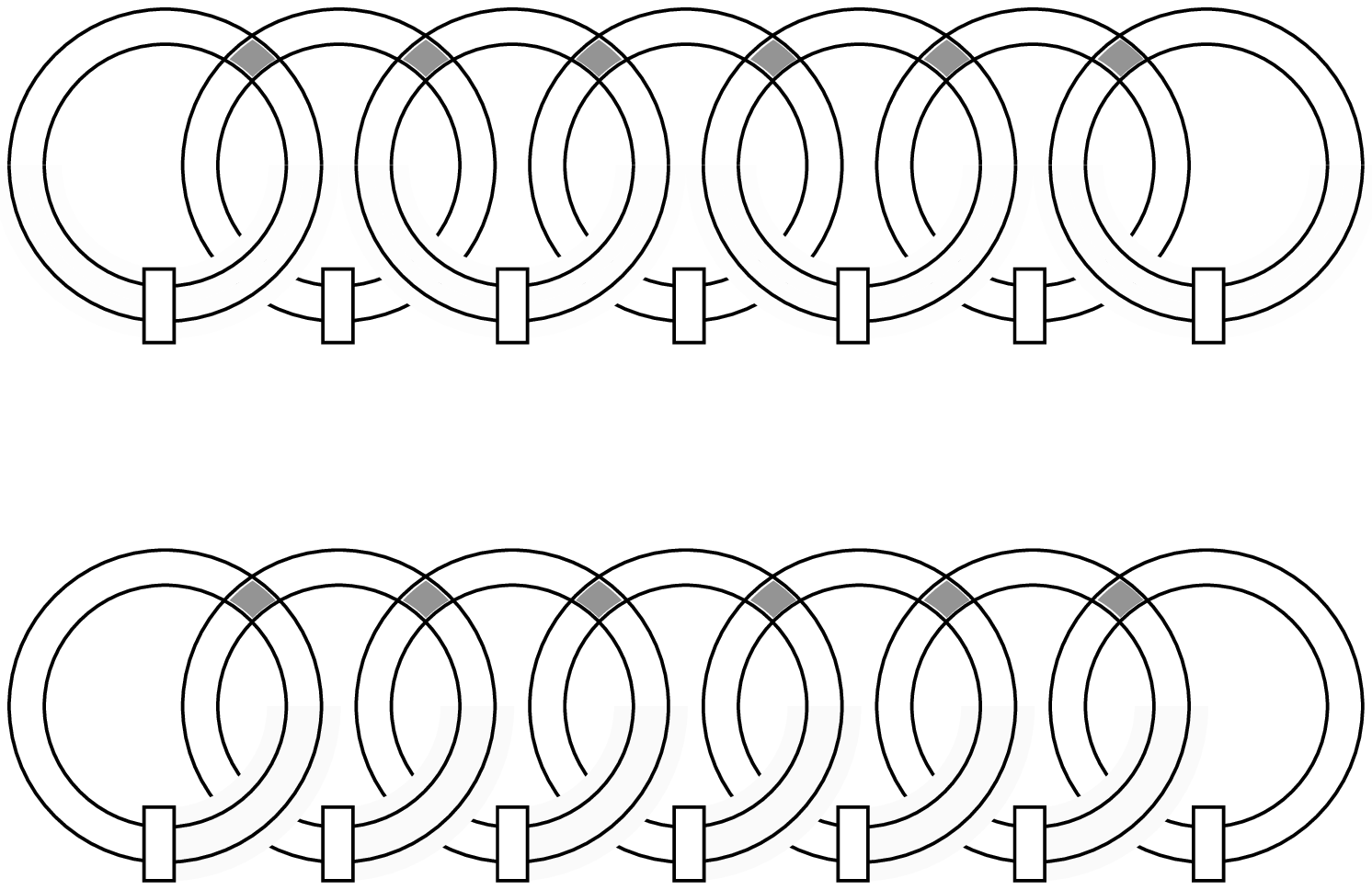}{6.0}{4.3}}
\def\rokuichi{\byouga{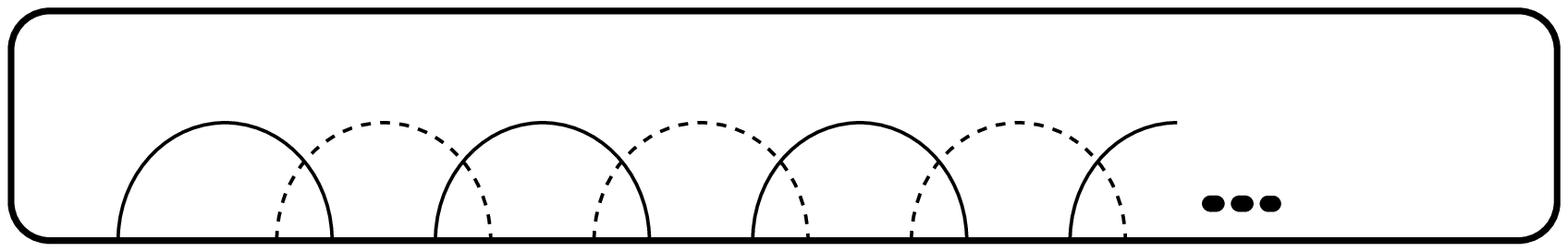}{7.5}{6.1: $r=[2a_1,2a_2,\dots]$}}
\def\rokuni{\byouga{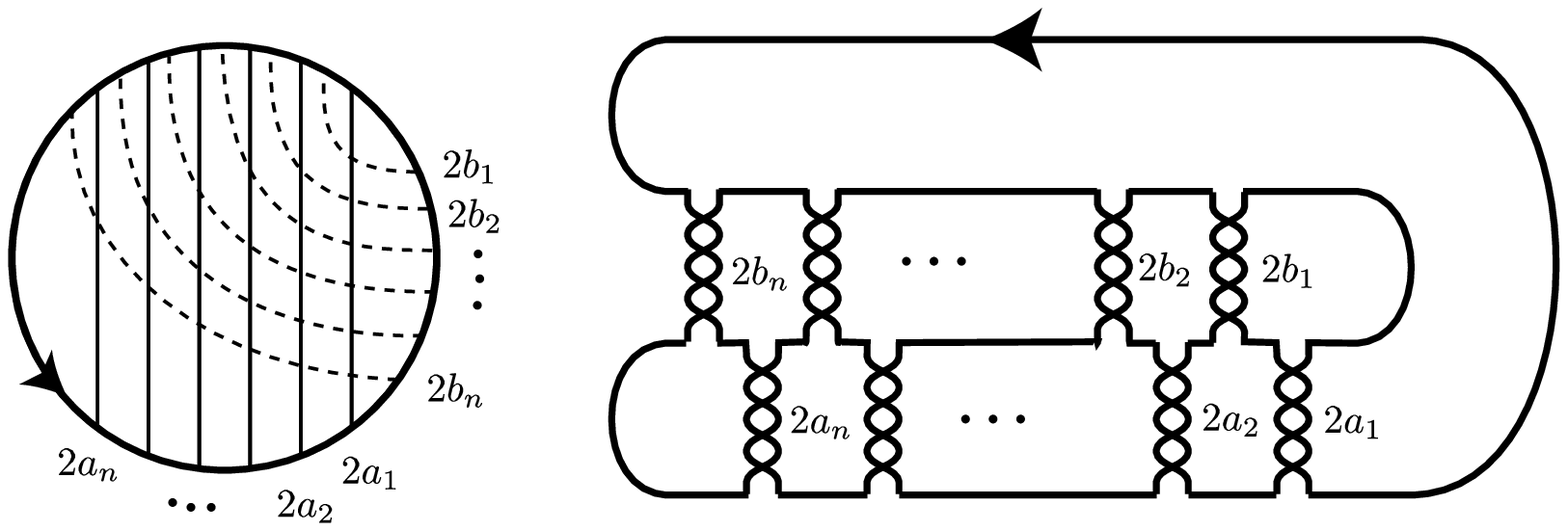}{9}{6.2}}
\def\rokusan{\byouga{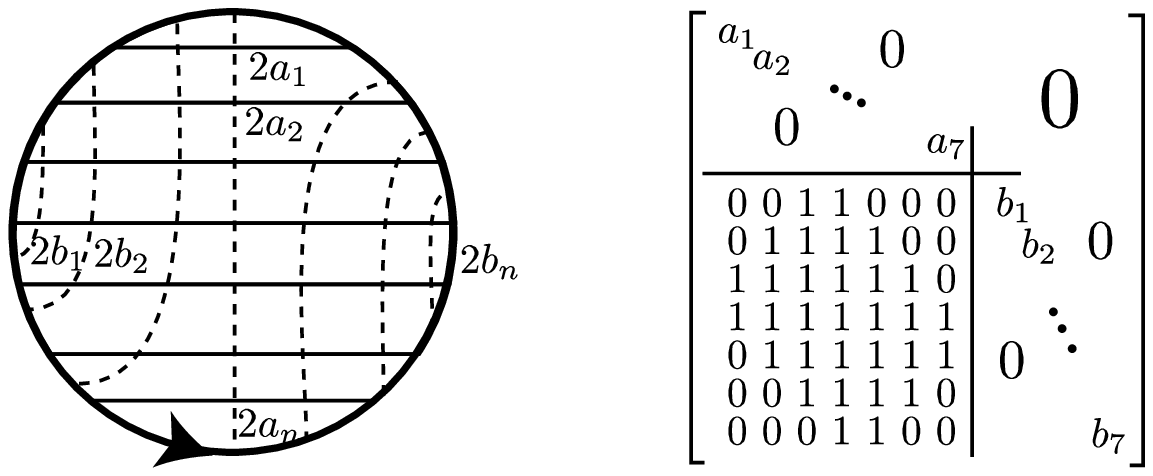}{7.5}{6.3}}
\def\nanaichi{\byouga{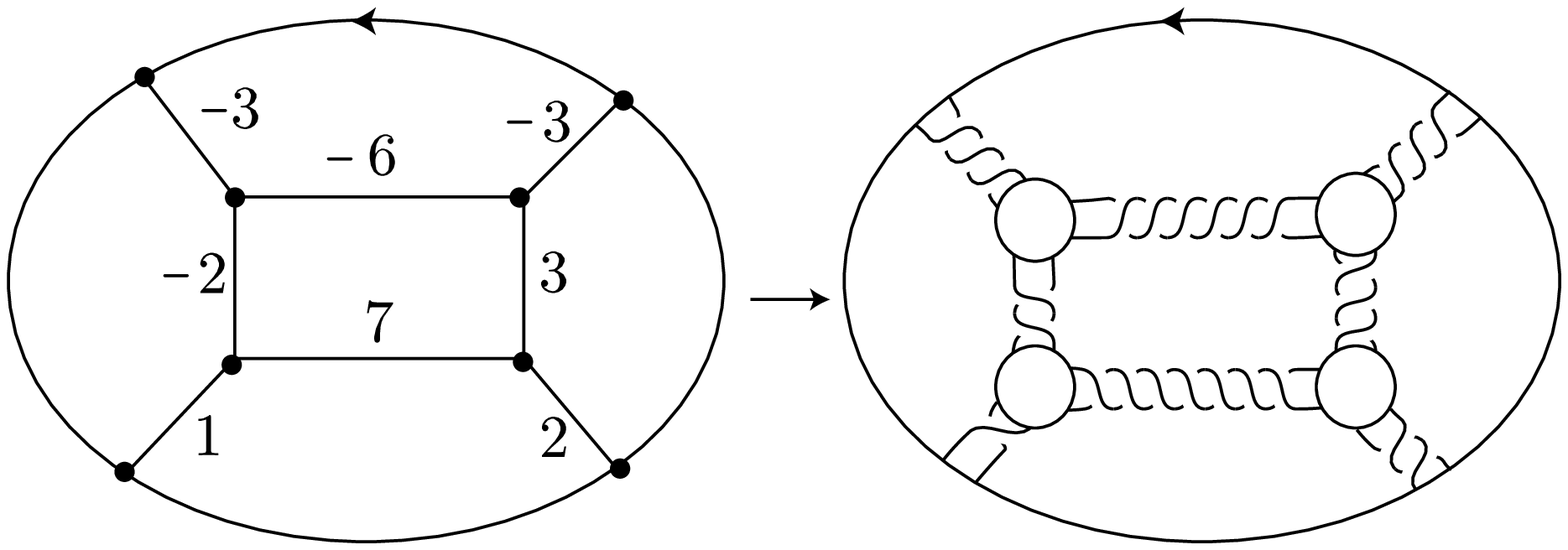}{8}{7.1}}
\def\nanani{\byouga{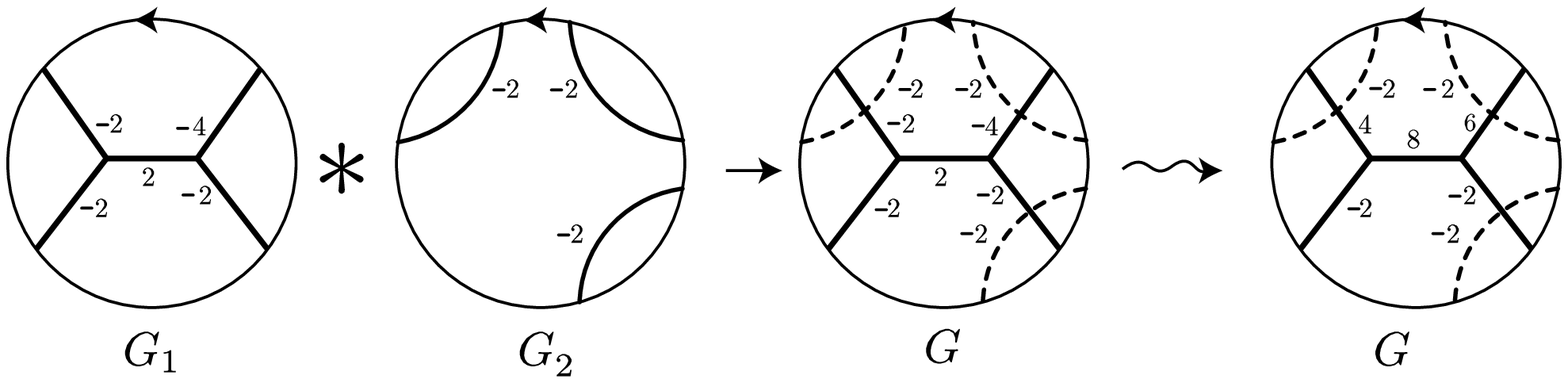}{12}{7.2}}
\def\kyuichi{\byouga{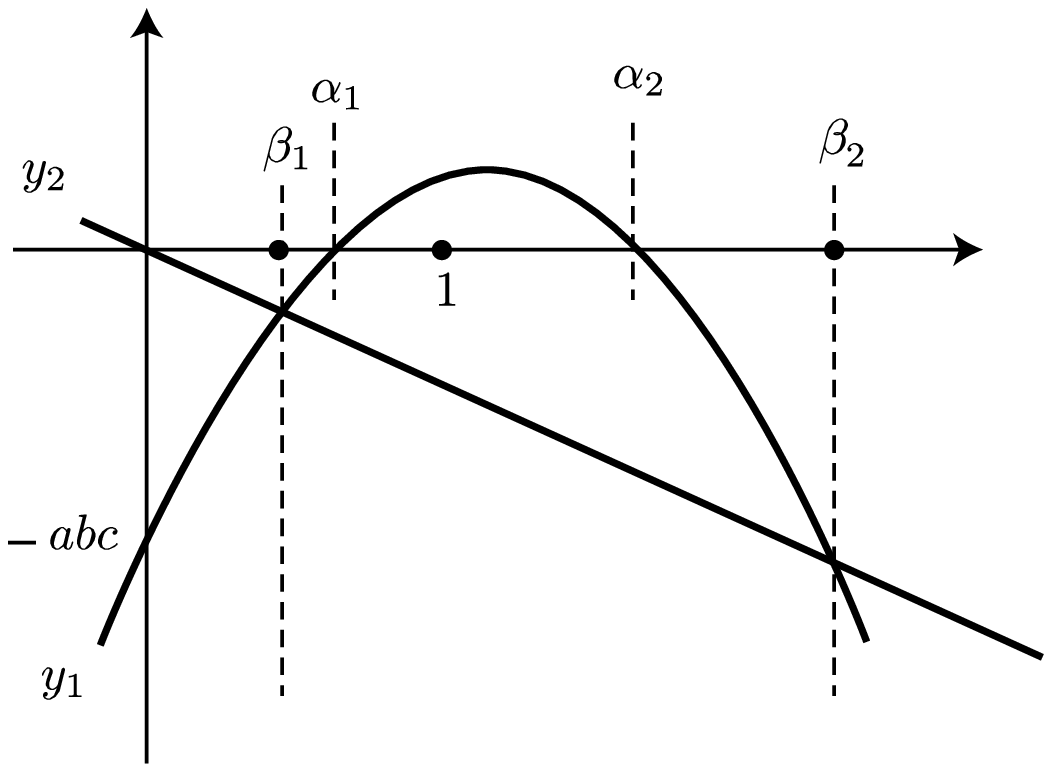}{6}{9.1}}
\def\kyuni{\byouga{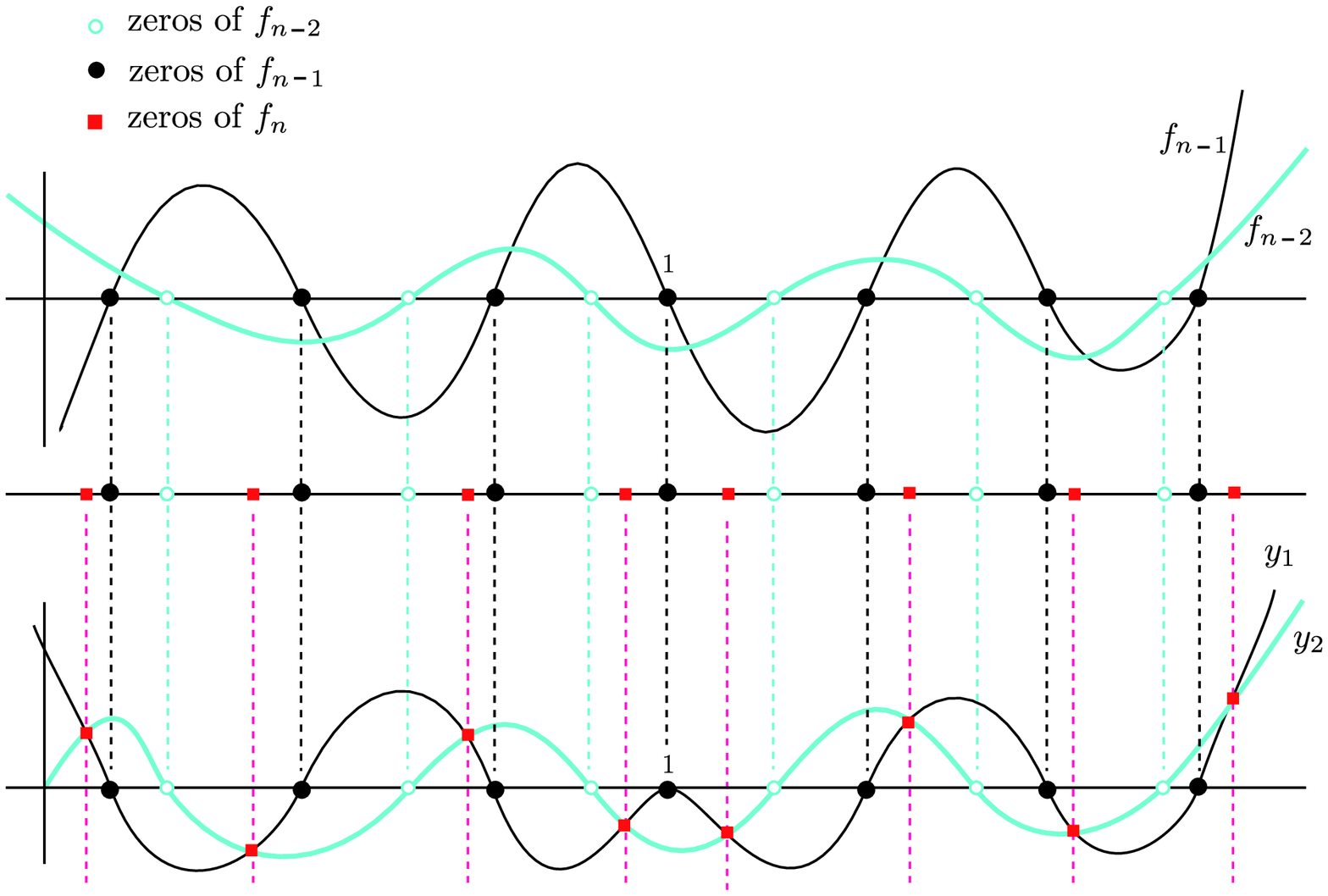}{11}{9.2}}
\def\kyusan{\byouga{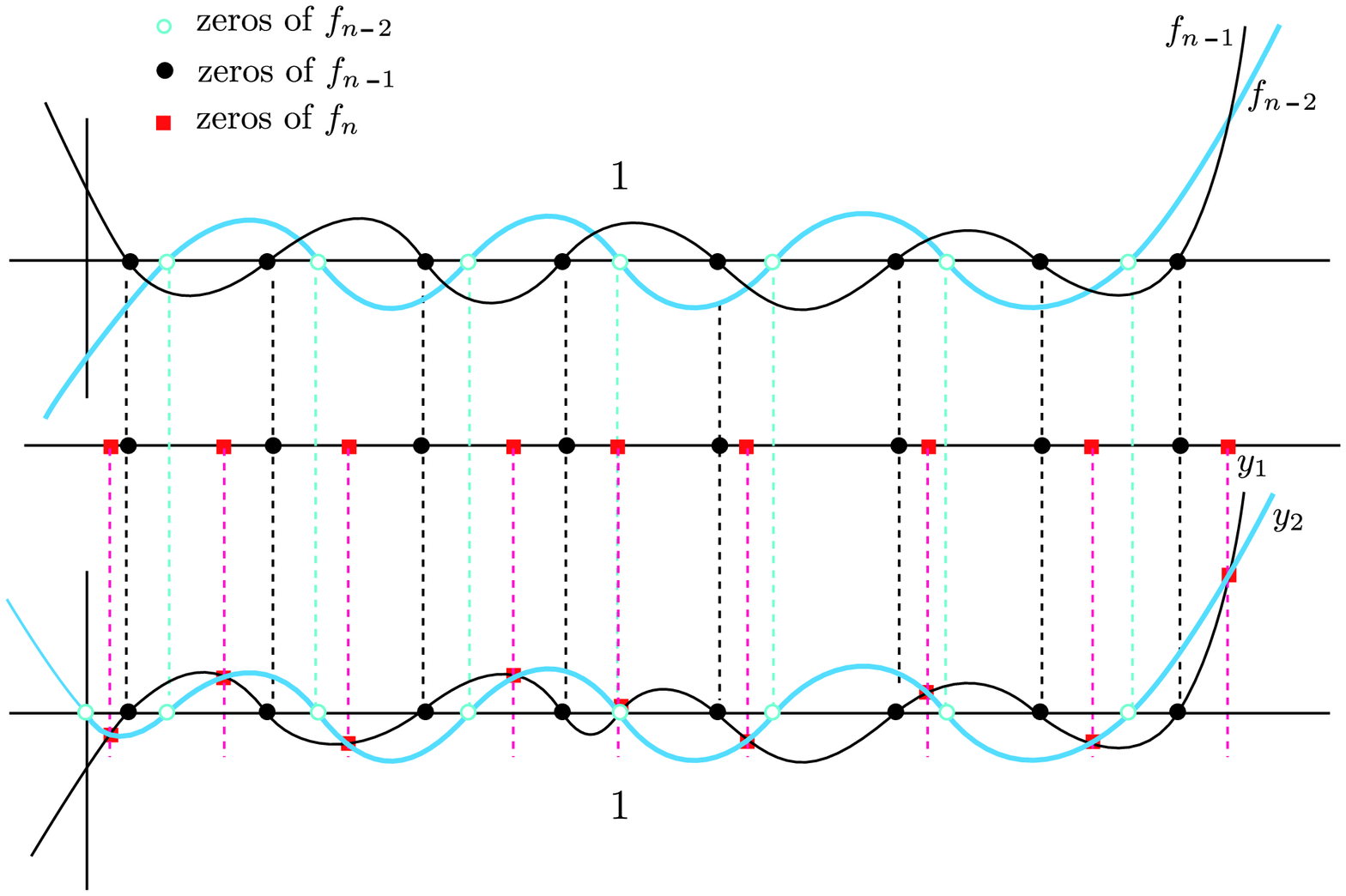}{11}{9.3}}
\def\kyuyon{\byouga{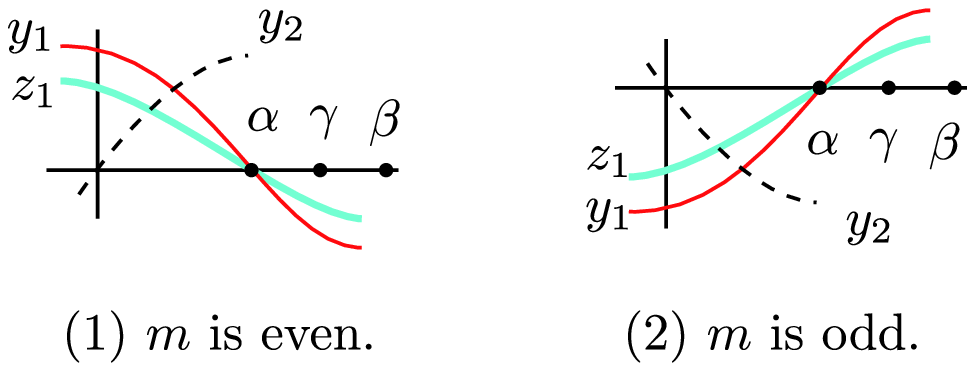}{6.5}{9.4}}
\def\kyugo{\byouga{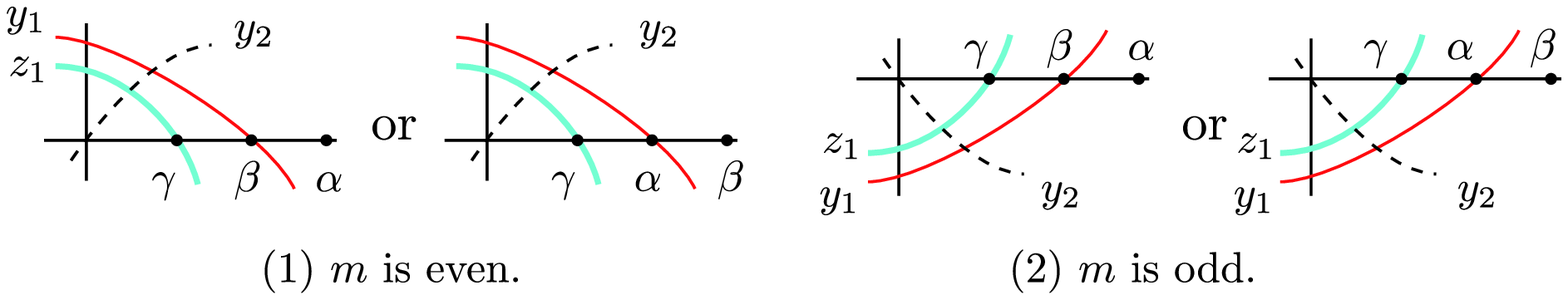}{12}{9.5}}

\def\jui{\byouga{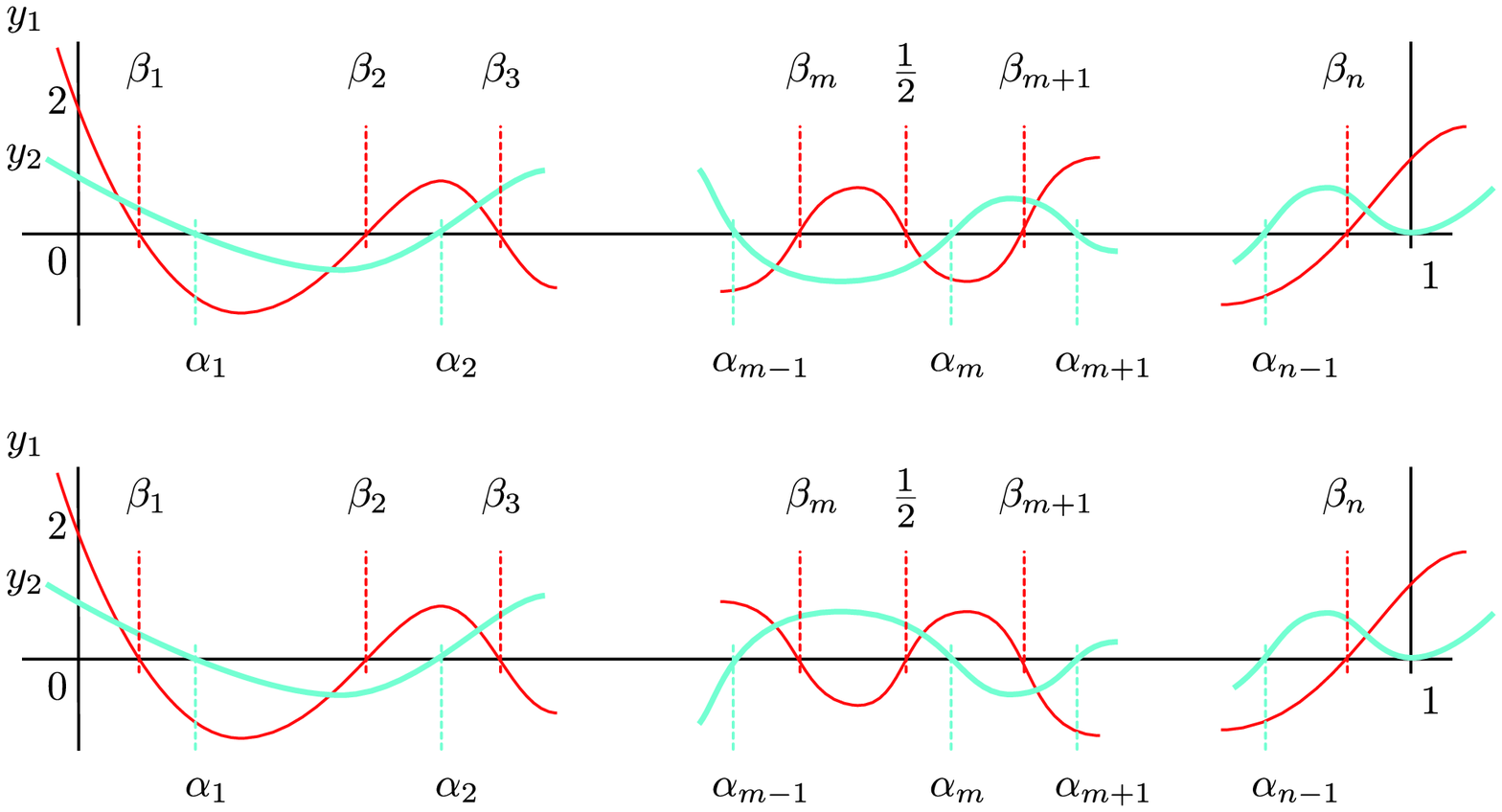}{10}{10.1}}
\def\juii{\byouga{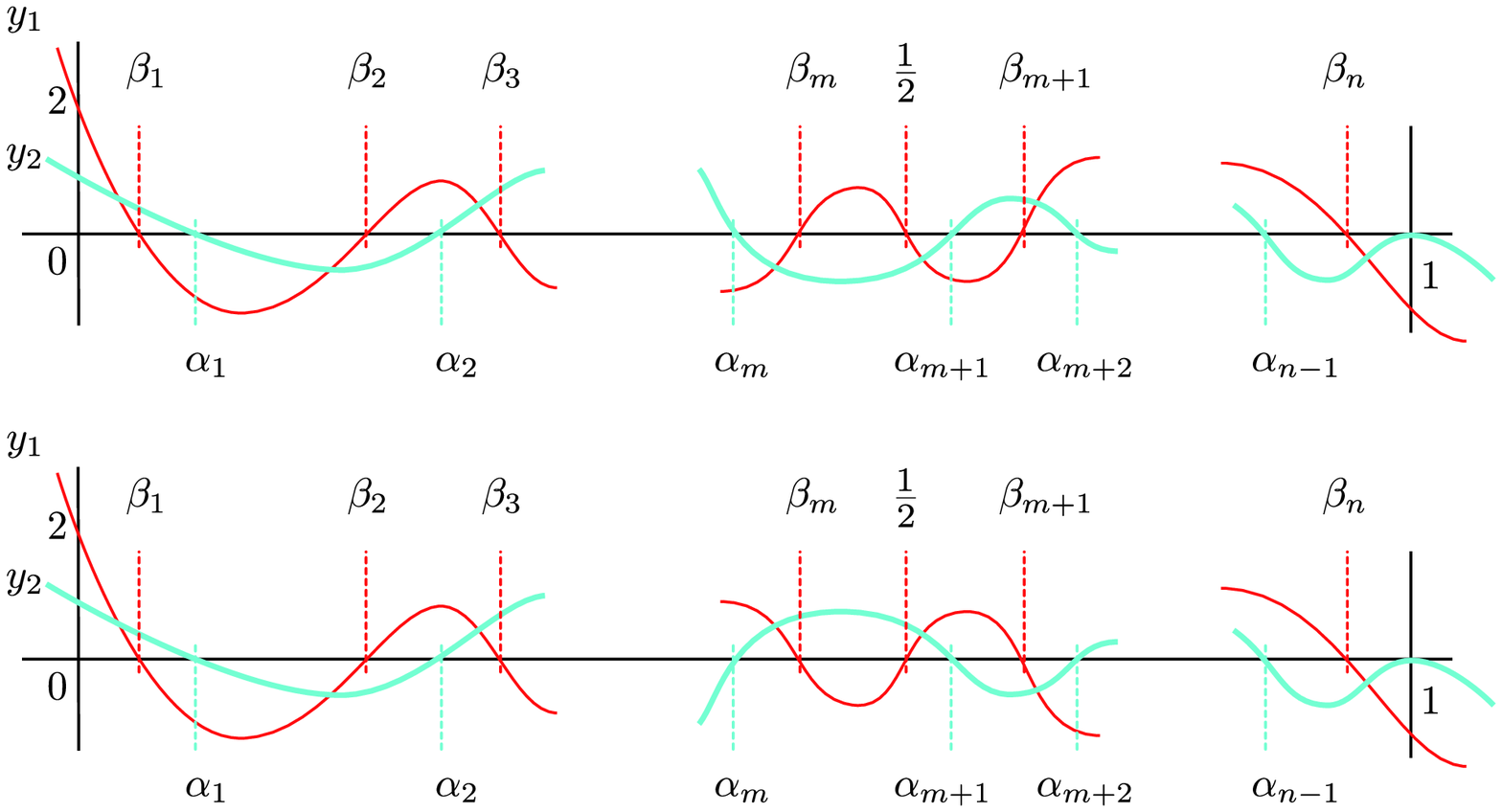}{10}{10.2}}

\def\dodecichi{\byouga{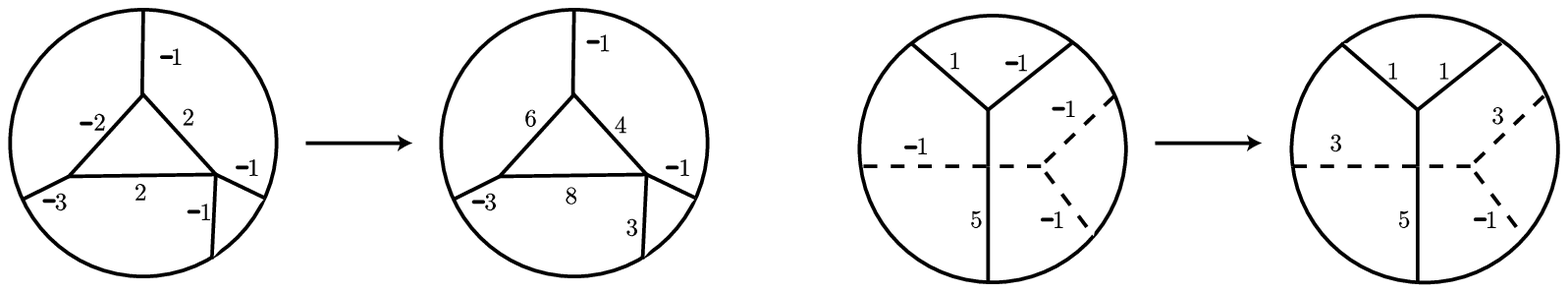}{12}{12.1}}

\def\kichi{\byouga{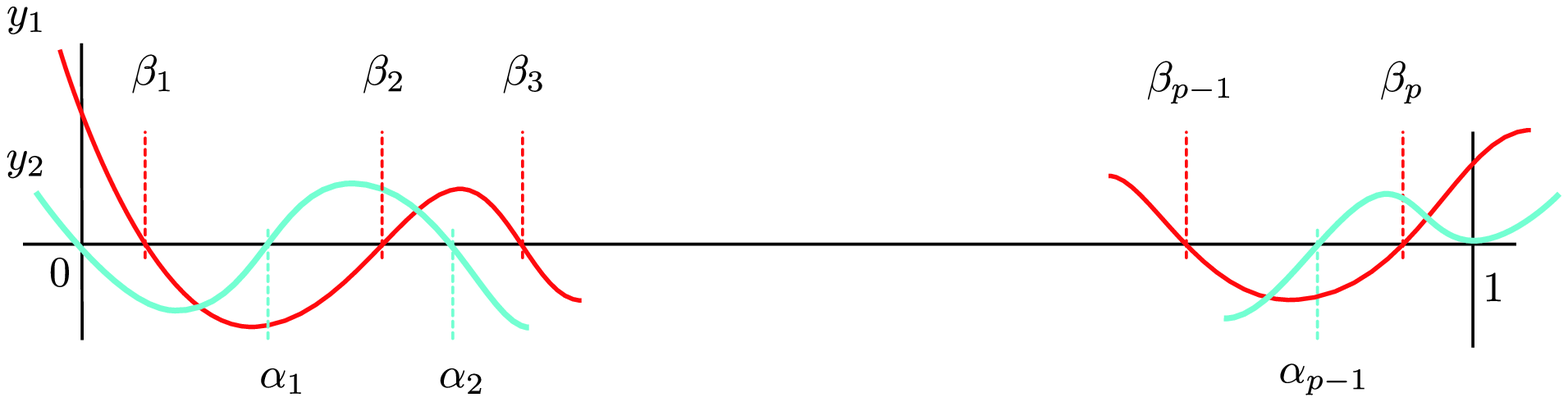}{10}{13.1}}
\def\kni{\byouga{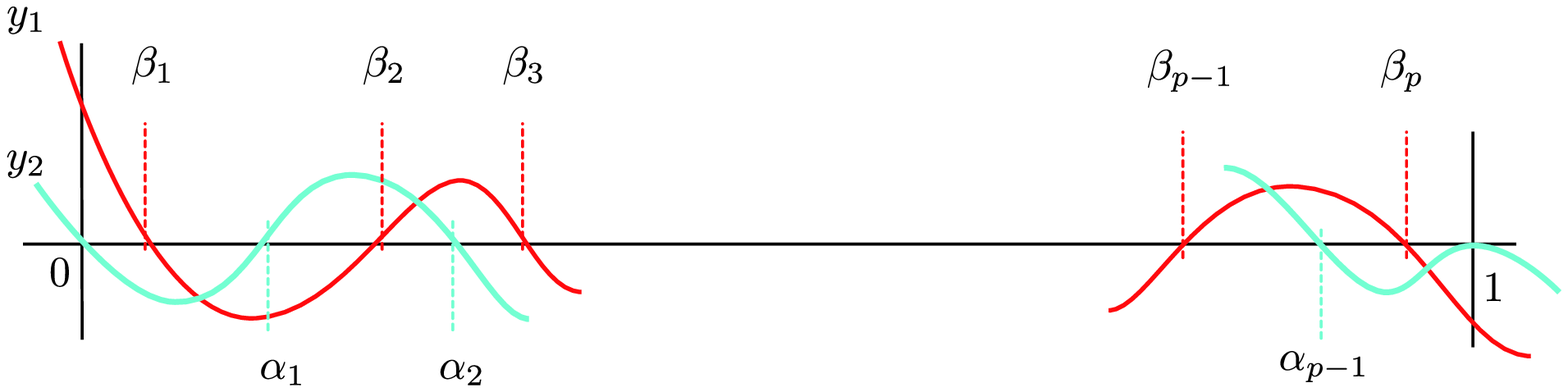}{10}{13.2}}
\def\ksan{\byouga{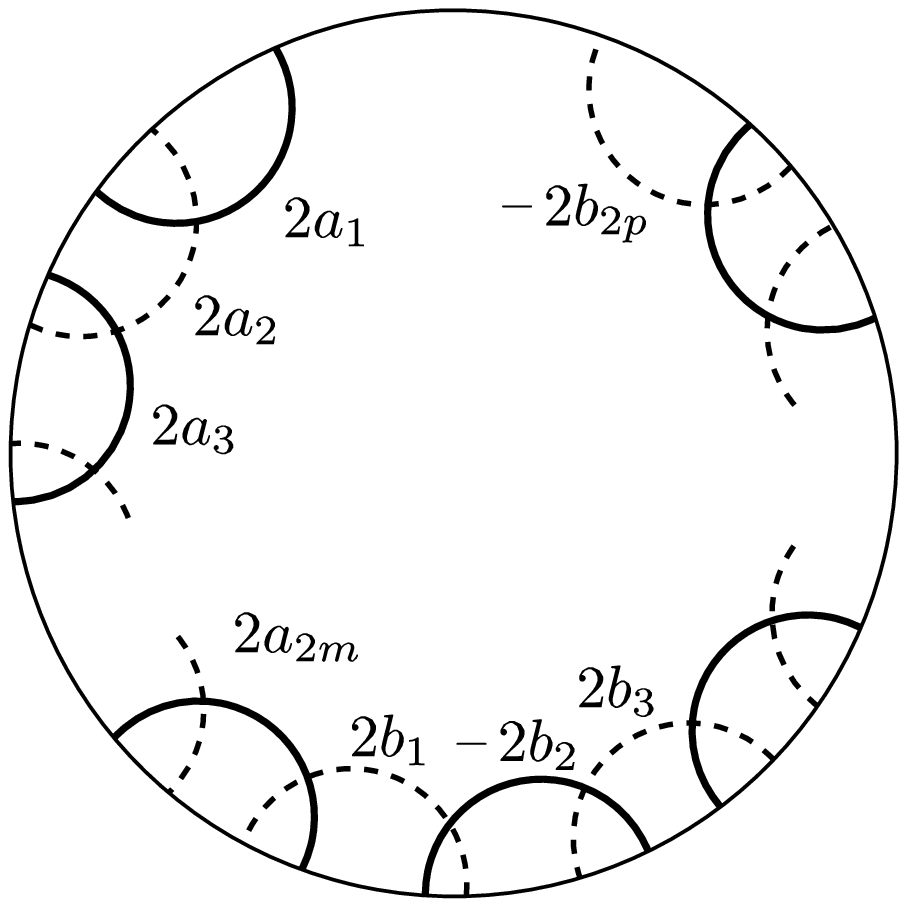}{4}{13.3}}
\def\kyon{\byouga{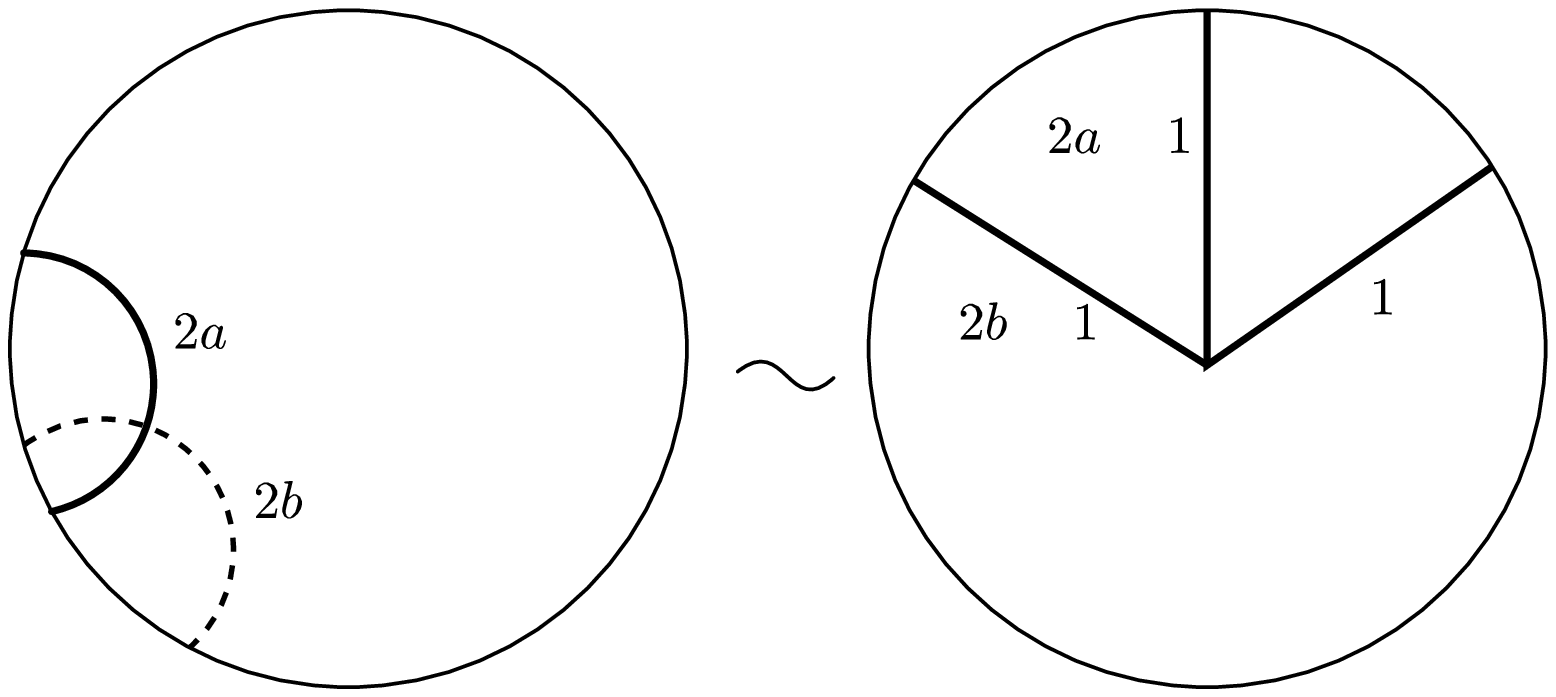}{6}{13.4}}
\def\kgo{\byouga{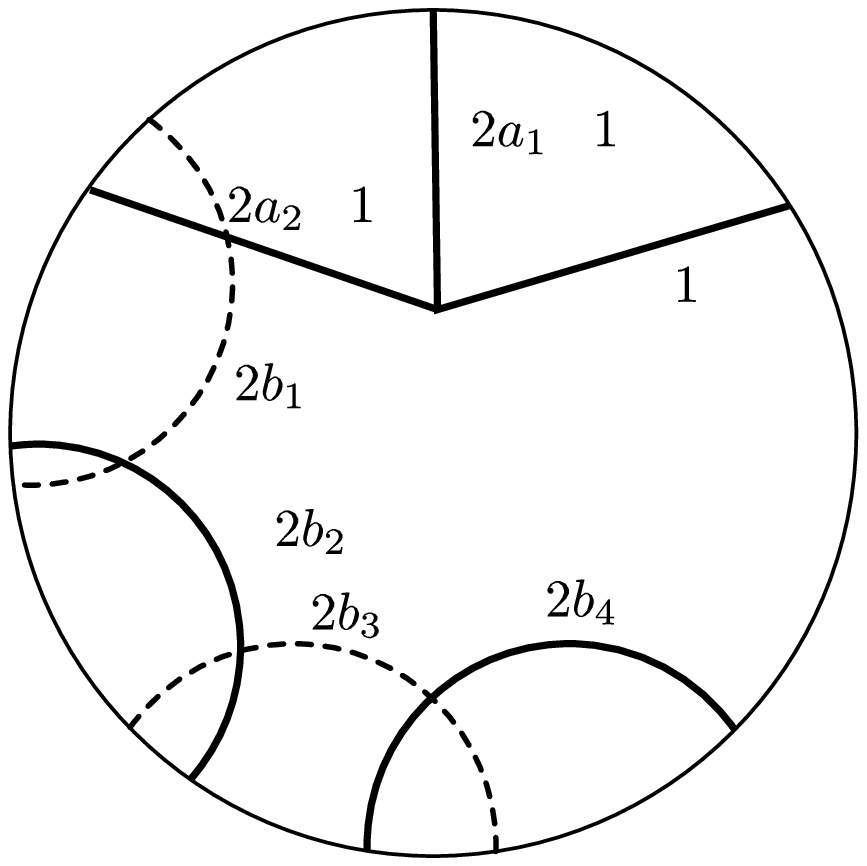}{3}{13.5}}
\def\kroku{\byouga{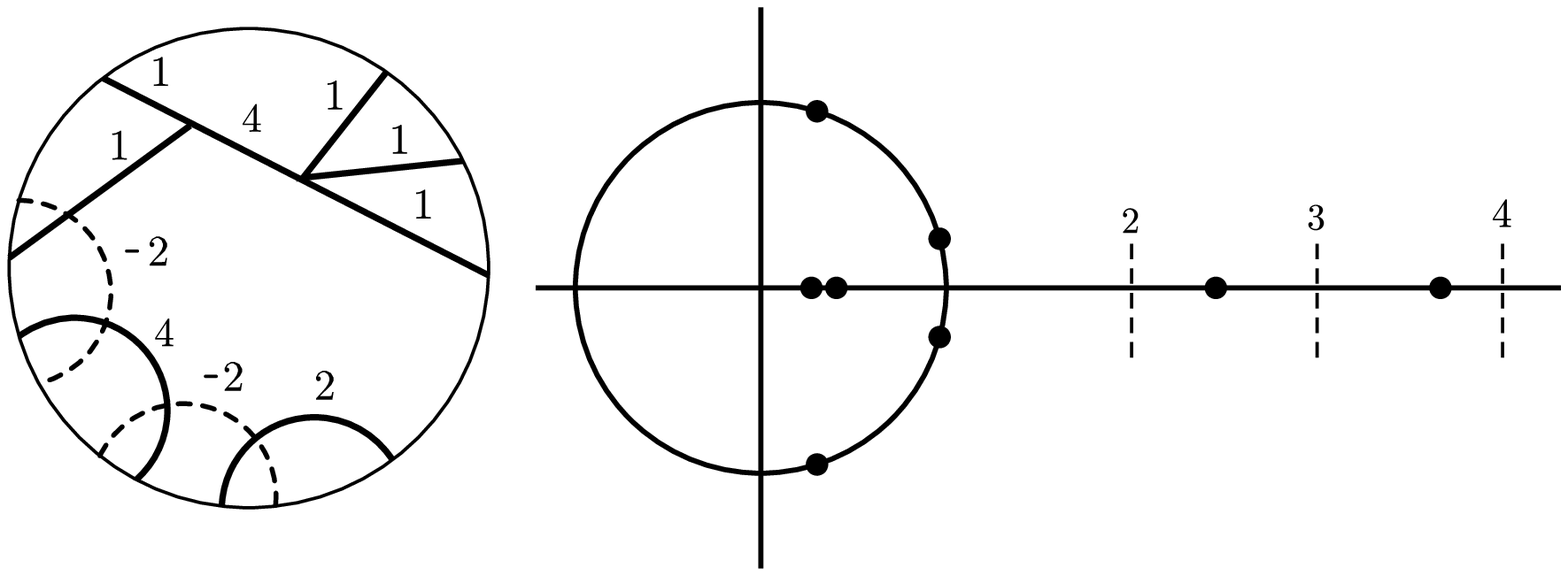}{11}{13.6}}
\def\knana{\byouga{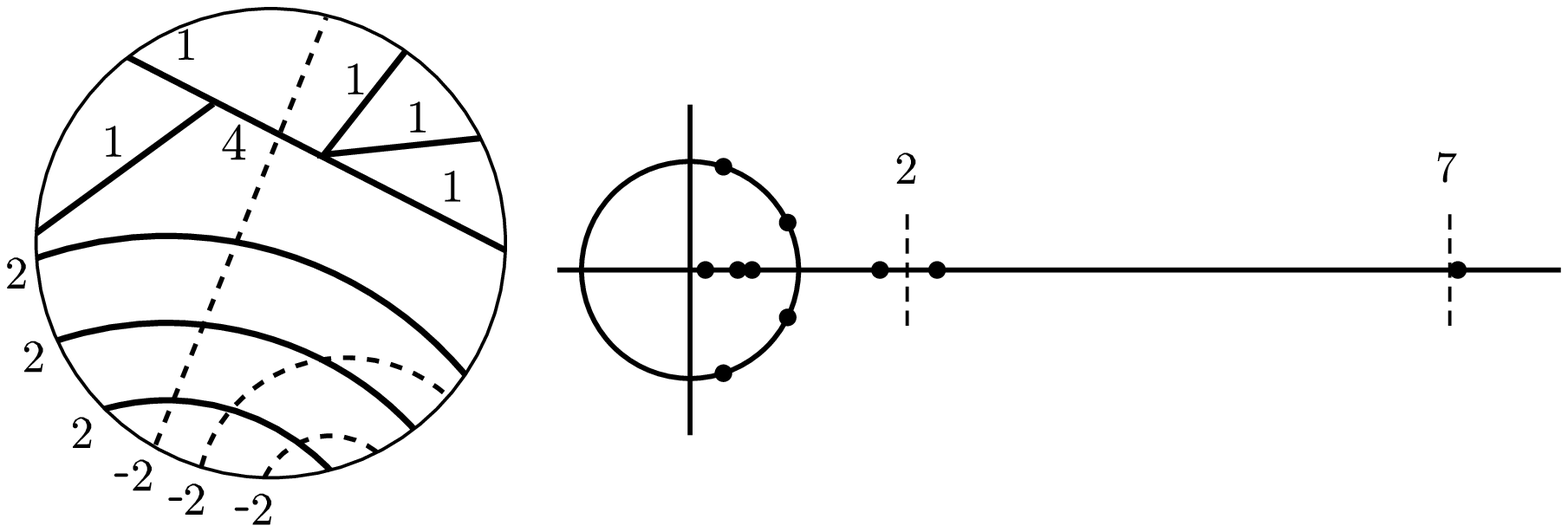}{11}{13.7}}

\def\jgi{\byouga{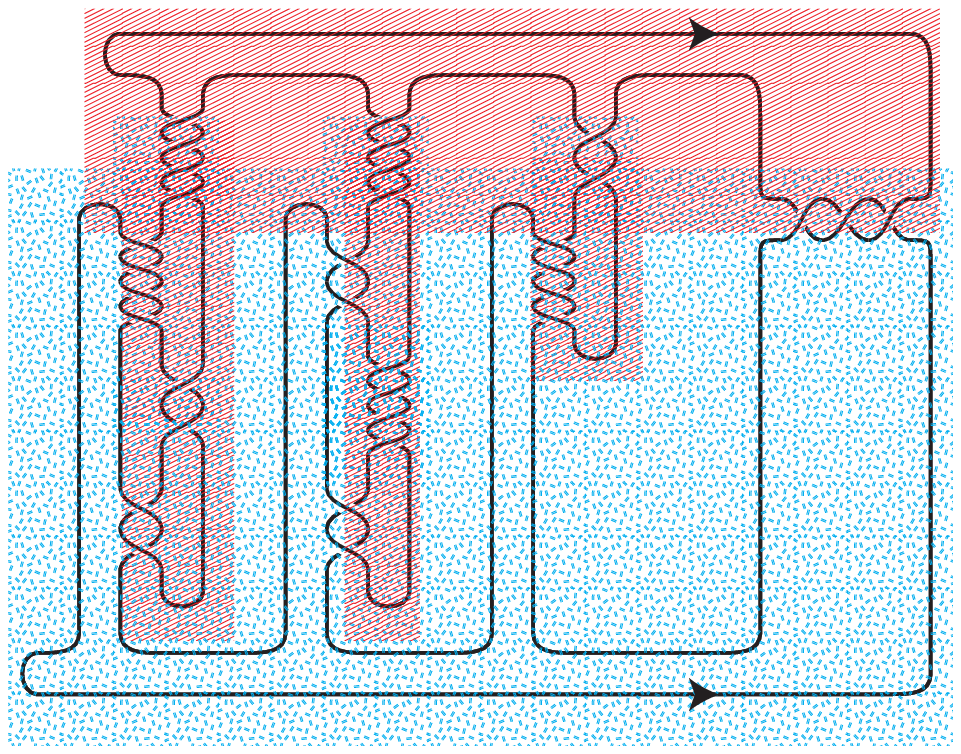}{8}{15.1}}
\def\jgii{\byouga{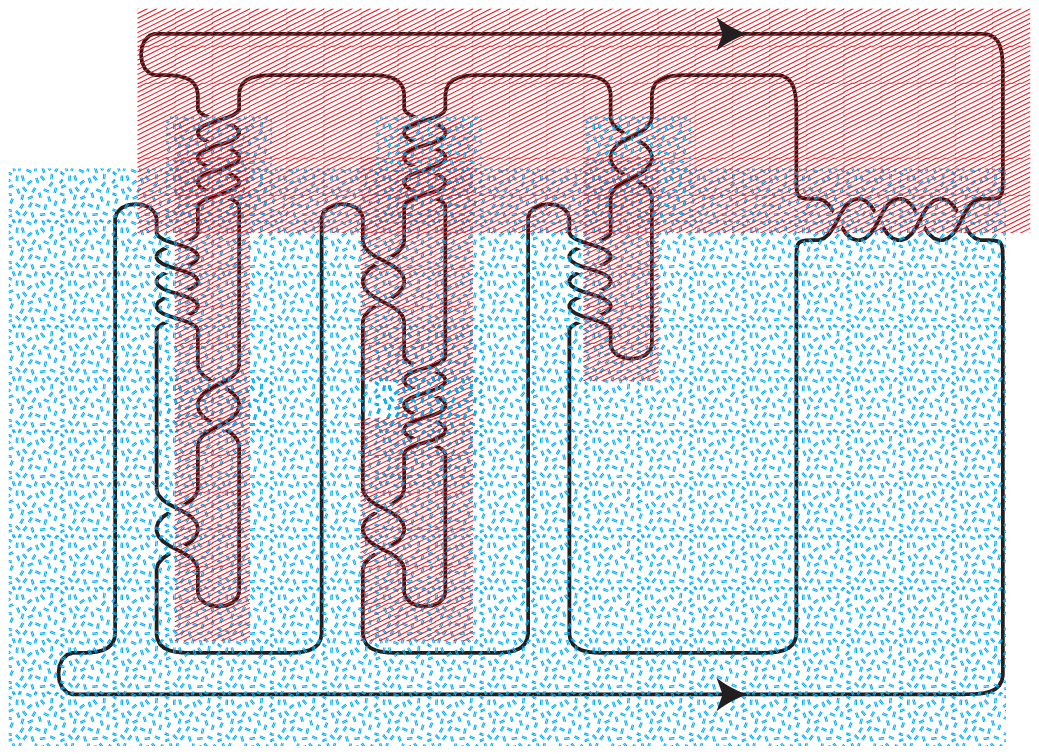}{8}{15.2}}
\def\jgiii{\byouga{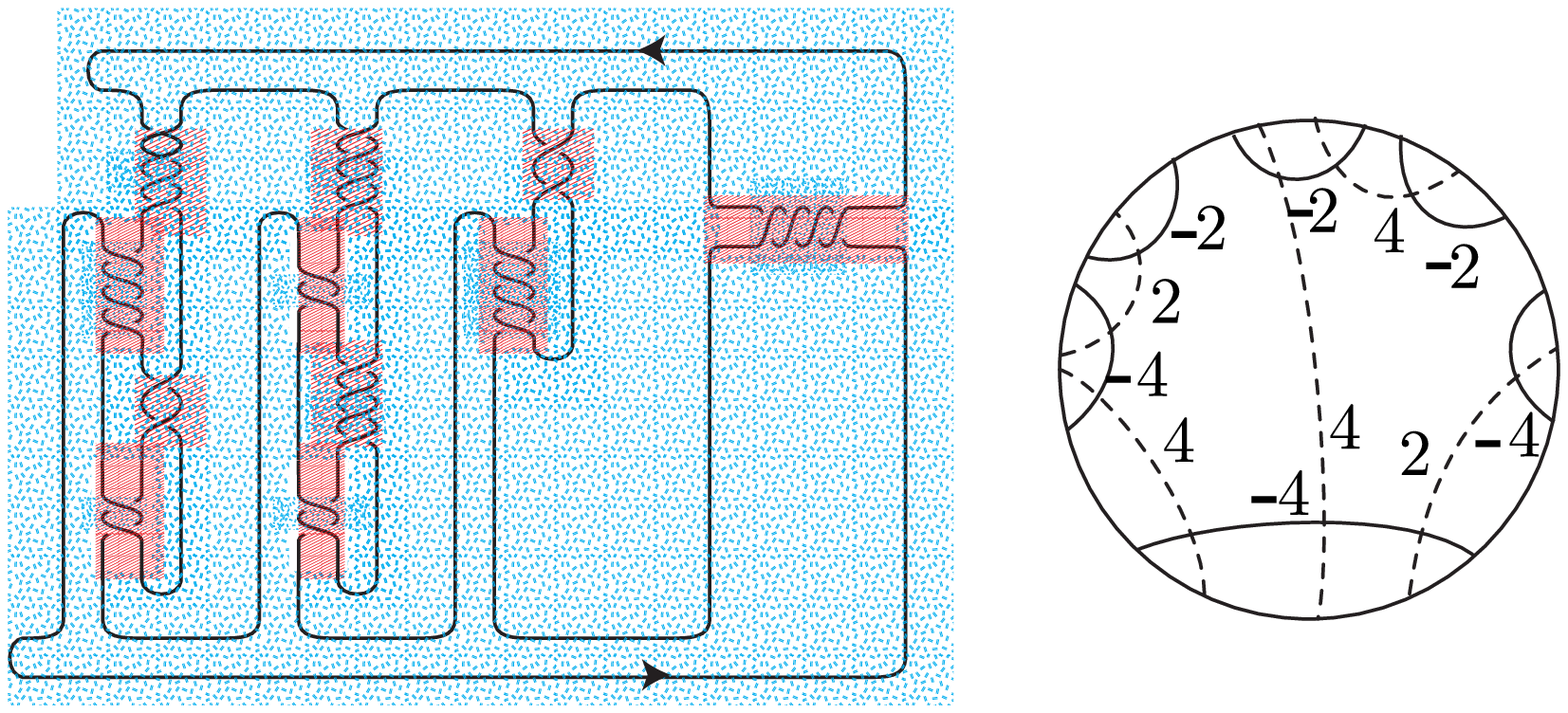}{12}{15.3}}
\def\jgiv{\byouga{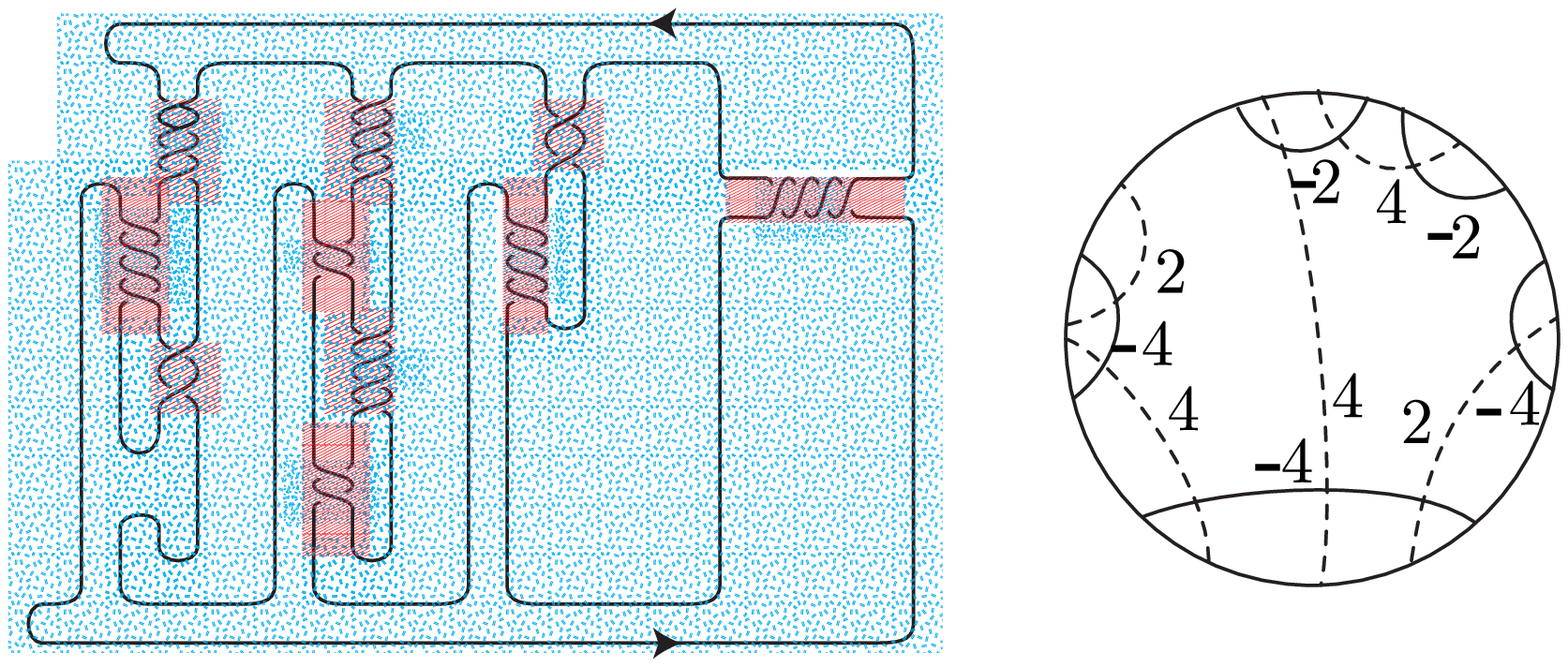}{12}{15.4}}
\def\jgv{\byouga{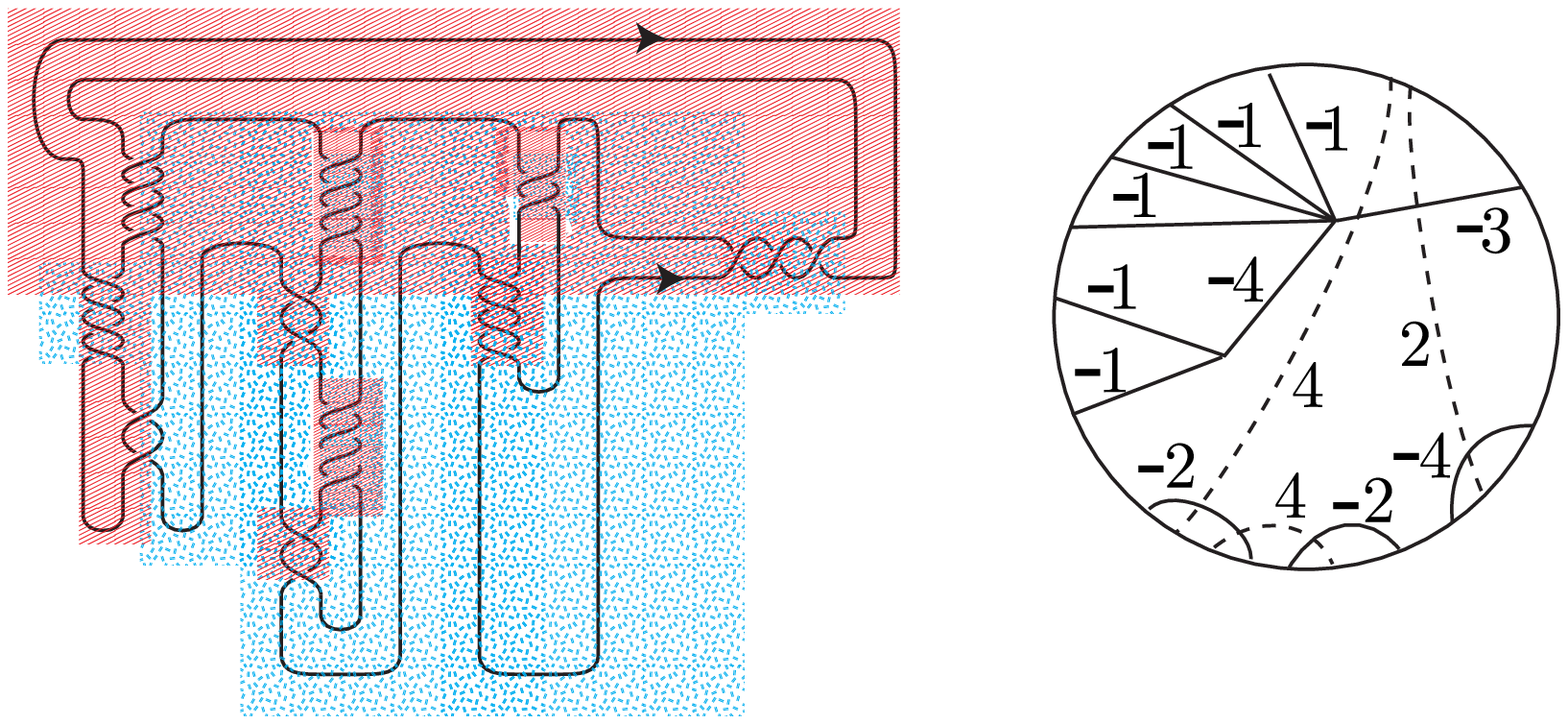}{12}{15.5}}

\def\bichi{\byouga{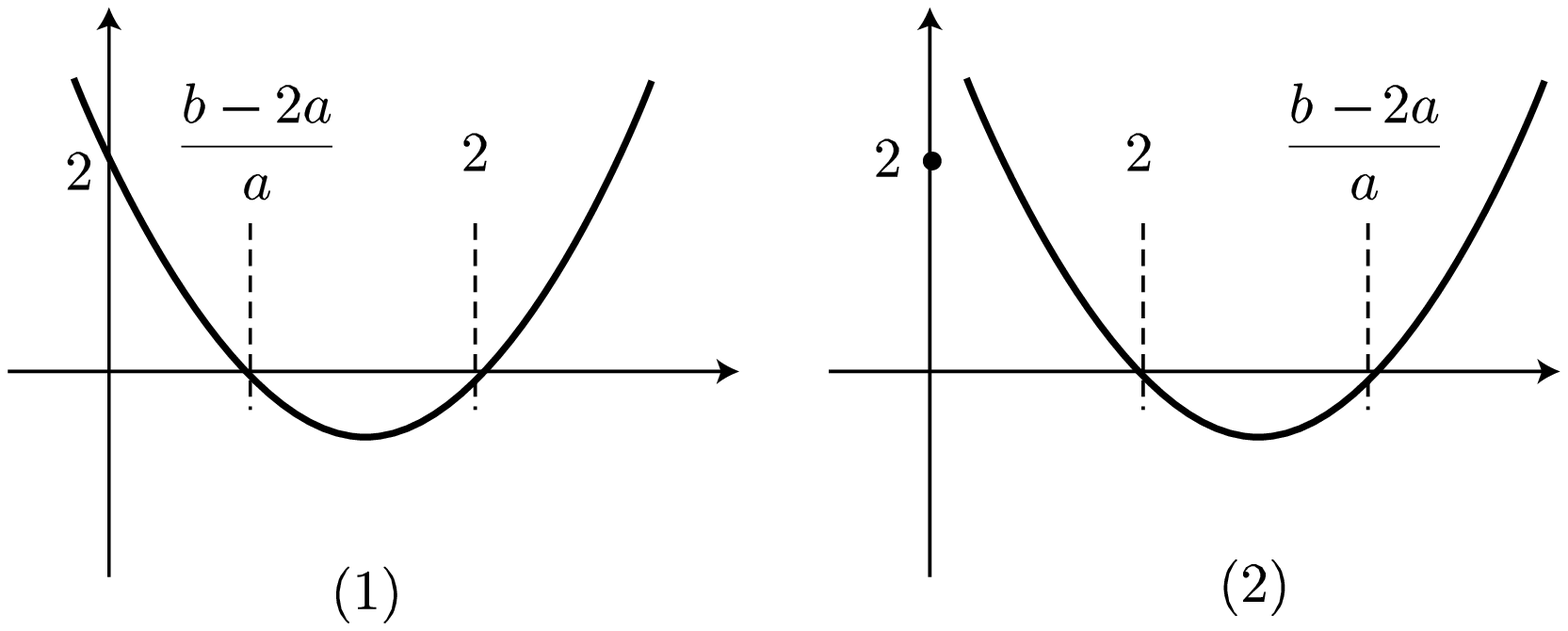}{7}{B.1}}

\def\dichi{\byouga{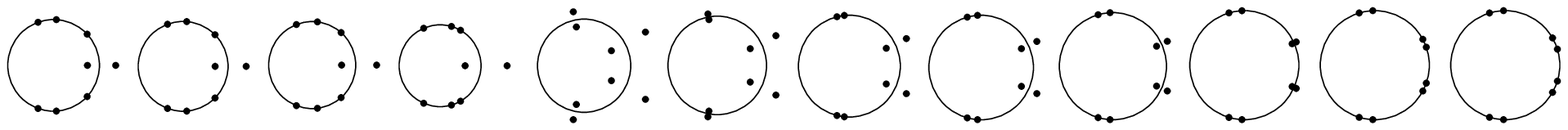}{12}{C.1}}
\def\dni{\byouga{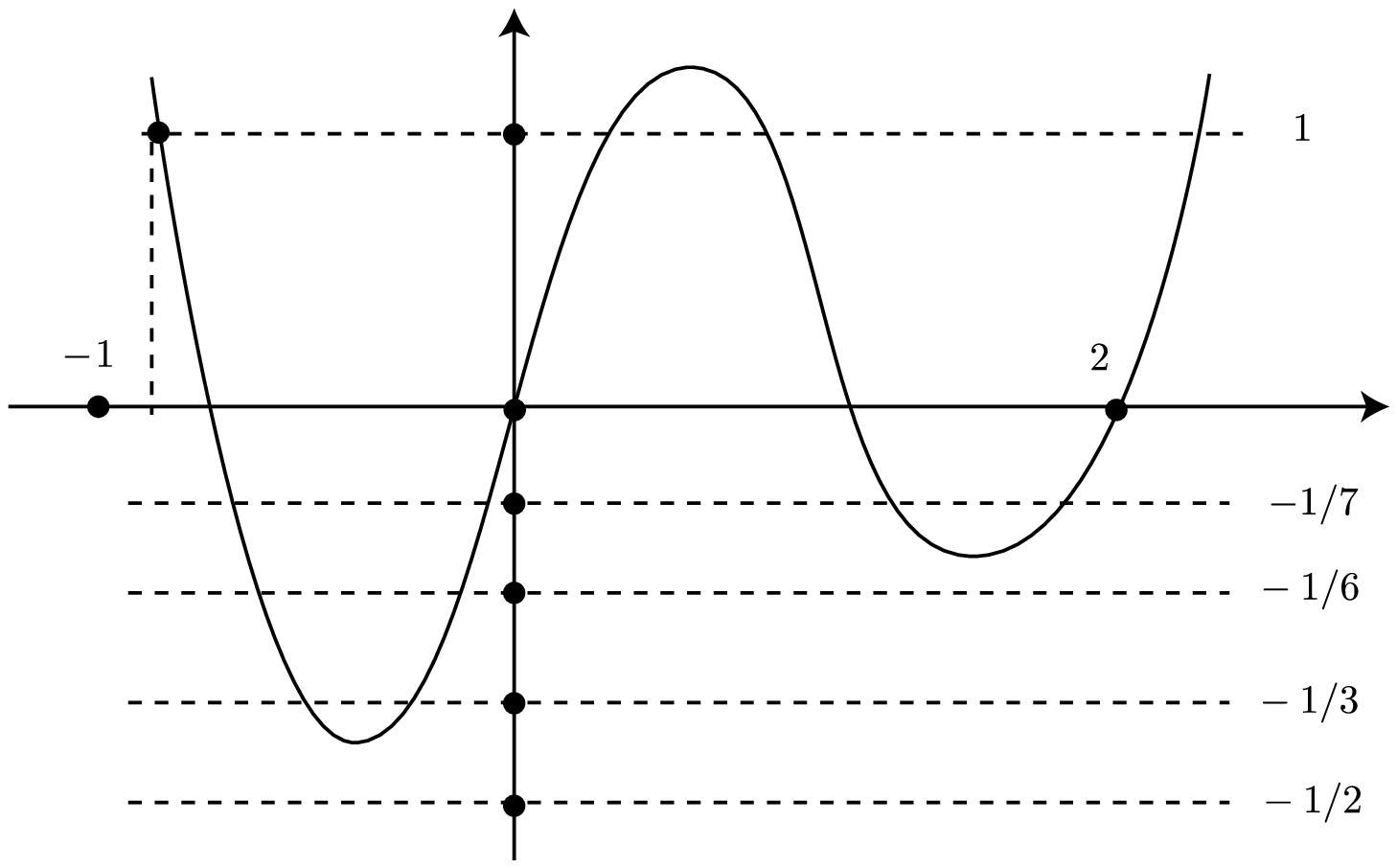}{7}{C.2}}
\def\dsan{\byouga{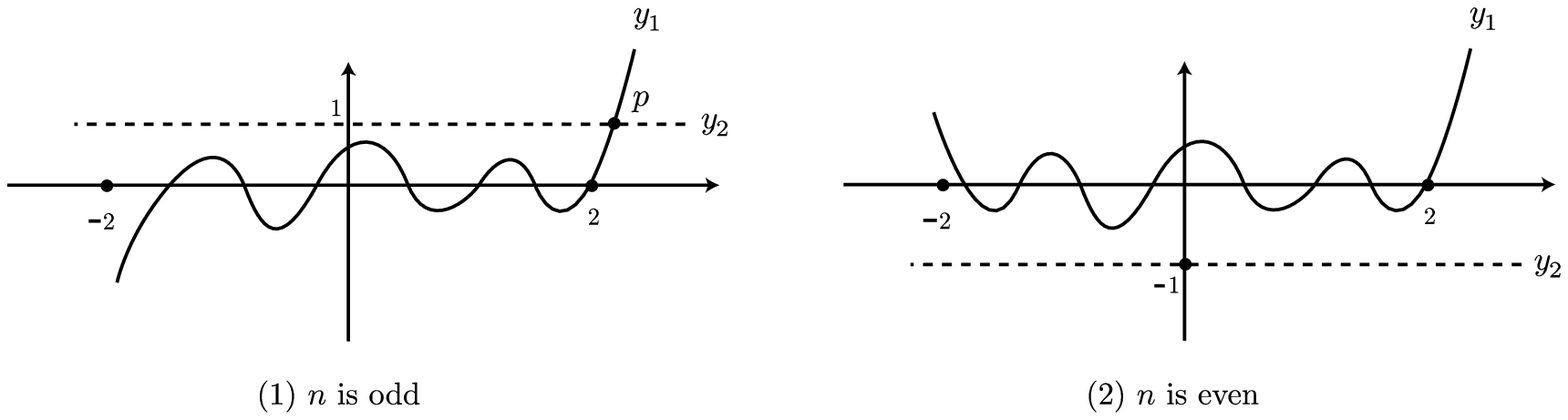}{12}{C.3}}
\def\dyon{\byouga{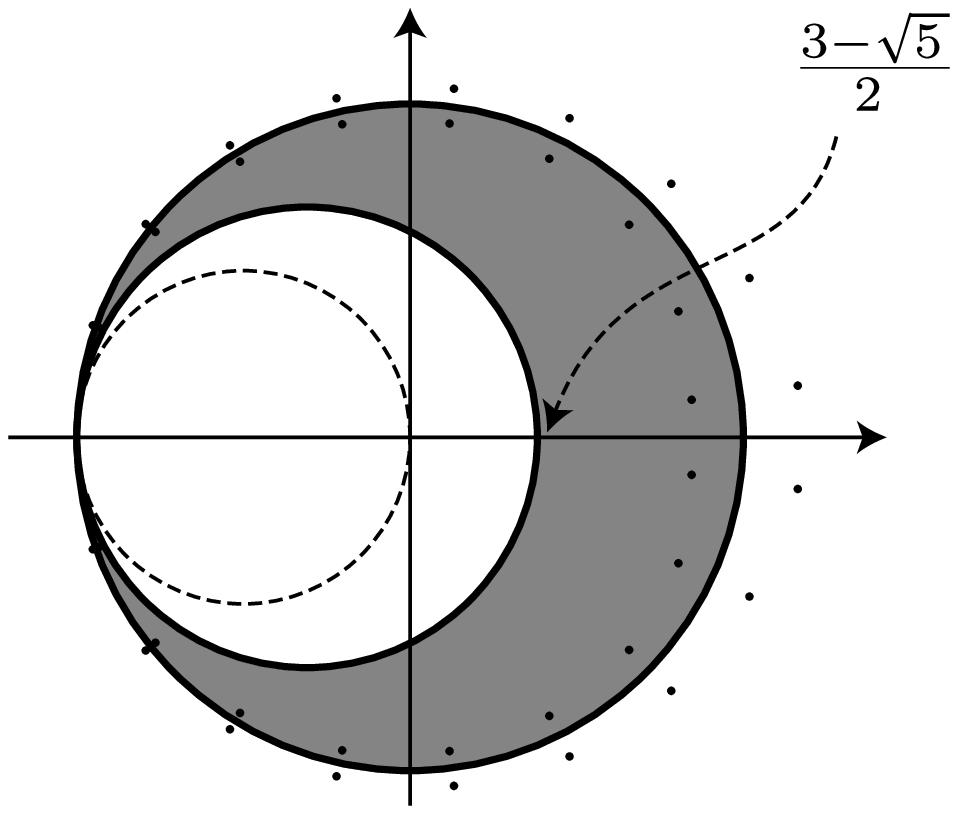}{5}{C.4}}
\def\dgo{\byouga{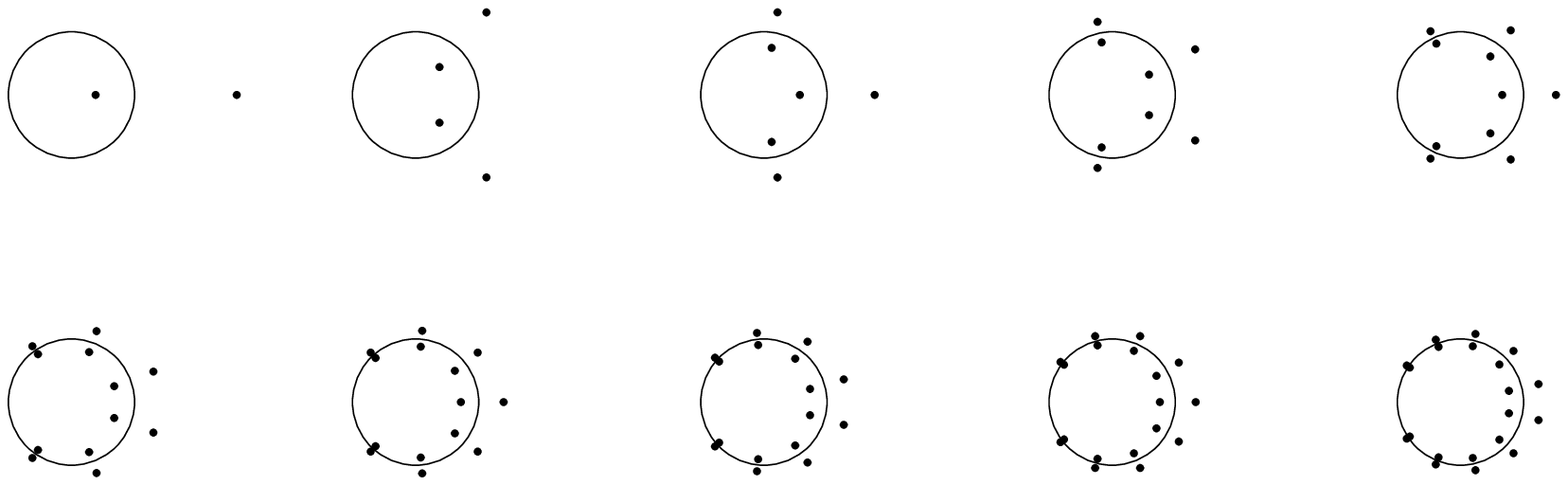}{12}{C.5}}

\def\draftnote{ }
\def\today{}
\arraycolsep=\labelsep

\def\diag{\mathrm diag}

\def\Hsymb#1{\mbox{\strut\rlap{\smash{\LARGE$#1$}}\quad}}

\usepackage{amssymb}\usepackage{amsmath}\usepackage{latexsym}

\def\H{\mathcal H}
\newcommand\ZZ{{\Bbb Z}}
\newcommand\RR{{\Bbb R}}
\newcommand\CC{{\Bbb C}}

\newcommand\qed{\hfill $\Box$ \hfill \\} 
\renewcommand\qed{\hfill $\Box$ \hfill \\} 
\newcommand{\bysame}{\mbox{\rule{3em}{.4pt}}\,}

\newtheorem{thm}{Theorem}[section]

\newtheorem{lemm}[thm]{Lemma}

\newtheorem{prop}[thm]{Proposition}
\newtheorem{cor}[thm]{Corollary}

\newtheorem{yosou}[thm]{Conjecture}

\newtheorem{dfn}[thm]{Definition}
\newtheorem{rem}[thm]{Remark}
\newtheorem{ex}[thm]{Example}

\newtheorem{qu}[thm]{Question}

\newtheorem{probl}[thm]{Problem}

 \newcounter{alphasect}
\def\alphainsection{0}

\let\oldsection=\section
\def\section{%
  \ifnum\alphainsection=1%
    \addtocounter{alphasect}{1}
  \fi%
\oldsection}%

\renewcommand\thesection{%
  \ifnum\alphainsection=1%
    \Alph{alphasect}%
  \else
    \arabic{section}%
  \fi%
}%

\newenvironment{alphasection}{%
  \ifnum\alphainsection=1%
    \errhelp={Let other blocks end at the beginning of the next block.}
    \errmessage{Nested Alpha section not allowed}
  \fi%
  \setcounter{alphasect}{0}
  \def\alphainsection{1}
}{%
  \setcounter{alphasect}{0}
  \def\alphainsection{0}
}%

\renewcommand{\diag}{\mbox{\rm diag}}

\newcommand{\tbt}[4]{
\mbox{$\left[\begin{array}{cc}
#1&#2\\
#3&#4
\end{array}\right]$}
}
\newcommand{\tv}[2]{
\mbox{$\left[\begin{array}{r}
#1\\
#2
\end{array}\right]$}
}

\def\hugesymbol#1{\mbox{\strut\rlap{\smash{\Huge$#1$}}\quad}}
\def\largesymbol#1{\mbox{\strut\rlap{\smash{\Large$#1$}}\quad}}
\def\smallsymbol#1{\mbox{\strut\rlap{\smash{\small$#1$}}\quad}}
\def\hsp{\hspace{-1mm}}

\def\smav{
\hspace{-3mm}\mbox{{\fontsize{5pt}{7pt}\selectfont
$\begin{array}{c}
1\\[-0.5mm]0\\[-1.2mm]\cdot\\[-2mm]\cdot\\[-2mm]\cdot\\[-1mm]0
\end{array}
$}}\hspace{-3mm}
}

\def\smavt{
\hspace{-2mm}\mbox{{\fontsize{7pt}{9pt}\selectfont
$\begin{array}{c}
t\\[-0.5mm]0\\[-1.2mm]\cdot\\[-2mm]\cdot\\[-2mm]\cdot\\[-1mm]0
\end{array}
$}}\hspace{-3mm}
}

\newcommand{\nodiv}{\mbox{$\not|\,$}}

\title{Various stabilities of the Alexander
polynomials of knots and links}

\author{Mikami Hirasawa}
\address{Department of Mathematics,
Nagoya Institute of Technology,\\
Nagoya Aichi 466-8555 Japan\\
{\it E-mail: hirasawa.mikami@nitech.ac.jp}
}

\author{Kunio Murasugi}
\address{Department of Mathematics,
University of Toronto,\\
Toronto, ON M5S2E4 Canada\\
{\it E-mail: murasugi@math.toronto.edu}
}

\begin{document}
\maketitle
 \pagestyle{myheadings}
 \markboth{\hiduke}
 {\hiduke}
 
\begin{abstract} 
In this paper, we study distribution of the zeros of
the Alexander polynomials of knots and links in $S^3$.
We call a knot or link {\it real stable} (resp. {\it circular stable})
if all the zeros of its Alexander polynomial are real 
(resp. unit complex).
We give a general construction of real stable and circular stable
knots and links. 
We also study pairs of real stable knots and links such that
the zeros of the Alexander polynomials are interlaced.
\end{abstract}



\newcommand{\sect}[3]{
#1:\ #2 \dotfill \pageref{#1}\linebreak}
\newcommand{\subsect}[3]{
\hspace*{5mm} #1:\ #2 \dotfill \pageref{#1}\linebreak}

\newpage

{\centerline {\bf Table of contents}}

\noindent
0: Introduction\dotfill\pageref{0}\\
\sect{1}{Stability Property}{2}
\subsect{1.1}{Half-plane property}{2}
\subsect{1.2}{$D$-stable polynomial}{3}
\sect{2}{Hurwitz-stability}{6}
\subsect{2.1}{Hurwitz-Routh Criterion}{7}
\subsect{2.2}{Lyapunov matrix}{7}
\sect{3}{Stable polynomial}{8}
\subsect{3.1}{Multivariate stable polynomials}{8}
\subsect{3.2}{Real stable univariate polynomials}{9}
\sect{4}{Preliminaries }{10}
\sect{5}{The Alexander polynomials of alternating knots}{18}
\subsect{5.1}{Hoste's Conjecture}{18}
\subsect{5.2}{Trapezoidal Conjecture}{19}
\sect{6}{Construction of real stable knots (I)}{20}
\subsect{6.1}{Quasi-rational knots}{21}
\subsect{6.2}{Stable quasi-rational knots}{21}
\subsect{6.3}{Examples}{22}
\sect{7}{Construction of real stable knots (II)}{23}
\subsect{7.1}{Positive or negative disks}{23}
\subsect{7.2}{Stable alternating knots and links}{24}
\subsect{7.3}{Pseudo-positive or pseudo-negative disk}{25}
\subsect{7.4}{Example}{25}
\sect{8}{Exceptional stable knots and links}{25}
\subsect{8.1}{Exceptional stable knots}{26}
\subsect{8.2}{Exceptional stable links}{26}
\sect{9}{Interlacing property (I) 2-bridge knots}{28}
\sect{10}{Interlacing property (II) Quasi-rational knots $X_n$}{33}
\sect{11}{Interlacing property (III) Quasi-rational knots $Y_{2n+1}$}{36}
\subsect{11.1}{Conway polynomials of $Y_n$}{37}
\subsect{11.2}{Alexander polynomials of $Y_n$}{40}
\sect{12}{$c$-stable knots and links }{42}
\subsect{12.1}{Regular and exceptional $c$-stable $2$-bridge knots and links}{42}
\subsect{12.2}{Construction of $c$-stable quasi-rational knots 
and links }{43}
\subsect{12.3}{General construction of $c$-stable knots and links }{44}
\subsect{12.4}{Interlacing property of zeros on the unit circle}{45}
\sect{13}{Bi-stable knots and links}{46}
\subsect{13.1}{Bi-stable 2-bridge knots and links}{46}
\subsect{13.2}{Exceptional bi-stable knots and Salem knots}{47}
\subsect{13.3}{General bi-stable knots and links}{48}
\sect{14}{Mobius Transformations}{49}
\sect{15}{Montesinos knots}{56}
\sect{16}{Multivariate stable link polynomials}{61}
\sect{17}{Inversive links}{65}
\subsect{17.1}{Standard inversive links}{65}
\subsect{17.2}{Exceptional inversive links}{67}
Appendix A: Representation polynomials\dotfill\pageref{A}\\
\subsect{A.1}{Parabolic representation}{68}
\subsect{A.2}{Dihedral representation}{70}
Appendix B: Determination of $\delta_4$\dotfill\pageref{B}\\
Appendix C: Distribution of the zeros\dotfill\pageref{C}\\
References\dotfill\pageref{ref}

\newpage

\setcounter{section}{-1}
\section{Introduction}\label{0}
Let $\H \subset \CC$ be an open right half-plane, 
i.e.,
$\{\alpha \in \CC | {\rm Re}(\alpha) > 0\}$, 
or an open upper half-plane,
i.e., $\{\alpha \in \CC | {\rm Im}(\alpha) > 0\}$. 
Let $f(z_1, \cdots, z_n) \in
\CC[z_1, \cdots, z_n]$ be a polynomial in $n$ variables,
$z_1, \cdots, z_n$. 
We say that $f(z_1, \cdots, z_n)$ is {\it $\H$-stable} if
for any values $\alpha_j \in \H, 1 \leq j \leq n,
f(\alpha_1, \cdots, \alpha_n) \neq 0$.
If $\H$ is an open right half-plane, 
then $f$ is called {\it Hurwitz stable}. 
If $\H$ is an open upper half-plane,
then $f$ is called a {\it stable polynomial}, and further,
if $f$ is a real polynomial, $f$ is sometimes called {\it real stable}.
The theory of stable polynomials has a long history, 
but the recent development of this theory is very impressive
and is summarized in a remarkable survey article \cite{wag}.

The purpose of this paper is to provide 
a recent study on various stabilities of 
the Alexander polynomials of knots or links 
in $S^3$.
The study was motivated by our desire to answer 
a question (later called conjecture) posed 
by Jim Hoste in 2002.
He asks if the real part of each zero of the Alexander polynomial
$\Delta_K(t)$ of an alternating knot $K$ is larger than $-1$.
It is exactly a question whether 
$\Delta_K (-(t+1))$ is 
(strongly) Hurwitz stable for an alternating knot $K$.
The question leads us to other problems on 
stabilities of the Alexander polynomial of a 
(not necessarily alternating) knot.
For example, since the sequence of the coefficients of 
a stable univariate real polynomial under a certain condition
is unimodal, 
we see immediately that the stable Alexander polynomial of an
alternating knot satisfies Trapezoidal Conjecture,
one of the outstanding conjectures that still remains open.
In \cite{LM}, 
it is shown that many 2-bridge knots or links satisfy 
Hostefs Conjecture.
Further, a few more subtle theorems on Hurwitz stability and real
stability of the Alexander polynomials
of 2-bridge knots or links are proven.

In this paper, knots or links are not necessarily alternating, 
and we discuss stabilities of the Alexander polynomials
of knots or links, and further, 
we discuss the third stability, called
{\it circular stability}, 
of the Alexander polynomials of knots or links.

This paper is organized as follows. 
It consists of two parts. 
The first part, consisting of Sections 1-3 is a quick review of various types
of stable polynomials. 
Almost all materials in this part are taken from
various known sources and hence proofs are entirely omitted.
The rest of the paper forms the second part of the paper.
In Section 4, first we introduce new notations and various terminologies
and then we prove a couple of propositions
on matrices.
We use these propositions as basic tools to prove many theorems in this
paper. For convenience, we say that a knot or link is {\it $\H$-stable}
if its Alexander polynomial is $\H$-stable.
In Section 5, we review some connections between stable Alexander
polynomials and various conjectures in knot theory.
In Sections 6 and 7, 
we study real stable Alexander polynomial of a knot or link.
The proto-type is a 2-bridge knot or link.
As is shown in \cite{LM}, 
a 2-bridge knot with an alternating continued
fraction expansion, i.e., $[2a_1, 2a_2, \cdots, 2a_m],
a_j a_{j+1} < 0$, is always stable.
By generalizing these knots, we construct more general
stable knots in these sections.
In Section 8, we discuss some exceptional stable 2-bridge knots or links
which have non-alternating continued fraction expansions.
We discuss the stability of the (2 variable) Alexander polynomials 
of such links in Section 16.
One of the important properties of stable polynomials is the
\lq\lq interlacing property\rq\rq\ of the zeros. 
In Sections 9 through 11,
we discuss this property, first, for 2-bridge knots or links (Section 9)
and then for a generalization of 2-bridge knots,
which we name {\it quasi-rational} knots  (Sections 10 and 11).
In Section 12, we study circular stable polynomials.
A real polynomial $f(z)$ is called {\it circular stable} (or $c$-stable)
if all the zeros of $f(z)$ lie on the unit circle, i.e.,
$|\alpha|$ = 1 for
any zero $\alpha$ of $f(z)$.
The Alexander polynomial of a special alternating knot is always $c$-stable.
But, many non-special alternating knots also have the 
$c$-stable Alexander polynomial.
In this section, we give a systematic way to construct 
$c$-stable knots or links. 
These knots or links are in general not alternating.
Now for convenience, we call a real polynomial $f(z)$ 
{\it bi-stable} if the zero
of $f(z)$ is either real or unit complex number.
In Section 13, we prove that 2-bridge knots of 
a certain special type are bi-stable. 
A proof of the bi-stability of a polynomial is rather
complicated.
To show the real stability, a Seifert matrix plays a key role, but to
show the bi-stability, the interlacing property of the zeros is crucial.
Bi-stable knots or links appear implicitly in \cite{hiro} and others.
Salem fibred knots are a special type of bi-stable knots.
We briefly discuss them in this section.
In Section 14, we study a Mobius transformation 
$\varphi: \CC \cup \{\infty
\} \longrightarrow \CC \cup \{\infty \}$. 
There is one special Mobius transformation 
$\varphi$ that relates the Alexander
polynomial $\Delta_K(t)$ of a knot $K$ to 
Hosokawa polynomial $\nabla_{L(K)}(t)$
of a link $L(K)$ with an 
arbitrary number of components in such a way that 
$\varphi$ maps all the zeros of
$\Delta_K(t)$ to all the zeros of $\nabla_{L(K)}(t)$. 
In particular,\\
(1) if $\Delta_K(t)$ is $c$-stable, then $\nabla_{L(K)}(t)$ is real stable,
and\\
(2) if $\Delta_K(t)$ is real stable, then $\nabla_{L(K)}(t)$ is $c$-stable
and further,\\
(3) if $\Delta_K(t)$ is bi-stable,
$\nabla_{L(K)}(t)$ is bi-stable.

Then we show that given $\Delta_K(t)$, we can express $\nabla_{L(K)}(t)$ in
terms of the coefficients of $\Delta_K(t)$.\\
The (reduced) Alexander polynomial of a link depends on orientation of
each component in a very delicate manner. 
In fact, there exists a 2-component link $L$ such that one orientation gives 
a stable link, but reversing the orientation of one component results in
a $c$-stable link. 
We call such a link {\it inversive}.
We have many inversive Montesinos links.
Therefore, in Section 15, we study the Alexander polynomials of
alternating Montesinos knots or links. 
We specify some class of
alternating Montesinos knots or links and prove that a knot 
$K$ or link in this class has the following property:\\
(a) If $K$ is a knot, then $K$ is $c$-stable, stable or bi-stable.\\
(b) If $K$ is a link, then $K$ is inversive.

In Section 16, we consider a 2-component link $L$.
Let $\Delta_L (x, y)$ be the Alexander polynomial of $L$.
The stability problem of $\Delta_L (x, y)$ is not an easy problem, unless
$\Delta_L (x, y)$ is multi-affine, i.e., each variable has
degree at most one in each term. 
If $\Delta_L (x, y)$ is stable, then so is
$\Delta_L (t, t)$, but in general $t^n \Delta_L (t, t^{-1})$ is not
stable, where $n =\deg_y\Delta_L (x, y)$.
Note that $t^n \Delta_L (t, t^{-1})$ is the Hosokawa polynomial of $L$ with
orientation of the second component reversed.
On the other hand, $t^n\Delta_L (t, -t^{-1})$ is always stable, if
$\Delta_L (x, y)$ is stable. 
In Section 16, we discuss mainly the stability problem 
of the Alexander polynomials $\Delta_{L}(x,y)$ of 
$2$-bridge links $L$.
In the last section, Section 17,
we study inversive 2-bridge links using 2-variable
Alexander polynomials $\Delta_{K}(x,y)$.\\
Appendix has three sections.
In Appendix A, we study the stability problem of  integer polynomials
considered in knot theory, particularly,
the stability of Riley polynomials associated to parabolic or dihedral
representation of the knot group. 
Riley studied these representations of the knot groups
$G(K(r))$ of 2-bridge knots $K(r)$. 
He defined an integer polynomial $\theta_K(z)$ associated
to the parabolic representation of $G(K(r))$ to $SL(2,\CC)$. 
It is known \cite{ri72} or \cite{swa}
that $\theta_{K(r)}(z)$ is real stable, if $r = 1/(2n+1)$.
However, if $r \ne 1/(2n+1)$, it is usually not stable.
The second polynomial Riley
studied is an integer polynomial
$\varphi_{2n+1}(z)$ associated to a trace-free representation of
$G(K(1/2n+1))$ onto a dihedral group 
$D_{2n+1} \subset GL(2, \CC)$ (see \cite{ri84}).
In this section, 
we prove that $\varphi_{2n+1}(z)$ is real stable. 
Since $\varphi_{2n+1}(z)$ is not reciprocal, we cannot apply the methods
we used in the previous sections and our approach here is quite different.
In Appendix B, we discuss the maximal values of the real parts of the zeros
of the Alexander polynomials of alternating knots.
Let $\delta (K)$ be the maximal value of the real parts of the zeros of
$\Delta_K (t)$. 
It is proved in Section 4 that even for
alternating knots $K$, $\delta (K)$ is not bounded, 
i.e., given any positive real number $\delta_0$, there exists an alternating knot $K_0$
such that $\delta(K_0) > \delta_0$. 
It should be noted that for a
2-bridge knot $K$, $\delta(K) < 6$ (\cite{LM}, \cite{St}). 
However, for alternating knots, we can modify this invariant as follows.
Let $\Gamma_{2n}$ be the set of all alternating knots $K$ with 
$\deg \Delta_K (t) = 2n$. 
We conjecture that $\delta(K)$ for $K$ in $\Gamma_{2n}$
is bounded, i.e., there exists a positive real number $\delta_{2n}$ such
that $\delta (K) \leq \delta_{2n}$ for $K \in \Gamma_{2n}$.
Further, we conjecture that $\delta_{2n}$ can be achieved by 
fibred alternating knots. 
This seems true for 2-bridge knots.
If the conjecture holds, 
then since the number of alternating fibred knots
in each $\Gamma_{2n}$ is finite, 
we can determine $\delta_{2n}$ for each $n$.
In this section, we prove that the conjecture holds for $n=1$ and $2$.
In the last section, Appendix C, we discuss the distribution of the
zeros of a series of some special type of 2-bridge knots.
It seems that these examples suggest many deep properties of
the distribution of the zeros of the Alexander polynomials of
alternating knots and links.

Finally, we note that some of the theorems in this paper
have been announced without proofs in the survey article \cite{HMsurvey}.

\section{Stability Property}\label{1}

\subsection{Half-plane property}\label{1.1}

Let $\H \subset \CC$ be an open half-plane such that $\partial \overline{\H}$
contains the origin.
Let $f(z_1, \cdots, z_n) \in 
\CC[z_1, \cdots, z_n]$ be a polynomial in $n$
variables.

\begin{dfn}\cite[p.303]{branden}\label{dfn:1.1}
$f \in \CC[z_1, \cdots, z_n]$ is said to be 
{\it $\H$-stable} if $f\equiv 0$ identically, or
for any values $\alpha_j \in \H, 1 \leq j \leq n,
f(\alpha_1, \cdots, \alpha_n) \ne 0$.
If $f(z_1, \cdots, z_n) \in \CC[z_1, \cdots, z_n]$ is 
$\H$-stable for some open half-plane, 
we say $f$ has a {\it half-plane property}.
\end{dfn}

There are two special cases.

\begin{dfn}\cite[p.303]{branden}\label{dfn:1.2}
(1) Let $\H$ be the right-half plane, i.e., 
$\H = \{\alpha \in \CC|{\rm \ Re}(\alpha) >0\}$. 
Then an $\H$-stable polynomial 
$f \in \CC[z_1, \cdots, z_n]$ is called
{\it Hurwitz-stable}. 
In other wards, $f$ is Hurwitz-stable if for any $\alpha_j
\in \CC, 1\leq j \leq n$, such that ${\rm Re}(\alpha_j) > 0$,
$f(\alpha_1, \cdots, \alpha_n) \ne 0$.
(2) Let $\H$ be the upper-half plane, i.e., 
$\H = \{\alpha \in \CC|Im(\alpha)>0\}$. 
Then an $\H$-stable polynomial $f \in \CC[z_1, \cdots, z_n]$ 
is called a {\it stable polynomial}.
\end{dfn}

\begin{rem}\label{rem:1.3}
If a real polynomial $f \in \RR[z_1, \cdots, z_n]$ is stable, 
$f$ is sometimes called {\it real stable}.
\end{rem}

From definitions we see immediately

\begin{prop}\label{prop:1.4}
Let $f(z) \in \RR[z]$ be a real univariate polynomial. 
Then
(1) $f(z)$ is real stable if and only if $f(z)$ has only real zeros.
(2) $f(z)$ is Hurwitz-stable if and only if for any zero 
$\alpha$ of $f(z)$,
${\rm Re}(\alpha) \leq 0$.
\end{prop}

\begin{ex}\label{ex:1.5}
(1) $f(t) = t^4 + 7t^3 +13t^2 +7t + 1$ is real stable and
also Hurwitz-stable.
(2) $f(t) = t^4 +2t^3 -5t^2 + 2t + 1$ is neither real stable, nor
Hurwitz-stable.
\end{ex}

The theorem below is elementary, but useful.

\begin{thm}\cite[Lemma 2.4]{wag}\label{thm:1.6}
The following operations preserve 
$\H$-stability in $\CC[z_1,\cdots, z_n ]$.

\noindent
(a) Permutation: For any permutation $\sigma \in S_n$,\\
\hspace*{5mm}
$f \longrightarrow f (z_{\sigma(1)}, \cdots, z_{\sigma(n)})$

\noindent
(b) Scaling: For any $c \in \CC,$ and $(a_1, \cdots, a_n) \in R_{+}^n$ 
(i.e., $a_j> 0, 1 \leq j \leq n$),\\
\hspace*{5mm}
$ f \longrightarrow cf(a_1 z_1, \cdots, a_n z_n)$

\noindent
(c) Diagonalization: For $\{i,j\}, 1 \leq i, j \leq n$,\\
\hspace*{5mm}
$f \longrightarrow f (z_1, \cdots, z_n) \mid_{z_i = z_j}$

\end{thm}

\subsection{$D$-stable polynomial}\label{1.2}

There is another type of stability.

\begin{dfn}\cite{bb09}\label{dfn:1.6} 
Let $D$ be the unit open disk in $\CC$. 
A polynomial $f(z) \in \CC[z]$ is called
{\it $D$-stable} if for any $\alpha \in D, f(\alpha) \ne 0$.
\end{dfn}

\begin{prop}\label{prop:1.7}
Suppose $f(z) \in \RR[z]$ is reciprocal, i.e., 
$f(z) = z^n f(z^{-1})$ for some $n$. 
If $f(z)$ is $D$-stable,
then all zeros 
$\alpha$ of $f(z)$ are on the unit circle, i.e., $|\alpha |= 1$.
\end{prop}

\begin{dfn}\label{dfn:1.8}
We say that $f(z) \in \CC[z]$ is {\it $c$-stable} if
for each zero $\alpha$ of $f$, $|\alpha|=1$. 
\end{dfn}

\begin{ex}\label{ex:1.8}
(1) $f(t) = 2t^6 - 4t^5 +6t^4 -7t^3 +6t^2 -4t +2$ is
$c$-stable, but not Hurwitz-stable. 
(2) $f(t) = 2t^6 - 4t^5 +6t^4 -9t^3 +6t^2 -4t +2$ is
is neither $c$-stable, nor Hurwitz-stable.\\
(3) $f(t) = 3t^4 -12t^3 + 17t^2 -12t + 3$ is neither $c$-stable, nor
Hurwitz-stable.
\end{ex}

\section{Hurwitz-stability}\label{2}

There are two basic tools to show Hurwitz-stability of 
a real univariate polynomial.

\subsection{Hurwitz-Routh Criterion}\label{2.1}

Let $f(z) = a_0 z^n + a_1 z^{n-1}+\dots + a_n \in \RR[z]$ 
be a real polynomial, where $a_0>0, a_j\in  
\RR, 0 \leq j \leq n$.
Define an $n\times n$ matrix $H_n$ as follows:

\begin{equation}
H_n=
\left[
\begin{array}{cccccc}
a_1&a_0&0&0&\cdots&0\\[11pt]
a_3&a_2&a_1&a_0&\cdots&0\\[11pt]
& &\ddots  & & &\\\
\vdots& &  & &  &\vdots\\
& &  & & \ddots &\\
a_{2n-1}&a_{2n-2}& \cdots& & a_{n+1}&a_n
\end{array}
\right],
\end{equation} 
where we define $a_j=0$ if $j > n$.

For $1 \leq k \leq n$, let $H_k$ be the first $k\times k$ principal submatrix of $H_n$. 
Namely, $H_k$ is the $k \times k$ submatrix consisting of the first 
$k$ rows and columns of $H_n$.

For example, $H_1=[a_1]$ and $H_2=\left[\begin{array}{cc}
a_1&a_0\\
a_3&a_2
\end{array}\right]$.

We say that $f(z)$ is {\it strongly Hurwitz-stable}
(or simply {\it s-Hurwitz-stable})
if any zero of $f(z)$ has a negative real part.

\begin{thm}[Hurwitz-Routh Criterion]\cite[Theorem 8.8.1]{lanc}
\label{thm:2.1}
A real polynomial 
$f(z) = \sum_{j=0}^{n} a_j z^{n-j}, a_0> 0, a_j\in \RR, 
1 \leq j \leq n$,  is strongly Hurwitz-stable  if and only if 
$\det H_k> 0$ for $1 \leq k \leq n$.
\end{thm}

Using Theorem \ref{thm:2.1}, 
we can characterize strongly Hurwitz-stable polynomials with small degrees.

\begin{ex}\label{ex:2.2}
(1) $f(z)=a_0 z + a_1, a_0>0$, is s-Hurwitz-stable if and only if 
$a_1> 0$.\\
(2) $f(z)= a_0 z^2 + a_1 z  +  a_2, a_0> 0$, is s-Hurwitz-stable  if 
and only if  $a_1, a_2> 0$.\\
(3) $f(z)=a_0 z^3+a_1z^2+a_2z+a_3,a_0>0$, is s-Hurwitz-stable if and only if
$a_1, a_2, a_3>0$ and  $a_1a_2>a_0a_3$.\\
(4) $f(z)=a_0z^4+a_1z^3+a_2z^2+a_3z+a_4, a_0>0$ is s-Hurwitz-stable if and only if
(i) $a_1, a_2, a_3, a_4>0$, (ii) $a_1a_2>a_0 a_3$, and (iii) $a_3(a_1a_2-a_0a_3)>a_1^2a_4$.
\end{ex}

\subsection{Lyapunov matrix}\label{2.2}

There is another important tool to study Hurwitz-stability of a real univariate polynomial given by Lyapunov.
Let $f(z)$ be a real polynomial of degree $n$.
Let $M$ be a companion matrix of $f(z)$.

\begin{thm}[Lyapunov, {\cite[Theorem 8.7.2]{lanc}}]\label{thm:2.3}
$f(z)$ is strongly Hurwitz-stable 
if and only if there exist two real positive definite (symmetric) 
matrices $V$ and $W$ such that 

\begin{equation}
VM + M^T V = - W.
\end{equation}
\end{thm}

For convenience, we call $V$ a {\it Lyapunov matrix} associated to $M$. 
It is often quite difficult to find a Lyapunov matrix 
even if $f(z)$ is known to be Hurwitz-stable. 

\begin{ex}
(1) $f(z) = z+a_1$.  Then $M=[-a_1]$.  
If $a_1<0$, Lyapunov matrix does not exist, since $M$ is positive definite.  
If $a_1>0$, then $V = E$ is a Lyapunov matrix associated to $M$ and 
$f(z)$ is s-Hurwitz-stable.\\
(2) Let  $f(z) =  z^2 + a_1 z  +  a_2$.  
If $a_1, a_2 > 0$, then we know $f(z)$ is s-Hurwitz-stable, 
see Example \ref{ex:2.2} (2).  
For example, if $a_1= 3$ and $a_2=4$, i.e., 
$M=\left[\begin{array}{cc}
0&-4\\
1&-3
\end{array}\right]$, then
$V = 
\left[\begin{array}{cc}
7/12&-1/2\\
-1/2&5/6
\end{array}\right]$ 
is a Lyapunov matrix and $W=E$.
\end{ex}

In graph theory, this concept appears in literatures. 
We mention one example.

\begin{ex}[{\cite[Theorem 1.1]{cosw}} and {\cite[p.208]{bb}}]
The spanning-tree polynomial of a connected finite graph is 
Hurwitz-stable and also stable.
\end{ex}

\section{Stable polynomial}\label{3}

\subsection{Multivariate stable polynomials}\label{3.1}

First, we state two basic properties of stable polynomials.

\begin{thm}[{\cite[Lemma 2.4]{wag}}] \label{thm:3.1}
The following operations preserve stability in 
$\CC[z_1, \dots, z_n]$.\\
(a) Specialization: For any $a\in \CC$ with $\mathrm{Im}(a) \ge 0$,   
     $f\rightarrow f(a, z_2, \dots, z_n)$\\
(b) Inversion: If $\deg_{z_1}(f) = d$, 
      $f\rightarrow z_1^d f(-z_1^{-1}, z_2, \dots, z_n)$.\\
(c) Differentiation (or contraction)
$f \longrightarrow \frac{\partial}{\partial z_1}f (z_1, \cdots, z_n)$

\end{thm}

Next, the following theorems give us systematic ways to construct stable polynomials.

\begin{thm}[{\cite[Proposition 2.4]{bb}}] \label{thm:3.2} 
Let $A_i, 1 \le i \le n$, be complex, semi-positive definite 
$m\times m$ matrices and 
$B$ be an $m \times m$ Hermitian matrix. Then,
$f(z_1, \dots, z_n) = \det[z_1 A_1+\dots+ z_n A_n + B]$ is stable.
\end{thm}

As a consequence of Theorem \ref{thm:3.2}, 
we have:

\begin{thm}[{\cite[p.308]{branden}}]
\label{thm:3.3} 
Let $Z ={\rm diag}(z_1, \dots, z_n)$ be a diagonal matrix. 
If $A$ is an $n\times n$ Hermitian matrix, 
then both $\det(Z+A)$ and $\det(E+AZ)$ are stable.
\end{thm}

If $n = 2$, then the converse of Theorem \ref{thm:3.2} holds for a real stable polynomial.

\begin{thm}[{\cite[Theorem 1.13]{bb2}}]\label{thm:3.4}
 (Characterization of real stable polynomials with two variables)
Let $f(x, y) \in \RR[x, y]$.  
Then $f$ is real stable if and only if $f$ is written as 
\begin{equation}
f(x, y) = \pm \det[xA + yB + C],
\end{equation}
\noindent
where $A$ and $B$ are positive semi-definite matrices and 
$C$ is a symmetric matrix of the same order.
\end{thm}

The following theorem claims that the stability of multivariate polynomials 
can be reduced to the stability of univariate polynomials.

\begin{thm}[{\cite[Lemma 2.3]{wag}}]\label{thm:3.5}
A polynomial $f  \in \CC[z_1, \dots, z_n]$ is stable if and only if 
for any $(a_1, \dots, a_n)\in\RR^n$  and $(b_1, \dots, b_n)\in\RR_+^n$ 
(i.e., $b_j>~0, 1 \le j \le n$),
$f(a_1+ b_1 t, \dots, a_n+b_n t) \in \CC[t]$ is stable.
\end{thm}

If a polynomial is of special type, the stability problem could be slightly simpler.

\begin{thm}[{\cite[Theorem 5.6]{branden}}]\label{thm:3.6}
Let $f\in \RR[z_1, \dots, z_n]$ be a multi-affine polynomial, 
(i.e., each variable $z_j$ has degree at most $1$ in each term).
Then $f$ is stable if and only if for all 
$(x_1, \dots, x_n)\in \RR^n$ and for $1 \le i, j \le n$, 
$\Delta_{i j}(f)(x_1, \dots, x_n) \ge 0$,
where $\Delta_{i j}(f) =\dfrac{\partial f}{\partial z_i}\dfrac{\partial f}{\partial z_j}-\dfrac{\partial^2 f}{\partial z_i \partial z_j} f$                       
\end{thm}

\begin{rem} 
If $f$ is not multi-affine, then in Theorem \ref{thm:3.4},
the \lq\lq only if\rq\rq\ part holds, but the \lq\lq if\rq\rq\ part does not.
\end{rem}

\begin{ex}[{\cite[Example 5.7]{branden}}]\label{ex:3.8} 
Let $f = a_{00}+a_{01}y+a_{10}x + a_{11}xy, a_{ij}\in\RR$.
Then $\Delta_{12}(f) = -\left[\begin{array}{cc}
a_{00}&a_{01}\\
a_{10}&a_{11}
\end{array}
\right]$.
Therefore, $f$ is stable if and only if 
$\det [a_{ij}] \le 0$. 
\end{ex}

\begin{thm}[{\cite[p.1]{ww}}]\label{thm:3.9}
Suppose $f\in \CC[z_1, \dots, z_n]$ is homogeneous. 
Then, 
$f$ is Hurwitz-stable if and only if $f$ is stable.
\end{thm}

\subsection{
Real stable univariate polynomials}\label{3.2}
In this subsection, we discuss
real stable univariate polynomials. 
We are particularly interested in them,
since they have many deep properties.

\begin{thm}[{\cite[p.307]{branden}}]\label{thm:3.10}
Let $f(z) = a_0 z^n+ a_1z^{n-1}+ \dots+ a_n \in 
\RR[z]$, 
$a_0\neq 0, a_j  \ge 0, 0 \le j \le n$.
Suppose $f(z)$ is real stable.  
If $a_ia_k\neq0$ for $i < k$, then for any 
$j, i < j < k,  a_j\neq 0$.
Therefore, if $a_n\neq 0$, then all 
$a_j \neq 0$, for $1 \le j \le n$.
\end{thm}  

Theorem \ref{thm:3.10}
shows that it is worth studying a 
sequence of  the coefficients of a real stable polynomial.

\begin{dfn}[{\cite[p.126]{wilf}}]\label{dfn:3.11}
A sequence $\{c_0, c_1, \dots,  c_n\}$ of positive numbers 
is called {\it unimodal} 
if there exist indices $r$, $s$ such that
\begin{equation} \label{siki:3.2}
c_0 \le c_1 \le \dots \le c_r= c_{r+1}=\dots= c_{r+s}\ge c_{r+s+1}  
\ge\dots \ge c_n.
\end{equation}
Further, $\{c_0, c_1, \dots,  c_n\}$ is called 
{\it log-concave} if
\begin{equation}\label{siki:3.3}
c_{j-1}c_{j+1}\le c_j^2\  {\rm for}\   
j = 1, 2,\dots, n-1.
\end{equation}

If \lq\lq$\le$\rq\rq\  is replaced by \lq\lq$<$\rq\rq\ in (\ref{siki:3.3}), 
then it is called {\it strictly log-concave}.
\end{dfn}

For example, binomial coefficients $\left\{
\binom{m}{k}\right\}_{k=0}^m$
is unimodal.

The following theorem is well-known.

\begin{thm} [{\cite[Proposition p.127]{wilf}}]\label{thm:3.12}
If a positive sequence $\{c_0, c_1, \dots,  c_n\}$
is log-concave, then it is unimodal.
\end{thm}

Now we have an important result.

\begin{thm}[{\cite[p.127]{wilf}}]\label{thm:3.13}
Let $f(z) = a_0z^n+ a_1z^{n-1}+ \dots+ a_n\in
\RR[z], a_0\neq 0, 
a_n\neq 0$.  
Suppose $a_j\ge 0, 0 \le j \le n$. 
If $f$ is real stable (and hence $a_j> 0$ for all $j \ge 0$), 
then $\{a_0, a_1, \dots ,a_n\}$ is strictly log-concave, 
and hence it is unimodal.
\end{thm}

In this case, we have either $a_0<a_1<\dots<
a_r>a_{r+1}>\dots>a_n$ or
$a_0<a_1<\dots<a_r=a_{r+1}>a_{r+2}>\dots>a_n$.

In the rest of this paper, we study various stabilities of 
(mostly) the Alexander polynomials of knots and links in $S^3$.

\section{Preliminaries}\label{4}
From this section on, we study polynomials of knots and links from
the view point of stabilities.
In this paper, we make a strict distinction between a knot and a link.
Namely, by a link it means a disjoint union of two or more simple closed
curves in $S^3$. 
If the material can be applied on knots and links as well, we
always write such as \lq\lq knots (or links)"
Unless specified otherwise, we assume that a link is oriented, 
but the orientation is not usually mentioned.

A 2-bridge knot (or link) $K$ is always 
represented by a rational number 
$r=\beta / \alpha$, where we assume
$0\le |\beta|\le \alpha$ and $\gcd(\alpha,\beta)=1$.
When $K$ is a link, $\alpha$ is even. When $K$ is a knot,
then $\alpha$ is odd and $\beta$ is assumed to be even.
Then, $r$ has a unique even continued fraction expansion
of the following form,
where $a_i \neq 0, 1\le i\le m$.
\begin{align*}
r      =\cfrac{1}{2a_1
          -\cfrac{1}{2a_2
          -\cfrac{1}{2a_3
          -\cfrac{1}{\ddots
          -\cfrac{1}{2a_{m-1}
          -\cfrac{1}{2a_m}}}}}}
          \end{align*}
\noindent
This expansion is written as $r = [2a_1 , 2a_2 , \cdots, 2a_m]$.
We call $K$ a 2-bridge knot (or a link) of type 
$r = \beta / \alpha$ and denoted it by $K(r)$. 
For example, $2/3 =[2,2]$ represents a trefoil knot and 
$2/5 =[2,-2]$ represents a figure eight knot.
For $2$-bridge links, we assume 
they are oriented as in Figure 4.1.
For example,
$1/4 =[4]$ represents a non-fibred link
and
$3/4 =[2,2,2]$ represents a fibred torus link.

\figyonichi

Now, we discuss briefly a Seifert surface and Seifert matrix of a 
$2$-bridge
knot or link $K(r)$. 
Since we assume that 
$r = [2a_1 , 2a_2 , \cdots, 2a_m]$,
$K(r)$ has a natural Seifert surface $F$ depicted in 
Fig. 4.2 below.

\figyonni

Let $\alpha_1 , \alpha_2 ,\cdots, \alpha_m$ be (oriented) simple
closed curves on $F$ as are shown in Fig. 4.2. 
Then ${{\alpha_j },1 \leq j\leq m}$
forms a basis for $H_1 (F, \ZZ)$. 
Let $u_{i,j} = lk(\alpha_i^+ ,
\alpha_j )$, where $\alpha_i^+$ denotes the simple
closed curve in $S^3$ that is a slight lift of $\alpha_i$ toward the
positive normal direction. 
Then $M = [u_{i,j}]_{1 \leq i,j \leq m}$
is a Seifert matrix of $K(r)$. 
In this paper, we call $M$ a standard
Seifert matrix of $K(r)$.
It is easy to see from 
Fig. 4.2 that $M$ is as below left (resp. right) when 
$m$ is even (resp. odd).\\
\centerline{
$\left[
 \begin{array}{ccccccc}
a_1     \hsp&\hsp   0 \hsp&\hsp        \hsp&\hsp         \hsp&\hsp\cdots    \hsp&\hsp       \hsp&\hsp0 \\
-1      \hsp&\hsp a_2 \hsp&\hsp1       \hsp&\hsp         \hsp&\hsp          \hsp&\hsp       \hsp&\hsp\\
0       \hsp&\hsp   0 \hsp&\hsp a_3     \hsp&\hsp0        \hsp&\hsp          \hsp&\hsp       \hsp&\hsp\\  
        \hsp&\hsp     \hsp&\hsp-1      \hsp&\hsp a_4      \hsp&\hsp1         \hsp&\hsp       \hsp&\hsp\\
\vdots  \hsp&\hsp     \hsp&\hsp        \hsp&\hsp         \hsp&\hsp\ddots    \hsp&\hsp       \hsp&\hsp\vdots\\
        \hsp&\hsp     \hsp&\hsp        \hsp&\hsp         \hsp&\hsp0         \hsp&\hsp a_{m-1}\hsp&\hsp0\\
0       \hsp&\hsp   0 \hsp&\hsp\cdots  \hsp&\hsp         \hsp&\hsp0         \hsp&\hsp-1     \hsp&\hsp a_m
\end{array}
\right],$
$\left[
\begin{array}{ccccccc}
a_1     \hsp&\hsp   0 \hsp&\hsp        \hsp&\hsp         \hsp&\hsp\cdots    \hsp&\hsp       \hsp&\hsp0 \\
-1      \hsp&\hsp a_2 \hsp&\hsp1       \hsp&\hsp         \hsp&\hsp          \hsp&\hsp       \hsp&\hsp\\
0       \hsp&\hsp   0 \hsp&\hsp a_3     \hsp&\hsp0        \hsp&\hsp          \hsp&\hsp       \hsp&\hsp\\  
        \hsp&\hsp     \hsp&\hsp-1      \hsp&\hsp a_4      \hsp&\hsp1         \hsp&\hsp       \hsp&\hsp\\
\vdots  \hsp&\hsp     \hsp&\hsp        \hsp&\hsp         \hsp&\hsp\ddots    \hsp&\hsp       \hsp&\hsp\vdots\\
        \hsp&\hsp     \hsp&\hsp        \hsp&\hsp         \hsp&\hsp-1         \hsp&\hsp a_{m-1}\hsp&\hsp1\\
0       \hsp&\hsp   0 \hsp&\hsp\cdots  \hsp&\hsp         \hsp&\hsp0         \hsp&\hsp0     \hsp&\hsp a_m
\end{array}
\right]
$
}

Then the Alexander polynomial 
$\Delta_{K(r)}(t)$ of $K(r)$ is defined by
$\Delta_{K(r)} (t) = \det (tM - M^{T})$.

For 2-bridge knots and links, we have other particular
forms of Seifert surfaces depicted in Fig. 4.3, where the
boxes contain some full-twists.
A Seifert surface for a 2-bridge knot or link is obtained
by successively plumbing unknotted twisted annuli,
where the shaded squares indicate the glueing squares for plumbing.
The usual way is depicted in Fig. 4.3 top.
This surface is isotopic to those in Fig 4.2.
The surface depicted in Fig. 4.3 bottom is new.
Note that two surfaces are not in general isotopic, but bound
the same 2bridge knot or link if the corresponding boxes contain
the same number of full-twists.

\yonsan

\noindent
To be more precise, given $r=[2a_1,2a_2, \dots, 2a_n]$, 
let $A_1, A_2,\dots, A_n$ be unknotted annuli
such that $A_i$ has  $a_i$ full-twists.
In Fig 4.3 top, for, $2\le i\le n$,
$A_i$ is plumbed on the negative (resp.
positive) side of $A_{i-1}$ if $i$ is even (resp. odd).
This type of plumbed surface is said to be 
{\it of chain type}.
Meanwhile in Fig 4.3 bottom, every annulus is
plumbed on the negative side of the proceeding annulus.
This type of plumbed surface is said to be 
{\it of twisted chain type}.
A plumbed surface of twisted chain type has
a Seifert matrix of the following form,
which is also said to be of twisted type.\\
\begin{equation}\label{siki:4.1}
\left[\begin{array}{ccccc}
a_1 \hsp&\hsp 1 \hsp&\hsp 0 \hsp&\hsp\cdots \hsp&\hsp 0\\
0 \hsp&\hsp a_2 \hsp&\hsp 1 \hsp&\hsp 0 \hsp&\hsp \vdots\\
\vdots \hsp&\hsp \ddots\hsp&\hsp \ddots\hsp&\hsp \ddots  \hsp&\hsp 0\\
\hsp&\hsp \hsp&\hsp 0\hsp&\hsp a_{n-1} \hsp&\hsp 1\\
0  \hsp&\hsp \cdots\hsp&\hsp  \hsp&\hsp 0\hsp&\hsp a_n
\end{array}\right]
\end{equation}
If $K$ is a link (of $\mu$ components, 
$\mu \geq 2$), $\Delta_K (t_1 , t_2, \cdots, t_{\mu})$
denotes the (multivariate) Alexander polynomial of $K$. 
(\cite{tor})
Then
$\Delta_K (t, t, \cdots,t) / (t-1)^{\mu - 2} := \nabla_K (t)$ 
is called the {\it Hosokawa polynomial} of $K$. (\cite{hoso}) 
The degree of $\nabla_K (t)$ is even and 
$(t-1)\Delta_K (t, t, \cdots,t) = (t-1)^{\mu -1}\nabla_K (t)$ is called the
{\it reduced Alexander polynomial} of a link $K$. 
We denote by $\Delta_K (t)$ the reduced
Alexander polynomial of a link $K$, but we call it simply the Alexander
polynomial of a link $K$. 
Generally, if $M$ is a Seifert matrix of a knot 
(resp. a link) $K$,
then the Alexander polynomial 
(resp. reduced Alexander polynomial) is
defined as $\det(t M-M^T)$.
To study the zeros of a real stable polynomial 
$f(z) \in \RR[z]$, 
it is usually assumed that the leading coefficient is positive.
Since $\Delta_K (t)$ is defined up to $\pm t$, 
we denote by $D_K (t)$
the Alexander polynomial of a knot (or link) $K$ with a positive
leading coefficient and a non-zero constant term. 
We may call 
$D_K(t)$ the {\it normalized Alexander polynomial} of $\Delta_K (t)$,
or the {\it normalization} of $\Delta_K (t)$.
We should note that for a 2-bridge knot or link $K(r)$ with
$r=[2a_1,2a_2,\dots,2a_n]$, 
the normalization $D_{K(r)}(t)$ of $\Delta_{K(r)}(t)$
is given by $D_{K(r)}(t)=\varepsilon \det (tM-M^{T})$, where
$\varepsilon=\prod_{j=1}^{n}\frac{a_j}{|a_j|}$.
For a knot $K$, $\Delta_K (t)$ is reciprocal, namely, 
$\Delta_K (t) = t^n \Delta_K (t^{-1})$
for some integer $n$. 
However, for a link, $\Delta_K (t)$ is not 
necessarily reciprocal, but the Hosokawa polynomial 
$\nabla_K (t)$ is always
reciprocal.

Suppose that $f(t)$ is a reciprocal real polynomial of degree $2n$. 
Let $x = t + \frac{1}{t}$. Then $t^{-n}f(t)$ can be written uniquely
as a real polynomial $F(x)$ in $x$ of degree $n$. 
For convenience, we call $F(x)$
the {\it modified polynomial} or {\it modification} of $f(t)$. 
We should note that
if $f(t) = \Delta_K (t)$ for a knot $K$, 
then $F(x)$ is equivalent to the Conway polynomial of $K$. 
To be more precise,
let $C_K (z) = a_0 z^{2n} + a_1 z^{2n-2} + 
\cdots + a_{n-1} z^2 + a_n$ be the Conway polynomial of a knot $K$. 
Then $F(x) = a_0 (x-2)^n + a_1 (x-2)^{n-1} +
\cdots + a_{n-1} (x-2) + a_n$.
$F(x)$ is not necessarily reciprocal.

For convenience, we call a knot $K$ (or link) {\it (real) stable}, 
{\it $c$-stable} or
{\it bi-stable} if the Alexander polynomial of $K$ is, respectively,
(real) stable, $c$-stable or bi-stable.
Further, we call the complex zero $\alpha$ of $\Delta_K (t)$ with $|\alpha|=1$
a {\it unit complex zero} of $\Delta_K (t)$.
If the bi-stable Alexander polynomial has both real zeros and unit complex
zeros, we call it {\it strictly bi-stable} and such a knot or link
is called {\it strictly bi-stable}. 
If any zero of $\Delta_K (t)$ is neither
real nor unit complex, we say $\Delta_K (t)$ is {\it totally unstable}
and a knot (or link) $K$ is called {\it totally unstable}. 
Therefore, a knot (or link) is classified into five classes: 
stable, $c$-stable, strictly bi-stable, totally unstable and none of them. 
Note that in literature, a stable knot may have appeared in a different sense. 
However, in this paper, we use our terminologies.

A link is never totally unstable, since $\Delta_K (1)=0$.
If a knot $K$ is totally unstable, then $\deg\Delta_K (t)$ is divisible by $4$. 
This is because if $\alpha$ is a complex zero of $\Delta_K (t)$, then so are $\bar{\alpha},
\frac{1}{\alpha}, \frac{1}{\bar{\alpha}}$.
The Hosokawa polynomial $\nabla_K (t)$ of a link can be
totally unstable if $\deg\nabla_K (t)$ is a multiple of $4$.

The stability problem of $\Delta_K (t)$ can be checked by using modified
polynomial of $\Delta_K (t)$ as follows.

\begin{prop}\label{prop:4.1}
Let $F(x)$ be the modified polynomial of 
$\Delta_K(t)$ of a knot $K$. Then $K$ is bi-stable if and
only if $F(x)$ is real stable. 
(Therefore, if $F(x)$ does not have a real zero, $K$ is totally unstable.)
Further, if $K$ is bi-stable, the number of the real zeros of 
$\Delta_K (t)$ is exactly twice the number of the real zeros 
$\alpha$ of $F(x)$ (counting multiplicity)
such that $|\alpha| \geq 2$.
\end{prop}

{\it Proof.}
 We prove a slightly more general statement. 
Let $f(t)$ be a real
reciprocal polynomial of degree $2n$. 
We write

$f(t) = c_0 t^{2n} + c_1 t^{2n-1} + \cdots + c_{2n-1} t + c_{2n},
c_0 > 0, c_j = c_{2n-j},
0 \leq j \leq 2n$.
Express $f(t) = c_0 \prod_{j=1}^{2n} (t- \alpha_j)$ and $F(x) = c_0
\prod_{j=1}^n (x- \beta_j)$.
Suppose $\beta_j$ is real and $|\beta_j | \geq 2$. 
Then $\beta_j$ gives two real
zeros of $f(t)$, since $t - \beta_j + 1/t = 0$ has two real zeros.
However, if $\mid \beta_j | < 2$, then $t - \beta_j + 1/t = 0$ gives two
unit complex zeros. This proves the \lq\lq if\rq\rq\ part.

Conversely, suppose that $f(t)$ is bi-stable. 
If $\alpha_j$ is real, then
$\alpha_j + \frac{1}{\alpha_j}$ is real and hence the corresponding zero of
$F(x)$ is real 
and further, $|\alpha_j + \frac{1}{\alpha_j} | \geq 2$.
If $\alpha$ is unit complex, then $\alpha_j +
\frac{1}{\alpha_j} = \alpha_j + \overline{\alpha}_j$ is
real, and hence the corresponding zero of $F(x)$ is real and
$|\alpha_j + \frac{1}{\alpha_j}| < 2$.
\qed

\begin{rem}\label{rem:4.2}
In \cite{wu}, Wu proved a similar proposition using Conway polynomial.
\end{rem}

\begin{ex}\label{ex:4.n7}
Let $f(t) = t^6 - 3t^5 + 2t^4 - t^3 +2t^2 -3t +1$.
Then $F(x) = (x^3 -3x) - 3(x^2 -2) + 2x - 1 
= x^3 - 3x^2 - x + 5$. 
$F(x)$ has three real zeros, 
two of which are in the interval $(-2, 2)$.
Therefore, $f(t)$ has two real zeros and four unit 
complex zeros, and hence
$f(t)$ is strictly bi-stable.
\end{ex}

In the rest of this section, 
we show four elementary but useful propositions.
The first proposition is well-known.

\begin{prop}\label{prop:4.4} (Min-Max Theorem)
Let $M$ be a real symmetric matrix of order $n$. 
Let $\alpha_1 \leq \alpha_2 \leq \cdots \leq \alpha_n$ be eigenvalues of $M$. 
Then,
\begin{equation*}
\alpha_1 \leq \min\{diagonal\ entries\ of\ M\}{\rm \  and\ }
\alpha_n \geq \max\{diagonal\ entries\ of\ M\}.
\end{equation*}
\end{prop}

Since the second proposition is less known, we give a proof.

\begin{prop}\label{prop:4.5}
(Strong Positivity Lemma)
Let $M = [a_{i,j}]_{1 \leq i,j \leq n}$ be an $n \times n$ real matrix
such that
%
for $i = 1,2, \cdots,n$,
\begin{align}\label{siki:4.2}
&(1)\ a_{i,i} > 0,
 \nonumber\\
&(2)\ 
a_{i,i} > \mid a_{i,1} \mid + \mid a_{i,2} \mid + \cdots +\mid
a_{i,i-1} \mid + \mid a_{i,i+1} \mid + \cdots + \mid a_{i,n} \mid.
\end{align}

Then $\det M > 0$.
\end{prop}

For convenience, a row that satisfies (2)  is called {\it excessive}.

{\it Proof.} If $n=1$ or $2$, the proposition is trivially true. Suppose the
proposition holds for $(n-1)\times(n-1)$ matrices.
Let $P = [\lambda_{i,j}]$ be an $n \times n$ matrix of the form
%
\begin{align}\label{siki:4.3} 
&(1)\ \lambda_{i,i}= 1,  1 \leq i \leq n,\nonumber\\
&(2)\ \lambda_{1,k}= -\frac{a_{k,1}}{a_{1,1}}, 1 \leq k \leq n,\nonumber\\
&(3)\ \lambda_{i,j} = 0, {\rm \ for\ } i \ne j {\rm \ or}\ i\ne 1.
\end{align}

Then $PM =
\left[
\begin{array}{c|ccc}
a_{11}&a_{12}&\cdots&a_{1n}\\
\hline
0& & &\\
\vdots& & M'& \\
0& & & 
\end{array}
\right]
$, where
$M^{\prime} =[ a_{i,j}^{\prime}]_{2 \leq i,j \leq n}$ 
is an $(n-1) \times (n-1)$ matrix of the form:
For $i, j =2,\cdots,n$,
\begin{align}\label{siki:4.4}
a_{i,j}^{\prime} = a_{i,j} - a_{1,j}\frac{a_{i,1}}{a_{1,1}}.
\end{align}

We claim that $a_{i,i}^{\prime} > 0$ and every row of $M^{\prime}$ is
excessive.
If $a_{i,1}= 0$, then $a_{i,2}^{\prime} = a_{i,2}, \cdots,
a_{i,n}^{\prime} = a_{i,n}$ and hence $a_{i,i}^{\prime} = a_{i,i}>0$,
and since the $i^{\rm th}$ row of $M$ is excessive, 
the $i^{\rm th}$ row of $M^{\prime}$ is
also excessive.
Suppose $a_{i,1} \ne 0$. 
We show
%
\begin{equation}\label{siki:4.5}
a_{i,i}^{\prime} > | a_{i,2}^{\prime} |+ | a_{i,3}^{\prime}| +
\cdots+|a_{i,i-1}^{\prime}|+ | a_{i,i+1}^{\prime} | + \cdots+
|a_{i,n}^{\prime} |.
\end{equation}
First, for $j \ne i$ and $j \geq 2, | a_{i,j}^{\prime}|=| a_{i,j} -
a_{1,j}\frac{a_{i,1}}{a_{1,1}}| \leq | a_{i,j}|+
\frac{|a_{1,j}| |a_{i,1}|}{a_{1,1}}$
and hence
%
\begin{equation}\label{siki:4.6}
\sum_{j=2, j \ne i}^n | a_{i,j}^{\prime} | 
\leq \sum_{j=2, j\ne i}^n | a_{i,j}| + \frac{|a_{i,1}|}{a_{1,1}} \sum_{j=2, j\ne i}^ n
|a_{1,j}|.
\end{equation}
Next, for $i \geq 2$, since $a_{i,i} > |a_{i,1}|$ and 
$a_{1,1} > |a_{1,i}|$, we see
%
\begin{equation}\label{siki:4.7}
a_{i,i}^{\prime} =a_{i,i}- \frac{a_{1,i} a_{i,1}}{a_{1,1}} \geq
a_{i,i}- \frac{|a_{1,i}| |a_{i,1}|}{a_{1,1}} > 0.
\end{equation}

From (\ref{siki:4.6}) and (\ref{siki:4.7}), we have;
%
\begin{align}\label{siki:4.8} 
a_{i,i}^{\prime} - \sum_{j=2, j \ne i}^n | a_{i,j}^{\prime}| 
&\geq a_{i,i} - \frac{|a_{1,i}| |a_{i,1}|}{a_{1,1}} -
\sum_{j=2, j \ne i}^n | a_{i,j} | -
\frac{|a_{i,1}|}{a_{1,1}}\sum_{j=2, j\ne i} |a_{1,j}| 
\nonumber\\
&=  a_{i,i} - \sum_{j=2, j \ne i} |a_{i,j} | -
\frac{|a_{i,1}|}{a_{1,1}}\left( |a_{1,i}| + \sum_{j=2, j \ne i}^n 
|a_{1,j}|\right)
\nonumber\\
&= a_{i,i} - \sum_{j=2, j \ne i}^n | a_{i,j}| - 
\frac{|a_{i,1}|}{a_{1,1}} \sum_{j=2}^n |a_{1,j}|.
\end{align}

By assumption,
\begin{align}\label{siki:4.9} 
a_{i,i} - \sum_{j=2, j \ne i}^n |a_{i,j}|
> |a_{i,1} | {\rm \  and}
\nonumber\\
\delta = \sum_{j=2}^n \frac{| a_{1,j} |}{a_{1,1}} < 1.
\end{align}
Therefore, from (\ref{siki:4.9}), we have;
%
\begin{equation}\label{siki:4.10} 
a_{i,i}^{\prime} - \sum_{j=2, j\ne i}^n |a_{i,j}^{\prime} |
> |a_{i,1} | - | a_{i,1} | \delta > 0.
\end{equation}
This proves (\ref{siki:4.5}).

Since all rows of $M^{\prime}$ are excessive and $a_{i,i}^{\prime} > 0, 
2\leq i \leq n$, 
by (\ref{siki:4.7}),
it follows by induction that 
$\det M^{\prime} > 0$ and hence 
$\det M =\det(PM)= a_{1,1} \det M^{\prime} > 0$. 
\qed

\begin{rem}\label{rem:4.6}
In the above proof, even if the $i^{\rm th}$ row of $M$ is not excessive,
i.e., $a_{i,i} = | a_{i,1}| +| a_{i,2} |+ \cdots+ 
|a_{i,i-1} | + |a_{i,i+1}|+ \cdots + | a_{i,n}|$, 
a new $i^{\rm th}$ row of $M^{\prime}$ becomes excessive,
provided $a_{i,1} \ne 0$ and the first row is excessive. 
In fact,
the first inequality in 
(\ref{siki:4.9}) becomes an equality, 
but the second
inequality holds, and hence 
(\ref{siki:4.10})
holds.
\end{rem}
Let $P$ be a matrix representing an elementary matrix operation. 
Then we
say, for convenience, that a matrix 
$PMP^{T}$ is obtained from $M$ by
applying
an elementary matrix operation on the rows (and columns) 
of $M$ simultaneously. 
For example, we say something like that a new matrix
$M^{\prime}$ is
obtained by interchanging the $i^{\mathrm th}$ row (and column) 
and the $j^{\mathrm th}$ row(and column) of
$M$ simultaneously.

Using Remark 
\ref{rem:4.6}, 
we can obtain the same conclusion under a slightly
weaker assumption than that of 
Proposition 
\ref{prop:4.5}. 

\begin{prop}\label{prop:4.7}
(Positivity Lemma)
Let $M = [a_{i,j}]_{1 \leq i,j \leq n}$ be an 
$n\times n$ real matrix.
Assume that
\begin{align} \label{siki:4.11} 
(1) &\ M {\mathrm \ cannot\ be\ transformed\ into\ a\ form:\ } 
\left[\begin{array}{cc}
A&B\\
O&C
\end{array}\right]
{\mathrm \  or\ }
\left[\begin{array}{cc}
A&O\\
B&C
\end{array}\right]
\nonumber\\
 &{\mathrm  \ by\ a\ sequence\ of\ exchanges\ of\ the}\ 
i^{\mathrm th} {\rm \  row\ (and\ column)\ and\ the\ } 
j^{\mathrm th} {\rm\ row\ }
\nonumber\\
&{\rm  (and\ column)\ simultaneously,\ where\ } 
A {\rm \  and}\ 
C {\rm \ are\ square\ matrices}.
\nonumber\\
(2)\ &{\rm  For\ } i=1,2, \cdots,n,\ a_{i,i} > 0.
\nonumber\\
(3)\ &{\rm  For}\ i=1,2, \cdots, n,\
a_{i,i} \geq | a_{i,1}| + \cdots+ | a_{i,i-1} | + | a_{i,i+1}|
+ \cdots + | a_{i,n}|,
\nonumber\\
(4)\ &M {\rm \  has\ at\ least\ one\ excessive\ row}.
\end{align}

Then $\det M > 0$. 
Further, if $M$ is symmetric, then $M$ is positive definite.
\end{prop}

\begin{rem}\label{rem:4.8}
If (4) is dropped, then $\det M \geq 0$. 
Suppose (1) is dropped. 
Then if each of the block matrices $A$ and $C$ satisfies 
(\ref{siki:4.11}) 
(1) - (4),
then $\det M > 0$.
\end{rem}

{\it Proof of Proposition \ref{prop:4.7}. }
We may assume without loss of generality that $M$ has been arranged
by a sequence of exchanges of rows and columns, simultaneously,
in such a way that
%
\begin{align}\label{siki:4.12}
&(1)\ {\rm the\ first\ } k\ {\rm rows\ are\ excessive,\ 
but\ each\ of\ the\ remaining\ rows\ are\ not}
\nonumber\\
&(2)\  {\rm for\ each\ }
i=k+1, \cdots,n, {\rm \  at\ least\ one\ of\ }
a_{i,1}, a_{i,2}, \cdots, a_{i,i-1}
{\rm \  is\  not\ } 0.
\end{align}

Write $M =\tbt{A}{B}{C}{D}$, 
where $A$ is a $k \times k$ matrix. 
Starting from the first row, 
we can transform $A$ into an upper triangular matrix
$A^{\prime}$
by a sequence of row operations (without changing the value of the
determinant) to get $M^{\prime} = \tbt{A^{\prime}}{B^{\prime}}{C}{D}$,
where $A^{\prime} = 
\left[ \begin{array}{cccc}
\hspace{-1mm}a_{11}\hspace{-1mm} &     &        & \hspace{-1mm}* \\[-4pt]
   &\hspace{-2mm} a_{22}' \hspace{-2mm}&        &                \\[-4pt]
                   &     & \ddots&                \\[-4pt]
\hspace{-1mm}O\hspace{-1mm} &     &        &\hspace{-1mm} a_{kk}'\hspace{-1mm}
  \end{array} \right]
$.
By Remark \ref{rem:4.6},  
each row of the first $k$ rows is excessive.
Now by (\ref{siki:4.12}), we may assume without loss of generality that
one of $a_{k+1,1}, \cdots, a_{k+1,k}$ is not $0$.
Apply row operations repeatedly on the $(k+1)^{\rm st}$ row
so that the first 
$k$ entries of the $(k+1)^{\mathrm st}$ row of 
$M^{\prime}$ become $0$,
and also, the new $(k+1)^{\mathrm st}$ row is excessive by 
Remark \ref{rem:4.6}.  
We can repeat
this until all rows are excessive. 
Then apply Proposition  \ref{prop:4.5}, 
to show that $\det M>0$. Furthermore, the same argument
can be applied to show that all principal minors of
$M$ are positive, and hence if $M$ is symmetric, then $M$
is positive definite.
\qed

Finally, we prove the following proposition.

\begin{prop}\label{prop:4.9}
Let $M$ be a real matrix of the form:\\
$M =
\tbt{A}{O}{H}{B}$
or $\tbt{A}{H}{O}{B}$,
where $A$ is a $p\times p$ positive definite symmetric matrix and 
$B$ a $q \times q$
negative definite symmetric matrix and 
$H$ an arbitrary matrix. 
Then
$M^{-1} M^{T}$ is conjugate to a symmetric matrix in 
$GL(p+q, \RR)$.
Therefore, the characteristic polynomial $f(t)$ of the matrix 
$M^{-1} M^{T}$
is real stable.
\end{prop}

{\it Proof.} 
First, if $A=E_p$ and $B=-E_q$, then
$M^2 =E_{p+q}$, and hence, $M^{-1}M^{T} = M M^{T}$ 
is symmetric.
Now, consider the general case. 
Let $P_A$ and $P_B$ be, respectively,
matrices which diagonalize $A$ and $B$. 
Write $P_A A P_A^{T} =
\diag\{a_1, a_2, \cdots, a_p\}, a_j > 0$, 
for $1 \leq j \leq p$, and
$P_B B P_B^{T} = \diag \{ -b_1, -b_2, \cdots, -b_q\}$,
$b_j > 0$, for $1 \leq j \leq q$.
Let $D_a = \diag \{ 1/ \sqrt{a_1}, 1/ \sqrt{a_2}, \cdots, 1 / \sqrt{a_p}\}$
and $D_b = \diag \{ 1/ \sqrt{b_1}, 1/ \sqrt{b_2}, \cdots, 1 / \sqrt{b_q}\}$. 
Further, let $P = P_A \oplus P_B$ and 
$D = D_a \oplus D_b$.
Then a simple computation shows that
$
D P M P^{T} D^{T} = M_0, {\rm \  where\ } 
M_0 = 
\left[\begin{array}{cc}
E_p &O \\
H_0 & -E_q
\end{array}\right]$.
Now since $M_0^2 = E_{p+q}$, it follows that $M = P^{-1} D^{-1} M_0 (D^{T})^{-1}
(P^{T})^{-1}$ and $M^{-1}M^{T} = P^{T}D^{T} M_0 D P P^{-1}D^{-1}
M_0^{T} (D^{-1})^{T} (P^{-1})^{T}$ = $P^{T}D^{T} M_0 M_0^{T}
(D^{T})^{-1}(P^{T})^{-1}$, and hence, $M^{-1}M^{T}$
is conjugate of a symmetric matrix $M_0 M_0^{T}$. 
\qed

\section{The Alexander polynomials of alternating knots}\label{5}

Before we concentrate on the study of various stabilities
of knots or links, we discuss, in this section, some connection between the stability of alternating knots or links and various conjectures in Knot theory.

\subsection{Hoste's Conjecture}\label{5.1}

In 2002, based on his extensive calculations of the zeros of 
the Alexander polynomials, Hoste made the following conjecture.

\begin{yosou}[J. Hoste, 2002]\label{conj:5.1}
Let $K$ be an alternating 
knot and $\Delta_K(t)$ the Alexander polynomial of $K$. 
Then for any zero $\alpha$ of $\Delta_K(t)$, Re$(\alpha) > -1$.
\end{yosou}

One of the key observations is that Conjecture \ref{conj:5.1} 
is equivalent to the following 

\begin{yosou}\label{conj:5.2} 
Under the same assumption, 
$\Delta_K  (-(t+1))    \in\RR[t]$ is strongly Hurwitz-stable.
\end{yosou}

Using Lyapunov matrices, the following theorem is proved.

\begin{thm}[{\cite[Theorem 1]{LM}}]\label{thm:5.3}
 Let $K$ be a $2$-bridge knot (or link). 
 Then $\Delta_K (- (t+3))$ and 
 $\Delta_K  (t+6)$ are strongly Hurwitz-stable. 
Equivalently, any zero $\alpha$ of $\Delta_K(t)$ satisfies

\begin{equation}
-3 < Re(\alpha) < 6.
\end{equation}
\end{thm}

For other special results, see \cite[Theorems 3,4 and 5]{LM}.  

\begin{rem}\label{rem:5.4} 
A.Stoimenow proves in \cite{St} that for a 2-bridge knot (or link) $K$, any zero 
$\alpha$ of $\Delta_K (t)$ satisfies
\begin{equation}
\left| \sqrt{\alpha}   -\dfrac{1}{\sqrt{\alpha}}
\right| < 2
\end{equation}

This implies 
\begin{equation}
-1< {\rm Re}(\alpha)< 3+\sqrt{8}=5.8284...
\end{equation}
\end{rem}

It should be noted that for a non-alternating knot, 
neither a lower bound nor an upper bound of Re($\alpha$) exist 
\cite[Examples 1 and 2]{LM}.
Further, we think that an upper bound of 
Re($\alpha$) does exist only for a family of $2$-bridge knots or links. 
In fact, there exists an infinite sequence of alternating stable (Montesinos) knots 
$K_1, K_2, \dots,  K_m,\dots$ such that 
the maximal value of the zeros of 
$\Delta_{K_m}(t)$ is at least $m+1$.
(See Theorem \ref{thm:15.2} Case 3.)
Therefore, in general, an upper bound of Re($\alpha$) does not exist, 
even for alternating knots.
However, an upper bound may exist for some family of the Alexander polynomials.  
For example, let $\Gamma_n$ be the set of all Alexander polynomials 
(of degree $n$) of alternating knots.

\begin{yosou}\label{conj:5.6}
There exist a real number $\delta_n>0$   
such that for any zero $\alpha$ of $\Delta_K(t)$ in $\Gamma_n$ 
\begin{equation}
Re(\alpha) \leq  \delta_n
\end{equation}
\end{yosou}    

It is known that Conjecture \ref{conj:5.6} 
is false for non-alternating knots \cite[Example 2]{LM}. 
Since the Alexander polynomial of an alternating knot $K$ is of the form
$\Delta_K(t)=\sum_{j=0}^{2n}(-1)^j c_j t^{2n-j}, c_j>0,0\le j\le 2n$,
it follows that if $\Delta_K(t)$ is real stable, then all the zeros are positive
and hence Conjecture \ref{conj:5.1} holds. 
Therefore, we have:

\begin{thm}\label{thm:5.6}
Let $K$ be an alternating knot.
If $K$ is bi-stable, then Conjecture \ref{conj:5.1} holds for $K$.
\end{thm}

\subsection{Trapezoidal Conjecture}\label{5.2}
Let
$\Delta_K(t)=\sum_{j=0}^{2n}(-1)^j c_j t^{2n-j}$
be the Alexander polynomial of an alternating knot $K$.
Then the Trapezoidal conjecture claims:

\begin{yosou}\label{conj:5.7} \cite{fox62}
There is an integer $k$, $1\le k\le n$, such that
\begin{equation}\label{siki:5.5}
c_0<c_1<\dots<c_k=c_{k+1}=\dots=c_{2n-k}>c_{2n-k+1}\dots>c_{2n}
\end{equation}
\end{yosou}

This conjecture has been proven for several families of alternating knots,
but not proven in general,
\cite{hart}, \cite{mu85} etc.
Further, this conjecture does not hold for Hosokawa polynomials of 
alternating links.
For example, a 2-bridge link $K(r), r=11/14$ is
bi-stable, and satisfies Trapezoidal conjecture, since 
$\Delta_K(t)=t^5-3t^4+3t^3-3t^2+3t-1$, but its Hosokawa polynomial
$\nabla_K(t)=\Delta_K(t)/(1-t)=t^4-2t^3+t^2-2t+1$ does not.

The trapezoidal property of the coefficients is 
quite similar to the unimodality of a sequence
considered in Section 3.

For an alternating knot $K$, $\Delta_K(-t)=\sum_{j=0}^{2n} c_j t^{2n-j}$
satisfies all assumptions of Theorem \ref{thm:3.13} and hence the coefficient sequence
is strictly log-concave.
Therefore, we have
\begin{equation}\label{siki:5.6}
c_0<c_1<\dots<c_n>c_{n+1}>c_{n+2}>\dots>c_{2n}
\end{equation}
and hence we obtain the following:

\begin{thm}\label{thm:5.8}
For an alternating stable knot (or link) $K$,
Trapezoidal conjecture holds.
\end{thm}

We should note that if the number of components of an alternating stable link
$K$ is even, then $\deg\Delta_K(t)$ is odd, say $2n+1$, and
(\ref{siki:5.6}) should be replace by (\ref{siki:5.7}) below.
\begin{equation}\label{siki:5.7}
c_0<c_1<\dots<c_n=c_{n+1}>c_{n+2}>\dots>c_{2n+1}
\end{equation}

For a knot $K$, there is another necessary condition for $\Delta_K(t)$ 
to be stable.

\begin{prop} \label{prop:5.9}  
If $K$ is a stable knot, 
then the signature $\sigma(K)$ of $K$ is zero.
\end{prop}

In fact, if the signature is not zero, 
$\Delta_K(t)$ has at least two zeros on the unit circle \cite{mi}. 
However, the converse of Proposition \ref{prop:5.9} 
is false.
For example let $K(r)$ be the 2-bridge knot where $r=[2,2,-4,-2]$.
Then, $\sigma(K)=0$ but the zeros of $\Delta_K(t)=2 - 6 t + 9 t^2 - 6 t^3 + 2 t^4$ 
are $1\pm i$ and $\frac{1\pm i}{2}$.
Further, the following example 
shows that for links, 
Proposition \ref{prop:5.9} does not hold.

\begin{ex}\label{ex:5.10} 
Let $L$ be an alternating pretzel link $P(2,4,4)$, oriented so that $L$ is a 
special alternating $3$-component link.  
Then the reduced Alexander polynomial $\Delta_L(t) = 8(t-1)^2$ that is stable,
while $\sigma(L)=2$.
\end{ex}

We suspect that Hoste's conjecture and the Trapezoidal conjecture are independent.
However, for alternating knots the condition $\sigma(K)=0$ may imply (\ref{siki:5.6}).
Therefore, we propose the following conjecture.

\begin{yosou}\label{conj:5.11}
Let 
$\Delta_K(t)=\sum_{j=0}^{2n}(-1)^j c_j t^{2n-j}, c_j>0$
be the Alexander polynomial of an alternating knot $K$.
If $\sigma(K)=0$, then the coefficient sequence satisfies (\ref{siki:5.6}),
i.e.,
$c_0<c_1<\dots<c_n>c_{n+1}>c_{n+2}>\dots>c_{2n}
$.
More generally, if $\sigma(K)=2k$, then the coefficient sequence satisfies
(\ref{siki:5.8}) below:
\begin{equation}\label{siki:5.8}
c_0<c_1<\dots<c_{n-m-1}<c_{n-m}=\dots=c_{n+m}>c_{n+m+1}>\dots>c_{2n},
\end{equation}
\noindent
where $m\le k$.
\end{yosou}

This conjecture is quite likely true for $2$-bridge knots.
However, it is false for non-alternating knots.
In fact, the signature of a non-alternating knot $10_{132}$ is $0$,
but the Alexander polynomial is $t^4-t^3+t^2-t+1$.

\begin{rem}\label{rem:5.12}
In graph theory, the concept of unimodality has been used in \cite{gr}.
Very recently, the following long-standing conjecture was proved in \cite{huh}.
\end{rem}

\begin{yosou}\label{conj:5.12}\cite[p. 534]{gr}
The sequence of the coefficients of the 
chromatic polynomial is unimodal.
\end{yosou}

\section{Construction of real stable knots (I)}\label{6}

The first family of alternating stable knots or links was given in
Theorem \ref{thm:6.1} below.
In the proof, Seifert surfaces in Fig. 4.2 were used to show that
there exists a symmetric companion matrix $M$
of the Alexander polynomial of $K(r)$, 
and hence, all the eigenvalues of $M$ are real.

\begin{thm}[{\cite[Theorem 2]{LM}}] \label{thm:6.1} 
Let $r=[2a_1, 2a_2 , \cdots, 2a_m]$ be an
even continued fraction expansion of a rational number
$r = \beta/ \alpha$.
If the sequence $\{a_1, a_2, \cdots, a_m\}$ alternates in sign, i.e., 
$a_j a_{j+1} < 0$, $1\le j\le m-1$  then the 
$2$-bridge knot (or link) $K(r)$ is real stable.
\end{thm}

In this section, first, we construct a new surface that is a slight
generalization of a Seifert surface for a $2$-bridge knot (or link),
and then, using these surfaces, we define a knot (or link), called a
quasi-rational knot (or link), and generalize Theorem \ref{thm:6.1}.
 
\subsection{Quasi-rational knots} \label{6.1} 
In this subsection, we define the class of
\lq\lq quasi-rational links\rq\rq\ 
which is a generalization of 2-bridge links.

\begin{dfn}\label{dfn:6.2}
Let $D$ be a disk with two families 
$\Gamma_1=\{\alpha_1,\dots, \alpha_p\}$,
$\Gamma_2=\{\beta_1,\dots,\beta_q\}$ of properly embedded arcs in $D$, 
where no arcs share their end points on $\partial D$. 
In each family, the arcs are disjoint, but
arcs from different families may intersect one another.
Each arc, say $\gamma$, is assigned with a non-zero integer 
$w(\gamma)$ called a weight.
Push the interior of $\alpha_i$'s (resp. $\beta_j$'s) in the
positive (resp. negative) normal direction of $D$,
and along each pushed arc $\gamma'$, attach a band to $D$ with
$w(\gamma)$ half-twists.
The boundary of the resulting surface $F(\Gamma_1,\Gamma_2)$ is called
a {\it quasi-rational} knot or link.
Conventionally, the arcs in $\Gamma_2$ are depicted in dotted lines. 
(See Fig. 6.1.)
\end{dfn}

The surface $F$ is
orientable if and only if all weights are even.
In that case, $F(\Gamma_1,\Gamma_2)$ is a Murasugi sum of two Seifert surfaces 
$F(\Gamma_1,\emptyset)$ and $F(\emptyset,\Gamma_2)$
along $D$, where the summands are respectively
connected sums of elementary torus links. 
Hence by \cite{ga}, $F(\Gamma_1,\Gamma_2)$ is of minimal genus.

\begin{ex} 
As particular example, we see
in Fig. 4.2 or Fig. 6.1 that
2-bridge knots and links $K(r)$ are 
quasi-rational.\\
For weights $\{2w(\alpha_1),2w(\alpha_2), \dots, 2w(\alpha_p)\}$,
$\{2w(\beta_1),2w(\beta_2), \dots, 2w(\beta_q)\}$, a nice Seifert matrix
$M=
\left[\begin{array}{cc} A_1&O\\C&A_2\end{array}\right]$ 
is obtained from $F(\Gamma_1,\Gamma_2)$, where 
$A_1=\diag\{w(\alpha_1),w(\alpha_2), \dots, w(\alpha_p)\},
A_2=\diag\{w(\beta_1),w(\beta_2), \dots, w(\beta_q)\}$ and
in $C$, all diagonal entries are $-1$, $(k,k+1)$-entries 
are $1$ and 
the other entries are $0$.
For convenience, we say that $M$ is of {\it split type}.
\end{ex}

\rokuichi

\subsection{Stable quasi-rational knots}\label{6.2}
In this subsection we generalize Theorem \ref{thm:6.1} to some classes 
in quasi-rational knots or links.

\begin{thm}\label{thm:6.3}
Let $F(\Gamma_1,\Gamma_2)$ be a Seifert surface for a 
quasi-rational knot or link.
Suppose the weights for $\Gamma_1$ (resp. $\Gamma_2$)
are all positive (resp. negative) and even.
Then $L=\partial F(\Gamma_1,\Gamma_2)$ is alternating and stable.
\end{thm}

{\it Proof.}
By turning the arcs in $\Gamma_1$ outside of the disk $D$,
we have an alternating diagram for $L$.

Let $\{2w(\alpha_1),2w(\alpha_2), \dots, 2w(\alpha_p)\}$ 
and $\{2w(\beta_1),2w(\beta_2), \dots, 2w(\beta_q)\}$ 
be weights for $\Gamma_1$ and $\Gamma_2$,
where $w(\alpha_i)$'s are positive and $w(\beta_j)$'s are negative.
Take a natural bases of $H_1(F(\Gamma_1,\Gamma_2))$, where
each loop is a union of the core curve of an attached band and a curve
in $D$. The orientations of the loops are arbitrary.
Then we have a Seifert matrix
$M=
\left[\begin{array}{cc} A_1&O\\C&A_2\end{array}\right]$, where 
$A_1=\diag\{w(\alpha_1),w(\alpha_2), \dots, w(\alpha_p)\},
A_2=\diag\{w(\beta_1),w(\beta_2), \dots, w(\beta_q)\}$.
Since $\Delta_L(t)=\det(t M-M^T)=(\det M)\det(t E - M^{-1} M^T)$,
it follows from Proposition \ref{prop:4.9} that $\Delta_L(t)/\det M$ is
the characteristic polynomial of a symmetric matrix, and hence
$\Delta_L(t)$ is stable.
\qed

Finally, we see that the signature of $L$ is $p-q$. 
Since a stable knot has signature $0$, 
it follows that if $p \ne q$, $L$ is not a stable knot,
but a stable link.

\subsection{Examples}\label{6.3}

In this subsection, we construct two series of stable knots, both of which
are quasi-rational knots. 
These knots contain first two alternating non-$2$-bridge stable knots
whose Alexander polynomials have the real zeros larger than $6$.

\begin{ex}\label{ex:6.4}
The family of knots denoted by\\
\centerline{
$X_n(2a_1, 2a_2, \cdots, 2a_n \mid 
2b_1, 2b_2, \cdots, 2b_n)$
} is depicted in Fig.6.2.
For example, 
$X_1 (2 \mid-2)$ is $4_1$ and $X_2 (2,2 \mid-2,-2)$
is
$8_{12}$ and $X_3 (2,2,2 \mid-2,-2,-2)$ is $12_{a0125}$.
If all $a_j, 1 \leq j \leq n$, are positive and all $b_j, 1 \leq j
\leq n$, are negative, then $X_n$ is alternating and always stable.
A Seifert surface is obtained by applying Seifert's algorithm to
Fig.6.2 (right).
\end{ex}

\rokuni

\begin{ex}\label{ex:6.5}
The family of knots denoted by\\
\centerline{
$Y_{2n+1} (2a_1 , 2a_2, \cdots,
2a_{2n+1} \mid 2b_1, 2b_2, \cdots, 2b_{2n+1})$
} is depicted in Figure 6.3, together with a
Seifert matrix in the case of $n=3$.
For example, $Y_1 (2 \mid-2)$ is $4_1$ and $Y_3 (2,2,2 \mid-2,-2,-2)$
is $12_{a1124}$. As before, if all $a_j, 1 \leq j \leq 2n+1$,
are positive and all $b_j, 1 \leq j \leq 2n+1$, are negative, then
$Y_{2n+1}$ is alternating and always stable.
\end{ex}

\rokusan

We note that $12_{a0125}$ and $12_{a1124}$ are only alternating knots
with at most $12$ crossings such that the real part of the zeros is larger
than $6$.
In fact, they are $6.904\cdots$ for $12_{a0125}$ and $7.699\cdots$ for
$12_{a1124}$. Furthermore, these values increase as $n$ increases
unboundedly
if all $a_j =1$ and all $b_j = -1$, 
as is proved in Theorem \ref{thm:6.6} below. 

\begin{thm}\label{thm:6.6}
(1) Let $K_n^{(1)} = X_n (2,2, \cdots,2 \mid -2,-2,\cdots,-2)$. 
Then
$K_n^{(1)}$ is stable and the maximal value of 
the zeros is at least $n+1$.

(2) Let $K_{2n+1}^{(2)} = Y_{2n+1}(2,2, \cdots,2 \mid -2,-2, \cdots,-2)$.
Then $K_{2n+1}^{(2)}$
is stable and the maximal value of the zeros is at least $2n+1$.
\end{thm}

{\it Proof.} 
(1) In the proof of Theorem \ref{thm:6.3}, 
the matrix $C$ is a lower triangular
matrix such that every entry under the diagonal entries and diagonal entry
as well are $1$, but $0$ elsewhere. 
Since all $a_j$ are $1$ and all $b_j$ are $-1$, 
it follows from Proposition \ref{prop:4.9} that a companion matrix of the Alexander
polynomial of $K_n^{(1)}$ is $\tbt{E_n}{C^T}{C}{E_n + C C^T}=S$.
It is obvious  that the maximal value of the diagonal 
entries of $S$ is $n+1$.
This proves (1) by Min-Max Theorem.

(2) 
In this case $C$ is like the lower left matrix in Fig.6.3.
Then by the same argument, we see that the maximal value
of the diagonal entries of 
$\tbt{E_{2n+1}}{C^T}{C}{E_{2n+1} +C C^T}$
is $2n+1$. 
This proves (2). 
\qed

These two quasi-rational knots are quite unlikely Montesinos knots.
However, the infinite sequences $\{K_n^{(1)}\}$ and $\{K_{2n+1}^{(2)}\}$ give the
evidence
that the maximal value of the real parts of the zeros of
Alexander polynomials of alternating knots is unbounded.
Further, the maximal zero $7.699\cdots$ of the Alexander polynomial of
$12_{a1124}$ is quite likely the number $\delta_6$ defined in Section 5.
(See Appendix B.)

\begin{rem}\label{rem:6.8}
Since $K_n^{(1)}$ and $K_{2n+1}^{(2)}$ 
both are strongly negative amphicheiral, 
their Alexander polynomials $\Delta(t)$ have the property:
$\Delta(t^2) = f(t) f(-t)$ for some $f(t) \in \ZZ[t]$.
See \cite{hartleyKawauchi}.
Therefore, their Conway polynomials are of the form: 
$g(z) g(-z)$ for some $g(z) \in \ZZ[z]$.
Furthermore, it is easy to show that 
$X_n(2a_1, 2a_2, \cdots, 2a_n| 2b_1,2b_2, \cdots, 2b_n)$ 
is strongly negative amphicheiral
if $b_j = -a_{n+1-j}, 1 \leq j \leq n$. 
Similarly, 
$Y_{2n+1}(2a_1, 2a_2,\cdots, 2a_{2n+1}| 
2b_1, 2b_2, \cdots, 2b_{2n+1})$ is
strongly negative amphicheiral if 
$b_j = -a_j, 1 \leq j \leq 2n+1$.
Therefore, the Conway polynomials of these knots are
of the form: 
$g(z) g(-z)$ for some $g(z) \in \ZZ[z ]$ 
and generally these knots are neither stable nor 
$c$-stable unless all $a_j$'s have the same sign.
\end{rem}

\section{
Construction of real stable knots (II)}\label{7}

In this section, we show a more general construction 
of stable knots or links.

\subsection{Positive or negative disks}\label{7.1}

Let $D$ be a disk and divide $D$ 
into small domains by a 
(not necessarily connected) plane graph
$G$.
Let $\{v_1, v_2, \cdots, v_n, {v_1}^{\circ},{v_2}^{\circ}, \cdots ,
{v_k}^{\circ}\}$ and $\{e_1, e_2, \cdots , e_{\ell}\}$ be, respectively,
the set of vertices and edges of $G$, where 
${v_j}^{\circ}, 1 \leq j \leq k$, 
is a vertex on $\partial D$, and $e_j, 1 \leq j \leq {\ell}$,
is not a part of $\partial D$ and $e_j$ does not intersect
$\partial D$ except its ends. 
We call $v_j$ an interior vertex and ${v_j}^{\circ}$
a boundary vertex. 
Any part of the boundary of $D$ is not considered as an edge of $G$.

We assume
%
\begin{align}\label{siki:7.1} 
(1)\  &{\rm To\ every\ interior\ vertex\ }
v_j {\rm \  of}\  G, 
{\rm \ there\ is\ a\ path\ in\ }\ 
G {\rm \ that\  connects\ }
v_j {\rm \  to\ some\ }
\nonumber\\
& {\rm boundary\ vertex\ } 
{v_i}^{\circ}.
\nonumber\\
(2)\ &{\rm The\ valency\ of\ }
v_j {\rm \  is\ at\ least\ 2}.
\end{align}
Now each edge $e_j, 1 \leq j \leq {\ell}$, has a weight, $w(e_j) = m_j$, a
non-zero integer, as before.
If all weights are even, we call such a graph an even graph.
If some weights are odd, then we assume that $G$ satisfies condition (\ref{siki:7.2})
below.
Let $d_1, d_2, \cdots, d_m$ be the domains in which $G$ divides $D$ and
$\{e_{j_1}, e_{j_2}, \cdots, e_{j_s}\}$ be the set of all edges on the
boundary of
$d_j, 1 \leq j \leq m$. Then for $1\leq j \leq m$,

\begin{equation}\label{siki:7.2}
\sum_{i=1}^s w(e_{j_i}) \equiv 0 \ ({\rm mod}\ 2)
\end{equation}

If $G$ is an even graph, then (\ref{siki:7.2}) 
is satisfied automatically.
A weighted graph $G$ that satisfies (\ref{siki:7.2}) 
 is called {\it admissible}.
If every edge $G$ has positive (or negative) weight, $G$ is called a positive
(or negative) graph.
Now, first, we replace each interior vertex with a small disk and then, as
we did in the previous section, replace each edge 
$e_j, 1\leq j \leq {\ell}$, with a narrow band of $w(e_j) = m_j $ half-twists.
The resulting surface is denoted by $F(G)$ and is called the surface
representing a weighted graph $G$. The projection of the 1-skeleton of
$F(G) - int(D)$ on $D$ is $G$. 
$F(G)$ is orientable if and only if $G$ is admissible. 
See Fig. 7.1.

\nanaichi

It is well-known that any knot or link is obtained as the boundary of
a Murasugi sum of finitely many surfaces $F(G_j)$ representing admissible
weighted
graphs $G_j$.
In particular, any alternating knot or link is obtained as the boundary of
a Murasugi sum of finitely many surfaces $F(G_j)$ representing admissible
positive and/or negative weighted graphs $G_j$, where glueing is only
allowed between $F(G_p)$ with a positive graph $G_p$ 
and $F(G_q)$ with a negative graph $G_q$.

\subsection{Stable alternating knots and links}\label{7.2}

We prove the following theorem that is a generalization of 
Theorem \ref{thm:6.3}.

\begin{thm}\label{thm:7.1}
Let $\{F(G_1^+), \cdots, F(G_p^+)\}$ and $\{F(G_1^{-}),
\cdots, F(G_q^{-})\}$ be
respectively the set of surfaces representing even
positive and negative graphs $\{G_j^+\}$ and $\{G_j^{-}\}$.
Suppose $K$ is obtained as the boundary of a Murasugi sum
of these surfaces,
where glueing is only
allowed between surfaces of positive graphs 
and those of negative graphs.
Then $K$ is alternating and real stable.
\end{thm}

{\it Proof.} From a construction, 
a Seifert matrix $M$ of $K$ is of the form of (\ref{siki:7.3}) below,
where $A$ is the direct sum of positive definite symmetric matrices 
and $C$ is the direct sum of negative definite symmetric matrices,
and $B$ is obtained from the information of the gluing 
process of a Murasugi sum. 

\begin{equation}\label{siki:7.3}
M=\left[\begin{array}{cc}
A&O\\
B&C
\end{array}
\right]
\end{equation}

Then by Proposition \ref{prop:4.9}, 
$\Delta_K (t)/ \det(M)$ is the characteristic polynomial of a symmetric
matrix, and hence $\Delta_K (t)$ is stable. 
This proves Theorem \ref{thm:7.1}.   
\qed

\subsection{Pseudo-positive or pseudo-negative disk}\label{7.3}

Let $G$ be a weighted even plane graph on a disk $D$ that divides 
$D$ as before.
Suppose $G$ is neither positive nor negative.
Let $F(G)$ be the surface obtained from $G$. 
We call $G$ a {\it pseudo-positive} (or
{\it pseudo-negative}) if the Seifert matrix obtained from $F(G)$
is positive definite (or negative definite). 
In general, $F(G)$ is not alternating. 
However, we have the following theorem.

\begin{thm}\label{thm:7.2}
Let $\{G_i^+\}, 1 \leq i \leq p$ and 
$\{G_j^-\}, 1 \leq j \leq n$ be
the sets of even pseudo-positive and pseudo-negative graphs,
respectively. 
Let $\{F(G_1^+), \cdots, F(G_p^+)\}$ and
$\{F(G_1^{-}), \cdots,
F(G_n^{-})\}$ be respectively the sets of surfaces representing
$\{G_i^+\}$ and $\{G_j^{-}\}$. 
Suppose $K$ is obtained as the boundary of a
Murasugi sum of these surfaces,
where glueing is only
allowed between surfaces of pseudo-positive graphs 
and those of pseudo-negative graphs.
Then $K$ is real stable.
\end{thm}

Since a proof is essentially the same as that of Theorem \ref{thm:7.1}. 
we omit the details. We note that $K$ is not necessarily alternating.

\begin{rem}\label{rem:7.3}
If an even graph $G$ is neither positive nor negative, 
then a Seifert matrix $M$
of $F(G)$ may not be positive definite or negative definite.
However, we may change at most $s$ weights of $G$ so that $F(G)$ 
becomes
pseudo-positive or pseudo-negative, where $s$ is the first 
Betti number of
$H_1(F(G);\ZZ)$. 
For details, see Section 12.
\end{rem}

\subsection{Example}\label{7.4}

Let $G_1$ and $G_2$ be even weighted graphs on disks as in Fig.7.2.
Let $F(G)=F(G_1) \ast F(G_2)$, a Murasugi sum of $F(G_1)$ and
$F(G_2)$. 
Then $K= \partial F(G)$ is not stable, but bi-stable, since
$\Delta_K (t) = t^6 -4t^4 + 7 t^3 - 4t^2 +1$ 
has two real and four unit complex zeros.
Now change $G_1$ to $G_1^{\prime}$
by giving three new even weights to $G_1$ as shown in Fig 7.2. 
Then
$G_1^{\prime}$ is pseudo-positive and 
$K^{\prime} =\partial{F(G_1^{\prime})\ast F(G_2)}$ is stable. 
In fact, 
$\Delta_{K^{\prime}}(t)
= t^6 -15t^5 + 60t^4-93t^3 +60t^2 -15t +1$ is stable.

\nanani

\section{Exceptional stable knots and links}\label{8}

If the sequence of a continued fraction expansion of $r$ alternates in sign, 
then $K(r)$ is stable. (Theorem \ref{thm:6.1} or {\cite[Theorem 2]{LM}}). 
However, the converse is not necessarily true. 
There are many stable $2$-bridge knots 
with non-alternating sequences. 
For example, if $r = [2,-2,-8,2]$,
$K(r)$ is stable, since $\Delta_{K(r)}(t) = (2t^2 -5t+2)^2$. 
Further, for $2$-bridge links, it is possible to construct systematically
such exceptional stable links. 
In this section, we show some of these
knots and links.

\subsection{Exceptional stable knots}\label{8.1}

We begin with the following proposition.

\begin{prop}\label{prop:8.1}
Let $r=[2a,-2,-2b,2c]$, where $a,b,c > 0$. 
Then $K(r)$ is stable if $bc \geq 2a(c+1)$.
\end{prop}

{\it Proof.} 
First we see that $\Delta_{K(r)}(t)= At^4-Bt^3+Ct^2-Bt+A$, where
$A=abc$, $B=4abc+bc-ac+a$ and $C=6abc+2bc-2ac+2a+1$. 
Consider the modified polynomial $f(x)$ of $\Delta_{K(r)}(t)$, where
$f(x)= Ax^2-Bx+(C-2A)$.
Since the discriminant $d$ of $f(x)$ is $bc(bc-2a(c+1))+a^2(c-1)^2$, it follows
that $d \geq 0$ if $bc \geq 2a(c+1)$. Let $\alpha$ and $\beta$ be
two real zeros of $f(x)$. 
We claim that $|\alpha|, |\beta| > 2$. 
In fact,
$\alpha$ and $\beta$ are given as $\frac{B \pm \sqrt{d}}{2abc}
=2 + {\frac{bc-a(c-1) \pm \sqrt{d}}{2abc}}$. 
But since $bc \geq 2a(c+1)$,
we see that $bc-a(c-1) \geq 0$ and also $(bc-a(c-1))^2 > d$,
and hence $\alpha, \beta > 2$.
\qed

\begin{ex}\label{ex:8.2}
If $a=1, b=4$ and $c=1$, then $K(r)$ is stable.
\end{ex}

A similar argument can be applied to show the following proposition.

\begin{prop}\label{prop:8.3}
Let $r=[2a,2b,-2b,-2a], a,b > 0$. 
Then $K(r)$ is stable if
and only if $a \geq 4b$.
Here, $f(x)=
1 + 2 a^2 - 4 a b + 4 a^2 b^2 + (-a^2 + 2 a b - 4 a^2 b^2) x + a^2 b^2 x^2$.
\end{prop}

\begin{ex}\label{ex:8.4}
(1) For $r=[8,2,-2,-8]$, 
$\Delta_{K(r)}(t) = (4t^2-9t+4)^2$ is stable.\\
(2) Let $r=[10,2,-2,-10]$,
$\Delta_{K(r)}(t) =
25t^4-115t^3+181t^2-115t+25$ is stable.
\end{ex}

\subsection{Exceptional stable links}\label{8.2}

In this sub-section, we study exceptional stable $2$-bridge links.

\begin{prop}\label{prop:8.5}
Let $r=[2a,2b,-2c]$, $a,b,c > 0$. 
Then $K(r)$ is stable if
and only if $a \geq c$.
\end{prop}

{\it Proof.} 
We see that $\Delta_{K(r)} = (t-1)(At^2-Bt+A)$, where $A=abc$ and
$B=2abc+a-c$. 
Therefore $g(t) =At^2-Bt+A$ has two real zeros
if the discriminant $d = (a-c)(4abc+a-c) \geq 0$. 
Since $4abc+a-c > 0$, the proposition follows easily.
\qed

Now we construct exceptional stable links systematically.

\begin{dfn}\label{dfn:8x1} 
Let $r=[2a_1, 2a_2,\dots, 2a_n]$ be a sequence of non-zero integers.\\
(1) $N(r)$ denotes $t M(r) - M(r)^T$,
where $M(r)$ is the $n \times n$ matrix $\{m_{i,j}\}$
such that for all $k$'s,
$m_{k,k}=a_k$, $m_{k,k+1}=1$, 
and other entries are $0$. See (\ref{siki:4.1}).\\
(2)
we write
$-r=[-2a_1, -2a_2,\dots, -2a_n]$, and
$r^{-1}=[2a_n, 2a_{n-1},\dots, 2a_1]$.
\end{dfn}

\begin{lemm}\label{lem:8x2} 
For a given sequence $r=[2a_1, 2a_2,\dots, 2a_n]$,
we have;\\
$\det N(r)=\det N(r^{-1})=(-1)^n \det N(-r)=(-1)^n \det N(-r^{-1})$.
\end{lemm}

{\it Proof.}
This lemma can be simply proven by induction on $n$,
but understood better in terms of Alexander polynomials
of 2-bridge knots and links.
Actually,
$\det N(r)$ coincides with 
the Alexander polynomial $\Delta_{K(r)}(t)$.
And hence $\det N(r), \det N(r^{-1}),
\det N(-r)$ and $\det N(-r^{-1})$ are
equivalent up to multiplication of $\pm 1$ and powers of $t$.
In this case, the differences are detected by their constant terms,
which are respectively 
$(-1)^n \prod_{i=1}^n a_i,(-1)^n \prod_{i=1}^n a_i, (-1)^n\prod_{i=1}^n (-a_i)$,
and $(-1)^n \prod_{i=1}^n (-a_i)$.
\qed

For given two sequences $r=[2a_1,2a_2, \cdots, 2a_n]$ and 
$s=[2b_1,2b_2, \cdots, 2b_n]$, and an integer $k$,
let $[r,2k,s]$ denote
$[2a_1,2a_2, \cdots, 2a_n ,2k,2b_1,2b_2,
\cdots, 2b_n]$.

\begin{thm}\label{thm:8x3}  
Given $r = [2a_1,2a_2, \cdots, 2a_n]$, 
let $T_1 =[r,2k,r]$, 
$T_2=[r,2k,-r]$, 
$T_3=[r,2k,r^{-1}]$ and
$T_4=[r,2k,-r^{-1}]$. 
Then for any integer $k \ne 0$, we have;\\
(1)	$\Delta_{K(T_i)}(t)= \Delta_{K(r)}(t) f(t)$,
where f(t) is an integer polynomial. ($1\le i\le 4$)\\
(2) $\Delta_{K(T_4)}(t) = k(t-1) [\Delta_{K(r)}(t)]^2$.
\end{thm}

\begin{cor}\label{cor:8x4} 
$\Delta_{K(T_4)}(t)$ is stable if and only if
$\Delta_{K(r)}(t)$ is stable.
\end{cor}

For other cases, $T_1,T_2$ and $T_3$, generally $\Delta_{K(T_i)}(t)$ is
not stable unless the sequence of $T_j$ alternates in sign. 
For example,
if $m$ is odd and $r$ is an alternating sequence 
with $2a_1 > 0$ and $k < 0$,
then $K(T_1)$ and $K(T_3)$ are stable.

\begin{ex}\label{ex:8x6} 
(1) Let $s=[4,2,-2]$. 
By Proposition \ref{prop:8.5}, 
$K(s)$ is stable.
Thus, for $r=[4,2,-2,-4,2,-2,-4]$, $K(r)$ is stable by Theorem \ref{thm:8x3},
and further, for $r^{\prime} = [r,2k, -r^{-1}]$, $K(r^{\prime})$ is also
stable, if $k \ne 0$.

(2) Let $s=[2,-2,-8,2]$. Then $K(s)$ is stable by Proposition 
\ref{prop:8.1}.
Therefore,
for $r=[2,-2,-8,2,4,-2,8,2,-2]$, $K(r)$ is stable.
\end{ex}

Theorem \ref{thm:8x3} is a corollary of the following lemma.
For a given matrix $M$, $M_{ij}$ denotes the matrix obtained from $M$ by
deleting the $i^{\rm th}$ row and $j^{\rm th}$ column.

\begin{lemm}\label{lem:8x5} 
For sequences $r=[2a_1,2a_2, \cdots, 2a_n]$ and
$s=[2b_1,2b_2, \cdots, 2b_n]$,
let $A$ (resp. $B$)
denote $t M(r) -M(r)^{T}$ 
(resp. $t M(s) -M(s)^{T}$), where
$M(r)$ and $M(s)$ are as in Definition \ref{dfn:8x1}.
Then we have:
$\Delta_{K([r, 2k, s])}(t)
=\det N([r, 2k, s])=
k (t-1)\det A\det B+ 
t \left( \det A_{n,n}\det B+ 
\det A \det B_{1,1}
\right)
$.
\end{lemm}

{\it Proof of Theorem \ref{thm:8x3}.}
(1) If $s$ is equal to 
$r, r^{-1},-r$
 or 
$-r^{-1}$, then
$\det A=\varepsilon \det B$,
where $\varepsilon$ equals $1$ or $-1$ according to 
Lemma \ref{lem:8x2}.
Since $\Delta_{K(r)}(t)=\det A$, we
have the conclusion by Lemma \ref{lem:8x5}. 
(2) If $r=
-r^{-1}$, then 
(i) $\det A\det B
=(-1)^n \det A \det A$,
(ii)
$\det A_{n,n} \det B
=(-1)^n \det A_{n,n}\det A$,
and 
(iii)
$\det A \det B_{1,1}
=\det A (-1)^{n-1}\det A_{n,n}$. 
Therefore, by Lemma \ref{lem:8x5}, we have the conclusion.
\qed

The following formula is often used in this paper.
A proof is an exercise.
\begin{prop}\label{prop:LA}
Let $A$ and $B$ be square matrices of sizes $n$ and $m$.
Let $M$ be the matrix obtained from  
$A\oplus B$ by changing the $(\alpha,n+\beta)$-entry to $x$ and
$(n+\gamma,\delta)$-entry to $y$,
where $1\le \alpha, \delta\le n$ and $1\le \beta,\gamma\le m$.
Then we have:
\begin{equation}\label{siki:formula}
\det M=\det A\det B-(-1)^{\alpha+\beta+\gamma+\delta} xy
\det A_{\alpha,\delta}\det B_{\gamma,\beta}
\end{equation}
\end{prop}

{\it Proof of lemma \ref{lem:8x5}.}
By (\ref{siki:formula}), we have the following:\\
\begin{align*}
\det N([r,2k,s])
=&
\det \left[\begin{array}{cr|c|cc}
& & & &\\
\Hsymb{A} & & & & \\
&  & t & & \\\hline
& -1& k(t-1) & t & \\
\hline
& & -1 & &\\
& & & &\\
& & &  &\Hsymb{B}
 \end{array}
 \right]\\
=&
\det A \det 
\left[\begin{array}{c|cc}
k(t-1) & t & \\
\hline
-1 & &\\
 & &\\
&  &\Hsymb{B}
 \end{array}
 \right]+t \det A_{n,n}\det B\\
=& \det A \left\{(k(t-1) \det B +t  
\det B_{1,1}\right\}+t \det A_{n,n}\det B\\
=&k(t-1)\det A\det B+t(\det A_{n,n}\det B
+ \det A \det B_{1,1}) 
\end{align*}
\qed

\begin{qu}\label{qu:8.13}
To what extent does the stability property of the Alexander
polynomials of an alternating knot $K$ reflect the topological
properties of $K$?
\end{qu}

\begin{probl}\label{probl:8.14}
Characterize stable alternating knots and links.
\end{probl}

\section{Interlacing property (I) 2-bridge knots}\label{9}
For a series of stable real polynomials,
the interlacing property of two sets of zeros is an interesting and
important property.
First, in this section, we prove a simple, but useful basic theorem in
this paper (Theorem \ref{thm:9.4}).
We begin with a definition.

\begin{dfn}[{\cite[p.310]{branden}}]\label{dfn:9.1}
Let $f, g \in \RR[z]$ be univariate polynomials. 
Suppose $f, g$ are real stable.  
Let $\alpha_1\le \alpha_2\le \dots\le \alpha_n$ and   
$\beta_1\le\beta_2\le\dots\le \beta_m$
be the zeros of $f$ and $g$, respectively.  
Then we say that the zeros 
$\{\alpha_j\}$  and $\{\beta_k\}$ are {\it interlaced},
(or we simply say that $f$ and $g$ are {\it interlaced})
if the following conditions are satisfied.\\
(i)	$|m - n| \le 1$,\\
(ii) they can be ordered so that 

(a) if $n= m$, then  
\begin{center}
$\displaystyle{\ \ \alpha_1\le\beta_1\le\alpha_2\le\beta_2\le\dots\le\alpha_n\le\beta_n, {\mathrm or}}$\\
$\displaystyle{\beta_1\le\alpha_1\le\beta_2\le\alpha_2\le\dots\le\beta_n\le\alpha_n}$,
\end{center}

(b) if $n = m+1$, 
then \\
\centerline{
$\alpha_1\le\beta_1\le\alpha_2\le\beta_2\le\dots\le\alpha_m\le\beta_m\le\alpha_{m+1}(=\alpha_n)$,}

(c) if $m = n+1$, then \\
\centerline{
$\beta_1\le\alpha_1\le\beta_2\le\alpha_2\le\dots\le\alpha_n\le\beta_{n+1}(=\beta_m)$.}
\end{dfn}

\begin{dfn}[{\cite[p.56]{wag}}]\label{dfn:9.2}
For $f, g \in \CC[z]$, 
we define the {\it Wronskian} $W[f, g]$ as 
$W[f, g] = f'g - f g'$.  
For $f (\neq 0), g (\neq 0)\in \RR[z]$, 
we say that two real stable $f, g$ are 
in {\it proper position} (denoted by $f\ll g$) if $W[f, g] \le 0$ on all real values.
\end{dfn}

If the zeros of $f$ and $g$ are interlaced, 
then either $W[f,g] \le 0$ or $W[f, g] \ge 0$ 
on all real values and hence $f \ll g$ or $g \ll f$
\cite[p.56]{wag}.

\begin{thm}[{\cite[p.57]{wag}}]\label{thm:9.3}
(Hermite-Kakeya-Obreschkoff  Theorem)
Let $f, g\in    \RR[z]$. 
Then all non-zero polynomials in 
$\{af + bg | a,b \in\RR\}$ are real-rooted if and only if 
(1) $f , g$ are real stable and 
(2) $f \ll g,  g \ll f$ or $f = g = 0$.
\end{thm}

Now first we prove the following theorem. 
Although 
Theorem \ref{thm:9.3} is a  
strong tool, our proof is not a simple application of 
Theorem \ref{thm:9.3}.  

\begin{thm}\label{thm:9.4}
Let
$s=[2a_1, -2a_2, \cdots, (-1)^{k-1}2a_k, \cdots, (-1)^{n-1} 2a_n]$,
$s^{\prime}=[2a_1, -2a_2, \cdots, (-1)^{n-2}2a_{n-1}]$, 
where $a_j > 0$, for $1 \leq j \leq n$. Then
$\Delta_{K(s)} (t) \Delta_{K(s^{\prime})}(t)$ is simple, and
$\Delta_{K(s)}(t)$ and $\Delta_{K(s^{\prime})}(t)$ are interlaced.
\end{thm}

Note that by reversing the sequences,
the same conclusion of Theorem \ref{thm:9.4} holds for the case 
$s'=[ -2a_2, 2a_3, \cdots, (-1)^{n-1}2a_n]$.

{\it Proof of Theorem \ref{thm:9.4}.}
Our proof is by induction. We use Seifert matrices of twisted 
chain type (\ref{siki:4.1}).

Case 1. $n=2$. 
$\Delta_{K(s')}(t)=a_1(t-1)$ and
$\Delta_{K(s)} (t) 
= -a_{1}a_{2} + t + 2a_{1}a_{2}t -a_{1}a_{2}t^2$.
Since $\Delta_{K(s)}(1)=1$,
$\Delta_{K(s)}(t)\Delta_{K(s')}(t)$ is simple and 
$\Delta_{K(s)}(t)$ has two real zeros, 
$\alpha_1$ and
$\alpha_2$
with $\alpha_1 < 1 < \alpha_2$. 

Case 2. $n=3$. 
Write 
$s=[2a,-2b,2c], s'=[2a, -2b], s''=[2a], a,b,c > 0$. 
By a Seifert matrix 
$M = \left[\begin{array}{ccc}
a&1&0\\
0&-b&1\\
0&0&c
\end{array}\right]$
for $K(s)$, we have
\begin{align*}
\Delta_{K(s)}(t)=
&c(t-1) \Delta_{K(s')}(t)+ t \Delta_{K(S'')}(t)\\
=&(t-1)\left\{c \left( -ab+(1+2ab)t-a b t^2 \right) +at\right\}.
\end{align*}
Consider two curves
$y_1= c\left(-ab+(1+2ab)t-abt^2 \right)$ and $y_2 = -at$.
From the observation in case $n=2$,
we see that these two curves intersect in two points at
say, $t=\beta_1$ and $t=\beta_2$
such that $\beta_1<\alpha_1<1<\alpha_2<\beta_2$.
See Fig.9.1.
Since the zeros of $\Delta_{K(s)}(t)$ are 
$\beta_1,1$ and $\beta_2$,
we have the conclusion for $n=3$.

\kyuichi

Note that for alternating stable knots and links,
all the zeros of the Alexander polynomials are positive, because
the coefficients are non-zero and have alternating signs. 
Moreover, since the Alexander polynomials are reciprocal,
each zero less than $1$ has its counterpart greater than one,
and vice versa.

Case 3. General case.\\
Assume inductively that 
$\Delta_{K(s')}(t)\Delta_{K(s'')}(t)$ is simple
and $\Delta_{K(s')}(t)$ and $\Delta_{K(s'')}(t)$ 
are interlaced, where
for simplicity, we write the sequences as follows:
\begin{align*}
s&=[2a_1,-2a_2,\dots, (-1)^{n-1}2a_{n}],\\
s'&=[2a_1,-2a_2,\dots, (-1)^{n-2}2a_{n-1}],\\
s''&=[2a_1,-2a_2,\dots, (-1)^{n-3}2a_{n-2}].
\end{align*}
Use Seifert matrices $M_s,M_{s'},M_{s''}$
of twisted chain type and call 
\begin{align*}
f_n(t)&:=\det(t M_s -M_s^{T}),\\
f_{n-1}(t)&:=\det(t M_{s'}- M_{s'}^{T}),\\
f_{n-2}(t)&:=\det(t M_{s''}-M_{s''}^{T}).
\end{align*}
Expand $\det(t M_s -M_s^{T})$ at the last row and column, 
and we have:
\begin{equation}
f_n(t)=(-1)^{n-1} a_{n}(t-1) f_{n-1} + t f_{n-2}
\end{equation}
Consider two curves 
$y_1=(-1)^{n-1} a_{n}(t-1) f_{n-1}$
and 
$y_2=-t f_{n-2}$. 
We show that $y_1$ and $y_2$ intersect in $n$ points 
and 
(i) $f_n f_{n-1}$ is simple and 
(ii) $f_n(t)$ and $f_{n-1}(t)$ are interlaced. 
First, note that
the leading coefficients of
$y_1$ and $y_2$ have the same sign.

Case 3.1. $n$ is even, say $2m$.
Fig. 9.2 depicts the case $m$ is even (in, particular $m=4$).
The case $m$ is odd is similar, and the arguments are the same.
By induction hypothesis $f_{n-1}f_{n-2}$ is simple, and
since $K(s')$ is a link and $f_{n-1}$ is simple, 
$(t-1)$ divides $f_{n-1}$ exactly once.
Hence $t=1$ is the double zero of $y_1$.
Meanwhile, since $K(s'')$ is a knot, $(t-1)$ does not divide $f_{n-2}$.
Also, by induction hypothesis,
$f_{n-1}$ and $f_{n-2}$ have interlaced zeros.
If $m$ is even (resp. odd), 
then the leading coefficients of $y_1$ and $y_2$
are both  positive (resp. negative).
The sets of the zeros of $f_{n-1}$ and $y_1$ are the same and
the set of the zeros of $y_2$ coincides with that of $f_{n-2}$
with $0$ added.
Then two curves $y_1$ and $y_2$ intersect exactly $n$ times
and (i) $f_n f_{n-1}$ is simple and (ii)
$f_n$ and $f_{n-1}$ are interlaced, as shown 
in Fig. 9.2.
Since $\deg f_n=n$. We have the conclusion.

\kyuni

\kyusan
Case 3.2. $n$ is odd, say $2m+1$.
Fig. 9.3 depicts the case $m$ is even (in, particular $m=4$).
The case $m$ is odd is similar, and the arguments are the same.
$K(s')$ is a knot, and hence $f_{n-1}(1)\neq 0$.
By induction hypothesis, (i) $f_{n-1}f_{n-2}$ is simple and
(ii) $f_{n-1}$ and $f_{n-2}$ have interlaced zeros.
If $m$ is even (resp. odd), 
then the leading coefficients of $y_1$ and $y_2$
are both  positive (resp. negative).
$y_1$ and $y_2$ have interlacing zeros, except sharing $1$ in common.
Then two curves $y_1$ and $y_2$ intersect exactly $n$ times
and (i) $f_{n}f_{n-1}$ is simple and (ii) 
$f_n$ and $f_{n-1}$ are interlaced.
\qed

In proving Theorem \ref{thm:9.4},
we have the following theorem:

\begin{thm}\label{thm:9.5new}
Let
$s=[2a_1, -2a_2, \cdots, (-1)^{k-1}2a_k, \cdots, (-1)^{n-1} 2a_n]$
and
$r=[2a_1, -2a_2, \cdots, (-1)^{k-1}2a_k, \cdots, (-1)^{n-1} 2(a_n-1)]$,
(or $r=[2(a_1-1), -2a_2, \cdots, 
(-1)^{k-1}2a_k, \cdots, (-1)^{n-1} 2a_n], a_1>1$),
where $a_j > 0$, for $1 \leq j \leq n$ and $a_n>1$. 
Let $\{\alpha_j,1\le j\le [\frac{n}{2}]\}$ and 
$\{\beta_j,1\le j\le [\frac{n}{2}]\}$ be the zeros of
$\Delta_{K(r)}(t)$ and $\Delta_{K(s)}(t)$ respectively
in $(0,1)$. 
Then $\{\alpha_j\}$ and $\{\beta_j\}$ are disjoint and are 
interlaced, i.e.,
$0<\alpha_1<\beta_1<\dots<
\alpha_{[\frac{n}{2}]}<\beta_{[\frac{n}{2}]}<1$.
\end{thm}

{\it Proof.}
As in the proof of Theorem \ref{thm:9.4},
$f_n=(-1)^{n-1} a_n (t-1)f_{n-1}+t f_{n-2}$.
If we replace $a_n$ by $a_n -1$, the curve $y_1$ is squeezed toward
the $t$-axes, while the intersection points with the $t$-axes are
fixed.
Since $y_2$ is irrelevant to $a_n$,  each of the  zeros of $f_n$ 
less than (resp. more than) $1$ is moved toward 
(but never beyond) its left (resp. right) neighbour.
Therefore, we have the conclusion.
\qed

Theorem \ref{thm:9.5new} is
generalized as follows:

\begin{thm}\label{thm:9.a1}
Let $s=[2a_1,-2a_2,\dots, (-1)^{k-1}2a_k,
\dots, (-1)^{n-1} 2a_n]$
and $r=[2a_1,-2a_2,\dots, (-1)^{k-1}2(a_{k}-1),\dots,
(-1)^{n-1}2a_n]$, where $a_j>0,1\le j\le n$
and $a_k>1,1\le k\le n$.
Let $\delta(K)$ be the maximal value of the
zeros of  $\Delta_{K}(t)$. 
Then $\delta(K(s))<\delta(K(r))$.
\end{thm}

{\it Proof.} 
We prove the theorem for the case where
$n$ is even, say, $n=2m$ and $k$ is odd.
The same argument works for the other cases.

Let $s_1=[2a_1,-2a_2,\dots, (-1)^{k-2}2a_{k-1}]$,
$s_2=[(-1)^{k-1}2a_k,\dots,(-1)^{n-1}2a_n]$,
$r_2=[(-1)^{k-1}2(a_{k}-1), \dots, (-1)^{n-1}2a_n]$,
$s_1'=[2a_1,-2a_2,\dots, (-1)^{k-3}2a_{k-2}]$,
$s_2'=[(-1)^k 2a_{k+1},\dots,(-1)^{n-1}2a_n]$.
We use Seifert matrices 
$M(s)$, $M(r)$, $M(s_1)$, $M(s_2)$, $M(r_2)$,
$M(s_1')$ and $M(s_2')$ of twisted chain type.
Let $f_n=\det(t M(s)-M(s)^T)$, 
$\widehat{f_n}=\det(tM(r)-M(r)^T)$,
$f_{k-1}=\det(tM(s_1)-M(s_1)^T)$,
$g_{n-k+1}=\det(tM(s_2)-M(s_2)^T)$,
$h_{n-k+1}=\det(tM(r_2)-M(r_2)^T)$,
$f_{k-2}=\det(tM(s_1')-M(s_1')^T)$,
$g_{n-k}=\det(tM(s_2')-M(s_2')^T)$.
Now by Proposition \ref{prop:LA}, we have
$f_n=f_{k-1}g_{n-k+1}+t f_{k-2} g_{n-k}$
and
$\widehat{f_n}=f_{k-1}h_{n-k+1}+t f_{k-2}g_{n-k}$.
Let $\alpha, \beta$ and $\gamma$ be respectively the 
smallest zeros of
$f_{k-1},g_{n-k+1}$ and $h_{n-k+1}$.
Then $\gamma<\beta$ be Theorem \ref{thm:9.5new}.
Let $y_1=f_{k-1}g_{n-k+1},z_1=f_{k-1}h_{n-k+1}$
and $y_2=-tf_{k-2}g_{n-k}$.
We note that the signs of the leading coefficients
of $y_1$ and $z_1$ are both $(-1)^{m}$ and that of 
$f_{k-2}g_{n-k}$ is $(-1)^{m-1}$.
Also, the smallest positive zero of $y_2$ is 
larger than $\alpha$ or $\beta$ by 
Theorem \ref{thm:9.4}.
Consider the intersection of three curves $y_1,z_1$ and $y_2$.
The smallest value of the intersection will be seen from
the diagrams below.
See Fig 9.4 for the case $\alpha\le\gamma$, and Fig 9.5 for the case
 $\gamma<\alpha$.

\kyuyon

\kyugo

In each ease, $d(K(r))<d(K(s))$, where $d(K)$ denotes the
smallest (positive) zero of $\Delta_{K}(t)$,
and hence $\delta(K(r))>\delta(K(s))$.
\qed

As an immediate consequence of Theorem \ref{thm:9.a1},
we have the following theorem:

\begin{thm}\label{thm:9.a3}
Let $F_n$ be the set of stable 
2-bridge knots or links $K(r)$ with
$r=[2a_1,-2a_2,\dots, (-1)^{n-1}2a_n], a_j>0, 1\le j\le n$.
Then $\delta(K(r))$ is maximal if and only if 
$r=[2,-2,\dots, (-1)^{n-1} 2]$,
i.e., $K(r)$ is fibred.
\end{thm}

\begin{ex}\label{ex:9.9}
(1) Let $s=[4,-2,2,-6,4,-2], s'=[4,-2,2,-6,4]$.
Then the zeros of $\Delta_{K(s)}$ are approximately
$\{0.2866, 0.4550, 0.7654, 1.3065, 2.1976, 3.4888\}$
and those of $\Delta_{K(s')}$ are approximately
$\{0.2877, 0.6179, 1.0000, 1.6183, 3.4761\}$\\
(2) Let
$s=[4,-2,2,-6,4,-2,4,-4], s'=[4,-2,2,-6,4,-2,4]$.\\
Then the zeros of $\Delta_{K(s)}$ are approximately\\
$\{
0.2857, 0.3535, 0.6148, 0.8171, 1.2237, 1.6265, 2.8287, 3.4999\}$
and those of $\Delta_{K(s')}$ are approximately
$\{0.2859, 0.3716, 0.6772, 1., 1.4767, 2.6913, 3.4973\}$\\
\end{ex}

\begin{rem}\label{rem:9.10}
For exceptional stable 2-bridge knots or links, the
interlacing property may not hold. 
For example, let $s = [10,2,-2,-10]$ and
$s'=[10,2,-2]$. 
Then $K(r)$ and $K(s')$ are both stable,
but they are not interlaced. 
In fact, the zeros of $\Delta_{K(s)}(t)$
are approximately $\{0.4923, 0.7592, 1.3172, 2.0313\}$, but those of
$\Delta_{K(s')}(t)$ are approximately $\{0.4202, 1, 2.3797\}$.
Therefore, they are not interlaced.
\end{rem}

\section{Interlacing property (II) Quasi-rational knots ${X_n}$}\label{10}

In the following two sections, we prove the interlacing property for
two series of alternating stable knots ${X_n }$ and ${Y_{2n+1}}$ considered
in Section 6.3.
The idea of our proof is similar to the proof of 
Theorem \ref{thm:9.4},   
but we need
a lot of computations. The first series of knots ${X_n}$
shows us that the zeros of the Alexander polynomials of alternating knots
is unbounded. On the other hand, the Alexander polynomial of each knot
in the second series of knots ${Y_{2n+1}}$ is not irreducible.
Nevertheless, the maximal value of the zeros of a factor (of degree 4) of
the Alexander polynomial of $Y_3$ is quite likely equal to $\delta_6$ 
defined in Section 5.1.

Now consider a series of stable knots

$X_n (a,b) = X_n (2a_1 ,\cdots, 2a_n \mid -2b_1 , \cdots, -2b_n ), a_j 
,b_j > 0, 1 \leq j \leq n$.

We prove

\begin{thm}\label{thm:10.1}
(1) $X_n (a,b)$ is a stable alternating knot of genus $n$.

(2) Let $K_n = X_n (2,2, \cdots, 2|-2,-2, \cdots, -2)$. 
Then for $n \geq 2$, $\Delta_{K_n}(t)\Delta_{K_{n-1}}(t)$ is simple
and $\Delta_{K_n}(t)$ and $(t-1)\Delta_{K_{n-1}}(t)$
are interlaced.

(3) Let $\alpha_n$ be the maximal value of the zeros of
$\Delta_{K_n}(t)$. 
Then $\alpha_n \geq n+1$.
\end{thm}

We suspect that (2) holds for $X_n (a,b)$ with any 
$a_j > 0$ and $b_j>0, 1 \leq j \leq n$.
 Therefore we conjecture:

\begin{yosou}\label{conj:10.2}
Let $a_j > 0$ and $b_j > 0, 1 \leq j \leq n$. 
Then for $n\ge 2$, 
$\Delta_{X_n}(t)\Delta_{X_{n-1}}(t)$ is simple,
and $\Delta_{X_n}(t)$ and $(t-1)\Delta_{X_{n-1}}(t)$ are
interlaced.
\end{yosou}

Now since (1) and (3) are already proved in 
Theorem \ref{thm:6.6} (1),
we prove only
(2) in this section.
For simplicity, we denote by $G(n)$ the normalization of 
$\Delta_{K_n}(t)$.

As the first step, we prove

\begin{prop}\label{prop:10.3}
Let $\lambda (t) = 2t^2 -5t + 2$.

(1) For $n \geq 2, G(n) = \lambda (t) G(n-1) - (t-1)^4 G(n-2)$.

(2) For $n \geq 0, \lambda(t) \nodiv G(n)$.

(3) For $n \geq 0, (t-1) \nodiv G(n)$,
where we define $G(0) = 1$.
\end{prop}

{\it Proof.}  Since $K_n$ is a knot, (3) holds trivially. 
Now Proposition \ref{prop:10.3}  
holds for $n=1$ and $2$. In fact, $K_1 = 4_1$
and $K_2 = 8_{12}$ and 
$\Delta_{K_2}(t)$ = $t^4 -7t^3 + 13t^2 - 7t +1
=(2t^2 -5t+ 2) (t^2 -3t+1) -(t-1)^4 = 
\lambda(t) G(1) -(t-1)^4 G(0)$.
Suppose $n \geq 3$. 
A Seifert matrix $M$ of $K_n$ is given in a proof of
Theorem \ref{thm:6.6} (1). 
It is of the form: $M = \tbt{E_n}{O}{C}{-E_n}$,
where $C$ is a lower triangular matrix defined 
in the proof of Theorem \ref{thm:6.6} 
(1). Since $\Delta_{K_n} (t) = \det [Mt -M^{T} ]$,
we see that $G(n) = \det N$, where $N = 
\tbt{(t-1) E_n}{ C^{T}}{Ct}{ (t-1)E_n}$.
To compute $G(n)$, first expand $\det N$ along the 
$n^{\mathrm th}$ row and
we obtain by Proposition \ref{prop:LA} that
$G(n) = \det N = (t-1) \det N_1 - t G(n-1)$,
where $N_1 = \tbt{(t-1)E_{n-1}}{C_1^{T}}{C_1 t}{(t-1) E_n}$
and $C_1$ is
an $n \times (n-1)$ matrix of the form:
\begin{center}
$C_1=
\left[
\begin{array}{rrrcr}
1&\multicolumn{4}{c}{}\\
1&1&\multicolumn{3}{c}{\largesymbol{O}}\\
1&1&1&\multicolumn{2}{c}{}\\
\multicolumn{3}{c}{}&\ddots& \\
1&1&1&\cdots&1\\
1&1&1&\cdots&1
\end{array}
\right]
$
\end{center}

Next, subtract the next to last row of $\det N_1$ from the last row and
then subtract the next to last column of $\det N_1$
from the last column and then expand it along the last
row, and we obtain
$\det N_1 = 2(t-1)G(n-1) -(t-1)^3 G(n-2)$.
Therefore, $G(n) = (t-1) \det N_1 - t G(n-1)$
= $(t-1)\{2(t-1) G(n-1) -(t-1)^3 G(n-2)\}- t G(n-1)$
=$\{2(t-1)^2 -t\}G(n-1) -(t-1)^4 G(n-2)
= \lambda(t) G(n-1) - (t-1)^4 G(n-2)$.
This proves (1).
Using (1), (2) follows by induction. 
\qed

To prove Theorem \ref{thm:10.1} (2),  
let $\alpha_1 < \alpha_2 < \cdots <
\alpha_{n-1} < 1$ be the zeros of $(t-1) G(n-1)$ in 
$[0,1]$ and
let $\beta_1 < \beta_2 < \cdots < \beta_n$ be
the zeros of $G(n)$ in $[0,1]$. 
Then it suffices to prove the following
proposition.

\begin{prop}\label{prop:10.4}
(1) For $n\ge 1$, $G(n)G(n-1)$ is simple.
(2) $\{\alpha_j, 1 \leq j \leq n-1\}\cup \{1\}$ and $\{\beta_j,1 \leq j
\leq n\}$ are interlaced, namely,

\begin{equation}\label{siki:10.1}
\beta_1 < \alpha_1 < \beta_2 < \alpha_2 < \cdots <
\alpha_{n-1} < \beta_n < 1.
\end{equation}

\noindent
(3) (a) If $n$ is even, say $2m$, 
then $\alpha_m < \frac{1}{2} <
\beta_{m+1}$, and

(b) if $n$ is odd, say $2m+1$, then $\beta_{m+1} <\frac{1}{2} 
< \alpha_{m+1}$.
\end{prop}

We note that $\frac{1}{2}$ (and 2) are the zeros of 
$\lambda(t)$ and also
we need (3) to show (2) by induction.

{\it Proof.} We use induction.

Case $n=1$. The zeros of $G(1) = t^2 -3t +1$ in $[0,1]$ is 
$\beta_1 =
0.38\cdots$. Since $G(0) = 1$, we have $\beta_1 < 1$ and further
$\beta_1 < \frac{1}{2}$. 
This proves Proposition \ref{prop:10.4} 
when $n=1$.

Case $n=2$. 
The zeros of $G(2)$ in $[0,1]$ are $\beta_1 =0.228\cdots$ and
$\beta_2 = 0.5449\cdots$, while the zeros of $(t-1) G(1)$ in 
$[0,1]$ are
$\alpha_1 = 0.38\cdots$ and 1. 
Therefore, $\beta_1 < \alpha_1 <
\beta_2 < 1$ and $\alpha_1 < \frac{1}{2} < \beta_2$.

Case $n \geq 3$. Suppose $G(n)G(n-1)$ is simple.
Let $\{\alpha_j, 1 \leq j \leq n-1\}\cup \{1\}$ be the zeros of
$(t-1)G(n-1)$ in $[0,1]$ and $\{\beta_j ,1 \leq j \leq n\}$ be 
the zeros of
$G(n)$ in $[0,1]$.

Inductively, we assume that they are interlaced, namely,
%
\begin{equation}\label{siki:10.2}
\beta_1 < \alpha_1 < \beta_2 < \alpha_2 < \cdots <
\alpha_{n-1} < \beta_n < 1.
\end{equation}
Consider the zeros of $G(n+1)$ in $[0,1]$.
Since $G(n+1) = \lambda(t) G(n) - (t-1)^4 G(n-1)$, the zeros of 
$G(n+1)$ in $[0,1]$ are determined by the 
intersection of two curves 
$y_1 = \lambda(t) G(n)$ and
$y_2 = (t-1)^4 G(n-1)$. 
Note
that $\lambda(t) G(n)$ is simple and since $(t-1) \nodiv G(n-1)$, 
$1$ is the
zero of $y_2$ of exactly order 4.

Case (a) $n+1$ is even, say $2m$.
Then by induction, since $n$ is odd, we have

\begin{align}\label{siki:10.3}
\beta_1 &< \alpha_1 < \beta_2 < \alpha_2 < \cdots < \beta_m <
\frac{1}{2} < \alpha_m \nonumber\\
&< \beta_{m+1} < \cdots < \alpha_{n-1} 
< \beta_n <1.
\end{align}

Note that $y_1(0)= 2$ and $y_2 (0) =1$, since $G(k)$ is monic
for all $k \geq 1$.

When $m$ is  even (resp. odd),
two curves are depicted in Fig. 10.1 top (resp. bottom).%

\jui

There are exactly $(n+1)$ points of intersection 
$\{\gamma_j, 1 \leq j \leq n+1\}$
in $[0,1]$ which are the zeros of $G(n+1)$ in $[0,1]$:

\begin{align*}
\gamma_1 &< \beta_1 < \gamma_2 < \beta_2 < \gamma_3 < \cdots <
\gamma_m < \beta_m < \frac{1}{2} < \gamma_{m+1} \nonumber\\
&< \beta_{m+1}
< \cdots < \gamma_n < \beta_n < \gamma_{n+1} < 1.
\end{align*}

Thus $G(n+1)G(n)$ is simple and
$\{\gamma_j \}$ and $\{\beta_j \}\cup \{1\}$ are interlaced and
$\beta_m < \frac{1}{2} < \gamma_{m+1}$.
Proposition \ref{prop:10.4} 
is now proved for this case.


Case (b) $n+1$ is odd, say $2m+1$. 
Then by induction, since $n$ is even, we see
\begin{align*}
\beta_1 &< \alpha_1 < \beta_2 < \alpha_2 < \cdots < \alpha_m <
\frac{1}{2} < \beta_{m+1} \nonumber\\
&< \alpha_{m+1} < \cdots < \alpha_{n-1} <
\beta_n < 1.
\end{align*}
When $m$ is  even (resp. odd),
two curves are depicted in Fig. 10.2 top (resp. bottom).%

\juii

There are $(n+1)$ points of intersection 
$\{\gamma_j, 1 \leq j \leq n+1\}$ in $[0,1]$ and
\begin{align*}
\gamma_1 &< \beta_1 < \gamma_2 < \beta_2 < <\cdots < \beta_m <
\gamma_{m+1} < \frac{1}{2} < \beta_{m+1} 
\nonumber\\
&< \gamma_{m+2}
< \cdots < \beta_n < \gamma_{n+1} < 1.
\end{align*}

Now we covered all cases and a proof is completed. 

\qed

\begin{ex}
Given below is the list of the zeros of
the Alexander polynomials for
$G(k), k=1,2,\dots,8$.
The table satisfies Theorem \ref{thm:10.1}.
\end{ex}

$
\begin{array}{c}
 0.382 \\
 2.618
\end{array}
\begin{array}{c}
 0.228 \\
 0.544 \\
 1.838 \\
 4.390
\end{array}
\begin{array}{c}
 0.145 \\
 0.458 \\
 0.578 \\
 1.730 \\
 2.186 \\
 6.904
\end{array}
\begin{array}{c}
 0.098 \\
 0.382 \\
 0.526 \\
 0.591 \\
 1.692 \\
 1.900 \\
 2.618 \\
 10.193
\end{array}
\begin{array}{c}
 0.070 \\
 0.320 \\
 0.474 \\
 0.557 \\
 0.597 \\
 1.674 \\
 1.797 \\
 2.109 \\
 3.129 \\
 14.273
\end{array}
\begin{array}{c}
 0.052 \\
 0.269 \\
 0.426 \\
 0.519 \\
 0.573 \\
 0.601 \\
 1.664 \\
 1.746 \\
 1.927 \\
 2.349 \\
 3.719 \\
 19.155
\end{array}
\begin{array}{c}
 0.040 \\
 0.228 \\
 0.382 \\
 0.481 \\
 0.544 \\
 0.582 \\
 0.603 \\
 1.658 \\
 1.717 \\
 1.838 \\
 2.077 \\
 2.618 \\
 4.390 \\
 24.841
\end{array}
\begin{array}{c}
 0.032 \\
 0.194 \\
 0.343 \\
 0.446 \\
 0.515 \\
 0.560 \\
 0.589 \\
 0.605 \\
 1.654 \\
 1.699 \\
 1.786 \\
 1.943 \\
 2.242 \\
 2.915 \\
 5.144 \\
 31.333
\end{array}
$

\section{Interlacing property (III) Quasi-rational knots $\{Y_{2n+1}\}$}
\label{11}

In this section, we discuss a slightly different sense of interlacing
property of the second series of stable alternating knots.
The Alexander polynomials of knots in this series may 
have multiple zeros and thus, they are not irreducible. Nevertheless, some zero has
the largest value up to 12 crossing knots.

Let $Y_{2n+1}(2a_1, 2a_2, \cdots, 
2a_{2n+1}\mid 2b_1, 2b_2, \cdots, 2b_{2n+1})$ be a
quasi-rational knot obtained as the boundary of the surface 
constructed in Example \ref{ex:6.5}. 
Define a series of alternating quasi-rational knots 
$Y_{n}$ by 
$Y_{2n+1}(
\underbrace{-2, -2, \cdots, -2}_{2n-1} | 
\underbrace{2, 2, \cdots, 2}_{2n-1})$.
Then $Y_1$ is the knot $4_1$ and $Y_2$ is 
the 12 crossing alternating knot $12_{a1124}$. 
We note that the largest zero of $Y_2$ is $7.69853\cdots$, 
which attains the largest real part of all the zeros of Alexander polynomials of 
alternating knots up to 12 crossings.

\subsection{Conway polynomials of $Y_n$}\label{11.1}
In this subsection, we give inductive formulae for 
the Conway polynomials of $Y_n$. 
(For the Conway polynomial, see \cite{Li}.)
Denote by $W_n$ the link obtained from $Y_n$ 
by removing band corresponding the 
the middle vertical edge of the graph in Fig.6.3.

Denote by $c_n(z)$ (resp. $d_n(z)$) 
the Conway polynomial of $Y_n$ (resp. $W_n$).
For simplicity, we write a function $f_n(z)$ of $z$ as $f_n$.

\begin{prop}\label{prop:11.1}
Let $a(z)=1-z^2$, $b(z)=1+z^2$.
Then  we have:
\begin{align}\label{siki:11.1}
&c_1=a, \nonumber\\
&c_2=a(a^2-4z^2), \  and\ for \ n\ge 1, \nonumber\\
&c_{n+2} = a^2 (2 c_{n+1} - b^2 c_{n}).
\end{align}
\end{prop}

{\it Proof.}
In (\ref{siki:11.2}) below, we can write
$d_{n}$ (resp. $d_{n+1}$) in terms of 
$c_n$ and $c_{n+1}$ (resp. $c_{n+1}$ and $c_{n+2}$).
Substituting these into (\ref{siki:11.3}), we have the conclusion.
Using skein trees, we prove Lemma \ref{lem:11.2} 
at the end of this subsection.
\qed

\begin{lemm}\label{lem:11.2}
\begin{align}
& d_1=z,\ and\ for\ n\ge 1, \nonumber\\
& \label{siki:11.2} c_{n+1}=(3z^4-4z^2+1)c_{n}+2z(z^4-1)d_{n}\\
& \label{siki:11.3} d_{n+1}=2z(1-z^2) c_{n}+(1-z^4)d_{n}.
\end{align}
\end{lemm}

We know that the Conway polynomials $c_n$ have a factor $a=1-z^2$.
In the proposition below, we determine the exponent of $a$ in $c_n$.

\begin{prop}\label{prop:11.3}
Let $f_{2m-1}=\dfrac{c_{2m-1}}{a^{2m-1}}, 
f_{2m}=\dfrac{c_{2m}}{a^{2m-1}}$. 
Then for $m\ge1$, we have the following, where
$a(z)=1-z^2$, $b(z)=1+z^2$.
\begin{align}
\label{siki:11.4} & f_{2m-1}, f_{2m} \in \ZZ[z]\\
\label{siki:11.5} & a\not| f_{2m-1}\  
and\  a\not\vert f_{2m} ,\\
\label{siki:11.6} & f_{2m+1}=2f_{2m}-b^2 f_{2m-1}\\
\label{siki:11.7} & f_{2m+2}=2a^2 f_{2m+1}-b ^2 f_{2m}
\end{align}
\end{prop}

{\it Proof.}
First we prove (\ref{siki:11.4}) by induction.
For $c_1$ and $c_2$, the claim is trivial. 
Assume that the claim holds up to $2m$.
By (\ref{siki:11.1}) and induction hypothesis,
 we have 
 \begin{align}
c_{2m+1}&=a^2(2 c_{2m}-b^2 c_{2m-1})\nonumber\\
   &=a^2(2a^{2m-1}f_{2m}-b^2a^{2m-1}f_{2m-1})\nonumber\\
\label{siki:11.8}    &=a^{2m+1}(2 f_{2m}-b^2 f_{2m-1}),
\end{align}
and hence $c_{2m+1}$ has factor $a^{2m+1}$.
We also have 
 \begin{align}
 c_{2m+2}&=a^2(2 c_{2m+1}-b^2c_{2m})\nonumber\\
       &=a^2(2 a^{2m+1} f_{2m+1}-b^2 a^{2m-1} f_{2m})\nonumber\\ 
 \label{siki:11.9}      &=a^{2m+1}(2a^2 f_{2m+1}-b^2 f_{2m}),
       \end{align}
and hence $c_{2m+2}$ has factor $a^{2m+1}$.
Therefore, we have (\ref{siki:11.4}).
By (\ref{siki:11.8}) and (\ref{siki:11.9}),
we have
(\ref{siki:11.6}) and (\ref{siki:11.7}).
Finally, we prove (\ref{siki:11.5}), using the fact
that if a polynomial $f(t)$ is divided by 
$a=1-z^2$ then $f(1)=0$.
Let $e_n=f_n(1)$, and we prove that $e_k\neq 0$ for $k\ge 1$. 
By putting $z=1$ in (\ref{siki:11.6}) and (\ref{siki:11.7}),
we have
\begin{align}
\label{siki:11.10} e_{2m+1}&=2e_{2m}-4e_{2m-1}\\
\label{siki:11.11} e_{2m+2}&=-4e_{2m}
\end{align} 
Since $e_1=1, e_2=-4$, by (\ref{siki:11.11}) we have
$e_{2m}=(-4)^{m}$ and hence $e_{2m}\neq 0$.
Then by (\ref{siki:11.10}), $e_{2m+1}=2(-4)^{m}-4e_{2m-1}$
and hence $e_{2m+1}$ is positive (resp. negative) 
if $m$ is even (resp. odd), and in any case, $e_{2m+1}\neq 0$.
In fact, we see that $e_{2m-1}=(-4)^{m-1}(2m-1)$.
\qed

The following lemma is used to prove Lemma \ref{lem:11.2}.
The equation (\ref{siki:11.12}) means the relation of the
Conway polynomials of three knots or links that differ 
only locally as depicted in the diagrams.
In Lemma \ref{lem:11.4} and its proof, 
we adopt such a convention.
Lemma \ref{lem:11.4} below reduces the Conway polynomial of
a diagram with two or more parallel bands (of positive writhe)
to those with one and zero bands.

\begin{lemm}\label{lem:11.4}
\begin{align}
\label{siki:11.12}
\incskein{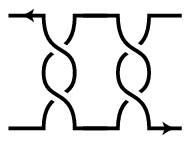}&=2z\incskein{1}-z^2\incskein{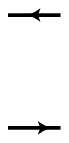}\\
\label{siki:11.13}
\incskein{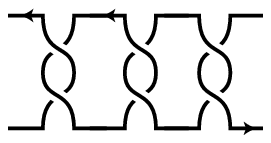}&=3z^2\incskein{1}-2z^3\incskein{none}
\end{align}
\end{lemm}

{\it Proof.}
First, we have $\incskein{1}=\incskein{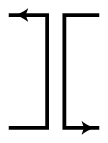}+z\incskein{none}$
and hence $\incskein{0}=\incskein{1}-z\incskein{none}$. \\

\noindent
 $\incskein{11}=\incskein{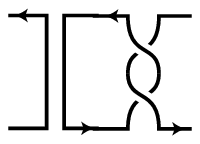}+z\incskein{1}
=\incskein{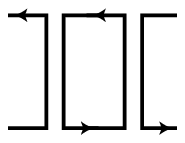}+z\incskein{0}+z\incskein{1}
=z\Bigl(\incskein{1}-z\incskein{none}\Bigr)+z\incskein{1}$.\\

\noindent
$\incskein{111}=2z\incskein{01}+z^2\incskein{1}
=2z\Bigl(2z \incskein{1}-z^2\incskein{none}\Bigr)
-z^2\incskein{1}$.\\
In fact, we can inductively prove 
$\underbrace{\incskein{1}\cdots\incskein{1}}_{n}=
nz^{n-1}\incskein{1}-(n-1)z^n\incskein{none}$.
\qed

To prove Lemma \ref{lem:11.2} using
skein trees, notice the following (see Figure 11.1).
Suppose that a band $b_1$ crosses over 
exactly one other band $b_2$, then untwisting $b_1$
results in cutting both $b_1$ and $b_2$.
Suppose a band $b_1$ crosses over exactly two bands
$b_2$ and $b_3$. Then untwisting $b_1$ results in
removing $b_1$ and merging the band $b_2$ and $b_3$.

\medskip
\begin{minipage}{1\hsize}
\centerline{
(1)\ \raisebox{-13pt}{\inclf{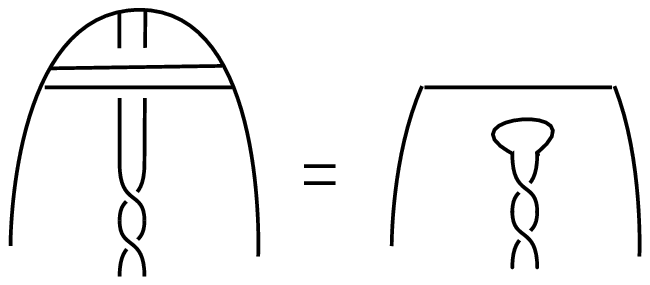}{2.5} } 
\hspace{0.8cm}
(2)\ \raisebox{-13pt}{\inclf{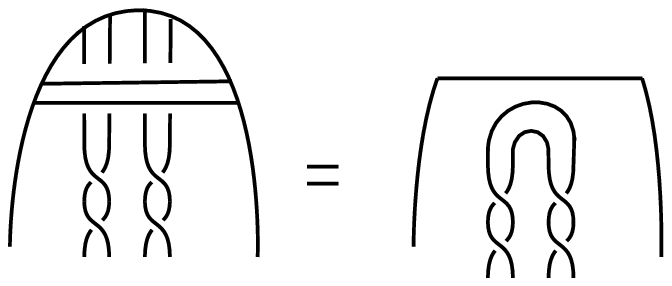}{2.5}}
}
\medskip
\centerline{Fig. 11.1}
\vspace*{2mm}
\end{minipage}

{\it Proof of Lemma \ref{lem:11.2}.}
We depict skein trees for $Y_n$'s and $W_n$'s.
The knots $Y_n$'s and links $W_n$'s are represented
by diagrams in Figures 11.2 through 11.5.
Each dotted arc indicates the site where an arc is removed,
and hence dotted arcs are not counted as arcs.
Horizontal arcs are assigned with weight $-2$ and
non-horizontal ones with weight $2$, except for the ones
with label $4$ in Fig. 11.4. 
Recall that a diagram represents a knot or link
on the boundary of a Seifert surface obtained
from a disk by attaching twisted bands along the arcs,
where  horizontal arcs contribute bands on the top side,
and non-horizontal arcs on the back side.

First, we prove (\ref{siki:11.3}). We have a skein tree as
in Fig. 11.2, where the sites of crossing changes and splicing 
are marked with $*$.
In the first step, we use Fig. 11.1 (1).

\begin{minipage}{1\hsize}
\vspace*{3mm}
\centerline{\inclf{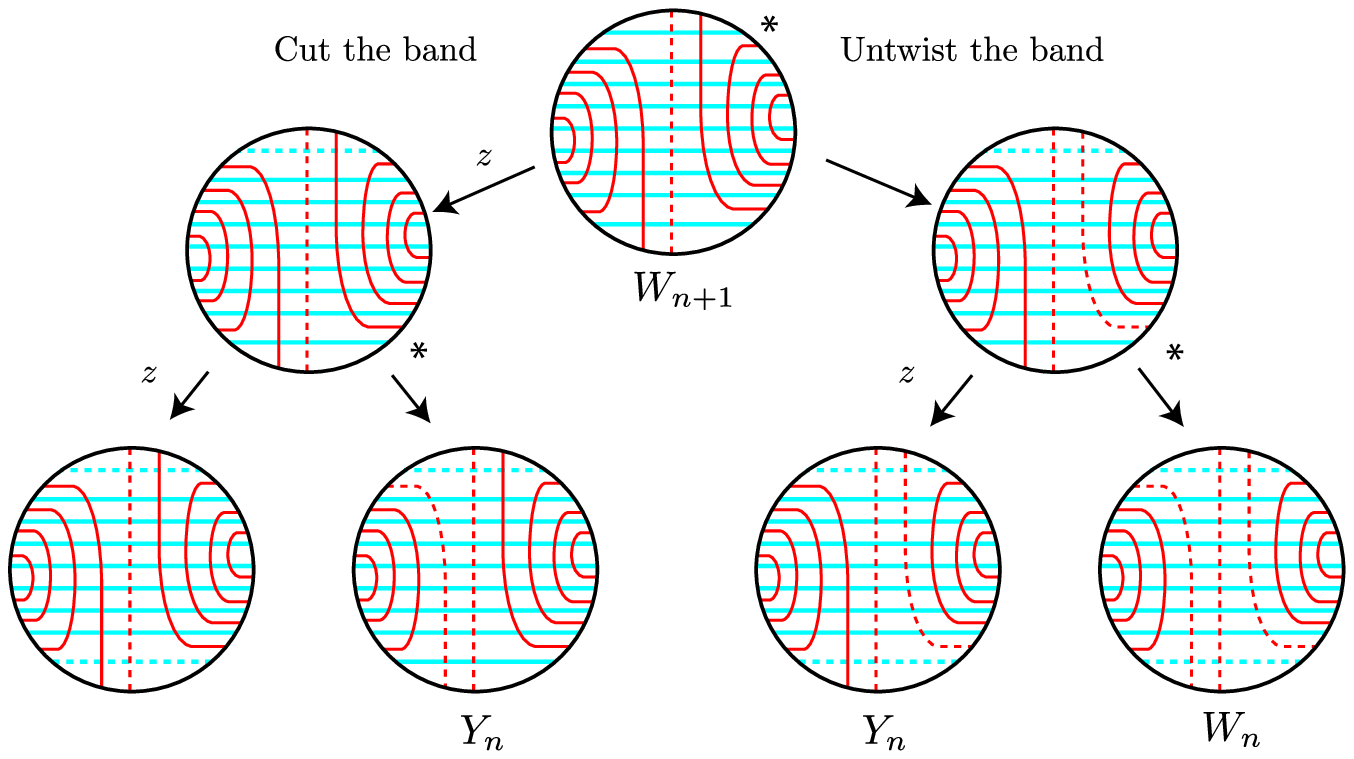}{9}}
\centerline{Fig. 11.2}
\vspace*{3mm}
\end{minipage}

In the left of the bottom row, 
we have two parallel bands.
We apply Lemma \ref{lem:11.4},
but since the writhe is negative, we replace $z$ by $-z$.
See Fig. 11.3.
Then we have (\ref{siki:11.3}). 

\begin{minipage}{1\hsize}
\vspace*{3mm}
\centerline{\inclf{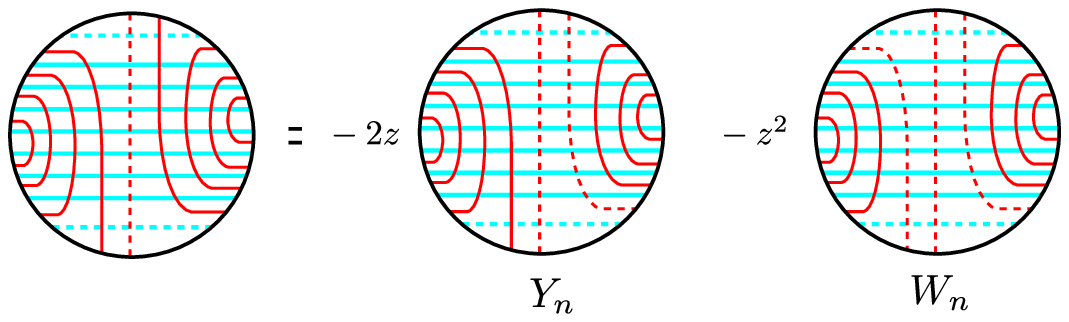}{7}}
\centerline{Fig. 11.3}
\vspace*{3mm}
\end{minipage}

Next, we prove (\ref{siki:11.2}).
We have a skein tree as
in Fig. 11.4.
In the first step, we use Fig. 11.1 (2) and
have a band with two full-twists.
On the way, we have the connected sum
of $Y_n$ and the 2-bridge link $K(\frac{1}{4})$, 
whose Conway polynomial  is equal to that of  $Y_n$ multiplied by
$-2z$.

\begin{minipage}{1\hsize}
\vspace*{3mm}
\centerline{\inclf{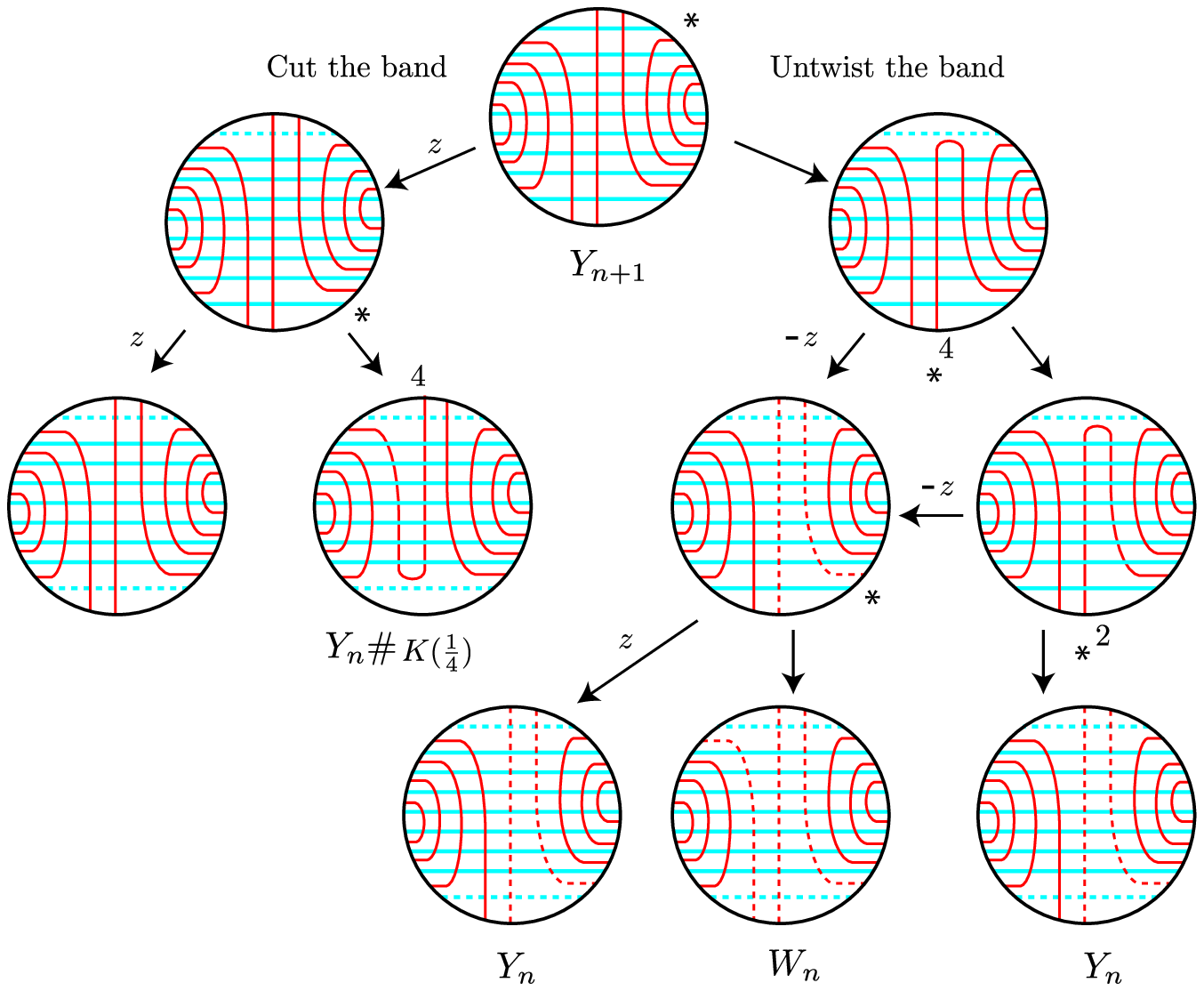}{9}}
\centerline{Fig. 11.4}
\vspace*{3mm}
\end{minipage}

In the left of the third row,
we have three parallel bands, and apply Lemma \ref{lem:11.4}
with $z$ replaced by $-z$. 
See Fig. 11.5.
Then we have (\ref{siki:11.2}).\\
Now we have completed a proof of Lemma \ref{lem:11.2}.
\qed

\begin{minipage}{1\hsize}
\vspace*{3mm}
\centerline{\inclf{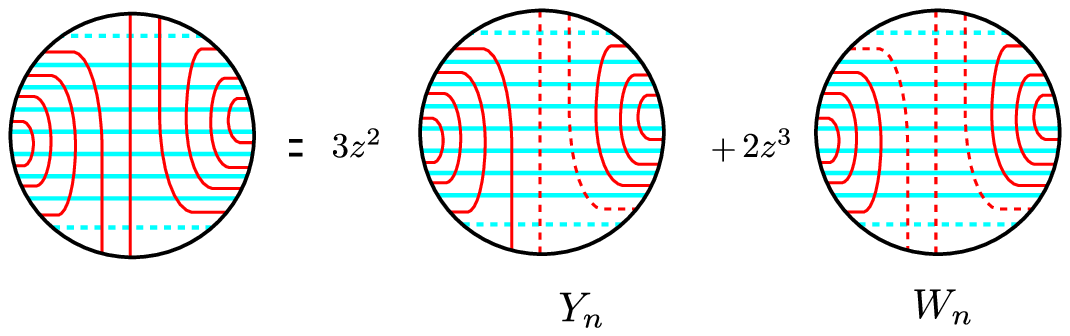}{7}}
\centerline{Fig. 11.5}
\vspace*{3mm}
\end{minipage}

\subsection{Alexander polynomials of $Y_n$}\label{11.2}
In the previous subsection we investigated the Conway polynomial 
$c_n(z)$ of $Y_n$. 
Now we translate $c_n(z)$ to the normalized Alexander polynomials 
$h_n(t)$ of $Y_n$.
By Proposition \ref{prop:11.3} we inductively have the following:

\begin{cor}\label{cor:11.5}
$\deg f_1(z)=1,\deg f_2(z)=4$. 
For $m>1$, $\deg f_{2m-1}(z)=4(m-1),\deg f_{2m}(z)=4m$.
The leading coefficient of $f_k$ is equal to 
$1$, for $k\ge 1$.
\end{cor}

The Alexander polynomial of $Y_n$ is obtained by
putting $z=\sqrt{t}-\frac{1}{\sqrt{t}}$ in $c_n(z)$.

\begin{prop}\label{prop:11.6}
Let $h_n(t)$ be the normalized Alexander polynomial of $Y_n$.
Let $\mu(t)=1-3t+t^2$ and $\rho(t)=1-t+t^2$.
Let $g_{2m-1}(t)=f_{2m-1}(\sqrt{t}-\frac{1}{\sqrt{t}})t^{2(m-1)}$
and $g_{2m}(t)=f_{2m}(\sqrt{t}-\frac{1}{\sqrt{t}})t^{2m}$

Then we have:
\begin{align}
h_{2m-1}(t)&=
\mu(t)^{2m-1} g_{2m-1}(t)\\
h_{2m}(t)&=
\mu(t)^{2m-1} g_{2m}(t) 
\end{align}
Moreover, $\mu(t)\not|g_{2m-1}(t), \mu(t)\not|g_{2m}(t)$.
\end{prop}

\begin{prop}\label{prop:11.7}
$\deg g_{2m-1}(t)=4(m-1), \deg g_{2m}=4m$,
$g_k(t)$'s are monic and reciprocal.
\begin{align}
\label{siki:11.16} & g_1=1, g_2=1-10t+19t^2-10t^3+t^4\\
\label{siki:11.17} & g_{2m+1}=2 g_{2m}-\rho(t)^2 g_{2m-1}\\
\label{siki:11.18} & g_{2m+2}=2\mu(t)^2 g_{2m+1}-\rho(t)^2 g_{2m}
\end{align}
\end{prop}

Now we study the interlacing property of $g_{k}$ and 
$g_{k+1}$. 
The following is the main theorem of this section.

\begin{thm}\label{thm:11.8}
For $m\ge1$, we have the following, where $\mu_0\doteq0.382$ is
the zero of $\mu(t)=1-3t+t^2$ in $[0,1]$.\\
(1) $g_{2m-1}g_{2m}$ is simple,\\
(2) $(t-1)\mu(t)g_{2m-1}$ and $g_{2m}$ are interlaced.\\
Namely, let 
$\alpha_1 < \alpha_2 < \cdots < \alpha_{2(m-1)}$ 
 be the zeros of
$g_{2m-1}$ in $[0,1]$,
and $\beta_1 < \beta_2 < \cdots < \beta_{2m}$ be the zeros of 
$g_{2m}$ in $[0,1]$.
Then,\\
$
\beta_1<\alpha_1<\beta_2
<\dots
< \alpha_{m-1} < \beta_m 
< \mu_0 <
 \beta_{m+1}< \alpha_{m} 
 <\dots< 
\alpha_{2(m-1)}<\beta_{2m}<1$.\\
(3) $g_{2m}g_{2m+1}$ is simple,\\
(4) $(t-1)g_{2m}$ and $\mu(t)g_{2m+1}$ are interlaced.
Namely, let
$\gamma_1 < \gamma_2 < \cdots < \gamma_{2m}$ 
be the zeros of 
$g_{2m+1}$ in $[0,1]$.
Then,\\
$
\gamma_1<\beta_1<\gamma_2
<\dots
< \gamma_{m} < \beta_m 
< \mu_0 <
 \beta_{m+1}< \gamma_{m+1} 
 <\dots< 
\beta_{2m}<\gamma_{2m}<1$.
\end{thm}

\begin{ex}\label{ex:11.9}
We list $g_k$'s with $k$ up to 6.
Note that in general, $g_k(0)=g_k(1)=1$.
\begin{align*}
g_1=&1,
g_2=1 - 10 t + 19 t^2 - 10 t^3 + t^4,
g_3=1-18 t+35 t^2-18 t^3+t^4\\
g_4=&1-36 t+266 t^2-784 t^3+1107 t^4-784 t^5+266 t^6-36 t^7+t^8\\
g_5=&1 - 52 t + 458 t^2 - 1424 t^3 + 2035 t^4 - 1424 t^5 + 458 t^6 - 52 t^7 + t^8\\
g_6=&(1 - 10 t + 19 t^2 - 10 t^3 + t^4) (1 - 68 t + 522 t^2 - 1552 t^3 + 
   2195 t^4 \\
   &- 1552 t^5 + 522 t^6 - 68 t^7 + t^8)
\end{align*}
Fig. 11.6 depicts adjacent pairs of $g_k$'s.
\begin{center}
\inclf{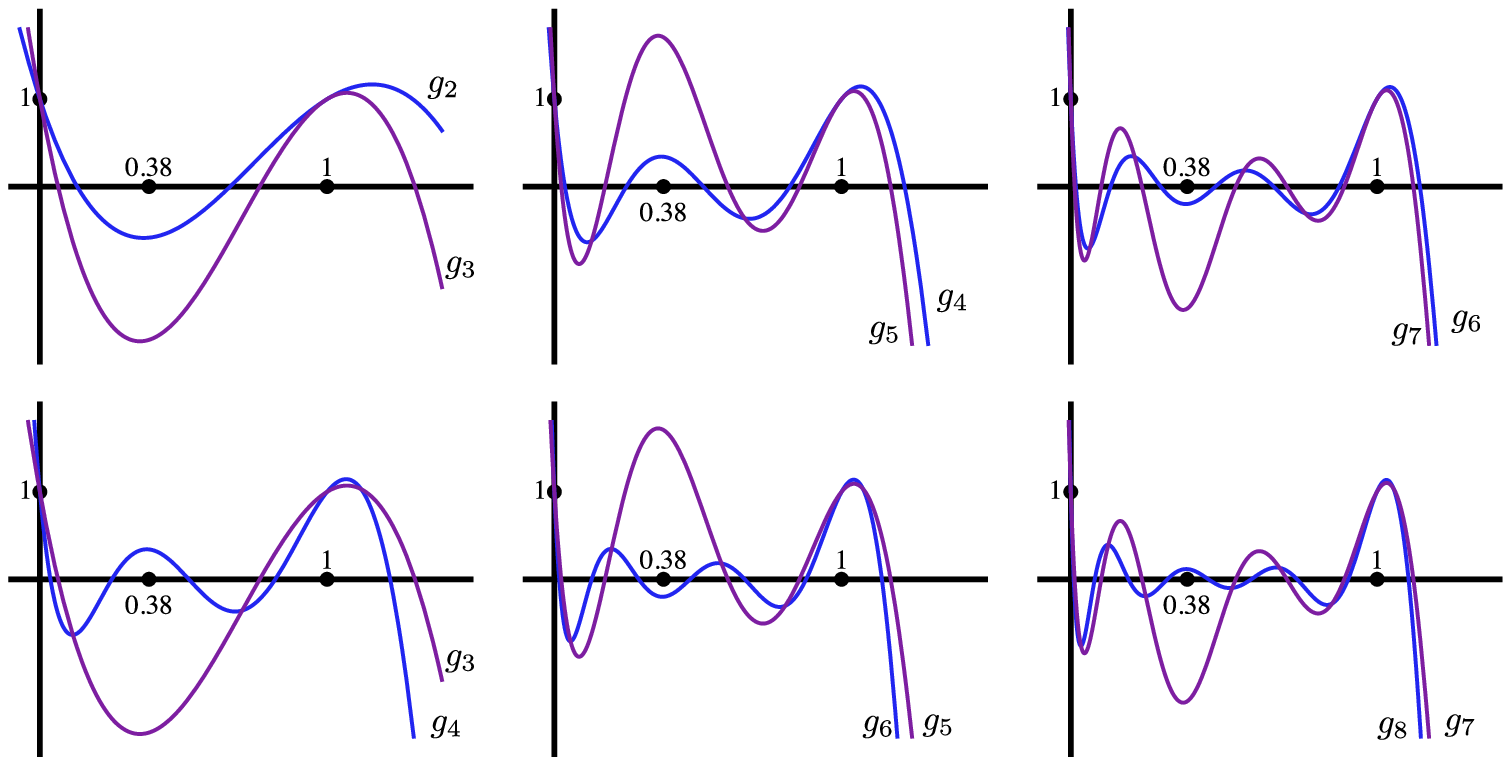}{12}
\centerline{Fig.11.6}
\end{center}
\end{ex}

{\it Proof of Theorem \ref{thm:11.8}.}
We prove the theorem by induction.
The zeros of $g_2$ in $[0,1]$ are approximately 
$0.129$ and $0.662$, and those of $g_3$ are
approximately $0.063$ and $0.765$.
See Fig. 11.7.
Since $\mu_0\doteq 0.38$, the claim is true for the pairs
$(g_1,g_2)$ and $(g_2,g_3)$,
i.e.,
$g_1 g_2$ is simple and $\mu(t)(t-1)g_1$ and $g_2$ are interlaced,
also, $g_2 g_3$ is simple and $(t-1)g_2$ and $\mu(t) g_3$ are interlaced.

\begin{center}
\inclf{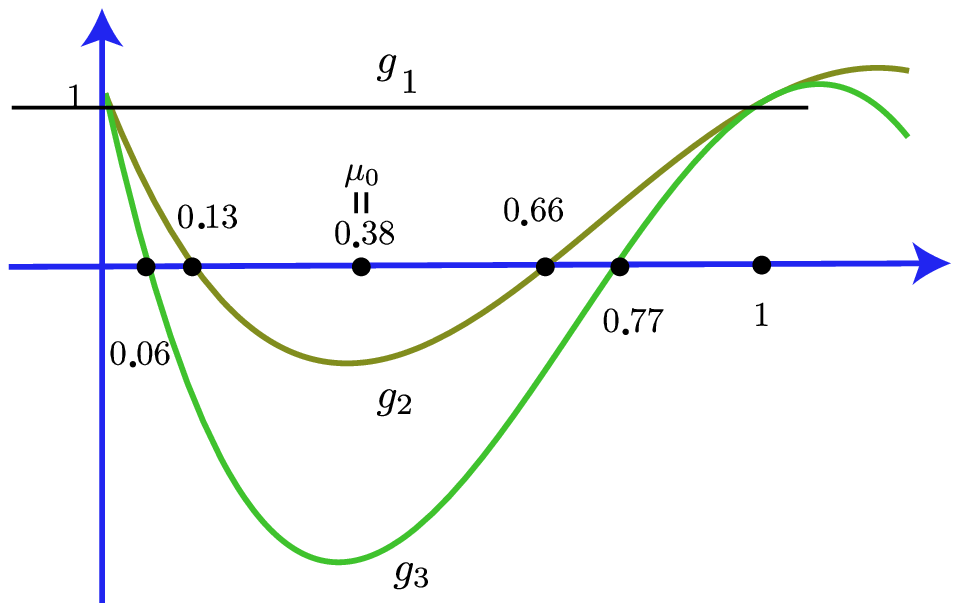}{5}
\centerline{Fig.11.7}
\end{center}

For general cases, we first fix $m$ and assume (1) and (2) of
the statement of Theorem \ref{thm:11.8}, prove (3) and (4) for that $m$,
and then prove (1) and (2) with $m$ replaced by $m+1$.

By Proposition \ref{prop:11.7}, 
$g_{2m+1}=2g_{2m}-(1-t+t^2)^2g_{2m-1}$.
So we examine intersection points of 
$y_1(t)=(1-t+t^2)^2g_{2m-1}$ and
$y_2(t)=2 g_{2m}$.
See Fig. 11.8 for $m$ odd (here $m=5$). 
For the case $m$ being even, the figure is similar.

Since $0<1-t+t^2\leq 1$ in $[0,1]$, $y_1$ and $g_{2m-1}$ have the
same real zeros, and the graph of $y_1$ is obtained from that of $g_{2m-1}$
by moving it toward the $x$-axis, 
fixing the points of intersection with the $x$- and $y$-axes.
By assumption, 
$\{$zeros of $y_1\}\cup\{\mu_0\}\cup\{1\}$ and 
$\{$zeros of $y_2\}$ are interlaced.
Since $g_{2m-1}(0)=g_{2m-1}(1)=g_{2m}(0)=g_{2m}(1)=1$,
we have $y_1(0)=y_1(1)=1, y_2(0)=y_2(1)=2$.
Let $\gamma_1<\cdots\gamma_{2m}$ be the zeros of $g_{2m+1}$ in $[0,1]$.
Then we see that $\gamma_1<\beta_1<\alpha_1<\dots<
\alpha_{m-1}<\gamma_{m}<\beta_{m}<\mu_0<
\beta_{m+1}<\gamma_{m+1}<\alpha_m<\beta_{m+2}<\gamma_{m+2}<
\dots<\alpha_{2m-2}<\beta_{2m}<\gamma_{2m}<1$.
Therefore, we have (3) and (4) of Theorem \ref{thm:11.8} for
the fixed $m$.

\begin{center}
\inclf{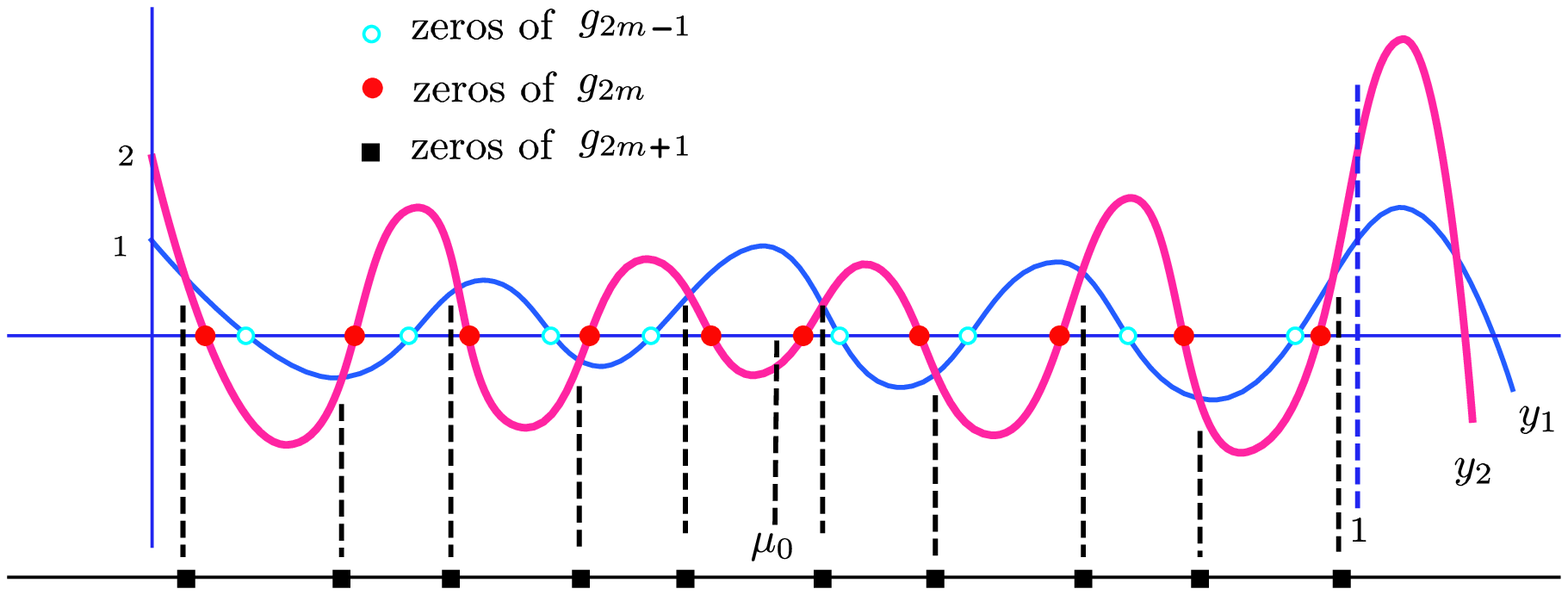}{11.5}
\centerline{Fig.11.8}
\end{center}

Next assume (3) and (4) holds for a fixed $m$.
Then we examine $g_{2m+1}$ and $g_{2m+2}$ and
prove (1) and (2) with $m$ replaced by $m+1$.
By proposition \ref{prop:11.7}, 
$g_{2m+2}=2(1-3t+t^2)^2g_{2m+1}-(1-t+t^2)^2g_{2m}$.
So we examine intersection points of 
$y_1(t)=(1-t+t^2)^2g_{2m}$ and
$y_2(t)=2 (1-3t+t^2)^2 g_{2m+1}$.
Note that $\mu_0$ is a zero of 
$y_2$ of order 2.
See Fig. 11.9 below  for $m$ odd (here $m=5$). 

\begin{center}
\inclf{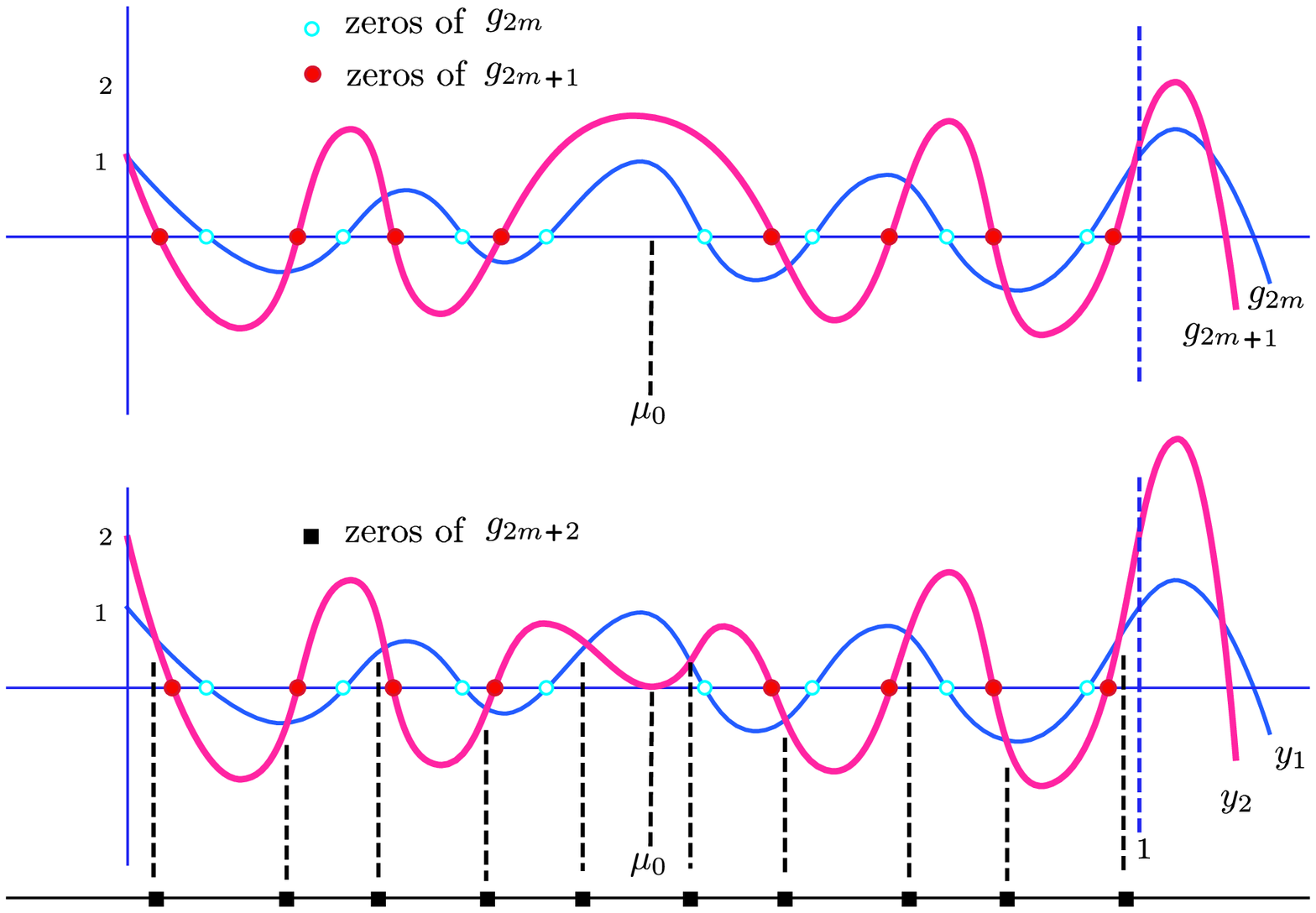}{11.5}
\centerline{Fig.11.9}
\end{center}

For the case $m$ being even, the figure is similar.
Let $\delta_1<\dots<\delta_{2m+2}$ be the zeros of $g_{2m+2}$
in $[0,1]$. 
Then we see that $\delta_1<\beta_1<\delta_2<\beta_2<\dots<
\delta_{m}<\beta_{m}<\delta_{m+1}
<\mu_0<
\delta_{m+2}<\beta_{m+1}<\dots
<\delta_{2m+1}<\beta_{2m}<
\delta_{2m+2}<1$.
Therefore, we have (1) and (2) with $m$ replaced by $m+1$.\\
Proof of Theorem \ref{thm:11.8} is now complete.
\qed

\begin{qu}\label{qu:11.5}
If the zeros of $\Delta_{K_1}(t)$ and 
$\Delta_{K_2}(t)$ are interlaced, how are $K_1$ and $K_2$
related geometrically ?
\end{qu}

\section{$c$-stable knots and links
}\label{12}

It is well-known \cite{mi} 
that if the absolute value of the signature of a
knot $K$ is equal to the degree of the Alexander polynomial,
then all the zeros of $\Delta_K(t)$ are on the unit circle and hence, 
$K$ is $c$-stable. 
However, the converse is not necessarily true,
even for $2$-bridge knots.
In this section, first we discuss $c$-stable
2-bridge knots and links, and then we show
some general construction of $c$-stable knots or links.

\subsection{Regular and exceptional $c$-stable $2$-bridge
knots and links}\label{12.1}

We begin with the following proposition:

\begin{prop}\label{prop:12.1}
Let $r=[2a_1,2a_2,\dots,2a_n]$. 
If all $a_i$'s have the same sign, then $K(r)$ is $c$-stable.
\end{prop}

{\it Proof.}
Suppose that all $a_j$'s are positive.
Let $M$ be a Seifert matrix of twisted chain type.
(See Fig 4.3). 
Then $M+M^{T}$ is positive definite by
Positivity lemma (Proposition \ref{prop:4.7}),
and hence $|\sigma(K(r))|=\deg \Delta_{K(r)}(t)$.
Therefore, $K(r)$ is $c$-stable.
\qed

The converse of Proposition \ref{prop:12.1} does not hold.
The simplest counter-example is $K(r)$, where
$r=[2,8,-2,-2]$.
In fact $\Delta_{K(r)}(t)=(2-3t+2t^2)^2$ and hence $K(r)$ is 
$c$-stable, but $\sigma(K(r))=0$.
In fact, 
we have the following proposition.

\begin{prop}\label{prop:12.2}
Let $r=[2,2k,-2,-2]$. Then we have:\\
(1) If $k<0$, then $K(r)$ is strictly bi-stable.\\
(2) If $k=1,2,3$, then $K(r)$ is totally unstable.\\
(3) If $k\ge 4$, then $(r)$ is $c$-stable.
\end{prop}

{\it Proof.}
First, we see that $\Delta_{K(r)}(t)=k(t-1)^2(t^2-t+1)+t^2$,
hence the modification of $\Delta_{K(r)}(t)$ is
$f(x)=k(x-1)(x-2)+1$. The conclusion follows by checking the
intersection of two curves $y_1=(x-1)(x-2)$ and
$y_2=-\frac{1}{k}$.
\qed

\begin{yosou}\label{conj:12.3}
Let $r_m=[
\underbrace{2,2, \dots, 2}_{m-1},2k,
\underbrace{-2,-2,\dots, -2}_{m}]$.
Then if $m$ is even, $K(r_m)$ is $c$-stable for
sufficiently large $k$.
To be more precise, there exists a positive
integer $N_m$ such that (1) if $m$ is even, then $K(r_m)$
is $c$-stable for $k\ge N_m$ and (2) if $m$ is odd, then 
$K(r_m)$ is $c$-stable for $k\le -N_m$. 
\end{yosou}
We can show that $N_2=4, N_3=3$ and $N_4=7$
(See Appendix C).
These knots $K(r_m)$ are exceptional $c$-stable knots.\\

Example \ref{ex:12.4} below gives an exceptional $c$-stable link.

\begin{ex}\label{ex:12.4}
Let $r=[2,2,2,-6,-2]$. 
Then $K(r)$ is a $c$-stable link.
In fact, $\Delta_{K(r)}(t)=(t-1)(3t^4-6t^3+7t^2-6t+3)$ is $c$-stable.
\end{ex}

\subsection{ Construction of $c$-stable quasi-rational knots 
and links
}\label{12.2}

Let $K$ be a quasi-rational knot or link such that
a Seifert matrix of $K$ is of the form :$
M=\tbt{A}{O}{B}{C}$, where
$A = \diag\{a_1, a_2, \cdots,
a_p\}, a_j > 0, 1 \leq j \leq p$ and 
$C=\diag\{c_1,c_2, \cdots, c_q\}, c_j>0,
1 \leq j \leq q$. 
Let $B = [b_{i,j}]_{1 \leq i \leq p, 1 \leq j \leq q}$.

\begin{prop}\label{prop:12.5}
Suppose $M$ satisfies the following conditions:\\
(1) 
$a_k > \frac{1}{2} \{|b_{1,k}|+|b_{2,k}|+ \cdots +|b_{q,k}|\}$ for
$k=1,2, \cdots, p$,\\
(2) $c_{\ell} > \frac{1}{2} \{|b_{\ell,1}|+|b_{\ell,2}|+ \cdots
+|b_{\ell,p}|\}$ for $\ell =1,2, \cdots, q$.\\
Then $K$ is $c$-stable.
\end{prop}

{\it Proof.}
 Since a symmetric matrix $\widehat{M}$ = $M + M^{T}$ is positive
definite, the signature of $\widehat{M}$ is equal to $p+q$ and also 
$p+q$ is
the degree of the Alexander polynomial of $K$. 
Hence, $K$ is $c$-stable.
\qed

Note that $K$ is generally non-alternating.

\begin{rem}\label{rem:12.6}
Conditions (1) and (2) in Proposition \ref{prop:12.5} 
are sufficient conditions 
for knots $X_n=X_n(2a_1,2a_2,\dots, 2a_n|2b_1,2b_2,\dots, 2b_n)$
defined in Example \ref{ex:6.4} to be $c$-stable. 
Suppose that all $a_j>0$ and $b_j>0$. 
Then if $X_n$ is $c$-stable, at least (1) or (2) is necessary.
In fact, $K_3=X_3(2, 2, 2|2,2,2)$ is not $c$-stable, but strictly bi-stable.
Here, $\Delta_{K_3}(t)=t^6-4t^4+7t^3-4t^2+1$ and $K_3$ is not alternating.
On the other hand, $X_3(4,2,2|2,2,2)$ is $c$-stable.
\end{rem}

\begin{prop}\label{prop:12.7}
Let $X_n = X_n (2a_1, 2a_2, \cdots, 2a_n \mid 2b_1, 2b_2, \cdots, 2b_n)$
be a quasi-rational knot defined in Example \ref{ex:6.4}. 
Suppose\\
(1) $a_1\geq n/2, a_2 \geq (n-1)/2, \cdots, a_k \geq (n-k+1)/2, \cdots,
a_n \geq 1/2$, and\\
(2) $b_1 \geq 1/2, b_2 \geq 2/2, \cdots, b_k \geq k/2, \cdots, b_n 
\geq n/2$.\\
Then $X_n$ is $c$-stable.
\end{prop}

{\it Proof.} Let $M$ be a Seifert matrix given in Section 6.3. Since $a_j$ and
$b_j$ are integers, it follows that $a_n \geq 1$
and $b_1 \geq 1$. Therefore, the $n^{\rm th}$ row and $(n+1)^{\rm st}$ row
are excessive. Apply the proof of Positivity Lemma on
the matrix $M + M^{T}$ to show that $M + M^{T}$ is positive definite.
\qed 

\begin{prop}\label{prop:12.8}
Let $Y_{2n+1} = Y_{2n+1}(2a_1,2a_2, \cdots, 2a_{2n+1} \mid
2b_1, 2b_2, \cdots, 2b_{2n+1})$ be a quasi-rational knot 
defined in Example \ref{ex:6.5}. 
Suppose\\
(1) $a_1$ and $a_{2n+1} \geq 1, a_2$ and $a_{2n} \geq 2, \cdots, a_k$ 
and $a_{2n+2-k} \geq
k, \cdots, a_{n+1} \geq n+1$,\\
(2) $b_1$ and $b_{2n+1} \geq 1, b_2$ and $b_{2n} \geq 2, \cdots,
b_k$ and
$b_{2n+2-k} \geq k, \cdots, b_{n+1} \geq n+1$.\\
Then $Y_{2n+1}$ is $c$-stable.
\end{prop}

{\it Proof.} $M + M^{T}$ satisfies all conditions of Positivity Lemma.
\qed

\subsection{General construction of $c$-stable knots and links
}\label{12.3}

The previous propositions show 
that given an arbitrary quasi-rational knot or
link, we can make it $c$-stable by changing the number of full twists
on some bands. 

In this subsection, we generalize this result to that
we can construct a $c$-stable knot or link from
a  given Seifert surface.

In case of Seifert surfaces specified by graphs as before,
we can construct a $c$-stable knot or link with the same
underlying graph.

\begin{thm}\label{thm:12.9}
Let $F$ be a Seifert surface for a knot or link $K$, with
${\rm rank}H_1(F,\mathbb{Z})=n$.
Suppose that a system of mutually disjoint $n$ arcs
$\alpha_1, \dots, \alpha_n$ properly embedded in $F$
is specified so that $F\setminus \cup_i\alpha_i$ is a disk.
Let $\widetilde{F}$ be a Seifert surface obtained by full-twisting 
$F$ along each arc $\alpha_i$, $k_i$ times.
Denote by $\widetilde{K}$ the knot or link $\partial \widetilde{F}$.
Then, there exist $N_i \in{\mathbb N}\ (i=1, 2\dots, n)$ such that 
if $k_i\geq N_i$ for each $i$, then $\widetilde{K}$ 
is $c$-stable.
\end{thm}

{\it Proof.}
Let $L=\{\ell_1,\dots,\ell_n\}$ be a set of embedded loops in $F$ such that
for each $i$, $\alpha_i \cap (\cup_i \ell_i)$ is a single transverse point in $\ell_i$. 
Note that such a system $L$ is unique up to isotopy
since $F\setminus \cup_i \alpha_i$ is a disk. 
Then $L$ with an arbitrary orientation
gives Seifert matrices $S$ for $F$
and $\widetilde{S}$ for $\widetilde{F}$.
Since twisting $F$
along $\alpha_i$ affects only the self-linking number of
$\ell_i$, $\widetilde{S}-S$ is a diagonal matrix whose $i$th diagonal entry is $k_i$.
Let $M$ be the symmetric matrix $\widetilde{S}+\widetilde{S}^T=(m_{i,j})$.
If each $k_i$ is large enough, we have $m_{i,i}>0$ and
$m_{i,i}>\Sigma_{j\neq i}\ |m_{i,j}|$ and hence by
Strong Positivity Lemma (Proposition \ref{prop:4.5}), $M$ is positive definite.
Then the signature $\sigma(M)$ is equal to $n$.
By \cite{mi}, $\Delta_{\widetilde{K}}(t)$ has at least $n$ of its zeros on
the unit circle, and hence $n\le \deg\Delta_{\widetilde{K}}(t)$.
Since $n=2g(\widetilde{F})$, we have ${\rm deg}\Delta_{\widetilde{K}}(t)\le n$.
Therefore, $\deg\Delta_{\widetilde{K}}(t)= n$ and the conclusion follows.
\qed

Note that if $M$ is positive definite, 
we have $n\le \deg\Delta_{\widetilde{K}}(t)\le 2g(\widetilde{K})\le n$,
and hence $\widetilde{F}$ is a minimal genus Seifert surface for $\widetilde{K}$.

Before we discuss some application of 
Theorem \ref{thm:12.9} we prove one proposition.

\begin{prop}\label{prop:12.10}
Let $G$ be a positive (or negatie) admissible connected 
planer graph on a disk $D$.
Suppose that $G$ satisfies (\ref{siki:7.1}).
Let $F(G)$ be the surface representing $G$.
Then $K=\partial F(G)$ is alternating and $c$-stable.
\end{prop}

We should note that $G$ is not necessarily an even graph.

{\it Proof.}
Since a diagram is special alternating, $K$ is special alternating.
Now let $M$ be a Seifert matrix obtained from $F(G)$.
Then $M+M^T$ is positive (or negative)
definite by Positivity lemma and hence $K$ is $c$-stable.
\qed

Now take finitely many disks $D_1,D_2,\dots, D_n$ each of which has
a positive admissible graph $G_j$ ($1\le j\le n$)
that satisfies (\ref{siki:7.1}).
Consider a Murasugi sum $F$ of surfaces $F(G_1), F(G_2), \dots F(G_n)$
glued by an arbitrary fashion.
Then the knot $K=\partial F$ is generally not $c$-stable,
but by Theorem \ref{thm:12.9}, we can make $K$ to be $c$-stable,
by changing at most $s$ weights in $\{G_1,G_2,\cdots, G_n\}$, where
$s=$rank$H_1(F;\ZZ)$.

\begin{ex}\label{ex:12.11}
The knot or link in the left is not $c$-stable, but by changing 
at most four weights, it becomes $c$-stable.
\end{ex}

\dodecichi

\subsection{Interlacing property of zeros on the unit circle}\label{12.4}

In this sub-section, 
we define the interlacing property for two $c$-stable
real polynomials.

\begin{dfn}\label{def:12.12}
Let $f(t)$ and $g(t)$ be $c$-stable real
polynomials, and let 
$\{\alpha_j, 1 \leq j \leq n\}$ and 
$\{\beta_k, 1 \leq k \leq m\}$
be, respectively, the unit complex zeros of 
$f(t)$ and $g(t)$ with
a property that ${\rm Im}(\alpha_j) \geq 0$ and 
${\rm Im} (\beta_k) \geq 0$. 
Then we say
that $f(t)$ and $g(t)$ are {\it interlaced} if 
$\{{\rm Re}(\alpha_j), 1 \leq j \leq
n\}$ and
$\{{\rm Re}(\beta_k), 1 \leq k \leq m\}$ are interlaced.
\end{dfn}

As a typical example, we prove the following proposition:

\begin{prop}\label{prop:12.13}
Let $r_n = [\underbrace{2,2, \cdots, 2}_{n}]$. Then
$\Delta_{K(r_n)}(t)$ and
$\Delta_{K(r_{n-1})} (t)$ are interlaced.
\end{prop}

{\it Proof.} 
The unit complex zeros 
$\{\alpha_k\}$of $\Delta_{K(r_n)} (t)$
with
${\rm Im}(\alpha_j) \geq 0$ are:
(1) if $n$ is even, say $2m$, then 
$\{\alpha_k\}=\{e^{\frac{(2k+1)\pi}{2m+1}}, 0 \leq k \leq
m-1\}$, and
(2) if $n$ is odd, say $2m+1$, then 
$\{\alpha_k\}=\{e^{\frac{2k\pi}{2m+2}}, 0 \leq k \leq m\}$.

Then, the proposition follows from inequalities below.

\begin{align}\label{siki:12.1}
(1)& \cos \frac{2k\pi}{2m+2} > \cos \frac{(2k+1)\pi}{2m+1} 
> \cos\frac{(2k+2)\pi}{2m+2}, 0 \leq k \leq m-1.\nonumber\\
(2)& \cos \frac{(2k+1)\pi}{2m+3} > \cos \frac{(2k+2)\pi}{2m+2} 
> \cos\frac{(2k+3)\pi}{2m+3},
0 \leq k \leq m-1
\end{align}
\qed

The following theorem is the $c$-stable version of
Theorem \ref{thm:9.4}, and is proved by using modified Alexander polynomials
instead of Alexander polynomials. Therefore, the details are omitted.

\begin{thm}\label{thm:12.14}
Let $r =[2a_1, 2a_2, \cdots, 2a_n], a_j > 0, 1 \leq j
\leq n$ and $s= [2a_1, 2a_2, \cdots, 2a_{n-1}]$.
Then are $\Delta_{K(r)}(t)$ and $\Delta_{K(s)}(t)$ interlaced.
\end{thm}

\begin{probl}\label{probl:12.15}
Characterize $c$-stable alternating knots and links.
\end{probl}

\section{Bi-stable knots and links}\label{13}

\subsection{Bi-stable 2-bridge knots and links}\label{13.1}

A bi-stable knot has not only real zeros, but also unit complex zeros. 
Therefore, we could say that it combines two parts, one is a $c$-stable part
and another is a real stable part. From this point of view, the following
theorem is not surprising, although a proof is not straightforward.

\begin{thm}\label{thm:13.1}
Let $r = [2a_1, 2a_2, \cdots, 2a_{2m}, 2b_1, -2b_2, 2b_3,-2b_4, 
\cdots, -2b_{2p}]$, 
(or $r = [ 2b_1, -2b_2, 2b_3, -2b_4, \cdots, 
-2b_{2p}, 2a_1, 2a_2, \cdots, 2a_{2m}]$),
where $a_j>0, 1\le j\le 2m$ and $b_k > 0, 1\le k\le 2p$, 
Then $K(r)$ is bi-stable. 
The number of the real zeros is $2p$ and that of the unit complex
zeros is $2m$.
\end{thm}

{\it Proof.} 
Since the signature of $K(r)$ is $2m$, it follows that the
number of the unit complex zeros is at least $2m$. 
Therefore, it suffices to show
that the number of the real zeros is (at least) $2p$.
First we prove the following lemma.
\begin{lemm}\label{lem:13.2}
(1) Let $r^{\prime} = [2a_1, 2a_2, \cdots, 2a_{2m-1}]$. 
Then $(t-1)$ divides $\Delta_{K(r^{\prime})}(t)$, but $(t-1)^2$ does not.
(2) Let $r^{\ast} = [-2b_2, 2b_3, -2b_4, \cdots, -2b_{2p}]$. 
Then $(t-1)$ divides $\Delta_{K(r^{\ast})}(t)$, but $(t-1)^2$ does not.
\end{lemm}

{\it Proof.} 
(1)
$\Delta_{K(r')} (t) /(t-1) = \Delta_{K(r')}(t,t)$, 
where
$\Delta_{K(r')}(x,y)$ denotes the 
2-variable Alexander polynomial of
a 2-component link $K(r')$. 
Then $|\Delta_{K(r')}(1,1)|$
is the absolute
value of the linking number $\ell$ between two components of
$K(r')$. \cite{tor} 
Since $|\ell | = |a_1 + a_3 + \cdots +a_{2m-1}| > 0$,
$\Delta_{K(r')}(t, t)$ is not divisible by $t-1$.
(2) $K(r^*)$ is stable, and all the zeros are simple.
\qed

\begin{lemm}\label{lem:13.3}
Let $D_K (t)$ be the normalization of $\Delta_K (t)$. 
Then we have
\begin{equation}
D_{K(r)}(t) = D_{K(r_1)} (t) D_{K(r_2)}(t) +t
D_{K(r^{\prime})} (t) D_{K(r^{\ast})}(t),
\end{equation}
\noindent
where $r_1 = [2a_1, 2a_2, \cdots, 2a_{2m}]$, $r_2 = [2b_1, -2b_2, 2b_3,
-2b_4, \cdots, -2b_{2p}]$, and $r^{\prime}$ and $r^{\ast}$ are given in
Lemma \ref{lem:13.2}.  
\end{lemm}

{\it Proof.}
Using a twisted chain type Seifert surface of $K(r)$, 
we have a Seifert matrix $M$ of the form: 
$M=\tbt{A}{B}{O}{C}$,
where $A$, $C$ are Seifert matrices
of $K(r_1)$ and $K(r_2)$, respectively, and $B$ has only $1$ at the
$(2m,2m+1)$-entry and $0$ elsewhere (see (\ref{siki:4.1})). 
Then it is easy to see that
$D_{K(r)}(t)$ has the required form.
\qed

We return to a proof of Theorem \ref{thm:13.1}. 
We know now
%
\begin{align}\label{siki:13.2}
(1)\  &f_m(t) = D_{K(r_1)}(t){\rm \  is}\  c{\rm -stable\ and\ hence\ } 
f_m > 0 {\rm \  for\ any\ real\ } t.
\nonumber\\
(2)\ &D_{K(r^{\prime})} (t){\rm \  is\ } c{\rm -stable,\  and\  has\  
only\ one\ real\ zero\ that\ is\  1},
\nonumber\\
&{\rm and\  hence\  we\  can\  write\ } D_{K(r^{\prime})}(t) = (t-1) g_m (t)
{\rm \  and\ } 
g_m (t) >0
{\rm \  for\  any\  real\  }t, \nonumber\\
(3)\ &D_{K(r_2)}D_{K(r^*)} {\rm \ is\  simple,}\nonumber \\
(4)\ &D_{K(r_2)}(t){\rm  \ is\ stable\ and\ has\ }
 2p {\rm \ positive\ real\ zeros,\ say,\ }
\beta_1 < \beta_2 < \cdots < \beta_p, {\rm \ in\  }[0,1],
\nonumber\\
(5)\ &D_{K(r^{\ast})}(t) = (t-1) h_p (t){\rm \  is\ stable\ and\ has\ }
 (2p-1) {\rm \ real\  zeros,\ say,\ } 
 \nonumber\\
 &\alpha_1 < \alpha_2 < \cdots < \alpha_{p-1} < \alpha_p (=
1){\rm \  in\ } [0,1], {\rm \ and\ }
h_p (1) \ne 0. 
\nonumber\\
&{\rm \ Further,\ } 
\{ \beta_j, 1 \leq j \leq p\} {\rm \  and\ } 
\{\alpha_j,1 \leq j \leq p\}{\rm \ are\ interlaced,\  i.e.,\ } 
\nonumber\\
& \beta_1 < \alpha_1 < \beta_2 < \alpha_2 <
\cdots < \alpha_{p-1} < \beta_p < \alpha_p = 1.
\end{align}
Using these notations, we can write
%
\begin{align}\label{siki:13.3}
D_{K(r)}(t) = f_m(t) D_{K(r_2)}(t) +t (t-1)^2 g_m (t) h_p
(t),{\rm \  and\ }
g_m(1) \ne 0 \ne h_p (1).
\end{align}
Now we calculate the number of real zeros of $D_{K(r)} (t)$.
From (\ref{siki:13.3}), 
we see that the real zeros of $D_{K(r)}(t)$ are
determined by the intersection 
$\{\gamma_j\}$ of two curves $y_1 = f_m (t) D_{K(r_2)}(t)$ 
and
$y_2 = - t (t-1)^2 g_m (t) h_p (t)$. Using the fact that $\{\alpha_j 
\}$ and $\{\beta_j \}$
are interlaced, we have a graph below. 
Note that $f_m(t), g_m(t) > 0$ for any real $t$, and $t=1$ is a double zero for $y_2$. Further, $y_1(0)>0$ and $y_2'(0)<0$.

(1) If $p$ is even, we have Fig. 13.1.

\kichi

From the graph, we see that 
there are (at least) $p$ points of
intersection
$\{\gamma_j, 1 \leq j \leq p\}$ in $[0,1]$ and
\begin{equation*}
\beta_1 < \gamma_1 < \beta_2 < \gamma_2 < \beta_3 < \cdots <
\beta_{p-1} < \gamma_{p-1} < \beta_p < \gamma_p < 1.
\end{equation*}

(2)	If $p$ is odd, then we have the following graph.

\kni

Therefore, $D_{K(r)}(t)$ has at least (and hence exactly) $2p$ real zeros. 
\qed

\begin{ex}\label{ex:13.4}
Let $r=[4,2,6,2,4,-6,2,-4]$.
$\Delta_{K(r)}(t)$ has hour real zeros and four unit complex zeros.
\end{ex}

\subsection{Exceptional bi-stable knots and Salem knots}\label{13.2}

A fibred knot (or link) $K$ is called a 
{\it Salem fibred knot ({\rm or} link)} (\cite{hiro}),
if $\Delta_K (t)$ is bi-stable and has exactly two real 
( $\neq 1$) zeros.
A typical example of a Salem fibred knot is a 2-bridge
knots $K(r_m)$ by Theorem \ref{thm:13.1}, 
where $r_m = [2, 2, \cdots, 2, -2]$, $m$ odd, 
(or abbreviated $[(2)^m,-2]$). 
Modifying $K(r_m)$, we obtain
a series of exceptional bi-stable knots given below.

\begin{prop}\label{prop:13.n5}
Let $r(m,n) = [(2)^m, -2, (2)^n], m \geq n \geq 0, m+n$
being odd. Then $K(r(m,n))$ is a Salem fibred knot.
\end{prop}

{\it Proof.} 
By induction on $m$ and $n$, 
it is shown easily that the
normalized Alexander polynomial 
$D_{m,n}(t)$ of $K(r(m,n))$ is given by
the following formula:

$D_{m,n}(t) = 
\sum_{k=0}^n (-1)^k (4k+1) t^k 
+ (4n+3) \sum_{k=n+1}^m(-1)^k t^k
+ \sum_{j=0}^n (-1)^{m+j+1} (4n+1-4j) t^{m+j+1}$

Since $D_{m,n}(1) = -1$, 
$D_{m,n}(t)$
has at least two real zeros.
Further, since $\sigma(K(r(m,n)) = m+n-1$, 
$D_{m,n}(t)$ has at least $m+n-1$ unit complex zeros and
hence $K(r(m,n))$ is a Salem fibred knot. 
\qed

Let $\mu (K)$ denote the maximal absolute value of 
the real zeros of $\Delta_K(t)$. 
We note that if $K$ is a Salem fibred knot,
$\mu(K)$ is equal to Mahler measure of $\Delta_K(t)$.

Now our computation suggests that for $m \geq 1$,

\begin{align}
(1)&\ \mu (K(r(m+2,0)) < \mu (K(r(m,0))\ {\rm and} \nonumber\\
(2)&\ \mu (K(r(m+2,1)) > \mu (K(r(m,1)).\ {\rm Further},\nonumber\\
(3)&\ \mu (K(r(m+n,0)) < \mu (K(r(m,n)) < \mu (K(r(m+n-1,1)).\ 
{\rm Finally} \nonumber\\
(4)&\  \lim_{m \to \infty} \mu (K(r(m,0)) = 2,\ {\rm and}\nonumber\\
(5)&\ \lim_{m \to \infty} \mu(K(r(m,1)) \doteq 3.41421
\end{align}

Beside these bi-stable knots, 
Hironaka showed two more Salem fibred
2-bridge knots \cite{hiro}. 

\begin{align*}
(1)&\  K_1 = K(s_1), s_1 = [(2)^5, (-2)^3],\ {\rm  and}\ 
\mu(K_1) \doteq 1.63557,\nonumber\\
(2)&\ K_2 = K(s_2), s_2= [(2)^9, (-2)^5]\ {\rm  and}\  
\mu(K_2) \doteq 1.42501.
\end{align*}

We find three more sporadic Salem fibred 2-bridge knots.
\begin{align}
(3)& \ K_3 = K(s_3), s_3 = [(2)^6, -2, 2, -2, -2],\ 
{\rm  and}\  \mu(K_3) \doteq 3.94748,\nonumber\\
(4)& \ K_4 = K(s_4), s_4 = [(2)^4, (-2)^3,2]\ {\rm  and}\  
\mu (K_4) \doteq 2.38215,\nonumber\\
(5)&\ K_5 = K(s_5), s_5 = [(2)^6, (-2)^5, (2)^3]\  {\rm and}\ 
\mu(K_5) \doteq 1.80017.
\end{align}

We suspect that there exist other Salem fibred 2-bridge knots. 
However, contrary to knots,
we find many Salem fibred 2-bridge links and we will study 
these links in a separate paper.

\subsection{General bi-stable knots and links}\label{13.3}
A 2-bridge knot in Theorem \ref{thm:13.1} is given as a 
quasi-rational knot shown in Fig 13.3.

\ksan

The first half part $[2a_1, 2a_2, \cdots, 2a_{2m}]$ is a special
alternating knot and it is also represented by an admissible
positive graph $G_0$. 
For example, the surface $F(G_0)$ of $K(r)$, $r=[2a,2b]$ is 
represented by $G_0$. See Fig.13.4.

\kyon

Therefore, for example, $F(G)$ of $K(r)$, 
$r=[2a_1,2a_2, 2b_1,-2b_2,2b_3,-2b_4]$
is represented by a graph $G$ below.

\kgo

This observation suggests us a construction of 
general bi-stable knots as in Proposition \ref{prop:13.5}
below. See Fig. 13.6 for example, where a bi-stable knot is
depicted by a graph, and the zeros are plotted.

\begin{prop}\label{prop:13.5}
Let $G_0$ be an admissible positive (or negative) graph
on a disk.
Attach $p$ mutually disjoint positive (or negative)
arcs to $\partial D$ and then
$p$ negative (or positive) arcs to $\partial D$ from the back side
in such a way that the first (or the last) arc 
crosses exactly one edge of $G_0$,
where the $2p$ arcs attached to $\partial D$ represent
a 2-bridge knot (or link) as is shown in Fig 13.6.
$G_0$ together with these $2p$ arcs forms a graph $G$.
Then the  knot $K=\partial F(G)$ is bi-stable.
\end{prop}
A proof is exactly the same as that of Theorem 
\ref{thm:13.1}
and we omit the details.
\qed

\kroku

A crucial point of this construction is the interlacing property of a 2-bridge knot $K(r)$, $r=[2b_1,-2b_2,2b_3,\dots, -2b_p]$.
Therefore, a knot $K(r)$ may be replaced by other stable knots
which have some kind of interlacing property.

In Fig 13.7 below, a 2-bridge knot $K(r)$ is 
replaced by a stable knot 
$K_3= X_3(2,2,2|-2,-2,-2)$. 
The knot $K$ thus obtained is bi-stable.
In fact, $\Delta_K(t) = -3+44 t-235 t^2+662 t^3-1161 t^4+1387 t^5-1161 t^6+662 t^7-235 t^8+44 t^9-3 t^{10}$ is bi-stable. The zeros are plotted in Fig 13.7 right.

\knana

However, $K(r)$ may not be replaced by an exceptional stable 
2-bridge knot $K'$. 
If the c-stable part $[2a_1, 2a_2 , \cdots, 2a_{2m}]$ is
replaced by an exceptional $c$-stable 2-bridge knot, 
then the knot is generally not bi-stable. 
For example, neither $[2,8,-2,-2,4,-2,6,-8]$ nor
$[2,8,-2,-2,10,2,-2,-10]$ is bi-stable, 
where $[2,8,-2,-2]$ is exceptionally
$c$-stable and $[10,2,-2,-10]$ is exceptionally stable.
However, it is interesting to see that 
if the second term 8 in both cases is replaced by 
a sufficiently large positive integer, 
the knots become bi-stable.

\section{Mobius Transformations}\label{14}

In this section, we study the image of the zeros of the Alexander
polynomial of a knot by a Mobius transformation $\varphi$.
We begin with the definition of a special Mobius transformation $\varphi$
that is used in this section.\\
Let $\varphi : C \cup \{\infty\} \longrightarrow C \cup \{\infty\}$ be a
Mobius transformation given by

\begin{equation}\label{siki:14.1}
\varphi (z) = \frac{1-z i}{z-i}.
\end{equation}

$\varphi$ has the following properties.

\begin{align}\label{siki:14.2}
&(1)\ \varphi{\rm \  is\ one\ to\ one\ and\ }
{\varphi}^{-1} {\rm \ is\ given\ by\ }
{\varphi}^{-1}(z) = \frac{1+zi}{z+i},
\nonumber\\
&(2)\  \varphi{\rm \  keeps\ two\ points\ }
z = \pm 1 {\rm \  fixed,}
\nonumber\\
&(3)\ 
\varphi(0) = i, \varphi (-i) = 0{\rm \  and\  }
\varphi (i) = \infty,{\rm \  and}
\nonumber\\
&(4)\  \varphi^2 (z) = 1/z.
\end{align}

We can easily check the following lemma.

\begin{lemm}\label{lem:14.1}
(1) $\varphi$ maps the interior of the unit circle centred
at 0 onto the upper half-plane, and the exterior of the unit circle
onto the lower half-plane.\\
(2) $\varphi$ maps the unit circle onto the real line and vice versa.
\end{lemm}
\qed

The following simple property of $\varphi$ is crucial to our purpose.
A proof follows from easy computations, and hence we omit the details.

\begin{prop}\label{prop:14.2}
For any $\alpha \in C$, $\alpha \ne 0, \pm i$,\\
(1) $\varphi (\alpha) + \varphi (\frac{1}{\alpha}) = 4 / (\alpha + \frac{1}{\alpha})$.\\
(2) $\varphi (\alpha) \varphi (\frac{1}{\alpha}) = 1$.\\
In particular, if $\alpha(\ne 0)$ is real or $|\alpha| = 1, \alpha \ne \pm
i$, then $\alpha + \frac{1}{\alpha}$
and $\varphi (\alpha) +\varphi (\frac{1}{\alpha})$ are both real.
\end{prop}
\qed

The main theorem in this section is the following:

\begin{thm}\label{thm:14.3}
Let $f(t)$ be a reciprocal real polynomial of even degree,
say $2n$. 
Assume that $0$ and $\pm i $ are not zeros of $f(t)$.
Then there exists a reciprocal real
polynomial $f^{\ast}(t)$ of the same degree $2n$ 
satisfying the following conditions.\\
(1) the zeros of $f^{\ast}(t)$ are exactly the image of the zeros of $f(t)$
under $\varphi$, namely,
if $\alpha_1 , \alpha_2 , \cdots, \alpha_{2n}$ are the zeros of $f(t)$,
then $\varphi (\alpha_1), \varphi(\alpha_2), \cdots, \varphi
(\alpha_{2n})$
are exactly the zeros of $f^{\ast}(t)$.\\
(2) If $f$ is an integer polynomial, then so is $f^*$.\\
(3) $|f^{\ast}(1) | = 2^n |f(1)|$.\\
(4) $|f^{\ast}(-1)| = 2^n |f(-1)|$.\\
(5) ${(f^{\ast})^{\ast}}(t) = 2^{2n} f(t)$.\\
Furthermore, such an $f^{\ast}(t)$ is unique up to $\pm 1$.
\end{thm}

Before we prove the theorem, we mention a couple of corollaries.

\begin{cor}\label{cor:14.4}
Let $\Delta_K (t)$ be the Alexander polynomial of a knot $K$
and degree $\Delta_K (t) = 2n$. Then we have the following:\\
(1) ${\Delta_K}^{\ast} (t)$ is a reciprocal integer polynomial of the
same degree, $2n$.
Therefore, ${\Delta_K}^{\ast} (t)$ is the Hosokawa polynomial of some link
$K^{\ast}$ (with an arbitrary number of components).\\
(2) $ | {\Delta_K}^{\ast}(1)| = 2^n$ and $|{\Delta_K}^{\ast} (-1)| =
2^n |\Delta_K (-1)|$.\\
(3) ${\Delta_K}^{{\ast}{\ast}}(t) = 2^{2n} \Delta_K (t)$.\\

\end{cor}

\begin{cor}\label{cor:14.5}
If $\Delta_K(t)$ is stable, then ${\Delta_K}^{\ast} (t)$
is $c$-stable. If $\Delta_K(t)$ is $c$-stable, then ${\Delta_K}^{\ast} (t)$
is stable. Further, if $\Delta_K(t)$ is bi-stable, so is
${\Delta_K}^{\ast} (t)$.
\end{cor}

Now we proceed to a proof of Theorem \ref{thm:14.3}. 

Write
%
\begin{equation}\label{siki:14.3} 
f(t) = c_0 t^{2n} + c_1 t^{2n-1} + \cdots + c_{2n},
\end{equation}
\noindent
where $c_0 > 0$ and $c_j = c_{2n-j}, 0 \leq j \leq 2n$.\\
Let $\alpha_1, 1/ \alpha_1, \alpha_2,1/ \alpha_2, \cdots, \alpha_n,1/ \alpha_n$ 
be all the zeros of $f(t)$. Then we can write
%
\begin{equation}\label{siki:14.4}
f(t) = c_0 \prod_{j=1}^n (t- \alpha_j) (t -\frac{1}{\alpha_j})
\end{equation}

\begin{lemm}\label{lem:14.6}
Let $A_j = \alpha_j +\frac{1}{\alpha_j}, 1 \leq j \leq n$. Then
$\lambda = \prod_{j=1}^n A_j$ is a real number.
\end{lemm}

{\it Proof.} Since $\alpha_j$ is a zero of $f(t)$, so is $\overline{\alpha_j}$ 
and hence
$(\alpha_j + \frac{1}{\alpha_j} )(\overline{\alpha_j} + \frac{1}{\overline{\alpha_j}})$
is a real number. \qed

Later we will see that $\lambda c_0$ is an integer and show
the following;

\begin{equation}\label{siki:14.5}
f^{\ast}(t) = \lambda c_0 \prod_{j=1}^n (t- \varphi(\alpha_j) )(
t - \varphi(\frac{1}{\alpha_j} ))
\end{equation}

Now consider $F(t) = f(t)/c_0 = t^{2n} + \frac{c_1}{c_0} t^{2n-1} +
\cdots+ \frac{c_{2n}}{c_0}$. Then

\begin{equation}\label{siki:14.6}
F(t) = \prod_{j=1}^n (t-\alpha_j) (t- \frac{1}{\alpha_j})
= \prod_{j=1}^n (t^2 - A_j t + 1).
\end{equation}

For $0\le k\le n$, 
define $X_k = \sum_{j_1, \cdots, j_k} A_{j_1} A_{j_2} \cdots A_{j_k}$,
where the summation runs over all $ j_1, j_2, \cdots, j_k$
such that $1 \leq j_1 < j_2 < \cdots < j_k \leq n$. 
In particular, 
$X _0= 1$ and $X_n = \lambda$.
By expanding the right hand side of (\ref{siki:14.6}), 
we have the following system of relations.

Case I. $n=2m$.

For $k = 0, 1, 2, \cdots, m$,

\begin{equation}\label{siki:14.7}
\frac{c_{2k}}{c_0} = \binom{n}{k} X_0 + \binom{n-2}{k-1} X_2 +
\cdots + \binom{n-2k}{k-k} X_{2k},
\end{equation}

and for $k = 1,2, \cdots, m$,

\begin{equation}\label{siki:14.8}
 -\frac{c_{2k-1}}{c_0} = \binom{n-1}{k-1} X_1 +
\binom{n-3}{k-2} X_3 + \cdots
+ \binom{n-(2k-1)}{k-k}X_{2k-1}
\end{equation}

For simplicity, let $M_0$ and $N_0$ be, respectively, the coefficient
matrices of the system of relations of 
(\ref{siki:14.7}) and (\ref{siki:14.8}).
Namely, $M_0$ and $N_0$ are lower triangular integer matrices 
of sizes respectively $m+1$ and $m$, and each
with determinant $1$.

\begin{align}\label{siki:14.11} 
&M_0 (X_0, X_2, \cdots, X_{2m})^T = \frac{1}{c_0}(c_0, c_2,
\cdots, c_{2m})^T\ {\rm and}
\nonumber\\
&N_0 (X_1, X_3, \cdots, X_{2m-1})^T = -\frac{1}{c_0}(c_1,
c_3, \cdots, c_{2m-1})^T,
\end{align}

and hence

\begin{align}\label{siki:14.12}
&(X_0, X_2, \cdots, X_{2m})^T = \frac{M_0^{-1}}{c_0}(c_0, c_2,
\cdots, c_{2m})^T {\rm \  and}
\nonumber\\
&(X_1, X_3, \cdots, X_{2m-1})^T = -\frac{N_0^{-1}}{c_0}(c_1, c_3, \cdots,
c_{2m-1})^T.
\end{align}

Case II. $n=2m+1$.

The same argument shows the following


For $k = 0, 1, 2, \cdots, m$,
\begin{equation}\label{siki:14.13}
\frac{c_{2k}}{c_0} = \binom{n}{k} X_0 + \binom{n-2}{k-1} X_2 +
\cdots + \binom{n-2k}{k-k} X_{2k},
\end{equation}

and

\begin{equation}\label{siki:14.14}
-\frac{c_{2k+1}}{c_0} = \binom{n-1}{k} X_1 +
\binom{n-3}{k-1}X_3 + \cdots
+ \binom{n-(2k+1)}{k-k}X_{2k+1}
\end{equation}

Using coefficient matrices $M_1$ and $N_1$ of these systems of
relations, we can write
\begin{align}\label{siki:14.15}
& M_1 (X_0,X_2, \cdots X_{2m})^T = \frac{1}{c_0} (c_0, c_2,
\cdots, c_{2m})^T{\rm  and}
\nonumber\\
& N_1 (X_1,X_3, \cdots X_{2m+1})^T = -\frac{1}{c_0} (c_1, c_3, \cdots,
c_{2m+1})^T.
\end{align}

Here $M_1$ and $N_1$ are $(m+1) \times (m+1)$ 
lower triangular integer matrices with determinant $1$ and hence%

\begin{align}\label{siki:14.18}
&(X_0, X_2, \cdots, X_{2m})^T = \frac{M_1^{-1}}{c_0}(c_0, c_2,
\cdots, c_{2m})^T, {\rm  and}
\nonumber\\
&(X_1, X_3, \cdots, X_{2m+1})^T = -\frac{N_1^{-1}}{c_0}(c_1, c_3, \cdots,
c_{2m+1})^T
\end{align}

Now we study $f^{\ast}(t)$:

\begin{align}\label{siki:14.19}
f^{\ast}(t) 
&= \lambda c_0 \prod_{j=1}^n (t -
\varphi(\alpha_j))(t - \varphi(\frac{1}{\alpha_j} ))
\nonumber\\
&= \lambda c_0 \prod_{j=1}^n [t^2 - (\varphi(\alpha_j) + 
\varphi(\frac{1}{\alpha_j} )) t +\varphi(\alpha_j )\varphi(\frac{1}{\alpha_j})]
\nonumber\\
&= \lambda c_0 \prod_{j=1}^n (t^2 - \frac{4}{A_j} t + 1).
\end{align}

We write it as

\begin{equation}\label{siki:14.20}
f^{\ast}(t) = \lambda c_0 (d_0 t^{2n} + d_1 t^{2n-1} + \cdots
+ d_{2n}), d_0 =1.
\end{equation}

If we compare (\ref{siki:14.19}) with (\ref{siki:14.6}), 
we see immediately the following relations.

Case (I) $n=2m$.

For $k = 0, 1, 2, \cdots, m$,
\begin{align}\label{siki:14.21} 
\frac{d_{2k}}{d_0} = 
&\binom{n}{k} X_0 +
\binom{n-2}{k-1} \sum_{j_1,j_2} \frac{4^2}{A_{j_1} A_{j_2}} + 
\binom{n-4}{k-2} \sum_{j_1, \cdots, j_4} \frac{4^4}{A_{j_1}\cdots A_{j_4}}
+ \cdots \nonumber\\
&+ \binom{n-2k}{k-k} \sum_{j_1, \cdots, j_{2k}}
\frac{4^{2k}}{A_{j_1} \cdots A_{j_{2k}}},
\end{align}

and hence,

\begin{align}\label{siki:14.22} 
\lambda d_{2k} 
&= \binom{n}{k} X_n + \binom{n-2}{k-1} 4^2
X_{n-2} + \binom{n-4}{k-2} 4^4 X_{n-4}
+ \cdots 
\nonumber\\
&+ \binom{n-2k}{k-k} 4^{2k} X_{n-2k}
\end{align}

Similarly, for $k = 1,2, \cdots,m$,

\begin{align}\label{siki:14.23}
- \lambda d_{2k-1} 
&=\binom{n-1}{k-1} 4X_{n-1} + 
\binom{n-3}{k-2} 4^3 X_{n-3}
+ \cdots \nonumber\\
&+ \binom{n-(2k-1)}{k-k} 4^{2k-1} X_{n-(2k-1)}.
\end{align}

Let $P_{\ell}$ be a diagonal matrix of order $\ell$ of the form:
\begin{equation*}
P_{\ell} = \diag\{1, 4^2 ,4^4 , \cdots , 4^{2(\ell -1)}\}.
\end{equation*} 
Then
(\ref{siki:14.22}) and (\ref{siki:14.23}) can be written as

\begin{align}\label{siki:14.24}
&\lambda (d_0, d_2, \cdots,d_{2m})^T = M_0 P_{m+1} (X_n,
X_{n-2}, \cdots,X_0)^T, {\rm \ and}
\nonumber\\
- &\lambda (d_1, d_3, \cdots, d_{2m-1})^T = N_0 \widehat{P_m}
(X_{n-1}, X_{n-3}, \cdots, X_1)^T,
\end{align}
\noindent
where $\widehat{P_{\ell}} = 4 P_{\ell}$.

Case II. $n=2m+1$.

The same argument shows

\begin{align}\label{siki:14.25}
&\lambda (d_0, d_2, \cdots,d_{2m})^T = M_1 P_{m+1}
(X_{2m+1}, X_{2m-1}, \cdots,X_1)^T,{\rm \  and}
\nonumber\\
-& \lambda (d_1, d_3, \cdots, d_{2m+1})^T = N_1 \widehat{P_{m+1}}
(X_{2m}, X_{2m-2}, \cdots, X_0)^T.
\end{align}

Let $Q_{\ell}$ be an $\ell \times \ell$ matrix (that is the mirror of
the identity matrix), namely,
$Q_{\ell} = [ q_{i,j}]_{1 \leq i,j \leq \ell}$, where
$q_{i,j} = 1$, if $i+j=\ell+1$, and $q_{i,j} = 0$, otherwise.
Using $Q_{\ell}$, we have the final result.

Case (a) $n=2m$.
$(X_n, X_{n-2}, \cdots, X_0)^T = Q_{m+1}(X_0, X_2, \cdots, X_n)^T$
and
$(X_{n-1}, X_{n-3}, \cdots, X_1)^T = Q_m (X_1, X_3, \cdots, X_{n-1})^T$,
and hence, combining
(\ref{siki:14.24}) and (\ref{siki:14.12}), 
 we have
 %
\begin{align}\label{siki:14.26}
&\lambda c_0 (d_0, d_2, \cdots, d_{2m})^T = M_0 P_{m+1}
Q_{m+1}M_0^{-1} (c_0, c_2, \cdots, c_{2m})^T {\rm \ and}\\
 -& \lambda c_0 (d_1, d_3, \cdots, d_{2m-1})^T = N_0 
\widehat{P_m}Q_m (-N_0^{-1}) (c_1, c_3, \cdots, c_{2m-1})^T.
\end{align}

Case (b) $n=2m+1$.

Similarly, we have

\begin{align}\label{siki:14.28} 
\lambda c_0 (d_0, d_2, \cdots, d_{2m})^T 
&= M_1 P_{m+1}c_0
(X_{2m+1}, X_{2m-1}, \cdots, X_1)^T \nonumber\\
 &=
M_1 P_{m+1}Q_{m+1}(-N_1^{-1}) (c_1, c_3, \cdots, c_{2m+1})^T
{\rm \ and}\\
-\lambda c_0 (d_1, d_3, \cdots, d_{2m+1})^T 
&= N_1 
\widehat{P_{m+1}}c_0 (X_{2m}, X_{2m-2}, \cdots, X_0)^T 
\nonumber\\
 &=
N_1 \widehat{P_{m+1}}Q_{m+1}M_1^{-1} (c_0, c_2, \cdots, c_{2m})^T.
\end{align}

To be more precise, let $f(t) = \sum_{j=0}^{2n} c_j t^{2n-j}, c_0 > 0$,
and $f^{\ast}(t) = \sum_{j=0}^{2n}
a_j t^{2n-j}$. 
Then $a_j, 0 \leq j \leq n$, is obtained by the following
formulas.

Case (a) $n=2m$.
%
\begin{align}\label{siki:14.30}
&(a_0, a_2, \cdots,a_{2m})^T = M_0 P_{m+1} Q_{m+1} M_0^{-1}
(c_0, c_2, \cdots, c_{2m})^T{\rm \  and}
\nonumber\\
&(a_1, a_3, \cdots,a_{2m-1})^T = N_0 \widehat{P_m} Q_m N_0^{-1}
(c_1, c_3, \cdots, c_{2m-1})^T.
\end{align}

Case (b) $n=2m+1$.
%
\begin{align}\label{siki:14.31}
&(a_0, a_2, \cdots,a_{2m})^T = -M_1 P_{m+1} Q_{m+1} N_1^{-1}
(c_1, c_3, \cdots, c_{2m+1})^T{\rm \  and}\nonumber\\
&(a_1, a_3, \cdots,a_{2m+1})^T = -N_1 \widehat{P_{m+1}} Q_{m+1}
M_1^{-1} (c_0, c_2, \cdots, c_{2m})^T.
\end{align}
This proves Theorem \ref{thm:14.3} (1). 

Since all the matrices involved in the proof are
integer matrices, it follows that 
if $f(t)$ is an integer polynomial, then
so is $f^*(t)$.
This proves (2).

To prove (3) and (4), we compute $f(\pm 1)$ and $f^{\ast}(\pm 1)$.
Since $f(t) = c_0 \prod_{j=1}^n (t^2 - A_j t +1)$, it follows that $f(1)
= c_0 \prod_{j=1}^n (2-A_j)$ and $f(-1) = c_0 \prod_{j=1}^n (2+A_j)$.
Meanwhile, $f^{\ast}(t) = \lambda c_0 \prod_{j=1}^n (t^2 - \frac{4}{A_j} t
+1)$, and hence
$f^{\ast}(1)$ = $\lambda c_0 \prod_{j=1}^n (2 - \frac{4}{A_j})$ =$ \lambda
c_0 \prod_{j=1}^n \frac{1}{A_j}(2A_j - 4)$
= $\lambda c_0 \frac{1}{\lambda} 2^n \prod_{j=1}^n (A_j - 2) = 2^n (-1)^n
f(1)$.
Similarly,
$f^{\ast}(-1)$ = $\lambda c_0 \prod_{j=1}^n (2 + \frac{4}{A_j})$ = $c_0
\prod_{j=1}^n (2A_j +4)$ =
$c_0 2^n \prod_{j=1}^2 (A_j +2) = 2^n f(-1)$.
This proves (3) and (4).

To show (5), first we note that $\varphi^2 (z) = 1/z$. Thus the set of
the zeros of $f^{{\ast}{\ast}}(t)$ and that of $f(t)$ are identical.
Therefore, $f(t)$ divides $f^{{\ast}{\ast}}(t)$ or $f^{{\ast}{\ast}}(t)$
divides $f(t)$. 
However, $f^{{\ast}{\ast}}(1) =
2^n f^{\ast}(1)= 2^{2n}f(1)$ and hence, $f^{{\ast}{\ast}}(t) = 2^{2n}f(t)$.
Finally, the uniqueness is evident.
A proof of Theorem \ref{thm:14.3} is now completed.
\qed

%

\begin{ex}\label{ex:14.7}
Let $f(t) = \sum_{j=0}^{2n}c_j t^{2n-j}, c_0> 0$ and
$f^{\ast}(t) = \sum_{j=0}^{2n}a_j t^{2n-j}$. \\
(1) (i) Let $n = 1$ and $m = 0$. 
Then $M_1 = N_1 = P_1 = Q_1 = [1]$, and
hence $a_0 = - c_1$ and $a_1 = -4c_0$.
For example, if $f(t) = t^2 - 3t + 1$, then $f^{\ast}(t) = 3t^2 -4t +3$.\\
(ii) Let $n=3$ and $m=1$. 
Then $M_1= \left[\begin{array}{cc}
1&0\\
3&1
\end{array}\right]$
and $N_1 = \tbt{1}{0}{2}{1}$, 
and hence\\
$\tv{a_0}{a_2}=\tbt{2}{-1}{-10}{-3}\tv{c_1}{c_3}$
 and $\tv{a_1}{a_3}=\tbt{12}{-4}{-40}{-8}
 \tv{c_0}{c_2}$.
For example, if $f(t) = t^6 - t^5 + t^3 - t + 1$, then 
$f^{\ast}(t) =
-(3t^6 -12t^5 - 7t^4 + 40 t^3 -7t^2 -12t+3)$.\\
(2) (i) Let $n=2$ and $m=1$. Then 
$M_0 = \tbt{1}{0}{2}{1}$ 
and $N_0= [1]$.
Therefore, $\tv{a_0}{a_2}=\tbt{-2}{1}{12}{2}
\tv{c_0}{c_2}$ and
$(a_1)=(4c_1)$.
For example, if $f(t) = t^4 -t^3 +t^2 - t +1$, then $f^{\ast}
(t) = - (t^4 + 4t^3 -14t^2 + 4t + 1)$.\\
(ii) Let $n=4$ and $m=2$. 
Then $M_0 = 
\left[\begin{array}{rrr}
1&0&0\\
4&1&0\\
6&2&1
\end{array}\right]$,
and $N_0 = \tbt{1}{0}{3}{1}$, 
and hence, 
$\left[\begin{array}{r}
a_0\\
a_2\\
a_4
\end{array}\right]=
\left[\begin{array}{rrr}
2&-2&1\\
-56&8&4\\
140&20&6
\end{array}\right]
\left[\begin{array}{r}
c_0\\
c_2\\
c_4
\end{array}\right]
$
and
$\tv{a_1}{a_3}=\tbt{-12}{4}{28}{12}\tv{c_1}{c_3}$.
For example, if 
$f(t) = \sum_{j=0}^8 (-1)^j t^{8-j}$, then
$f^{\ast}(t) = t^8 +8t^7 - 44t^6 - 40t^5 + 166t^4 - 40t^3 - 44t^2 + 8t +1$.
\end{ex}

\begin{rem}\label{rem:14.8}
Even if $f(t)$ is monic, $f^{\ast}(t)$
is not necessarily monic. Furthermore, even if $\{c_0 , c_1 , c_2 ,
\cdots, c_{2n}\}$
alternates in sign, $\{a_0 , a_1, a_2, \cdots, a_{2n}\}$ may not alternate
in sign.
\end{rem}

\begin{qu}\label{qu:14.9}
Let $K$ be a $c$-stable knot and $K^*$ be a stable link
obtained by $\varphi$.
What can we say on $K^*$?
Does there exist a geometric way to construct $K^*$ from $K$?
\end{qu}

\section{Montesinos knots}\label{15}

In this section, we study the various stabilities of alternating
Montesinos knots or links.
It is not surprising to see that many Montesinos knots or links are
quasi-rational, and hence, their stability properties can be determined
by our method discussed earlier.

Now we begin with a well-known result of a characterization of alternating
Montesinos knots or links.
Let $K= M(e \mid \beta_1/ \alpha_1,\beta_2/ \alpha_2, \cdots,\beta_n/
\alpha_n )$ be a Montesinos
knot or link. We assume that $n \geq 3$.
Montesinos knots or links have two classes.\\
Class I (1) $\beta_i/ \alpha_i> 0 $ for any $i$, $1\leq i \leq n$, and $e
\geq 0$, or\\
(2) $\beta_i/ \alpha_i< 0$ for any $i$, $1 \leq i \leq n$, and $e \leq 0$.\\
Class (II) $0 < \sharp \{\beta_i/ \alpha_i > 0, 1 \leq i \leq n\}< n$
and $e = 0$.

The following proposition is well-known. (See \cite{lt}.) 

\begin{prop}\label{prop:15.1}
A Montesinos knot (or link) $K$ is alternating if and only if $K$
belongs to Class (I).
\end{prop}

Since we are interested in various stabilities of alternating knots or
links, we study the special classes of Montesinos knots or links
described in the following theorem.

\begin{thm}\label{thm:15.2}
Let $K = M(e \mid \beta_1 / \alpha_1, \beta_2 / \alpha_2 ,
\cdots, \beta_n / \alpha_n), n \geq3$,
be a Montesinos knot or link. We assume the following conditions:\\
(1) $\beta_i / \alpha_i > 0$ for any $i, 1 \leq i \leq n$,\\
(2) $e \ge 0$.\\
(3) At most one $\alpha_i $ , say $\alpha_1$, can be even.\\
(4) $\beta_i \equiv 0$ (mod 2), $1 \leq i \leq n$, unless $\alpha_i$
is even.\\
(5) Let $[2a_1^{(i)}, 2a_2^{(i)}, \cdots, 2a_{m_i}^{(i)}]$ be the even
continued fraction expansion of $\beta_i / \alpha_i, 1 \leq i \leq n$. 
For each $i$, the
sequence $\{2a_1^{(i)}, 2a_2^{(i)}, \cdots, 2a_{m_i}^{(i)}\}$ alternates
in sign. 
In particular, $2a_1^{(i)} > 0$ for $1 \leq i \leq n $.\\
Then we have the following conclusion.\\
Case 1. If all $\alpha_i , 1 \leq i \leq n$ are odd and e is odd, then $K$
is a special alternating knot, and hence $K$ is $c$-stable.\\
Case 2. If all $\alpha_i , 1 \leq i \leq n$, are odd and e is even,
then $K$ is a 2-component link and $K$ is inversive.
(A 2-component link $K$ is said to be inversive if the original (oriented)
link is $c$-stable (or stable), but if the orientation
of one component is revered, the resulting (oriented) link becomes a
stable (or $c$-stable) link.)\\
Case 3. If $\alpha_1$ is even and others are odd and $e$ is even, 
then $K$ is a knot and is stable. 
If $e$ is $2$ and all $|a_j^{(i)}| = 1$
for $1 \leq i \leq n$ and $1 \leq j \leq m_i $, then the maximal value
of the zero of $\Delta_K(t)$ is at least $n+1$.\\
Case 4. If $\alpha_1$ is even and others are odd, and $e$ is odd, 
then $K$ is a knot and $K$ is bi-stable.
\end{thm}
%
\begin{rem}\label{rem:15.3}
We should note that our cases do not contain all alternating
Montesinos knots. Since we are interested in stability of
the Alexander polynomial, the assumption (5) is crucial in the theorem.
Any knot or link in our list has some stability properties.
\end{rem}
Now, proofs of the first three cases are easy. For the first case, $K$ has a special alternating as in
Fig. 15.1 and hence, $K$ is $c$-stable. 

\jgi

For the second case, 
one orientation gives us a special alternating
diagram, and hence it is $c$-stable. 
If we reverse orientation of one component,
the diagram shows that $K$ is 
a quasi-rational links discussed in Section 6,
and hence, it is stable. See Figuress 15.2 and 15.3.

\jgii

\jgiii

For the third case, 
$K$ is also quasi-rational knot discussed in Section 6,
and hence, it is stable. 
Since the second statement can be proven
by applying the same argument 
used in the proof of Theorem \ref{thm:6.6} (1), 
we omit the details.
See Fig. 15.4.

\jgiv

The last case is the most complicated case. (See Fig.15.5.)

\jgv

Let $r_j = \frac{\beta_j}{\alpha_j}, 1\leq j\leq n$, and write:\\
$
r_1 =
[2a_1^{(1)}, -2a_2^{(1)}, \cdots, (-1)^{k-1}2a_k^{(1)}, 
\cdots,
2a_{2m_1+1}^{(1)}],\ {\rm where\ }
a_k^{(1)} > 0, 1 \leq k \leq 2m_1 +1.
$
Let $r_0 = [- e, r_1 ] = [- e, 2a_1^{(1)}, -2a_2^{(1)}, \cdots,
(-1)^{k-1}2a_k^{(1)}, \cdots, 2a_{2m_1+1}^{(1)}]$.
Since $e$ is odd, we see from the diagram that $K(r_0)$ is a special
alternating knot. 
Further, we see that $K$ is a Murasugi sum of 
$K(r_0)$
and the connected sum of remaining $(n-1)$ 2-bridge knots $K(r_2) \sharp K(r_3)
\sharp \cdots \sharp K(r_n)$. 
Since $K(r_0)$ is $c$-stable and
$K(r_2) \sharp K(r_3) \sharp \cdots \sharp K(r_n)$ is stable, it is not
surprising that a Murasugi sum of these knots is bi-stable, but a proof is not
immediate. The rest of this section will be devoted to a proof of this
case.\\
First, from the diagram, we see that a Seifert matrix $M$ of $K$ is a direct
sum of $M_j, j=0, 2,3, \cdots,n$, except the first column,
where $M_j$ is a Seifert matrix of $K(r_j)$ of twisted chain type. \\
\begin{center}
$M = \left[
\begin{tabular}{c|c|c|c|c|c}
\multicolumn{2}{c|}{}&\multicolumn{4}{c}{}\\[-2mm]
\multicolumn{2}{c|}{$M_0$} & \multicolumn{3}{c}{}\\[2mm]\hline
\smav& &$M_2$&  \multicolumn{3}{c}{} \\\hline
\smav& &     &$M_3$& \multicolumn{2}{c}{}\\\hline
$\vdots$& \multicolumn{3}{c|}{} &$\hspace{-1mm}\ddots\hspace{-1mm}$ & \\ \hline
\smav&\multicolumn{4}{c|}{}&
$M_n$
\end{tabular}
\right]
$
\end{center}

In fact, $M$ is of the
form above and has the following properties.
\begin{align}\label{siki:15.1} 
(1)\ &{\rm \ The\ diagonal\ entries\ of\ }
M_0{\rm \  is\ quite\ different\ from\ those\ of\ }
M_j, j \geq 2. {\rm \ They\  are\ }\nonumber\\
& \{-\frac{e+1}{2}, \underbrace{-1, \cdots, -1}_{2a_1^{(1)}-1},
 -(a_2^{(1)}+1), 
 \underbrace{-1, \dots, -1}_{2a_3^{(1)}-1}, 
 -(a_4^{(1)}+1) , 
 \underbrace{-1, \cdots, -1}_{2a_5^{(1)}-1},
\nonumber\\
&-(a_6^{(1)}+1), \cdots, -(a_{2m_1}^{(1)}+1), 
\underbrace{-1, \cdots, -1}_{2a_{m_1+1}^{(1)}-1}
\},
\nonumber\\
(2)\ &{\rm  Other\ \mbox{non-zero}\  entries\  of\  }
M_0 {\rm \  are\ those\ of\ the\ one\ line\ above\ the\ 
diagonal,}\nonumber\\
&{\rm all\ of\ which\ are\ } 1.
\nonumber\\
(3)\ &{\rm  The\ diagonal\ entries\ of\ }
M_j , j \geq 2,{\rm \ are\ } 
\{a_1^{(j)}, -a_2^{(j)}, \cdots, -a_{2m_j}^{(j)}\}.
\nonumber\\
(4)\ &{\rm  The\ size\ }
\rho_0{\rm  \  of\ } 
M_0 {\rm \  is\ } 
\rho_0 = 1 + \sum_{j \equiv 1 (2)}
(2a_j^{(1)}-1) + m_1 = \sum_{j \equiv 1 (2)} 2a_j^{(1)},
\nonumber\\
&{\rm \  while\ the\ size\ } 
\rho_j{\rm \  of\ }
M_j{\rm \  is\ } 
\rho_j =  2m_j, 2 \leq j \leq n.
\nonumber\\
(5)\ &{\rm  The\ extra\ 1\ on\ the\ first\ column\ of\ }
M{\rm \  appears\ only\ on\ the\ first\ row\ of\ }
\nonumber\\
&{\rm each\  block\ matrix\ }
M_2, M_3, \cdots, M_n,
{\rm \  namely\, 1\ appears\
on\ the\ } 
\nonumber\\
&(\rho_0+ 1, 1)-, (\rho_0 + \rho_2 + 1, 1)-, \cdots,
(\sum_{j=0, j \ne 1}^n \rho_j + 1, 1)\mbox{\rm -entries\  of\ }
M.
\end{align}
Now to study the Alexander polynomial of K, we consider
the determinant of  $tM -M^T$ that is of the form

$\det
\left[
\begin{array}{l|c|c|c|c|c|}
tM_0-M_0^T&
-1\, 0 \cdots 0&
\multicolumn{1}{c}{\cdots}&
\multicolumn{1}{c}{\cdots}&
\multicolumn{1}{c}{\cdots}&
\multicolumn{1}{c}{\cdots}\\\cline{1-2}
\smavt&tM_2-M_2^T&
\multicolumn{4}{c}{\hugesymbol{O} }\\
\cline{1-3}
\vdots& \multicolumn{1}{c|}{}&\multicolumn{1}{r|}{\smallsymbol{\ddots}}&\multicolumn{3}{c}{}\\[1mm]
\cline{3-4}
\vdots&\multicolumn{2}{c|}{}&tM_2-M_2^T&
\multicolumn{2}{c}{}\\\cline{4-5}
\vdots&\multicolumn{3}{c|}{}&{\ddots}& \multicolumn{1}{c}{}\\\cline{5-6}
\vdots&\multicolumn{4}{c|}{\hugesymbol{O}}&
\multicolumn{1}{c}{tM_p-M_p^T}\\
\end{array}
\right]
$

By expanding it along the first row and then the first column, we can show
easily that
\begin{align}\label{siki:15.2}
\Delta_K(t) &= \det (Mt -M^T )\nonumber\\
&= 
\prod_{j=0, j\ne 1}^n\det (tM_j - M_j^T)
\nonumber\\
&\ \ + t \det \widehat{M_0} \sum_{j=2}^n  \det (tM_2 - M_2^T ) \cdots \det
\widehat{M_j} \cdots
\det(tM_n -M_n^T) \nonumber\\
&= 
\Delta_{K(r_0)}(t) \Delta_{K(r_2)}(t) \cdots \Delta_{K(r_n)}(t)
\nonumber\\
&\ \ + t \det \widehat{M_0}[\sum_{j=2}^{n}  \Delta_{K(r_2)}(t) \cdots
\Delta_{K(\widehat{r_j})} (t) \cdots \Delta_{K(r_n)}(t)],
\end{align}
\noindent
where $\widehat{M_j}, j=0,2,3, \cdots, n$, is the matrix obtained from
$tM_j - M_j^T$
by deleting the first row and column and $\widehat{r_j} = [-2a_2^{(j)},
2a_3^{(j)}, \cdots, -2a_{2m_j}^{(j)}], j \geq 2$.

For simplicity, we denote
\begin{equation*}
f_0 = \det (tM_0 - M_0^T){\rm \  and\ } 
\widehat{f_0} = - \det \widehat{M_0}, 
\end{equation*}
and for $j = 2,3, \cdots, n$,
\begin{equation*}
f_j = (-1)^{m_j} \det (tM_j - M_j^T){\rm \  and\ } 
\widehat{f_j} = (-1)^{m_j} \det \widehat{M_j}.
\end{equation*}
Since the leading coefficients of $f_0, \widehat{f_0}, f_j$ and
$\widehat{f_j}$ are
all positive, these polynomials are normalizations of 
$\Delta_{K(r_0)}(t),
\det \widehat{M_0}, \Delta_{K(r_j)}(t)$ and 
$\Delta_{K(\widehat{r_j})}(t)$,
respectively.
Using these polynomials, we rewrite (\ref{siki:15.2}) 
as follows.
%
\begin{align}\label{siki:15.3}
\Delta_{K}(t) 
&= f_0 (t) f_2 (t) \cdots f_n (t) (-1)^{m_2 + \cdots +m_n}\nonumber\\
&\ \ +  t \widehat{f_0} (t) \widehat{f_2} (t) f_3(t) \cdots f_n (t)
(-1) (-1)^{m_2 + \cdots +m_n}
\nonumber\\
& \ \  t \widehat{f_0} (t) f_2(t) \widehat{f_3} (t) \cdots f_n (t)
(-1) (-1)^{m_2 + \cdots +m_n}
+ \cdots
\nonumber\\
& \ \  t \widehat{f_0} (t) f_2 (t) f_3(t) \cdots f_{n-1}(t)
\widehat{f_n} (t) (-1) (-1)^{m_2 + \cdots +m_n}.
\end{align}
Therefore, the normalization $F$ of $\Delta_{K}(t)$ is
\begin{align*}
F &=  f_0 (t) f_2(t) \cdots f_n (t)
- t \widehat{f_0}(t) \{  \widehat{f_2}(t) f_3 (t)
\cdots f_n (t)\\
& \ \ +  f_2(t) \widehat{f_3}(t) f_4 (t) \cdots f_n(t)
+ \cdots +  f_2 (t) f_3(t) \cdots f_{n-1}(t) \widehat{f_n}(t) \}.
\end{align*}
Further, $\widehat{f_0}, \widehat{f_2}, \cdots, \widehat{f_n}$ are
Alexander polynomials of links
and hence, they are divisible by $t-1$.
Let $f_0^{\ast} = \frac{\widehat{f_0}}{t-1}$ and $f_j^{\ast}
= \frac{\widehat{f_j}}{t-1}, 2 \leq j \leq n$. 
Then 
$f_0^{\ast}$ and $f_j^{\ast}$ are respectively 
reciprocal and $f_j^{\ast} (t)$ is also stable.
We can write
\begin{equation*}
F = f_0 f_2 \cdots f_n - t(t-1)^2 f_0^{\ast} \{ \sum_{j=2}^n  f_2
\cdots f_j^{\ast} \cdots f_n \}.
\end{equation*}
Let $F_1 = f_0 f_2 \cdots f_n$ and $F_2 = t(t-1)^2 f_0^{\ast} \{
{\displaystyle \sum_{j=2}^n}  f_2 \cdots f_j^{\ast} \cdots f_n \}$,
and further, for $2\le k \le n$, $F_{2,k}
=t(t-1)^2f_0^{*}f_2\cdots f_k^{*}\cdots f_n$.
Since $\deg F=\deg{\displaystyle \sum_{j=0,j\neq1}^{n}} f_j$ 
and $|\sigma(K)|=\rho_0=\deg f_0$,
it suffices to show that $F$ has at least $2q=\sum_{j=2}^{n} \deg f_j$ real zeros.
The proof will be divided into two parts.
In the first part we consider the case where no zeros of $f_2 (t),
\cdots, f_n (t)$ are in common, i.e., $f_2 (t) \cdots f_n (t)$ is simple.
Note that each $f_j (t), j \ne 0$, has no multiple zeros.
In the second part, we consider the case where these Alexander polynomials
have some zeros in common.

(1) Case (1).
$f_2 f_3 \cdots f_n$ is simple.

We will show that two curves $y_1=F_1$ and $y_2=F_2$ intersect at least 
$q =\sum_{j= 2}^n
\rho_j/2$ points in $[0,1]$.
Let $\gamma_1< \gamma_2< \cdots< \gamma_q$ be all the (real) zeros of
$f_2 f_3 \cdots f_n$ in $[0,1]$. 
Note that $\deg F_1=\deg F_2+1$ and $F_1(1)\neq 0$.

First we see that
(1)$F_1(0) > 0$ and
(2)$F_2(0) =0$ and $F_2^{\prime}(0) > 0$.\\
In fact, (1) $F_1(0) = \prod_{j \ne 1}f_j (0) > 0$ , since the leading
coefficient of $f_j, j \ne 1$, is positive and $f_j$ is reciprocal. 
(2) follows, since $f_j^{\ast}(0) > 0$ for $j \ne 1$.

Now, suppose $\gamma_1$ is the zeros of $f_k$. 
Then $f_k (\gamma_1) = 0$,
but $f_k^{\ast}(\gamma_1) \ne 0$, 
since $f_k f_k^{\ast}$ is simple. 
Therefore,
$F_{2,k}(\gamma_1) \ne 0$, but $F_{2,j}(\gamma_1) = 0$ for $j \ne k$.
Further, since $F_{2,k}^{\prime}(0) > 0$,
it follows that $F_{2,k}(\gamma_1) > 0$, and hence
$F_1$ and $F_2$ intersect in $[0, \gamma_1]$.

Next, we prove inductively the following lemma.

%
\begin{lemm}\label{lem:15.4}
(1) $F_1 (t) \leq 0$ in $[\gamma_{2k-1}, \gamma_{2k}], 
1 \leq k \leq [\frac{q}{2}]$,
and $F_1 (t) \geq 0$ in $[\gamma_{2k}, \gamma_{2k+1}], 
1 \leq k \leq [\frac{q-1}{2}]$.\\
(2) $F_2 (\gamma_{2k+1}) > 0, 0 \leq k \leq [\frac{q-1}{2}]$ 
and $F_2(\gamma_{2k}) < 0, 1 \leq k \leq [\frac{q}{2}]$.\\
Therefore, $F_1$ and $F_2$ intersect in $[\gamma_i , \gamma_{i+1}], 1 \leq i \leq
q-1$, and
hence $F_1$ and $F_2$ intersect at least $q$ points in $[0, 1]$.
\end{lemm}
Note that $F_1$ and $F_2$ do not intersect in $[\gamma_q, \gamma_{q+1}] \ni 1$.

{\it Proof.} Since (1) is obvious, we prove only (2) by induction on $k$.
Since we already showed that $F_2 (\gamma_1) > 0$, we consider
$F_2(\gamma_2)$.
Suppose $\gamma_2$ is the zero of $f_s (t)$. 
Then $F_{2,j}(\gamma_2) = 0$ for $j \ne s$.
We know that $\gamma_1$ is the zero of $f_k$.\\
(1) If $s=k$, 
then $f_k^{\ast}(\gamma_2) \ne 0$ and $F_{2,k}(\gamma_1) >0$.
Since $\gamma_1$ and $\gamma_2$ are the zeros of $f_k (=f_s)$, 
and $f_k$ and $(t-1) f_k^{\ast}$
are interlaced, we see that $f_k^{\ast}$ has the zero 
$\beta_0$ in $[\gamma_1, \gamma_2]$
and hence $F_{2,k}$ crosses the $t$-axis at $\beta_0$. 
Therefore,
$F_{2,k}(\gamma_2) < 0$.
Since $F_{2,j}(\gamma_2) = 0$ for $j \ne k$, it follows that 
$F_2(\gamma_2) < 0$.\\
(2) If $s \ne k$, then $F_{2,j}(\gamma_2)= 0$ for $j \ne s$ and further
$F_{2,s}(\gamma_2) < 0$,
since $F_{2,s}(\gamma_2) \ne 0$ and $F_{2,s}^{\prime} (\gamma_1) < 0$.
Therefore, $F_2 (\gamma_2) < 0$.

Now consider $F_2 (\gamma_m), 1 \leq m \leq q$.
Suppose $\gamma_m$ is the zero of $f_p$.

Case (1) $\gamma_m$ is the smallest zeros of $f_p$.
Then $f_p$ is not 0 at $\gamma_1, \gamma_2, \cdots, \gamma_{m-1}$ and it is
obvious that
(i) if $m$ is even, then $F_{2,p}(\gamma_m) < 0$, 
and $F_{2,\ell}(\gamma_m) = 0$ for $\ell \ne p$ 
and hence $F_2(\gamma_m) < 0$, and
(ii) if $m$ is odd, 
then $F_{2,p}(\gamma_m) > 0$, and
$F_{2,\ell}(\gamma_m) = 0$ for $\ell \ne p$ 
and hence $F_{2}(\gamma_m) > 0$.

Case (2) There exists $\gamma_h, h < m$, that is the closest zero of
$f_p$ to $\gamma_m$, i.e., $f_p(\gamma_h) = 0$, 
but $f_p (\gamma_{h+1}) \ne 0,\cdots, f_p(\gamma_{m-1})\ne 0$.
If $h$ is even, by induction assumption, $F_{2,p}(\gamma_h) < 0$.
Further, if $m$ is even, there are an odd number of zeros 
$\gamma_{h+1},\cdots, \gamma_{m-1}$
between $\gamma_h$ and $\gamma_m,$ and $F_{2,p}$ crosses the $t$-axis 
at these points.
However, since $f_p$ and $(t-1)f_p^*$
are interlaced, there is exactly one zero of $f_p^*$, say
$\beta_1$, in $[\gamma_h, \gamma_m],$
i.e., 
$\gamma_h < \beta_1 < \gamma_m$, and $F_{2,p}$ must cross the $t$-axis
at $\beta_1$ as well. 
Thus $F_{2,p}(\gamma_m) < 0$. 
Since $F_{2,\ell}(\gamma_m)= 0$ for $\ell \ne p$, 
we have $F_2 (\gamma_m)< 0$.

The same argument works for other cases where 
(a) $h$ is even and $m$ is odd,
and (b) $h$ is odd and $m$ is even or odd.
This proves Lemma \ref{lem:15.4} and 
Theorem \ref{thm:15.2} for non-multiple zero case.

Case (II) $f_2 f_3 \cdots f_n$ is not simple.

Let $\gamma_1 < \gamma_2 < \cdots < \gamma_{2d}$ be
all distinct real zeros of $f_2 f_3 \cdots f_n$. 
Let $p_k$ be the multiplicity of 
$\gamma_k, 1 \leq k \leq 2d$ 
and $a_j$ the leading coefficient of $f_j$. 
Note $a_j > 0, 2\leq j \leq n$. 
Since $f_2 f_3 \cdots f_n = a_2 a_3 \cdots a_n 
\prod_{k=1}^{2d} (t -\gamma_k)^{p_k}$, we can write
\begin{align*}
F = &f_0 f_2 f_3 \cdots f_n - t(t-1)^2 f_0^{\ast}
\sum_{j=2}^{n} f_2 \cdots f_j^{\ast}\cdots f_n\\
= & a_2 a_3 \cdots a_n f_0 \prod_{k=1}^{2d} (t- \gamma_k)^{p_k}
- t(t-1)^2 f_0^{\ast} \sum_{j=2}^n a_2 a_3 \cdots a_n
f_j^{\ast}[\prod_{k=1}^{2d}(t- \gamma_k)^{p_k} / f_j]\\
=& a_2 a_3 \cdots a_n \prod_{k=1}^{2d}(t- \gamma_k)^{p_k-1}
[f_0 \prod_{k=1}^{2d}(t- \gamma_k) - t(t-1)^2 f_0^{\ast}
\sum_{j=2}^n f_j^{\ast} \prod_{k=1}^{2d}(t- \gamma_k) / f_j].
\end{align*}

Since $f_j$ is simple and $\{\gamma_k\}$ is the set of all
distinct zeros of $f_2 f_3 \cdots f_n$, we see that 
$\prod_{k=1}^{2d}(t-\gamma_k)/ f_j$
is a (real) polynomial that is denoted by $g_j$.
Therefore to prove Theorem \ref{thm:15.2}, 
it suffices to show that
$G = f_0 \prod_{k=1}^{2d}(t- \gamma_k)
- t(t-1)^2 f_0^{\ast} \sum_{j=2}^n f_j^{\ast} g_j $
has $2d$ real zeros, or equivalently, two curves 
$y_1 = G_1 (t) = f_0 \prod_{k=1}^{2d}(t- \gamma_k)$
and $y_2 = G_2 (t) = t(t-1)^2 f_0^{\ast}\sum_{j=1}^n f_j^{\ast} g_j$ have
$d$ points of
intersection in $[0,1]$.

Let $G_{2,j} = t(t-1)^2 f_0^{\ast} f_j^{\ast} g_j$ and
hence $G_2 = \sum_{j=1}^n G_{2,j}$.
Note that $f_0$ and $f_0^{\ast}$ do not have any real zeros.
Suppose that $\gamma_j$ is the zero of $f_{j_1}, f_{j_2}, \cdots,
f_{j_{p_j}}$, i.e., $\gamma_j$ is the zeros of
$f_2 \cdots f_n$ of multiplicity $p_j$. 
Since the zeros of $g_k$ consist
of all real zeros of $G_1$ except
those of $f_k$, 
it follows that $G_{2,k}(\gamma_j)\neq 0$ if and only if
$\gamma_j$ is the zero of $f_k$.
Therefore, 
as it was proved in Lemma \ref{lem:15.4}, 
we can prove inductively that 
to each $\lambda, 1 \leq \lambda\leq p_j,
G_{2, j_{\lambda}}(\gamma_j)> 0$ or $< 0$, according as $j$ is odd or even.
Since $G_{2, \ell}(\gamma_j) = 0$,
if $\ell \ne j_{\lambda}, 1 \leq \lambda \leq p_j$,
it follows that  $G_2 (\gamma_j) = \sum_{\lambda = 1}^{p_j} G_{2,j_{\lambda}}
(\gamma_j)> 0$ or $< 0$, according as $j$ is odd or even.
Therefore, $y_1= G_1$ and
$y_2 = G_2$ intersect in $[\gamma_{j-1},\gamma_j]$ and
finally, two curves $y_1 = G_1$ and $y_2 = G_2$ intersect
at least $d$ points in $[0,1]$

A proof of Theorem \ref{thm:15.2} is now complete.

\qed

\section{Multivariate stable link polynomials}\label{16}

In this section, we study the stability of the 2-variable Alexander
polynomial $\Delta_{K(r)}(x,y)$ of a 2-bridge link $K(r)$,
$r= \beta/ 2\alpha, 0 < \beta < 2\alpha$.
Let $\H$ be the upper half-plane of $\CC$.
If the reduced Alexander polynomial of $K(r)$, i.e.,
$\Delta_{K(r)}(t) = (t-1)\Delta_{K(r)}(t,t)$ is not real stable,
then $\Delta_{K(r)}(x,y)$ is not $\H$-stable. 
Therefore, we only need to consider a real stable link.
If the continued fraction expansion of $r$ gives an alternating sequence, 
then $K(r)$ is real stable. 
However, for such links, $\H$-stability problem is completely solved
(Proposition \ref{prop:16.5}).
Therefore, we should study exceptional stable links.
In this section, we solve the $\H$-stability problem 
for the simplest
exceptional stable 2-bridge links.

Now we begin with a definition.

\begin{dfn}\label{dfn:16.1}
For a positive integer $n$, we define 
$G_n = \frac{x^n - y^n}{x-y}$ and 
$G_{-n}= \frac{-1}{(xy)^n} G_n$.
In particular, we define $G_0 = 0$.
\end{dfn}

It is easy to see that $G_n$ is $\H$-stable if and only if $|n| \leq 2$.

The following proposition is well-known.

\begin{prop}\label{prop:16.2}
Let $r = [2a_1, 2b_1, 2a_2,2b_2, \cdots, 2a_n, 2b_n,
2a_{n+1}]$. 
Then $\Delta_{K(r)}(x, y)$ is give by
\begin{equation}\label{siki:16.1}
\Delta_{K(r)}(x,y) = \sum_{0 \leq m \leq n} b_{j_1}b_{j_2}\cdots
b_{j_m} (x-1)^m (y-1)^m G_{\mu_1}G_{\mu_2} \cdots G_{\mu_{m+1}},
\end{equation}
where the summation is taken over all indices 
$j_k$ such that $1 \leq j_1 < j_2 < \cdots < j_m \leq n$ and 
$\mu_1 = a_1 + a_2 + \cdots + a_{j_1}$,
$\mu_2 = a_{j_1+1} + \cdots + a_{j_2}, 
\cdots, \mu_k =
a_{j_{k-1}+1} + a_{j_{k-1}+2} + \cdots + a_{j_k}$, $\cdots$,
$\mu_{m+1} = a_{j_m+1} + a_{j_m+2} + \cdots + a_n$.
\end{prop}

We should note that (\ref{siki:16.1})
is slightly different from the original formula
given in \cite{kane},  
since the orientation of one component of $K(r)$ in \cite{kane} 
is different form ours.

\begin{ex}\label{ex:16.3}
(1) If $r = [2a_1]$, then $\Delta_{K(r)}(x,y) = G_{a_1}$\\
(2) If $r = [2a_1, 2b_1, 2a_2]$, then $\Delta_{K(r)}(x,y) = b_1 (x-1)
(y-1) G_{a_1}G_{a_2} + G_{a_1+a_2}$\\
(3) If $r = [2a_1 , 2b_1, 2a_2, 2b_2, 2a_3]$, 
then \begin{align*}
\Delta_{K(r)}(x,y) &=
b_1 b_2 (x-1)^2 (y-1)^2 G_{a_1} G_{a_2} G_{a_3}\\
&\ \ +
(x-1)(y-1) \{b_1 G_{a_1}G_{a_2 +a_3} + b_2 G_{a_1+a_2} G_{a_3}\} + G_{a_1
+a_2+a_3}
\end{align*}
\end{ex}

Now the following simple proposition gives us a 
strong restriction on
$\H$-stability for 2-component links.

\begin{prop}\label{prop:16.4}
Let $K = K_1 \cup K_2$ be a 2-component link.
Suppose that $\Delta_K (x,y)$ is not a constant. 
(If $\Delta_K (x,y)$ is a constant, 
then $K$ is always $\H$-stable.)
If $K$ is $\H$-stable, then both $K_1$ and $K_2$ are stable and further,
$|lk(K_1, K_2)| \leq 2$.
\end{prop}

{\it Proof.} 
By Theorem \ref{thm:3.1}(a), 
if $\Delta_K (x,y)$ is $\H$-stable, then 
$\Delta_K (x,1)$ and $\Delta_K (1,y)$ are $\H$-stable.
Further, $\Delta_K (x,1) = \frac{x^{\ell} -1}{x-1} \Delta_{K_1} (x)$ 
and $\Delta_{K}(1,y)=\frac{y^\ell -1}{y-1}\Delta_{K_2}(y_2)$
\cite{tor}, 
where $\ell = {\it lk}(K_1, K_2)$, and hence,
$\Delta_{K_1}(x), \Delta_{K_2} (y)$ and $\frac{x^{\ell}-1}{x -1}$ are
$\H$-stable. The proposition
follows immediately.
\qed

\begin{prop}\label{prop:16.5}
Let $r = [2a_1, -2b_1, 2a_2, -2b_2 \cdots, 2a_n,
-2b_n, 2a_{n+1}] , a_j, b_j > 0$. 
Then $K(r)$ is $\H$-stable if and only
if
(1) $n = 0$ and $a_1 = 1$ or $2$, or
(2) $n=1$ and $a_1 = a_2 = 1$.
\end{prop}

{\it Proof.}
(a) Suppose $K(r)$ is $\H$-stable. 
Then $|\ell| = |{\it lk}(K_1, K_2)| \leq 2$. 
Since $\ell = \sum_{j =1}^{n+1}a_j$, we have (1) or (2).
(b) Suppose (1) or (2) holds. 
If (1) holds, then $\Delta_K (x,y) = 1$ or
$x+y$, and both are $\H$-stable. 
Suppose (2) holds.
Then $r =[2, -2b_1, 2]$ and 
$\Delta_K (x,y)$ = $-b_1 (x-1)(y-1) +(x+y)$
= $-b_1 + (b_1 +1) (x+y) - b_1 xy$.
Since $\Delta_K (x,y)$ is multi-affine, 
it is $\H$-stable if (and only if )
$\det \tbt{-b_1}{b_1+1}{b_1+1}{- b_1} < 0$
(Example \ref{ex:3.8}).
Since $b_1 > 0$, obviously, the determinant is negative. 
This proves Proposition \ref{prop:16.5}. 
\qed

We now consider exceptional stable links.

 \begin{prop}\label{prop:16.6}
 Let $r = [2a,2b, -2c], a,b,c > 0$ and $a > c$.
(1)Suppose $a=c$. Then $K(r)$ is $\H$-stable if and only if 
$a = 1$ or $2$.
(2) Suppose $a>c$. Then $K(r)$ is not $\H$-stable, unless 
$(a,c) = (2,1)$.
\end{prop}

{\it Proof.} 
(1) Suppose $a=c$. Then $\Delta_K (x,y)= b(x-1)(y-1)G_a^2$.
If $a \geq 3$, then $G_a$ is not $\H$-stable and 
hence $K(r)$ is not $\H$-stable.
However, if $a = c = 1$ or $2$, then each factor is $\H$-stable
and hence $K(r)$ is $\H$-stable.
(2) Suppose $a > c$. Since $lk(K_1,K_2) = a-c$, $c = a-1$ or $c = a-2$.
From (\ref{siki:16.1}), 
we have
$\Delta_{K(r)} (x,y)$ = $b (x-1) (y-1) G_a G_{-c}+ G_{a-c}
= b (x-1)(y-1) \frac{x^a - y^a}{x-y} \frac{x^c -
y^c}{x-y} \frac{-1}{(xy)^c} + \frac{x^{a-c} -
y^{a-c}}{x-y}$.
Let $F(x,y) = (xy)^c \Delta_{K(r)} (x,y)
=-b (x-1)(y-1) \frac{x^a - y^a}{x-y} \frac{x^c - y^c}{x-y} +
\frac{x^{a-c} - y^{a-c}}{x-y}(xy)^c$.
If $\Delta_{K(r)} (x,y)$ is $\H$-stable, 
so is $\Delta_{K(r)} (x,-1)$ by
Theorem \ref{thm:3.1} (a) 
and hence $f(x) = F(x, -1) = 2b (x-1) \frac{x^a - (-1)^a}{x+1} \frac{x^c -
(-1)^c}{x+1} + \frac{x^{a-c} - (-1)^{a-c}}{x+1}(-1)^c x^c$
must be $\H$-stable.

Case(I) $c=a-1$ and $a \geq 3$. Then
$f(x) =2b (x-1) \frac{x^a - (-1)^a}{x+1} \frac{x^{a-1}-(-1)^{a-1}}{x+1} +
(-1)^{a-1} x^{a-1}$.
We show that if $a \geq 3$, $f(x)$ is not (real) stable.
Since $f(x)$ is reciprocal, consider the modified polynomial $g(z)$ of
$f(x)$.

Case(a) $a$ is even, say $2m \geq 4$.

Then $f(x) =2b (x-1) \frac{x^{2m} - 1}{x+1} 
\frac{x^{2m-1} +1}{x+1} -x^{2m-1}$
and hence $g(z)$ is written as
$g(z) = 2b g_1 (z) g_2(z) -1$,
where $g_1(z)$ is the modification of $f_1(x) = (x-1)
\frac{x^{2m}-1}{x+1}$ and
$g_2(z)$ is that of $f_2(x) = \frac{x^{2m-1}+1}{x+1}$.
Since both $f_1$ and $f_2$ are $c$-stable, all
zeros of $g_1 (z)$ and $g_2 (z)$ are in $[-2,2]$ and further, 
since $m \geq 2$, at least one zero of 
$g_1 (z)$ and $g_2 (z)$ are in $(-2, 2)$.
Therefore it is impossible that all points of intersection of
two curves $z_1 = g_1(z) g_2(z)$ and $z_2= \frac{1}{2b}, b > 0$ are
outside of $(-2, 2)$, and hence $f(x)$ cannot be stable.

Case(b) $a$ is odd, say $2m+1, m \geq 1$.

Then $f(x) = 2b (x-1) \frac{x^{2m+1} + 1}{x+1}
\frac{x^{2m} -1}{x+1} + x^{2m}$ and the
modification is: 
$g (z) = 2bg_1(z) g_2 (z) +1$, where $g_1 (z)$ is the
modification of $\frac{x^{2m+1}+1}{x+1}$ and
$g_2(z)$ is that of $(x-1) \frac{x^{2m} - 1}{x+1}$. 
Therefore, all the
zeros of $g_1 (z)$ are in $(-2,2)$. 
As is proved in Case (a),
$f(x)$ cannot be stable.

Case(II). $c=a-2$ and $a \geq 3$. Then
\begin{align*}
\Delta_{K(r)}(x,y)
   &= b (x-1) (y-1) G_a G_{-(a-2)} + G_2\\
   &= b (x-1)(y-1) \frac{x^{a} - y^{a}}{x-y} 
   \frac{x^{a-2}- y^{a-2}}{x-y} \frac{-1}{(xy)^{a-2}} + x + y 
   \end{align*}

and 
\begin{align*}
f(x) &= (-1)^{a-2} x^{a-2} \Delta_{K(r)}(x, -1)\\
& = 2b (x-1) \frac{x^a -
(-1)^a}{x+1} \frac{x^{a-2} -
(-1)^{a-2}}{x+1} + (x -1) (-1)^{a-2} x^{a-2}.
\end{align*}

Case(a) $a$ is even, say $2m \geq 4$.

Then $f(x) = 2b (x-1) \frac{x^{2m} - 1}{x+1} \frac{x^{2m-2} - 1}{x+1} +
(x-1) x^{2m-2}$
and hence
$h(x)= \frac{f(x)}{x-1} = 2b (x-1)^2 \frac{x^{2m} - 1}{x^2-1}
\frac{x^{2m-2} - 1}{x^2-1} + x^{2m-2}$ is reciprocal.
The modification $\lambda (z)$ of $h(x)$ is
$\lambda (z) = 2b(z-2) \lambda_1(z) \lambda_2(z) +1$,
where $\lambda_1(z)$ and $\lambda_2 (z)$ are the modifications of
$\frac{x^{2m}-1}{x^2-1}$ and $\frac{x^{2m-2}-1}{x^2-1}$,
respectively. Since all the zeros of $\lambda_1(z)$ and $\lambda_2(z)$ are
in $(-2,2)$, $h (x)$ cannot be real stable.

Case(b) $a$ is odd, say $2m+1$.

The same argument works to show that $h(x)$ is not real
stable and hence $\Delta_{K(r)}(x,y)$ is not $\H$-stable. 
\qed

If $(a,c) = (2,1)$, we have the following proposition.

\begin{prop}\label{prop:16.7}
 Let $r = [4,2,-2]$. Then $K(r)$ is $\H$-stable.
\end{prop}

{\it Proof.} 
$xy\Delta_{K(r)} (x,y) = -(x-1)(y-1) G_2 G_1 +xy G_1 =
-(x-1)(y-1) (x+y) +xy = -(xy-x-y)(x+y-1)$.
Each factor is $\H$-stable by Example \ref{ex:3.8} 
and hence $K(r)$ is $\H$-stable. 
\qed

\begin{qu}\label{qu:16.8}
For $r=[4,2b,-2]$ and $b \geq 2$, is $K(r)$ $\H$-stable ?
\end{qu}

Finally, we prove that the 
2-variable Alexander polynomial of a 2-bridge
link has the same property as Theorem \ref{thm:8x3} (2). 
Using this, we can  systematically obtain exceptional 
$\H$-stable 2-component 2-bridge links.

\begin{thm}\label{thm:16.9}
Let $s = [2a_1, 2b_1, 2a_2,2b_2, \cdots, 2a_n,2b_n,
2a_{n+1}], a_j \neq 0 \neq b_j, 1 \leq j \leq n+1$.
Let $r = [s, 2k, -s^{-1}], k \neq 0$. 
Then,
\begin{equation}\label{siki:16.2}
\Delta_{K(r)}(x,y) = k(x-1)(y-1) [\Delta_{K(s)} (x,y)]^2,
\end{equation}
and hence, $K(r)$ is $\H$-stable if and only if $K(s)$ is $\H$-stable.
\end{thm}

{\it Proof.} 
Consider a sequence of integers\\
$A = \{a_1, b_1, a_2, b_2, \cdots, a_n, b_n,
\cdots,
a_{2n+1}, b_{2n+1}, a_{2n+2}\}$.
Take an ordered subset $C=\{b_1, b_2, \cdots, b_{2n+1}\}$. 
Let
$\widehat{C}$ be the set of all ordered subset of $C$, i.e.,
$\widehat{C} \ni U = \{b_{j_1}, b_{j_2},
\cdots, b_{j_k}\}$, where 
$1 \leq j_1 < j_2 < \cdots < j_k \leq 2n+1$.
To each set $U$ in $C$, 
we define a mapping $\rho_{2n+1} : \widehat{C}
\longrightarrow \ZZ[x^{\pm 1}, y^{\pm 1}]$ as
follows.

\begin{equation}\label{siki:16.n1}
\rho_{2n+1} (U) = b_{j_1}b_{j_2}\cdots
b_{j_k} (x-1)^k (y-1)^k G_{\mu_1}G_{\mu_2}\cdots G_{\mu_{k+1}},
\end{equation}
where 
$\mu_1 = a_1 + a_2 + \cdots + 
a_{j_1}$, $\mu_2 = a_{j_1+1} + \cdots
+ a_{\mu_{j_2}}, \cdots$,
$\mu_{k+1} = a_{j_k +1}+ a_{j_k+2} + \cdots + a_{2n+2}$.

For example,
$\rho_{2n+1} (\phi) = G_{a_1+ a_2 + \cdots + a_{2n+2}}$ and
$\rho_{2n+1} (C) = b_1 b_2 \cdots b_{2n+1}
(x-1)^{2n+1} (y-1)^{2n+1} G_{a_1}G_{a_2} \cdots G_{a_{2n+2}}$.

Now to each $U$, we call $U^{\ast}= \{b_{2n+2-j_k},
b_{2n+2-j_{k-1}},
\cdots, b_{2n+2-j_1}\}$ the dual of $U$.
We use these concepts to prove the theorem.

In the following, we assume that

\begin{align}\label{siki:16.n2} 
&{\rm for}\  
1 \leq j \leq n+1, a_{n+1+j} = - a_{n+2-j},\ {\rm and}\\
&{\rm for}\ 1 \leq j \leq n,b_{n+1+j} = - b_{n+1-j}.
\end{align}

Therefore, 
$A$ 
becomes 
\begin{align*}
\{a_1, b_1, a_2, b_2, \cdots,
a_{n+1},b_{n+1}, -a_{n+1}, -b_n, -a_n, \cdots,
-b_2, -a_2, -b_1, -a_1\}
\end{align*}
 and we can write

\begin{equation}\label{siki:16.n3}
\rho_{2n+1} (U^{\ast}) = b_{2n+2-j_k}\cdots b_{2n+2-j_1}
(x-1)^k (y-1)^k G_{-\mu_{k+1}}G_{-\mu_k}\cdots G_{-\mu_1}.
\end{equation}

Then we prove

\noindent
{\bf Claim 1.} {\it
$\rho_{2n+1}(U) + \rho_{2n+1}(U^{\ast}) = 0$, and hence
$\sum_{U}\rho_{2n+1}(U) = 0$, where the summation is taken over all $U$
that does not contain $b_{n+1}$.
}

{\it Proof.} 
If some $\mu_i = 0$, then $G_{\mu_1}=G_{-\mu_i} = 0$ and hence
$\rho_{2n+1}(U) = \rho_{2n+1}(U^{\ast}) = 0$. 
Therefore we may assume
that none of
$\mu_i$ is $0$.
Let $m$ be the number of negative elements in $U$. 
Then that number in $U^{\ast}$ is $k-m$. 
Let $q$ be the number of negative integers
in the set $\{\mu_{j_1}, \cdots, \mu_{j_{k+1}}\}$.
Then that number in
$U^{\ast}$
is $k+1-q$.
Therefore the number of occurrence of $(-1)$
in
$\rho_{2n+1}(U)$ is $m+q$, 
while that in $\rho_{2n+1} (U^{\ast})$ is
$k-m+k+1-q \equiv m+q+1$ (mod 2).
Next we count the exponent of the factor \mbox{$\frac{1}{xy}$} in $\rho (U)$ and
$\rho (U^{\ast})$.
Suppose $\mu_{\ell_1}, \mu_{\ell_2}, \cdots, \mu_{\ell_q}$ 
are negative.
Then in
$\rho_{2n+1}(U)$, the exponent of $\frac{1}{xy}$ is 
$ |\mu_{\ell_1}| +
| \mu_{\ell_2}| + \cdots +| \mu_{\ell_q}|$,
while in $\rho_{2n+1}(U^{\ast})$ that is 
$\sum_{\lambda \neq \ell_j}\mu_{\lambda}$.
Since $\sum_{\lambda \neq \ell_j}
\mu_{\lambda} - \sum_{j=1}^q
|\mu_{\ell_j}| = \sum_{j=1}^{2n+2} a_j= 0$, 
we see that $\rho_{2n+1}(U) + \rho_{2n+1}(U^{\ast}) = 0$. 
\qed

Now consider a short sequence
$A_0 = \{ a_1, b_1, a_2, b_2, \cdots, a_n, b_n, a_{n+1}\}$, 
the first half part of $A$. 
Let $B = \{ b_1, b_2, \cdots, b_n\}$ be an ordered set and
$\widehat{B}$
be the set of all ordered subset of $B$.
Then we have a mapping
$\rho_n: \widehat{B} \longrightarrow Z[x^{\pm 1}, y^{\pm 1}]$.

Given $U = \{b_{j_1}, b_{j_2}, \cdots, b_{j_p}, b_{n+1}, b_{j_{p+2}},
\cdots, b_{j_k}\} \in
\widehat{C}$, we can define two sets $W_{+}$
and $W_{-}$ in $\widehat{B}$
as follows.

$W_{+} = \{b_{j_1}, b_{j_2},\cdots, b_{j_p}\}$ and 
$W_{-} =\{b_{2n+2-j_k}, b_{2n+2-j_{k-1}}, \cdots,
b_{2n+2-j_{p+2}}\}$.

Then we claim

\noindent
{\bf Claim 2.}
{\it 
$\rho_{2n+1}(U) = b_{n+1}(x-1) (y-1) \rho_n (W_{+}) \rho_n
(W_{-}) \frac{-1}{(xy)^\alpha}$,
where
$\alpha = a_1+ a_2 + \cdots + a_{n+1}$.
}%

{\it Proof.}
Since $U \ni b_{n+1}$, we can write
\begin{align*}
\rho_{2n+1}(U) = &\\
&b_{n+1}(x-1)(y-1) \left(b_{j_1} \cdots b_{j_p}(x-1)^p
(y-1)^p G_{\mu_1} \cdots G_{\mu_{p+1}}\right)\\
&\left(b_{j_{p+2}} \cdots b_{j_k}
(x-1)^{k-p-1}(y-1)^{k-p-1} G_{\mu_{p+2}} \cdots G_{\mu_{k+1}}\right).
\end{align*}

Therefore, it suffices to show that

\begin{equation}\label{siki:16.n4}
(x-1)^{k-p-1}(y-1)^{k-p-1}(b_{j_{p+2}} \cdots b_{j_k}) 
G_{\mu_{p+2}}\cdots G_{\mu_{k+1}} =
\rho_n (W_{-}) \frac{-1}{(xy)^\alpha}.
\end{equation}

Since $\rho_n (W_{-}) = 
(x-1)^{k-p-1}(y-1)^{k-p-1}(b_{2n+2-j_k} \cdots
b_{2n+2-j_{p+2}})
G_{-\mu_{k+1}} \cdots G_{-\mu_{p+2}}$,
as done before, we compare the number of occurrences of 
$-1$ and the exponent
of $\frac{1}{xy}$, in LHS and RHS of 
(\ref{siki:16.n4}). 
First, let $d$ be the number of negative $b_i$ in LHS. 
Then
that number in RHS is $k-p-1-d$. 
Let $q$ be the number of negative
$\mu_{\lambda}$ in the set
$\{\mu_{p+2}, \cdots, \mu_{k+1}\}$ in LHS. 
Then
that number in RHS is $k-p-q$, 
and hence the sign of RHS is opposite to
that of LHS.

Next, we count the exponent of
$\frac{1}{xy}$.
Let $-\nu_1, -\nu_2, \cdots, -\nu_q$
be all negative members in $\{\mu_{p+2}, \cdots, \mu_{k+1}\}$, and
$\nu_{q+1}, \cdots, \nu_{k-p}$
be all positive members.
Then the exponent of $\frac{1}{xy}$
in LHS
is exactly $\nu_1 + \nu_2 + \cdots + \nu_q$, while that in RHS is
$\nu_{q+1} + \cdots + \nu_{k-p}$.
Since $\nu_{q+1}+ \cdots + \nu_{k-p} -(\nu_1 +
\cdots + \nu_q) = \sum_{j=1}^{n+1} a_j
= \alpha$, Claim 2 follows. \qed

Claim 2 implies easily the following

\noindent
{\bf Claim 3}. {\it
If
$U \in \widehat{C}$ contains $b_{n+1}$,
then $U^{\ast} \ni
b_{n+1}$ and $\rho_{2n+1}(U) = \rho_{2n+1}(U^{\ast})$.
}%

From Claims 1-3, we have,

\noindent{\bf Claim 4.}
{\it
$\sum_{U \in \widehat{C}}
\rho_{2n+1}(U)$ = $b_{n+1}(x-1)
(y-1) [\sum_{V \in \widehat{B}} 
\rho_n (V)]^2(\frac{-1}{(xy)^{\alpha}})$.
}%

Hence Theorem \ref{thm:16.9} follows.
\qed

\begin{ex}\label{ex:16.10}
Let $s=[4,2,-2]$. Then $K(s)$ is $\H$-stable by
Proposition \ref{prop:16.7}, 
and hence for $r=[s,2k,-s^{-1}], k \neq 0$, $K(r)$ is
$\H$-stable.
\end{ex}

If $K(s)$ is a 2-bridge knot, Theorem \ref{thm:16.9} 
does not hold, but
$\Delta_{K(r)}(x,y)$ will be of a nice form.
The following theorem is proven by applying a similar argument 
used in the proof of Theorem \ref{thm:16.9}.
The detail will appear in a separate paper.

\begin{thm}\label{thm:16.11}
Let $s = [2a_1, 2b_1, 2a_2, 2b_2, \cdots, 2a_n, 2b_n],
a_j \neq 0 \neq b_j, 1 \leq j \leq n,$
and $r = [s, 2k, -s^{-1}], k \neq 0$. Then

\begin{equation}\label{siki:16.4}
\Delta_{K(r)}(x,y) = G_k f(x,y) f(y,x),
\end{equation}
where $f(x,y) \in \ZZ[x,y]$ and $f(t,t) = \Delta_{K(s)}(t)$.
\end{thm}

We should note that if $|k| \geq 3$, then $G_k$ is not $\H$-stable.
Therefore, for $K(r)$ to be $\H$-stable, $k$ must be $\pm 1$ or $\pm 2$.

\section{Inversive links}\label{17}

A 2-component link $K$ is called 
{\it inversive} if the original link is stable
(or $c$-stable), 
but reversing the orientation of one component results
in a $c$-stable (or stable) link.
We see in Section 15 that some Montesinos links are 
inversive (Theorem \ref{thm:15.2}, Case 2).  
In this section, we study these links using 
2-variable Alexander polynomial $\Delta_K(x,y)$.

\subsection{Standard inversive links}\label{17.1}

From the definition, 
the following proposition is immediate.

\begin{prop}\label{prop:17.1} 
Let $K$ be a 2-component link and $\Delta_K (x,y)$ the
Alexander polynomial. 
Then $K$ is inversive if and only if
\begin{align}\label{siki:17.1}
(1)&\  \Delta_K (t,t) 
{\it \  is\  stable\ (or\ \mbox{$c$-stable})\ and}\nonumber\\
(2)&\  t^n \Delta_K (t, t^{-1}) {\it \  is\ \mbox{$c$-stable}\ 
 (or\ stable)}
\end{align}
\end{prop} 

\begin{rem}\label{rem:17.2}
Note that (2) in  (\ref{siki:17.1}) is equivalent to (2') below,
since $\Delta_K(x,y) = x^m y^n \Delta_K (x^{-1},y^{-1})$
for some integers $n$ and $m$.\\
\centerline{
(2')  $t^m \Delta_K (t^{-1},t)$ is $c$-stable (or stable).
}
For convenience, 
we call $\Delta_K (x,y)$ {\it inversive} if $\Delta_K (x,y)$
satisfies (1) and (2) in (\ref{siki:17.1}).
\end{rem}

Proposition \ref{prop:17.1} and 
Theorem \ref{thm:16.9} 
imply the following:

\begin{prop}\label{prop:17.3}
If a 2-bridge link $K(s)$ is inversive, then
$K(r)$ is inversive, where $r=[s, 2k, -s^{-1}], k \neq 0$.
\end{prop}

The simplest inversive 2-bridge link is 
$K(s), s=[2a], a \neq 0$.
Therefore we have the following corollary.

\begin{cor}\label{cor:17.4}
Let $r = [2a, 2k, -2a]$, where $k \neq 0$ and $a > 0$.
Then $K(r)$ is inversive.
\end{cor}

If $K(s)$ is a 2-bridge knot, 
then $K(r), r=[s,2k, -s^{-1}]$, may not be
inversive.

\begin{ex}\label{ex:17.5}
Let $s=[2,-2]$ and $r=[2,-2,2,2,-2]$. 
Then
$K(r)$ is not inversive. 
In fact, $\Delta_{K(r)}(x,y) =
(1-(2x+y)+xy)(1-(x+2y)+xy)$ and
$\Delta_{K(r)}(t,t) =(1-3t+t^2)^2$ is stable, 
but
$\Delta_{K(r)}(t,t^{-1}) = (1-2t+2t^2)(2-2t+t^2)$ 
is not $c$-stable.
\end{ex}

Now consider the general case. 
For convenience, we denote by $K^{\ast}$
the 2-component link
obtained from $K$ 
by reversing the orientation of one component of $K$.

\begin{prop}\label{prop:17.6}
Let $r = [2a_1, -2a_2, \cdots, 
(-1)^{2n} 2a_{2n+1}], a_j> 0$ 
for $1 \leq j \leq 2n+1$. 
Then $K(r)$ is inversive.
\end{prop}

{\it Proof.} 
First, by Theorem \ref{thm:6.1} 
$\Delta_{K(r)}(t)$ =$(t-1)\Delta_{K(r)}(t,t)$ is
stable. 
Further, we see that a diagram of $K(r)$ is 
alternating
and the diagram of $K^{\ast}(r)$ is 
special alternating, and hence
$\Delta_{K^{\ast}(r)}(t)$ is $c$-stable. 
Therefore, $K(r)$ is inversive.
\qed

If $r=[2a_1,2a_2, \cdots, 2a_{2n+1}], a_j > 0$ 
for $1 \leq j \leq 2n+1$,
then $K(r)$ is $c$-stable, 
by Proposition \ref{prop:12.1}. 
But, generally, $K^{\ast}(r)$ is not stable, 
and hence $K(r)$ is not inversive. 
The
following proposition, however, 
gives one sufficient condition for
$K(r)$ to be inversive.

\begin{prop}\label{prop:17.7}
Let $r=[2a_1,2a_2, \cdots, 2a_{2n+1}], a_j > 0$ for 
$1\leq j \leq 2n+1$. 
If $a_{2k+1}= 1$ for all $k$, $1 \leq k \leq n$, then
$K^{\ast}(r)$ is
stable and hence $K(r)$ is inversive.
\end{prop}

{\it Proof.} 
Write $r = \frac{\beta}{2\alpha}$. 
We may assume without loss of generality 
that $0 <\beta < 2\alpha$. 
Denote $r^{\ast} = \frac{2\alpha -\beta}{2\alpha}$.
Then it is known \cite[Proposition 3.17]{torti} 
that $K(r^{\ast})$ is equivalent to the mirror
image of $K^{\ast}(r)$. 
Therefore, $K^{\ast}(r)$ is stable if 
(and only if)
$K(r^{\ast})$ is stable.

Now, the even continued fraction expansion of 
$r^{\ast}$ is called the {\it dual} of 
(the even continued fraction expansion of) $r$ in 
\cite[Theorem 3.5]{torti}  
and there is
an algorithm to find the expression of $r^\ast$ 
\cite[p.7]{torti}.  
Using this algorithm, we can show that if all 
$a_{2k+1} = 1, 1 \leq k \leq n$, 
then $r^{\ast}=[2, -(2a_2 -2), 2, -(2a_4 -2), \cdots, 2,
-(2a_{2n} -2), 2]$.
Therefore, $K(r^{\ast})$ is stable and $K(r)$ is inversive. 
\qed

\begin{rem}\label{rem:17.8}
In Proposition \ref{prop:17.7}, 
if we assume that all 
$a_{2k}= 1$, for $1 \leq k \leq n$, 
then both $K(r)$ and $K(r^{\ast})$ are
$c$-stable.
In fact, it is easy to show that all entries of $r^{\ast}$ are positive.
\end{rem}

\begin{ex}\label{ex:17.9}
(1) Let $r = [2,4,2,6,2] = \frac{69}{118}$. 
Then $r^{\ast} =\frac{49}{118} 
= [2,-2,2,-4,2]$ and hence 
$K(r^{\ast})$ is stable. 
Since $K(r)$ is $c$-stable,
$K(r)$ is inversive.\\
(2) Let $r=[4,4,2] = \frac{7}{26}$. 
Then $r^{\ast} = \frac{19}{26}$ = [2,2,2,-2,2].
$K(r^{\ast})$ is not $c$-stable. 
In fact, $\Delta_{K(r^{\ast})} (x,y) =
x^2 y^2 - (2x^2 y +2xy^2) + 3xy - (2x+2y) + 1$, 
and hence $\Delta_{K(r^{\ast})} (t,t) 
= t^4-4t^3 +3t^2 -4t+2$. 
Then the modified polynomial $f(x)$ of 
$\Delta_{K(r^{\ast})}(t,t)$
has two real zeros, one of which is in $(-2,2)$ 
and another is larger than 2, 
and hence $\Delta_{K(r^{\ast})}(t, t)$ is strictly
bi-stable, and
$K(r)$ is not inversive.\\
(3) Let $r=[4,2,2,2,4] = \frac{13}{42}$. 
Then $r^{\ast} =\frac{29}{42}=[2,2,6,2,2]$
and hence $K(r^{\ast})$ is $c$-stable.
\end{ex}

\subsection{Exceptional inversive links}\label{17.2}

If the original 2-bridge link $K(r)$ is 
an exceptional stable link, 
then $K(r^{\ast})$ (and hence $K^{\ast}(r)$) may not be 
$c$-stable.
However, the following proposition shows that for some
exceptional 2-bridge link, 
$K(r^{\ast})$ is $c$-stable.

\begin{prop}\label{prop:17.10}
let $r=[4,2k,-2], k > 0$. Then $K(r)$ is inversive.
\end{prop}

{\it Proof.} 
Since $\Delta_{K(r)}(x, y) = k(x-1)(y-1)(x+y)-xy$,
we see that $\Delta_{K(r^{\ast})}(t,t) 
= t^2 \Delta_{K(r)}(t, t^{-1})
=kt^4 -2kt^3 +(2k+1)t^2 -2kt +k$. 
The modified polynomial $f(x)$ of 
$\Delta_{K(r^{\ast})}(t,t)$ is
$f(x) =kx^2 -2kx +1$, 
and both zeros of $f(x)$ are real and are in $(0,2)$
and hence $\Delta_{K(r^{\ast})}(t,t)$ is $c$-stable and 
$K(r)$ is inversive.
\qed

This is a rather exceptional case. 
For example, for $r$= [6,2,-2], $K(r)$
is not inversive. 
However, in general, we can prove
the following theorem.

\begin{thm}\label{thm:17.11}
Let $r = [2a,2k,-2c], a > c, k > 0$. 
If $k$ is sufficiently large, then $K(r)$ is inversive.
More precisely, there exists a positive integer 
$N(a,c)$ such that if $k \geq N(a,c)$, 
then $K(r)$ is inversive.
\end{thm}

{\it Proof.} 
Since $\Delta_{K(r)}(x,y) 
= k (x-1)(y-1) G_a G_{-c} + G_{a-c}$,
a simple calculation shows that
$\Delta_{K(r^{\ast})}(t, t) 
= k(t-1) \frac{t^{2a} -1}{t^2 -1}
\frac{t^{2c} -1}{t^2 -1}
+ \frac{t^{2(a-c)} -1}{t^2 -1} t^{2c}$,
and the modification $f(x)$ of $\Delta_{K(r)}(t,t)$ is
of the form: 
$f(x) = k(x-2) f_1 (x) f_2(x) + g(x)$, 
where $f_1, f_2$ and $g$ are, 
respectively, 
the modifications of $\frac{t^{2a} -1}{t^2-1}$, 
$\frac{t^{2c} -1}{t^2 -1}$
and $\frac{t^{2(a-c)} -1}{t^2 -1}$.
We note that all the zeros of $f_1, f_2$ and $g$ are 
in $(-2,2)$.
Consider two graphs $z_1 = (x-2)f_1 f_2$ and 
$z_2 =-\frac{g(x)}{k}$. 
If $k \longrightarrow \infty$, then 
$z_2\longrightarrow 0$
and hence if $k$ is sufficiently large, 
the points of intersection of two
curves
are almost the zeros (not 2) of $z_1$ and hence,
$\Delta_{K(r^{\ast})}(t,t)$ is $c$-stable
when $k$ is sufficiently large. 
Therefore $K(r)$ is inversive if $k$
is sufficiently large. 
\qed

\begin{probl}\label{probl:17.12}
Determine $N(a,c)$.
\end{probl}

We should note that if $a=c$ then $N(a,a) =1$.

\begin{ex}\label{ex:17.13}
It is easy to show that $N(3,1) = 3$ and
$N(3,2) = 2$.
\end{ex}

\begin{qu}\label{qu:17.14}
Can a 2-component inversive link $K$ be
characterized by the Alexander polynomial $\Delta_{K}(x,y)$?
\end{qu}

\medskip
\begin{alphasection}
\addtocounter{alphasect}{1}
{\bf Appendix A:  Representation polynomials}\label{A}
\setcounter{thm}{0}
\setcounter{subsection}{0}
\setcounter{equation}{0}

There are various integer polynomials associated to representations of
the knot group into $GL(2,\CC)$.
In this section, we discuss two particular representations of $G(K)$,
namely, a parabolic representation of $G(K(r))$,
the group of a $2$-bridge knot $K(r)$, to $SL(2,\CC)$ and a trace-free
representation of $G(K(r))$ to a dihedral group $D_{2n+1}\subset GL(2,\CC)$.

\subsection{Parabolic representation}\label{A.1}
Let $\theta_r (z)$ be the parabolic representation polynomial (Riley
polynomial) of $G(K(r))$ to $SL(2,\CC)$. (See \cite{ri72}.)
%
Suppose $r=\frac{1}{2n+1}$, and hence $K(r)$ is a torus knot of type
$(2,2n+1)$. 
Then $\theta_r (z) = \sum_{k=0}^n \binom{n+k}{2k} z^k$.

\begin{thm}\label{thm:A.1} \cite{ri72},\cite{swa}
If $r = \frac{1}{2n+1}, \theta_r (z)$ is real stable.
In fact, all the zeros of $\theta_r (z)$ are simple and they are
\begin{equation}\label{siki:A.1}
\alpha_k = - 4 \sin^2 \frac{(2k-1)\pi}{2(2k+1)}, 1 \leq k \leq n
\end{equation}
\end{thm}

\begin{rem}\label{rem:A.2}
For a rational number $r = \beta/ \alpha, 0 < \beta < \alpha,
\alpha$ odd,
$\theta_r (z)$ is an integer polynomial of degree $\frac{\alpha-1}{2}$ and
generally, $\theta_r(z)$ is not reciprocal.
\end{rem}

\begin{ex}\label{ex:A.3}
(1) Let $r = 2/5$. Then $\theta_r(z)= 1 -z + z^2$ is not stable, but
$c$-stable.
(2) Let $r = 5/7$. Then $\theta_r(z) = 1 + 2z +z^2 + z^3$ is not stable.
\end{ex}

\begin{probl}\label{probl:A.4}
Characterize $r$ so that $\theta_r(z)$ is stable.
\end{probl}

For a 2-bridge link $K(r), r= q/2n$, 
Riley polynomial is defined in a
slightly different manner.
Let $G(K(r)) = \langle x,y| Wy=yW\rangle$ be a presentation of 
the group of $K(r)$, 
where $x$ and $y$ are (oriented) meridian generators. 
Then $W$ is of the form:

\begin{equation}\label{siki:Ax.2} 
W = x^{\varepsilon_1} y^{\eta_1} x^{\varepsilon_2} y^{\eta_2}
\cdots x^{\varepsilon_{n-1}} y^{\eta_{n-1}} x^{\varepsilon_n},
\end{equation}
\noindent
where (1) $| \varepsilon_j | = | \eta_j | = 1$ for all $j$, and
(2) $\varepsilon_j = \varepsilon_{n-j+1}$ for $1 \leq j \leq n$, 
and $\eta_j = \eta_{n-j} $ for $1 \leq j \leq n-1$.

Let $\varphi : x \longrightarrow 
\left[\begin{array}{cc}
1& 1\\
0&1
\end{array}\right]$ %
and 
$y \longrightarrow 
\left[\begin{array}{cc}
1&0\\
z&1
\end{array}\right]$ %
be a parabolic representation of a free group $F(x,y)$, generated by 
$x$ and $y$, in $SL(2,\CC)$, and $z$ is a complex number that is determined later.
Then $\varphi$ defines a parabolic representation 
$\varphi_r$ of $G(K(r))$
in $SL(2,\CC)$ if and only if

\begin{equation}\label{siki:Ax.3}
\varphi (Wy) = \varphi (yW).
\end{equation}

Let $\varphi (W) = 
\tbt{a_r (z)}{ b_r (z)}{ c_r(z)}{ d_r(z)}$. 
Then a simple computation shows that (\ref{siki:Ax.3}) 
is equivalent to

\begin{equation}\label{siki:Ax.4}
z = 0{\rm \  or\ } b_r (z) = 0{\rm \  and\ } a_r(z) = d_r(z).
\end{equation}

We prove first that always $a_r (z) = d_r (z)$. To prove this, we need
the following simple lemma.
For convenience, we call a matrix $M = \tbt{a}{ b}{c}{d} \in
GL(2,\CC)$ is of $D$-type if $a=d$.

\begin{lemm}\label{lem:Ax.5}
If each of $M$ and $N$ in $GL(2,\CC)$ is of $D$-type, then $NMN$ 
is also of $D$-type.
\end{lemm}
Now $\varphi (x)$ and $\varphi (y)$ are of $D$-type, and so are $\varphi
(x^{-1})$ and $\varphi (y^{-1})$. 
Since the length of $W$ is $2n-1$,
$W$ has the central element. 
If $n$ is even, say $2m$, it is $y^{\eta_m}$. 
Since $\varphi (y^{\eta_m})$ is of $D$-type
and $\varepsilon_m = \varepsilon_{m+1}$, we see by Lemma \ref{lem:Ax.5}
that 
$\varphi (x^{\varepsilon_m} y^{\eta_m} x^{\varepsilon_{m+1}})$ is of
$D$-type.
Further, $\varphi (y^{\eta_{m-1}}(x^{\varepsilon_m} y^{\eta_m}
x^{\varepsilon_{m+1}}) y^{\eta_{m+1}})$ is also of
$D$-type. 
By repeating this process, we see that $\varphi (W)$ = $\varphi
(x^{\varepsilon_1} y^{\eta_1} \cdots y^{\eta_m}\cdots y^{\eta_{2m-1}}
x^{\varepsilon_{2m}})$ is of $D$-type.
If $n$ is odd, say $2m+1$, then the central element is 
$x^{\varepsilon_{m+1}}$, and the same argument works well.
If $z=0$, then $\varphi$ is an abelian representation, and hence we
ignore it. 
But each zero of $b_r (z)$ gives a (non-abelian)
parabolic representation $\varphi_r$. $b_r (z)$ is called 
{\it Riley polynomial}
$\theta_r (z)$ of $K(r)$. 
The degree of $\theta_r (z)$ is $n-1$.

\begin{ex}\label{ex:Ax.6}
Let $r=1/2n$, $n\ge 2$. Then $K(r)$ is an elementary
torus link, and it follows from (\ref{siki:A.k6}) and
\cite[Prop. 2.4]{eval} that the
Riley polynomial $\theta_r (z)$ of $K(r)$ is given 
by (\ref{siki:Ax.5}) below.

\begin{equation}\label{siki:Ax.5}
\theta_r (z) = \sum_{j=0}^{n-1} \binom{n+j}{2j+1} z^j.
\end{equation}

\end{ex}

It is known that $\theta_r (z)$ is real stable. 
In fact, the zeros of $\theta_r(z), r=\frac{1}{2n}$ are 
$-4\sin^2\frac{r\pi}{2n}, r=1,2\dots,n-1$
 \cite{swa}.
The following proposition confirms a conjecture by Dan Silver \cite{silv}.

\begin{prop}\label{prop:Ax.6}
Let $r = q/2n, 0 < q < 2n, \gcd(q, 2n) = 1$.
Then $|\theta_r (0)| = |lk(K(r))|$, where $lk(K(r))$ denotes the linking
number between two components of a 2-bridge link $K(r)$.
\end{prop}

{\it Proof.}
 $\theta_r (0)$ is determined by $\varphi_r (W)$ evaluated at $z=0$
which is
$\prod_{j=1}^n \varphi_r (x^{\varepsilon_j}) = 
\tbt{1}{\sum_{j=1}^n
\varepsilon_j }{0}{1}$.
Since $\sum_{j=1}^n \varepsilon_j $ is equal
to $lk(K(r))$, Proposition \ref{prop:Ax.6} follows. 
\qed

If $lk(K(r)) = 0$, then Dan Silver also conjectures that the
absolute value of the coefficient $c_1$ of $z$ of $\theta_r (z)$ is the
wrapping number of $K(r)$. However, examples below show that it is not
correct.
Note that if $lk(K(r)) = 0$, then $n$ is a multiple of $4$.

\begin{ex}\label{ex:Ax.7}
(1) Let $r=9/16 = [2,4,-2]$. 
Then $\theta_r (z) =
z(2+z^2)(2-4z+4z^2 -2z^3 +z^4)$, but the wrapping number is $2$.
(2) Let $r=11/24 =[2,-6,-2]$. Then $\theta_r (z)$ = $6z +18z^2 +35z^3
+48z^4 +56z^5 +44z^6 +36z^7 +16z^8 +10z^9 +2z^{10} +z^{11}$,
but the wrapping number is $2$.
\end{ex}

\subsection{Dihedral representation}\label{A.2}

It is well-known that there is a trace-free representation $\xi$ of a
dihedral group $D_p = \langle x,y| x^2 = y^2 = (xy)^p = 1 \rangle$,
$p=2n+1$, in $GL(2,\CC)$. 
$\xi$ is given by
$\xi(x)=\tbt{-1}{1}{0}{1}$ and
$\xi(y)=\tbt{-1}{0}{\omega}{1}$,
where $\omega$ is a zero of the integer polynomial
$\varphi_p (z) = \sum_{k=0}^n \frac{2n+1}{2k+1}\binom{n+k}{2k} z^k$.
See \cite{ri72}.

\begin{ex}\label{ex:A.5}
$\varphi_3 (z) = z +3$,
$\varphi_5 (z) = z^2 + 5z + 5$,
$\varphi_7 (z) = z^3 + 7z^2 + 14z +7$,
$\varphi_9 (z) = z^4 + 9z^3 + 27z^2 +30 z + 9 =(z+3)(z^3 +6z^2 +9z+3)$.
\end{ex}

Using $\xi$, we prove in \cite{hm09}  
some properties of the twisted Alexander
polynomial of a 2-bridge knot associated to a dihedral representation.
In this subsection, we prove for any odd number $2n + 1$,
$\varphi_{2n+1}(z)$ is stable.
Our proof is different from those of other parts in this paper.

\begin{thm}\label{thm:A.6}
Let $p = 2n+1, n \geq 1$, and
$\varphi_p (z) = \sum_{k=0}^n \frac{2n+1}{2k+1}\binom{n+k}{2k} z^k$.
Let $\zeta= e^{\frac{2\pi i}{p}}$. 
Then $z_k = \zeta^k + \zeta^{-k} -2, 1 \leq k \leq n$, 
are the zeros of $\varphi_p (z)$.
\end{thm}

{\it Proof.} 
Since $\zeta^k + \zeta^{-k} - 2 = (\sqrt{\zeta} - 
\frac{1}{\sqrt{\zeta}})^{2k}$, it suffices to show that
\begin{equation}\label{siki:A.k6}
\varphi_{2n+1} (z_1) = \sum_{k=0}^n
\frac{2n+1}{2k+1}\binom{n+k}{2k}(\sqrt{\zeta} - 
\frac{1}{\sqrt{\zeta}})^{2k}=0.
\end{equation}

Now we expand the right side of (\ref{siki:A.k6}) 
and let $A_k^{(n)}$ denote
the term of $\zeta^k, k=0, \pm 1, \pm 2, \cdots, \pm n$. 
Namely,
\begin{equation*}
\varphi_{2n+1} (z_1)= A_{-n}^{(n)} \zeta^{-n} + \cdots + A_{-1}^{(n)}
\zeta^{-1} + A_{0}^{(n)} + A_{1}^{(n)} \zeta
+ A_{2}^{(n)} \zeta^{2} + \cdots
+ A_{n}^{(n)} \zeta^{n}.
\end{equation*}

Then we see that
\begin{align}\label{siki:A.k7}
A_0^{(n)} &= \sum_{k=0}^n\frac{2n+1}{2k+1} \binom{n+k}{2k} (-1)^k
\binom{2k}{k}.\nonumber\\
A_1^{(n)} &= A_{-1}^{(n)} = \sum_{k=1}^n\frac{2n+1}{2k+1}
\binom{n+k}{2k} (-1)^k \binom{2k}{k+1}\nonumber\\
&\vdots\nonumber\\
A_m^{(n)} &= A_{-m}^{(n)} = \sum_{k=m}^n\frac{2n+1}{2k+1}
\binom{n+k}{2k)} (-1)^{k+m} \binom{2k}{k+m}\nonumber\\
&\vdots\nonumber\\
A_n^{(n)} &= A_{-n}^{(n)} = \frac{2n+1}{2n+1} \binom{2n}{2k} (-1)^2n
\binom{2n}{2n}=1
\end{align}

Therefore, to prove $\varphi_{2n+1}(z)= 0$, it suffices to show
\begin{equation*}
A_0^{(n)} = A_1^{(n)} = \cdots = A_n^{(n)} = 1.
\end{equation*}

We prove these equalities by applying generating function theory.

First we show that $A_0^{(n)} = 1$.
Consider the generating function $F_0(x)$ of $A_0^{(n)}$:
%
\begin{align}\label{siki:A.3}
F_0 (x) &= \sum_{n\geq 0} A_0^{(n)} x^n
\nonumber\\
&= \sum_{n \geq 0}\{x^n \sum_{k \geq 0}\binom{n+k}{2k} \binom{2k}{k}
(-1)^k \frac{2n+1}{2k+1}\}.
\end{align}
We prove that 
$F_0 (x) = \frac{1}{1-x} = 1 + x + x^2 + \cdots+ x^n +\cdots$.\\
Now, by interchanging the order of summations in (\ref{siki:A.3}), 
we have
\begin{equation*}
F_0 (x) = \sum_{k \geq 0}[\binom{2k}{k} \frac{(-1)^k}{2k+1} x^{-k}
\{\sum_{n \geq 0} (2n+1) \binom{n+k}{2k} x^{n+k}\}].
\end{equation*}
Note $\sum_{n \geq 0} (2n+1) \binom{n+k}{2k} x^{n+k}$ = $\sum_{n \geq
k} (2n+1) \binom{n+k}{2k} x^{n+k}$.
Let $n+k = r$. 
Then $2n+1=2(r-k)+1=2r-2k+1$ and

\begin{align}\label{siki:A.4}
&\sum_{n \geq k} (2n+1) \binom{n+k}{2k} x^{n+k} \nonumber\\
&\ \ = \sum_{r\geq 2k}(2r-2k+1) \binom{r}{2k} x^r
\nonumber\\
&\ \ = x^{2k}\{(2k+1) + (2k+3) \binom{2k+1}{1} x 
\nonumber\\
& \ \ \  \ + (2k+5) \binom{2k+2}{2}x^2
+ \cdots + (2k+2m+1) \binom{2k+m}{m} x^m + \cdots \}
\nonumber\\
&\ \ = x^{2k} \sum_{m \geq 0} (2k+2m+1) \binom{2k+m}{m} x^m.
\end{align}
We need the following lemma.

\begin{lemm}\label{lem:A.7}
(1) [W, p.50, (2.5.7)]
$\sum_{m \geq 0}\binom{2k+m}{m} x^m = \frac{1}{(1-x)^{2k+1}}$.\\
(2) [W, p.50, (2.5.11)] $\sum_{k \geq 0}\binom{2k}{k} x^k =
\frac{1}{\sqrt{1-4k}}$.\\
(3) [W p.51, (2.5.15)]
$\sum_{k \geq 0}\binom{2k+m}{k} x^k = \frac{1}{\sqrt{1-4x}} (\frac{1 -
\sqrt{1-4x}}{2x})^m$.\\
(4) [W p.32] Let $P(y)$ be a polynomial and $f= \sum_{n \ \geq 0} a_n x^n$. 
Then
$\sum_{n \geq 0}P(n) a_n x^n = P(x\frac{d}{dx})f$.
\end{lemm}

\begin{ex}\label{ex:A.8}
If $P(y) = 2y+2k+1$ and $f= \sum_{m \geq 0} 
\binom{2k+m}{m} x^m = (1-x)^{-(2k+1)}$.
Then $P(m) = 2m+2k+1$ and
\begin{align*}
&\sum_{m \geq 0}(2m+2k+1) \binom{2k+m}{m} x^m\\
&\ \ \ = 2x\frac{\;df}{dx}+ (2k+1)f \\
&\ \ \ = 2(-(2k+1))(-1)(1-x)^{-(2k+2)} x
+ (2k+1) (1-x)^{-(2k+1)}\\
&\ \ \ = 
\frac{(2k+1)(x+1)}{(1-x)^{2k+2}}.
\end{align*}
\end{ex}

Using Example \ref{ex:A.8}, 
we have from (\ref{siki:A.4}) 
%
\begin{equation}\label{siki:A.5}
\sum_{n \geq 0} (2n+1) \binom{n+k}{2k} x^{n+k} = x^{2k}
\frac{(2k+1)(x+1)}{(1-x)^{2k+2}}.
\end{equation}
Therefore,
\begin{align*}
F_0 (x) 
&= \sum_{k \geq 0}\binom{2k}{k}\frac{(-1)^k}{2k+1}
x^{-k}\left\{\frac{x^{2k}(2k+1) (x+1)}{(1-x)^{2k+2}}\right\}
\\
&= \frac{x+1}{(1-x)^2} \sum_{k \geq 0} \binom{2k}{k} (-1)^k
\left\{\frac{x}{(1-x)^2}\right\}^k.
\end{align*}
Let $y = \frac{x}{(1-x)^2}$. 
Then
$F_0 (x) = 
\frac{x+1}{(1-x)^2} \sum_{k \geq 0} (-1)^k \binom{2k}{k}y^k$.
By Lemma \ref{lem:A.7} (2), we see that 
$\sum_{k \geq 0} (-1)^k \binom{2k}{k} y^k = \frac{1}{\sqrt{1+4y}}$.
Since $1+4y = 1 + \frac{4x}{(1-x)^2} = \frac{(1+x)^2}{(1-x)^2}$, we have
$\sqrt{1+4y} = \frac{1+x}{1-x}$, and hence,
\begin{equation*}
F_0 (x) = \frac{x+1}{(1-x)^2} \frac{1-x}{1+x} = \frac{1}{1-x} = 1 + x +
x^2 + \cdots.
\end{equation*}
Therefore, for any $n$, 
the coefficient of $F_0 (x)$ is 1, i.e.,
$A_0^{(n)} = 1$.

Next, for $m \geq 1$, we show $A_m^{(n)}= 1$. 
Our approach is almost
the same, but we need a slight change in the process. Now,
$A_m^{(n)} = \sum_{k \geq m} (-1)^{k+m} \frac{2n+1}{2k+1} 
\binom{n+k}{2k} \binom{2k}{k+m}$.
As before, we interchange the order of summations of the generating
function $F_m (x)$ of $A_m^{(n)}$:
\begin{align*}
F_m (x)& = \sum_{n \geq 0} x^n \left\{\sum_{k \geq m} (-1)^{k+m}
\frac{2n+1}{2k+1} \binom{n+k}{2k} \binom{2k}{k+m}\right\} \\
&= \sum_{k \geq m}\frac{(-1)^{k+m}}{2k+1} \binom{2k}{k+m}
x^{-k}\left\{\sum_{n \geq 0}(2n+1) \binom{n+k}{2k} x^{n+k}\right\}\\
&= \sum_{k \geq m} \frac{(-1)^{k+m}}{2k+1} \binom{2k}{k+m} x^{-k}
\frac{x^{2k}(x+1) (2k+1)}{(1-x)^{2k+2}}\ \  {\rm \ by\ } (\ref{siki:A.5})\\
&= \frac{x+1}{(x-1)^2} \sum_{k \geq m} (-1)^{k+m} \binom{2k}{k+m}
\frac{x^k}{(1-x)^{2k}}.
\end{align*}

We show $F_m (x) = \frac{x^m}{1-x} = x^m + x^{m+1} + x^{m+2} + \cdots$.
We note that
\begin{align*}
&\sum_{k \geq m}(-1)^{k+m}\binom{2k}{k+m} \frac{x^k}{(1-x)^{2k}}\\
& \ \ \ = \binom{2m}{2m} \frac{x^m}{(1-x)^{2m}} - \binom{2m+2}{2m+1}
\frac{x^{m+1}}{(1-x)^{2(m+1)}} + \cdots\\
& \ \ \ = \frac{x^m}{(1-x)^{2m}}\left\{ \binom{2m}{2m} - \binom{2m+2}{2m+1}
\frac{x}{(1-x)^2} + \cdots \right\}\\
& \ \ \ = \frac{x^m}{(1-x)^{2m}} \sum_{k \geq 0} (-1)^k 
\binom{2k+2m}{k}\left\{\frac{x}{(1-x)^2}\right\}^k.
\end{align*}

Again, let $y= \frac{x}{(1-x)^2}$. Then $\sqrt{1+4y} = 
\frac{1+x}{1-x}$ and
$ 1 - \sqrt{1+4y} = 1 - \frac{1+x}{1-x} = \frac{-2x}{1-x}$ and $\frac{1
- \sqrt{1+4y}}{-2y} = 1 - x$.
Then by Lemma \ref{lem:A.7} (3), we have 
\begin{align*}
&\sum_{k \geq m} (-1)^{k+m} \binom{2k}{k+m} \frac{x^k}{(1-x)^{2k}}\\
&\ \ \ = \frac{x^m}{(1-x)^{2m}} \frac{1}{\sqrt{1+4y}} \left\{\frac{1 -
\sqrt{1+4y}}{-2y}\right\}^{2m}\\
& \ \ \ = \frac{x^m}{(1-x)^{2m}} \frac{1-x}{1+x} (1-x)^{2m}\\
& \ \ \ = \frac{x^m (1-x)}{1+x}.
\end{align*}

Therefore $F_m (x)$ = $\frac{x+1}{(1-x)^2} \frac{x^m (1-x)}{1+x} = 
\frac{x^m}{1-x}
= x^m (1+x+x^2 + \cdots)$.
Thus if $n < m$, then $A_m^{(n)} = 0$, and if $n \geq m$, then 
$A_m^{(n)}=1$,
i.e., for any $n \geq m, A_m^{(n)} = 1$.
\qed

\newpage

\addtocounter{alphasect}{1}
{\bf Appendix B:  Determination of $\delta_4$}\label{B}
\setcounter{thm}{0}
\setcounter{equation}{0}

Let $\Gamma_{2n}$ be the set of all Alexander polynomials $\Delta_K(t)$ of
alternating knots $K$ of genus $n$, i.e., $\deg\Delta_K(t) = 2n$.
Let $\delta_{2n}(K)$ be the maximal value of Re$(\alpha)$ of the zero
$\alpha$ of
$\Delta_K(t)$ and
$\delta_{2n} = \max_{\Delta_K(t) \in \Gamma_{2n}} \delta_{2n}(K)$.

\begin{yosou}\label{conj:A2.1}
$\delta_{2n}$ exists for any $n \geq 1$, and further,
there is a fibred stable alternating knot $K_n$
such that $\delta_{2n}(K_n) = \delta_{2n}$.
\end{yosou}

In this section, we prove Conjecture \ref{conj:A2.1} 
for $n=2$. 
Conjecture \ref{conj:A2.1} 
is trivially true for $n=1$. 
In fact, $\delta_2 = 2.618 \cdots$
that is a zero of $\Delta_{K_1}(t) = t^2 -3t+1$, where $K_1$ is $4_1$.

\begin{thm}\label{thm:A2.2}
Let $K_2 = 8_{12}$. Then $\delta_4 = \delta_4 (K_2) =
4.3902 \cdots$.
\end{thm}

Note that $K_2$ is a fibred stable alternating knot with 
$\Delta_{K_2}(t) = t^4 -7t^3 +13t^2 -7t +1$.

{\it Proof of Theorem \ref{thm:A2.2}.} 
Let $K$ be an alternating knot. Since $\delta_4 (K_2 )$ is the maximal
value among those of fibred knots of genus 2, we may assume
that $K$ is a non-fibred alternating knot.
Write $\Delta_K(t) = a t^4 - bt^3 + ct^2 - bt +a$, where $a,b,c >
0$ and further, $a \geq 2$.
To prove Theorem \ref{thm:A2.2}, 
we need the following theorem due to Jong.

\begin{thm}\label{thm:A2.3} \cite{jong} 
Let $\Delta_K(t) = a t^4 - bt^3 + ct^2 - bt + a$,
where $a,b,c > 0$, be the Alexander polynomial of
an alternating knot $K$ of genus 2. Then the following holds.\\
(1) if $\sigma(K) = 0$, then $3a-1 \leq b \leq 6a+1$,\\
(2) if $|\sigma(K)| = 2$, then $2a+1 \leq b \leq 6a -1$,\\
(3) if $|\sigma(K)| = 4$, then $2a-1 \leq b \leq 4a- 2$.
\end{thm}

Now there are three cases.\\
Case 1. $\Delta_K(t)$ has four complex zeros, none of which is 
a unit complex. 
Then, $\sigma(K) = 0$.
Let $\alpha, \overline{\alpha}, \beta, \overline{\beta}$ be all zeros
of$\Delta_K(t)$, where $\alpha \beta = 1$ and
$\overline{\alpha} \overline{\beta} = 1$.
First, the real part of each zero is positive. 
In fact, if $Re(\alpha) < 0$, 
then the real parts of all zeros are negative,
since $\alpha \beta$ = $\overline{\alpha} \overline{\beta} = 1$.
Therefore, $\alpha +\overline{\alpha} + \beta + \overline{\beta} =
2{\rm Re}(\alpha) + 2{\rm Re}(\beta)< 0$, 
but $\alpha +\overline{\alpha} + \beta + \overline{\beta} = b/a >
0$, a contradiction.
Now suppose ${\rm Re}(\alpha) \geq \delta_4$. 
Then $b/a =\alpha +\overline{\alpha} +\beta + \overline{\beta} 
= 2{\rm Re}(\alpha) + 2{\rm Re}(\beta) > 2\delta_4 > 8$,
but by Theorem \ref{thm:A2.3}, 
$b/a \leq 6 + \frac{1}{a} < 7$, a contradiction. 
Therefore, $\delta_4(K) < \delta_4$.\\
Case 2. $K$ is $c$-stable. Trivially, $\delta_{4}(K) < 1$ 
and hence $\delta_4 (K) < \delta_4$.
If $|\sigma(K)|=4$, then $K$ is $c$-stable and hence we may assume hereafter
that $|\sigma(K)| \leq 2$, and further,
$\Delta_K(t)$ has at least two real zeros. Therefore, the last case is the
following:\\
Case 3. $K$ is bi-stable, but not $c$-stable.
From the above remark, we see that $\delta _4 (K)$ is the maximal real
zero of $\Delta_K(t)$.
To show that $\delta_4 (K) < \delta_4$, first we consider the modified
polynomial $f(x)$ of $\Delta_K(t)$.
Write $\Delta_K(t) = a t^4 - bt^3 + (2b - 2a - 
\varepsilon )t^2 - bt + a$, where $\varepsilon= \pm 1$. Then
\begin{equation*}
f(x) = ax^2 - bx +(2b - 4a - \varepsilon)
= (x-2)(ax - (b-2a)) -\varepsilon.
\end{equation*}

Since $\delta_4(K)$ is less than the maximal real zero of $f(x)$, 
we compute the real zeros of $f(x)$.
Now the real zeros of $f(x)$ are determined by the 
intersection of two curves
$y_1 = (x-2)(ax - (b-2a))$ and $y_2 = \varepsilon$.
Since $y_1 (0) = 2(b-2a) \geq 2$, by Theorem \ref{thm:A2.3}, 
we have the following graphs:
(1) $\dfrac{b-2a}{a} \le 2$,\ \  (2) $\dfrac{b-2a}{a} \ge 2$

\bichi

The maximal real zero $\gamma$ of $f(x)$ (if exists) is given by
\begin{equation*}
\gamma = \frac{b}{2a} + \sqrt {\frac{d}{4a^2}},{\rm \  where\ } 
d = (b-4a)^2 +4a \varepsilon.
\end{equation*}

From this formula, we should note that when $a$ is fixed, $\gamma$ gets
larger as $b$ gets larger. 
Therefore, to obtain the maximal real zero of $f(x)$,
$b$ should be the maximal possible value.

Subcase (a) $\varepsilon = 1$.
From Fig. B.1, we see that $f(x)$ has two zeros, one is larger than 2, but the
other is less than 2. 
Therefore, $\Delta_K(t)$ has two unit complex zeros
and two real zeros, and hence $|\sigma(K)| = 2$. 
By Theorem \ref{thm:A2.3} (2), we 
see that $2a + 1 \leq b \leq 6a - 1$. When $b=6a-1$, $d= 4a^2 + 1$ and
hence
\begin{equation*}
\gamma = \frac{b}{2a} + \sqrt{\frac{d}{4a^2}} = \frac{b}{2a} + \sqrt{1+
\frac{1}{4a^2}}.
\end{equation*}

Since $a \geq 2$, we have $\gamma = \frac{6a -1}{2a} + \sqrt{1 +
\frac{1}{4a^2}} \leq 3 + \sqrt{1.0625}
= 4.03077\cdot < \delta_4$, and hence $\delta_4(K) < \delta_4$.

Subcase (b) $\varepsilon = -1$.
Since $f(x)$ has a real zero larger than 2, Fig B.1 (1) cannot occur. 
Therefore, from Fig B.1 (2), $f(x)$ has two real zeros greater than 2
and hence $\Delta_K(t)$ has four real zeros, and $|\sigma(K)| = 0$. 
Then by Theorem \ref{thm:A2.3} (1), we 
have that $3a-1 \leq b \leq 6a + 1$. 
When $b=6a+1$, $d=4a^2 +1$ and
$\gamma = \frac{6a+1}{2a} + \sqrt{1+ \frac{1}{4a^2}}$.
Since $a \geq 2$, it follows that $\gamma \leq 3 + \frac{1}{4} +
\sqrt{1.0625} = 4.281\cdots < \delta_4$.\\
A proof of Theorem \ref{thm:A2.2} 
is now complete.
\qed

\addtocounter{alphasect}{1}
{\bf Appendix C: Distribution of the zeros.}\label{C}
\setcounter{thm}{0}
\setcounter{equation}{0}

In this section, we discuss distribution of the zeros of the Alexander
polynomials of two infinite sequences of 2-bridge knots.
Namely, they are vertical and horizontal extensions of 
the 2-bridge knot $[2,2,2,2,-2,-2,-2,-2]$.

Let $r(k) = [2,2,2,2k, -2,-2,-2,-2], k \ne 0$. For simplicity, $K(r(k))$
will be denoted by $K(k)$.
The type of the zeros of the Alexander polynomial of $K(k)$ depends on $k$.
More precisely, we prove the following theorem.

\begin{thm}\label{thm:D.1}
Case 1. $k > 0$.\\
(1) If $k=1$ or $2$, then $\Delta_{K(k)}(t)$ 
is totally unstable, i.e., every
zero is a non-unit complex number.\\
(2) If $k=3,4,5$ or $6$, then $\Delta_{K(k)}(t)$ has 
four unit complex zeros and four
non-unit complex zeros, and hence $\Delta_{K(k)}(t)$ has no real zeros.\\
(3) If $k \geq 7$, then $\Delta_{K(k)}(t)$ has eight unit complex zeros, 
and hence
$K(k)$ is $c$-stable.

Therefore, in this case, $\Delta_{K(k)}(t)$ does not have real zeros.

Case 2. $k < 0$.\\
For all $k$, $\Delta_{K(k)}(t)$ has two real zeros and six unit complex
zeros, and hence $K(k)$ is strictly bi-stable.
\end{thm}

\begin{ex}
For any $k\neq 0$, 
$\Delta_{K(k)}(t)=
k -3kt +5kt^2 -7kt^3 +(8 k +1)t^4
-7kt^5 +5kt^6 -3kt^7 +kt^8$.
We plot the zeros around the unit circle
in Fig. C.1 for
$k=-4,-3,-2,-1,1,2,\dots,8$.
\end{ex}
\dichi

\begin{rem}\label{rem:D.2}
If $k > 0$, then $\sigma(K(k)) = 0$. But, if $k \geq 7$, 
$K(k)$
is $c$-stable. Therefore, $c$-stable alternating knots
are not necessarily special alternating. 
If $k < 0$, then $|\sigma(K(k))|= 2$, 
but $\Delta_{K(k)}(t)$ has always more than two unit complex zeros.
\end{rem}

{\it Proof of Theorem \ref{thm:D.1}}. 
First, using a standard Seifert matrix of $K(k)$, we can show
that $\Delta_{K(k)}(t) = k f(t) + t^4$, where $f(t) = (t -1)^2 (t^2 +1)
( t^4 -t^3 +t^2 -t+1)$.
Consider the modification $F(x)$ of $\Delta_{K(k)}(t): 
F(t) = kx(x-2)(x^2-x -1) +1$.
To prove the theorem, we study 
the intersection of two
curves, $y_1 = g(x) = x(x-2)(x^2 -x -1)$ and $y_2 = - 1/k$, $k \ne 0$.
By simple calculations, $y_1$ is depicted in Fig. C.2.
Then we see
(1) two curves $y_1$ and $y_2 = -1/k, k \leq -1$, intersect in
exactly four points, only one of which has $x$-coordinate greater than 2
and others in $(-2,2)$.
This proves Case 2.\\
(2) Suppose $k > 0$. 
If $k = 1$ or $2$, 
two curves do not intersect and hence,
$K(k)$ is totally unstable.
If $k = 3,4,5$ or $6$, then 
two curves
intersect in two points with $x$-coordinate in $(-2,2)$. 
If $k \geq 7$, two curves intersect in four points 
with $x$-coordinate between $-2$ and $2$,
and hence all the zeros are unit complex.
This proves Theorem \ref{thm:D.1}.
\qed

\dni

The previous sequence can be considered as a vertical extension of the
original 2-bridge knot $K(1) = [2,2,2,2,-2,-2,-2,-2]$.
The next sequence is a horizontal extension of $K(1)$. 
Consider the sequence $r[n] = [2,\dots,2, -2, \dots, -2]$,
$n$ consecutive $2$'s followed by $n$ consecutive $-2$'s, with $n\ge 1$.
$K(r[n])$ will be
denoted by $K[n]$. 
We prove the following theorem.

\begin{thm}\label{thm:D.4} 
(1) If $n$ is odd, then $\Delta_{K[n]}(t)$ has two real zeros, and other are
non-unit complex zeros.\\
(2) If $n$ is even, then $\Delta_{K[n]}(t)$ is totally unstable.
\end{thm}

{\it Proof.} First, using a standard Seifert matrix, it is easy to show by
induction on $n$ that
\begin{equation}\label{siki:D.2} 
\Delta_{K[n]}(t) = \sum_{k=0}^{n-1}(-1)^k (2k+1)\{t^k +
t^{2n-k}\} + (-1)^n (2n+1) t^n.
\end{equation}

To prove the theorem, we need the following lemma.

\begin{lemm}\label{lem:D.5} 
$\Delta_{K[n]}(t)$ does not have a unit complex zeros.
\end{lemm}

{\it Proof.}
 Obviously, $\pm 1$ is not the zeros of $\Delta_{K[n]}(t)$. Now,
we express $\Delta_{K[n]}(t)$ in a different form:

\begin{equation}\label{siki:D.3} 
(t+1)^2 \Delta_{K[n]}(t) = (t^{2n+1}-1)(t-1) + (-1)^n 4
t^{n+1}.
\end{equation}

Suppose $\Delta_{K[n]}(t)$ has a unit complex zero $\alpha = e^{i\theta},
\theta \ne 0, \pi$.
Then $(e^{(2n+1)i\theta} -1)(e^{i\theta}-1) = (-1)^{n+1}4
e^{(n+1)i\theta}$ and hence,
$|e^{(2n+1)i\theta} -1| |e^{i\theta}-1| = 4$. This is impossible, since
$|e^{(2n+1)i\theta} -1| \leq 2$ and $|e^{i\theta}-1| < 2$.
\qed

We return to a proof of the theorem.\\
Case 1. $n$ is odd, say $2m+1$. Then\\
$\Delta_{K[n]}(t) = (t-1)^2 (t^{2m}+t^{2m-2} + \cdots +t^2 +1)
(t^{2m}-t^{2m-1} +\cdots +
t^2 -t +1) -t^{2m+1}$.
Consider the modified polynomial $F(x)$ of $\Delta_{K[n]}(t)$. 
Then
$F(x) =(x-2)g(x)h(x) - 1$, where $g(x)$ and $h(x)$ are,
respectively, the modified polynomials of $f_1(t) =t^{2m} + t^{2m-2}+
\cdots + t^2 + 1$ and
$f_2 (t) = t^{2m} -t^{2m-1} +\cdots + t^2 -t +1$. 
Since
$f_1(t)$ and $f_2(t)$ have only unit complex zeros, all the zeros of $g(x)$
and $h(x)$ are real in $(-2,2)$. 
Then
the zeros of $\Delta_{K[n]}(t)$ are determined by 
intersection of
two curves $y_1 = (x-2)g(x)h(x)$ and $y_2 = 1$. 
See Fig. C.3 (1).
We see that $y_1$ and $y_2$ intersect in one point $P(p_1,p_2)$, with
$p_1 > 2$. 
If $y_2$ intersects $y_1$ at another point, say 
$Q (q_1,q_2)$, then $-2 < q_1 < 2$ and
the corresponding zeros of $\Delta_{K[n]}(t)$ are unit complex. 
This is
impossible by Lemma \ref{lem:D.5}. 

\dsan

Case 2. $n$ is even, say $2m$. 
Then again, we have
$\Delta_{K[n]}(t)= (t-1)^2 (t^{2m-2} +t^{2m-4}+\cdots +t^2+1)
(t^{2m}-t^{2m-1}+\cdots + t^2-t +1) +t^{2m}$.\\
From this form, we see easily that $\Delta_{K[n]}(t)$ has no real zeros,
since $\Delta_{K[n]}(t)> 0 $ if $t$ is real.
Now we did as before, consider the modified polynomial $F(x)$ of
$\Delta_{K[n]}(t)$ and see that 
$F(x) =(x-2) g(x)h(x) +1$, where $g(x)$ and $h(x)$ are the modifications
of $\frac{t^{2m}-1}{t^2 -1}$ and $\frac{t^{2m+1}+1}{t +1}$, respectively.
Both have only real zeros in $(-2,2)$.
Consider the intersection of $y_1 = (x-2)g(x)h(x)$ and $y_2 = -1$. 
Since $\deg y = 2m$, the graph appears as in Fig. C.3 (2).
As we proved earlier, $y_1$ and $y_2$ do not 
intersect, 
otherwise $\Delta_{K[n]}(t)$ would 
have a unit complex zero.
\qed

Now, none of the zeros of $\Delta_{K[n]}(t)$ is unit complex, but these zeros
seem to be distributed in a narrow strip containing the unit circle.
More precisely, we propose the following conjecture.

\begin{yosou}\label{conj:D.6} 
Let $\alpha= \frac{3- \sqrt{5}}{2}$. 
$\alpha$ is one of the real zeros of
$\Delta_{K[1]}(t)=t^2 -3t +1$.
Let $C$ be the circle with centre at $(\frac{\alpha -1}{2}, 0)$
and radius $\frac{\alpha+1}{2}$.
Then, for $n \geq 1$, all the zeros of
$\Delta_{K[n]}(t)$ with length $< 1$ lie in a narrow lunar domain bounded
by the unit circle and $C$.
\end{yosou}

\dyon

\begin{rem}\label{siki:D.7} 
If Hoste\rq s conjecture is true, then none of the zeros of the
Alexander polynomial of an alternating knot
is in the interior of the circle with centre at 
$(-1/2, 0)$ and radius $1/2$.
\end{rem}

\begin{ex}
Fig. C.5 below depicts the zeros for the cases $n=1,\dots, 10$ around
the unit circle.
\end{ex}

\dgo

\end{alphasection}

\end{document}